\documentclass[12pt]{article}

\usepackage{indentfirst}
\usepackage{amsfonts}
\usepackage{amsmath}
\usepackage{amssymb}
\usepackage{amsbsy, amsthm}
\usepackage[margin=1in]{geometry}

\newtheorem{dfn}{Definition} [section]
\newtheorem{theorem}[dfn]{Theorem}
\newtheorem{lemma}[dfn]{Lemma}

\newtheorem{corollary}[dfn]{Corollary}
\newtheorem{conjecture}[dfn]{Conjecture}
\newenvironment{pf}{\noindent{\bf Proof.}}
{\enspace\vrule height5pt depth0pt width5pt}

\def\X {{\mathcal X}}

\def\T {{\mathcal T}}

\def\Se {{\mathcal S}}

\def\F {{\mathcal F}}
\def\R {{\mathcal R}}
\def\P {{\mathcal P}}

\def\mf {{\rm mf}}
\def\ins {{\rm ins}}
\def\C {{\mathcal C}}
\def\Q {{\mathcal Q}}
\def\Z {{\mathbb Z}}
\def\M {{\mathcal M}}

\def\I {{\mathcal I}}
\def\W {{\mathcal W}}
\def\J {{\mathcal J}}
\def\A {{\mathcal A}}
\def\Y {{\mathcal Y}}
\def\HH {{\mathcal H}}

\begin{document}

\title{Packing Topological Minors Half-Integrally\footnote{This material is based upon work supported by the National Science Foundation under Grant No.~DMS-1664593, DMS-1929851 and DMS-1954054.}}
\author{Chun-Hung Liu\footnote{Email:chliu@math.tamu.edu. Partially supported by NSF under Grant No.~DMS-1664593, DMS-1929851 and DMS-1954054.} \\
\small Department of Mathematics, \\
\small Texas A\&M University,\\
\small College Station, TX 77843-3368, USA}
\maketitle

\begin{abstract}
The packing problem and the covering problem are two of the most general questions in graph theory.
The Erd\H{o}s-P\'{o}sa property characterizes the cases when the optimal solutions of these two problems are bounded by functions of each other.
Robertson and Seymour proved that when packing and covering $H$-minors for any fixed graph $H$, the planarity of $H$ is equivalent to the Erd\H{o}s-P\'{o}sa property.
Thomas conjectured that the planarity is no longer required if the solution of the packing problem is allowed to be half-integral. 

In this paper, we prove that this half-integral version of Erd\H{o}s-P\'{o}sa property holds for packing and covering $H$-topological minors, for any fixed graph $H$, which easily implies Thomas' conjecture. 
In fact, we prove an even stronger statement in which those topological minors are rooted at any choice of prescribed subsets of vertices.

A number of results on $H$-topological minor free or $H$-minor free graphs have conclusions or requirements tied to properties of $H$.
Classes of graphs that can half-integrally pack only a bounded number of $H$-topological minors or $H$-minors are more general topological minor-closed or minor-closed families whose minimal obstructions are more complicated than $H$.
Our theorem provides a general machinery to extend those results to those more general classes of graphs without losing their tight connections to $H$. 
\end{abstract}

\section{Introduction}

Numerous problems in graph theory and combinatorial optimization can be formulated as a packing problem or a covering problem.
Fix a family $\F$ of graphs, the packing problem asks for the maximum value $p$ such that the input graph $G$ contains $p$ pairwise disjoint subgraphs each isomorphic to a member of $\F$; the covering problem asks for the minimum value $c$ such that the input graph $G$ contains a set of at most $c$ vertices such that this set intersects every subgraph of $G$ isomorphic to some member of $\F$.
For example, when $\F=\{K_2\}$, the optimal value $p$ for the packing problem is the maximum size of a matching, and the optimal value $c$ of the covering problem is the minimum size of a vertex-cover.
The packing problem and the covering problems form a pair of dual integer programming problems.
It is clear that $p \leq c$, but $p \neq c$ in general.
One direction in approximation algorithm design is to find a relationship between $p$ and $c$ so that if any of $p$ or $c$ can be estimated efficiently, then the other one can be approximated efficiently.
A classical example is that when $\F=\{K_2\}$, $p' \leq p \leq c \leq 2p'$, where $p'$ is the size of any maximal matching, which can be found in linear time by a very simple greedy algorithm, so the minimum size of a vertex-cover can be approximated in linear time with factor 2 even though finding the minimum size of a vertex-cover is NP-hard.

A family $\F$ of graphs has the {\it Erd\H{o}s-P\'{o}sa property} if there exists a function $f$ with domain ${\mathbb N}$ such that for every graph $G$, either $G$ contains $k$ pairwise disjoint subgraphs where each of them is isomorphic to some member of $\F$, or there exists a subset $Z$ of $V(G)$ with $\lvert Z \rvert \leq f(k)$ such that $Z$ intersects every subgraph of $G$ isomorphic to some member of $\F$.
In other words, $p \leq c \leq f(p)$ if $\F$ has the Erd\H{o}s-P\'{o}sa property.
Therefore, one can approximate $p$ or $c$ by computing the optimal value of the LP-relaxation of the packing or covering problem.

Erd\H{o}s and P\'{o}sa \cite{ep} proved that the set of cycles has the Erd\H{o}s-P\'{o}sa property, and hence showed a relationship between the maximum number of disjoint cycles and the minimum size of a feedback vertex set.
Robertson and Seymour \cite{rs V} generalizes Erd\H{o}s and P\'{o}sa's theorem in terms of graph minors.
We say a graph $G$ contains another graph $H$ as a {\it minor} if $H$ can be obtained from a subgraph of $G$ by contracting edges.
Robertson and Seymour \cite{rs V} proved that $H$ is a graph such that the set of graphs containing $H$ as a minor has the Erd\H{o}s-P\'{o}sa property if and only if $H$ is a planar graph.
So Erd\H{o}s and P\'{o}sa's result follows from the case that $H$ is the one-vertex graph with one loop.
It also implies that some family of graphs does not have the Erd\H{o}s-P\'{o}sa property.

As shown above, for some family $\F$, the optimal value $c$ of the covering problem cannot be upper bounded by a function of the optimal value $p$ of the packing problem.
One attempt to remedy this situation is to relax the integral requirement. 
This leads to the following notion of half-integral packing. 

We say that a graph $G$ {\it half-integrally packs} $k$ (not necessarily pairwise non-isomorphic) graphs $H_1,H_2,...,H_k$ (for some positive integer $k$) if there exist subgraphs $G_1,G_2,...,G_k$ of $G$ such that 
	\begin{itemize}
		\item $G_i$ is isomorphic to $H_i$ for every $1 \leq i \leq k$, and
		\item for every $v \in V(G)$, there exist at most two indices $i$ with $1 \leq i \leq k$ such that $v \in V(G_i)$.
	\end{itemize}
The maximum value $k$ such that $G$ half-integrally packs $k$ members of $\F$ is indeed $2p_{1/2}$, where $p_{1/2}$ is the optimal value of the integer programming formulation of the packing problem with the integral condition relaxed to be half-integral.
Therefore, one can still approximate $c$ by computing the LP-relaxation of the packing problem or the covering problem, if the family of graphs has the Erd\H{o}s-P\'{o}sa property with respect to half-integral packing. 
	
Another motivation for half-integral packing arises from the following example showing that for every non-planar graph $H$, the set of graphs consisting of $H$-minors has no Erd\H{o}s-P\'{o}sa property.
Let $\Sigma$ be a surface with minimum genus in which $H$ can be embedded. 
For every positive integer $t$, consider $2t+1$ drawings of $H$ in $\Sigma$ such that every point of $\Sigma$ belongs to at most two drawings, and vertices of different drawings are disjoint.
Define $G_t$ to be the graph obtained from the union of those $2t+1$ drawings of $H$ by replacing each crossing generated by different drawings with a vertex of degree four.
By the minimality of the genus of $H$, every $H$-minor in $G_t$ partitions $\Sigma$ into disks.
So $G_t$ does not contain two disjoint $H$-minors.
On the other hand, every vertex of $G_t$ intersects at most two of those $2t+1$ drawings of $H$.
So every set of at most $t$ vertices of $G_t$ is disjoint from a drawing of $H$ and hence an $H$-minor in $G_t$.
This shows that the packing number in $G_t$ is one but the covering number in $G_t$ is greater that $t$, so the Erd\H{o}s-P\'{o}sa property does not hold.
However, $G_t$ half-integrally packs $2t+1$ $H$-minors.
So the Erd\H{o}s-P\'{o}sa property could still hold for sets of $H$-minors with non-planar $H$, if half-integral packing is considered.

Indeed, Thomas (see \cite{k}) conjectured that the planarity is not a necessary condition for having Erd\H{o}s-P\'{o}sa property with respect to the minor containment anymore if one relaxes packing to be half-integral packing.

\begin{conjecture}[See \cite{k}] \label{minor conj}
For every graph $H$, there exists a function $f$ such that for every graph $G$ and for every positive integer $k$, either $G$ half-integrally packs $k$ graphs where each of them contains $H$ as a minor, or there exists a set $Z \subseteq V(G)$ with $\lvert Z \rvert \leq f(k)$ such that $G-Z$ does not contain $H$ as a minor.
\end{conjecture}

Kawarabayashi \cite{k} proved Conjecture \ref{minor conj} when $H \in \{K_6,K_7\}$.
Norin \cite{n} announced a proof of Conjecture \ref{minor conj} in full but his proof was not published.

The main theorem of this paper is a strengthening of Conjecture \ref{minor conj} to topological minors.
A {\it subdivision} of a graph $H$ is a graph that can be obtained from $H$ by repeatedly subdividing edges.
Any vertex of a subdivision of $H$ not obtained by subdividing edges is called a {\it branch vertex}.
A {\it subdivision of $H$ in a graph $G$} is a subgraph of $G$ isomorphic to a subdivision of $H$.
For every collection $\R=(R_v: v \in V(H))$ of (not necessarily distinct) subsets of $V(G)$, an {\it $\R$-compatible subdivision of $H$ in $G$} is a subdivision of $H$ in $G$ whose branch vertex corresponding to $v$ belongs to $R_v$ for each $v \in V(H)$.
We say that a graph $G$ contains another graph $H$ as a {\it topological minor} (or an {\it $\R$-compatible topological minor}, respectively) if some subgraph of $G$ is a subdivision (or an $\R$-compatible subdivision, respectively) of $H$.

We remark that $\R$-compatible subdivisions are far-reaching generalizations of witnesses of positive instances of the classical $k$-Disjoint Path problem.
The $k$-Disjoint Path problem is a decision problem that asks that given an input graph $G$ and $k$ pairs $(s_1,t_1), (s_2,t_2),..., (s_k,t_k)$ of vertices of $G$, whether there exist $k$ pairwise disjoint paths $P_1,P_2,...,P_k$ in $G$ such that $P_i$ connects $s_i$ and $t_i$ for all $1 \leq i \leq k$.
Consider $kK_2$, which is the union of $k$ disjoint copies of $K_2$, with vertex-set $\{u_i,v_i: 1 \leq i \leq k\}$ and edge-set $\{u_iv_i: 1 \leq i \leq k\}$.
By choosing $R_{u_i}=\{s_i\}$ and $R_{v_i}=\{t_i\}$ for each $1 \leq i \leq k$, the union of any $k$ disjoint paths for the $k$-Disjoint Path problem is exactly an $\R$-compatible subdivision of $kK_2$.

The following is the main theorem of this paper.
We remark that our proof of this theorem does not rely on Norin's unpublished proof of Conjecture \ref{minor conj}.

\begin{theorem} \label{half EP}
For every graph $H$, there exists a function $f$ such that for every graph $G$, collection $\R=(R_v: v \in V(H))$ of subsets of $V(G)$, and positive integer $k$, either $G$ half-integrally packs $k$ graphs where each of them is an $\R$-compatible subdivision of $H$, or there exists a set $Z \subseteq V(G)$ with $\lvert Z \rvert \leq f(k)$ such that $G-Z$ does not contain any $\R$-compatible subdivision of $H$.
\end{theorem}

The following corollary immediately follows from Theorem \ref{half EP} if we choose $R_v=V(G)$ for every $v \in V(H)$.

\begin{corollary} \label{half EP non-rooted}
For every graph $H$, there exists a function $f$ such that for every graph $G$ and positive integer $k$, either $G$ half-integrally packs $k$ graphs where each of them is a subdivision of $H$, or there exists a set $Z \subseteq V(G)$ with $\lvert Z \rvert \leq f(k)$ such that $G-Z$ does not contain any subdivision of $H$.
\end{corollary}

We remark that Corollary \ref{half EP non-rooted} implies Conjecture \ref{minor conj}.

\bigskip

\noindent{\bf Proof of Conjecture \ref{minor conj} (assuming Corollary \ref{half EP non-rooted}):}
Let $\F$ be the family of graphs that can be obtained from $H$ by repeatedly splitting vertices with degree at least four.
Then $\F$ is a finite set of graphs such that any graph $G$ contains $H$ as a minor if and only if $G$ contains a subdivision of some graph in $\F$.
If $G$ does not half-integrally pack $k$ subgraphs where each of them contains $H$ as a minor, then for every graph $F \in \F$, $G$ does not half-integrally pack $k$ subdivisions of $F$, so Corollary \ref{half EP non-rooted} implies that for every $F \in \F$, there exist a function $f_F$ and a set $Z_F \subseteq V(G)$ with $\lvert Z_F \rvert \leq f_F(k)$ such that $G-Z_F$ does not contain any subdivision of $F$.
Therefore, $\bigcup_{F \in \F}Z_F$ is a set of size at most $\sum_{F \in \F}f_F(k)$ such that if $G$ does not half-integrally pack $k$ subgraphs where each of them contains $H$ as a minor, then $G-\bigcup_{F \in \F}Z_F$ does not contain any subdivision of $F$ for any $F \in \F$, and hence does not contain $H$ as a minor.
This proves that $\sum_{F \in \F}f_F$ is a function satisfying Conjecture \ref{minor conj}. 
$\Box$

\bigskip

We remark that unlike half-integral packing, the original packing problem with respect to topological minor containment does not behave nicely.
Answering a question raised by Robertson and Seymour \cite{rs V}, the author, Postle and Wollan \cite{lpw} provided a complete (but complicated) characterization for the graphs $H$ in which the set of graphs containing $H$ as a topological minor has the Erd\H{o}s-P\'{o}sa property, and proved that such a complicated characterization is unlikely to be avoidable in the sense that it is NP-hard to test whether an input graph $H$ has the property that the set of graphs containing $H$ as a topological minor has the Erd\H{o}s-P\'{o}sa property.
Hence the Erd\H{o}s-P\'{o}sa properties for minors respect to packing and half-integral packing only differ by the planarity, which is a property that can be test in linear time, but the gap between packing and half-integral packing for topological minors is larger, unless NP=P.

\subsection{Applications of Theorem \ref{half EP}} \label{sec: applications}

One consequence of Theorem \ref{half EP} or Corollary \ref{half EP non-rooted} is that one can modify a graph that does not half-integrally pack $k$ subdivisions of $H$ (or $H$-minors) into a graph containing no subdivision of $H$ (or $H$-minor) by deleting a bounded number of vertices, for any fixed graph $H$ and integer $k$.
Hence theorems on classes of graphs excluding a fixed graph $H$ as a topological minor or minor sensitive to the property of $H$ can be easily extended to classes of graphs that do not half-integrally pack $k$ subdivisions of $H$ or $H$-minors without losing its tight connection to $H$.
Note that the graphs with no subdivisions of $H$ (or $H$-minors) are exactly the graphs that do not half-integrally pack one subdivision of $H$ (or $H$-minor).

Note that any graph that does not half-integrally pack $k$ subdivisions of $H$ or $H$-minors does not contain a subdivision of $kH$ or a $kH$-minor, respectively, where $kH$ is a disjoint union of $k$ copies of $H$.
Hence results known to be true for any $H'$-subdivision-free or $H'$-minor-free graphs for any fixed graph $H'$ already imply results on graphs that do not half-integrally pack $k$ $H$-subdivisions or $H$-minors by simply taking $H'=kH$.
However, taking $H'=kH$ does not give effective conclusions in many situations since the conclusions of those theorems depend on parameters of $H'=kH$ which are significantly different from the corresponding parameters of $H$.

We provide some examples in this subsection to show how the aforementioned simple modification leads to new results without losing the its tight connection to $H$, and we expect much more similar applications of Theorem \ref{half EP} can be derived in a similar way.

\subsubsection{Partitioning $K_t$-topological minor free graphs into induced subgraphs}

Haj\'{o}s in 1940s conjectured that for every positive integer $t$, every graph with no subdivision of $K_{t+1}$ is properly $t$-colorable.
Equivalently, he conjectured that every graph with no subdivision of $K_{t+1}$ can be partitioned into $t$ induced subraphs whose every component has only one vertex.
This conjecture is stronger than the Four Color Theorem.
However, it is false for all sufficiently large $t$, since Erd\H{o}s and Fajtlowicz \cite{ef_hajos} showed that the number of required induced subgraphs is $\Omega(t^2/\log t)$.
On the other hand, the author and Wood \cite{lw} proved that $O(t)$ induced subgraphs are sufficient if one allows the components of each induced subgraph to have more than one vertex but still have bounded size, matching the order of a known lower bound $t$ for the number of required induced subgraphs.
This result can be strengthened by using Theorem \ref{half EP} as follows.

\begin{theorem}
For any positive integers $t \geq 2$ and $k$, there exists an integer $N$ such that every graph that does not half-integrally pack $k$ subdivisions of $K_{t+1}$ can be partitioned into $4t-4$ induced subgraphs whose every component contains at most $N$ vertices.
\end{theorem}

\begin{pf}
By Corollary \ref{half EP non-rooted}, there exists an integer $N_1$ (only depending on $t$ and $k$) such that for every graph $G$ that does not half-integrally pack $k$ subdivisions of $K_{t+1}$, there exists $Z \subseteq V(G)$ with $\lvert Z \rvert \leq N_1$ such that $G-Z$ does not contain a subdivision of $K_{t+1}$, so $G-Z$ can be partitioned into $4t-5$ induced subgraphs whose every component contains at most $N_2$ vertices by \cite[Theorem 6]{lw}, where $N_2$ only depends on $t$.
Together with the subgraph induced on $Z$, we obtain a partition of $G$ into $4t-4$ induced subgraphs whose every component contains at most $\max\{N_1,N_2\}$ vertices.
\end{pf}

\subsubsection{Structure theorems for excluding minors or topological minors}

A number of powerful structure theorems state that graphs that do not contain $H$ as a certain substructure are ``simpler'' than $H$ and hence provide a rough characterization of such graphs.

Classical results in this form were proved by Robertson and Seymour.
Their structure theorems for excluding a fixed graph as a minor are cornerstones of their Graph Minors series and lead to numerous applications in graph theory and theoretical computer science.
There are two versions of such structure theorems, where one \cite[Theorem (3.1)]{rs XVI} is stated in terms of tangles and the other \cite[Theorem (1.3)]{rs XVI} is stated in terms of tree-decompositions, and the former implies the latter.
(Tangles and tree-decompositions will be defined in Sections \ref{sec: using clique minors} and \ref{sec: surfaces}, respectively.)
In either version, it roughly says that for every graph $H$, every graph that does not contain $H$ as a minor can be obtained by gluing graphs that can be ``nearly embedded'' in surfaces in which $H$ cannot be embedded in a tree-like fashion.

Structure theorems for excluding a fixed graph as a topological minor were extensively explored recently.
Grohe and Marx \cite{gm} proved that if a graph $G$ does not contain a subdivision of $H$, then $G$ has a tree-decomposition such that every torso is either of ``nearly bounded maximum degree'' or ``nearly embeddable'' in a surface of bounded genus.
Though Grohe and Marx's theorem is sufficient to derive several nice algorithmic results, it does not provide enough structure information in the sense that the bound for the maximum degree and the genus in their theorem are much larger than the maximum degree and genus of $H$.
Dvo\v{r}\'{a}k \cite{d} improved their structure theorem by improving the genus part to an optimal statement.
The author and Thomas \cite{lt} further reduced the bound for the nearly maximum degree to be the maximum degree of $H$ minus one, which is optimal, though the theorem is stated in terms of tangles instead of tree-decompositions. 
In fact, the result in \cite{lt} works for $\R$-compatible subdivisions for a special kind of $\R$.

The improvement of the author and Thomas \cite{lt} on Dvo\v{r}\'{a}k's result allows us to conclude that if a graph $G$ does not contain a subdivision of $H$, then some parameter of $G$ is ``nearly'' less than the corresponding parameter of $H$, so an inductive argument can be applied. 
This theorem was also applied in their proof of a conjecture of Robertson on well-quasi-ordering \cite{l_thesis,lt_wqo} and the aforementioned result in \cite{lw}.

Theorem \ref{half EP} further strengthens the aforementioned results of Robertson and Seymour \cite{rs XVI}, Dvo\v{r}\'{a}k \cite{d} and the author and Thomas \cite{lt} by showing that the same structure holds for graphs that do not half-integrally pack $k$ $H$-minors or $k$ subdivisions of $H$.
We remark that such a strengthening does not immediately follow from the structure theorems for excluding $H'$-minors or $H'$-subdivisions by taking $H'=kH$, since the Euler genus of $kH$ is larger than $H$ when $H$ is non-planar.

We include a proof of the strengthening of the result of the author and Thomas \cite[Theorem 6.8]{lt} in this paper.
Results of Robertson and Seymour (\cite[Theorems (1.3) and (3.1)]{rs XVI}) and of Dvo\v{r}\'{a}k (\cite[Theorem 3]{d}) can be strengthened to graphs that do not half-integrally packing $k$ $H$-minors and $k$ subdivisions of $H$, respectively, in the same way.
Since formally stating those theorems requires definitions of several notions that will not be used in this paper, we omit the details.

For a graph $H$ and a surface $\Sigma$ in which $H$ can be drawn, $\mf(H,\Sigma)$ is defined to be the minimum size of a set $S$ of faces of a drawing of $H$ in $\Sigma$ such that every vertex of degree at least four in $H$ is incident with some face in $S$.
Other undefined notions mentioned in the following corollary can be found in Section \ref{sec: surfaces}.

\begin{corollary} \label{cor strengthened structure theorem}
Let $d \geq 4,h,k$ be positive integers.
Then there exist $\theta, \kappa, \rho, \xi, g \geq 0$ such that the following holds. 
If $H$ is a graph of maximum degree at most $d$ on $h$ vertices, $G$ is a graph and $X$ is a subset of $V(G)$ such that $G$ does not half-integrally pack $k$ subdivisions of $H$ whose branch vertices corresponding to vertices of degree at least four in $H$ are contained in $X$, then for every tangle $\T$ in $G$ of order at least $\theta$, there exists $Z \subseteq V(G)$ with $\lvert Z \rvert \leq \xi$ such that either
	\begin{enumerate}
		\item for every vertex $v \in V(G-Z) \cap X$, there exists $(A,B) \in \T-Z$ of order less than $d$ such that $v \in V(A)-V(B)$, or 
		\item there exists a $(\T-Z)$-central segregation $\Se = \Se_1 \cup \Se_2$ of $G-Z$ with $\lvert \Se_2 \rvert \leq \kappa$, having a proper arrangement in some surface $\Sigma$ of genus at most $g$ such that every society $(S_1,\Omega_1)$ in $\Se_1$ satisfies that $\lvert \overline{\Omega_1} \rvert \leq 3$, and every society $(S_2, \Omega_2)$ in $\Se_2$ is a $\rho$-vortex, and satisfies the following property: either
		\begin{enumerate}
			\item $H$ cannot be drawn in $\Sigma$, or
			\item $H$ can be drawn in $\Sigma$, $\mf(H,\Sigma) \geq 2$, and there exists $\Se'_2 \subseteq \Se_2$ with $\lvert \Se'_2 \rvert \leq \mf(H,\Sigma)-1$ such that for every $v \in V(G-Z) \cap X$, if there exists no $(A,B) \in \T-Z$ of order less than $d$ such that $v \in V(A)-V(B)$, then $v \in V(S)-\overline{\Omega}$ for some $(S,\Omega) \in \Se'_2$.
		\end{enumerate}
	\end{enumerate}
\end{corollary}

\begin{pf}
Let $\R=\{R_v: v \in V(H)\}$ be the collection such that $R_v=X$ for every vertex $v$ of $H$ with degree at least four, and $R_v=V(G)$ for other vertices $v$ of $H$.
So $G$ does not half-integrally pack $k$ $\R$-compatible subdivisions of $H$.
By Theorem \ref{half EP}, there exists $Z' \subseteq V(G)$ with $\lvert Z' \rvert \leq f(k)$ (for some function $f$ only depending on $H$) such that $G-Z'$ does not contain a subdivision of $H$ whose branch vertices corresponding to vertices of degree at least four in $H$ are contained in $X$.
Let $Z''$ and $\xi''$ be the set $Z$ and number $\xi$, respectively, mentioned in the statement of \cite[Theorem 6.8]{lt} when taking $G=G-Z'$, $X=X-Z'$ and $\T=\T-Z'$.
Then this theorem follows if we take $Z=Z' \cup Z''$.
Note that $\lvert Z \rvert \leq f(k)+\xi''$, and $f(k)+\xi''$ only depends on $H$ and $k$, and hence only depends on $d,h,k$.
\end{pf}

\subsubsection{Sparse universal graphs for apex-minor-free graphs}

In this subsection we generalize a result about apex-minor-free graphs to graphs that do not half-integrally pack many apex-minors.
We say that a graph $H$ is an {\it apex graph} if there exists $v \in V(H)$ such that $H-v$ is planar.
Note that $kH$ is not an apex graph if $H$ is a non-planar apex graph and $k \geq 2$.

For a class $\F$ of graphs and a positive integer $n$, an {\it $n$-universal graph for $\F$} is a graph that contains every $n$-vertex simple graph in $\F$ as a subgraph.
Universal graphs are motivated by applications in VLSI design and parallel computing.
Clearly, $K_n$ is an $n$-universal graph for any graph class.
So the question is how sparse an $n$-universal graph can be, for a fixed class $\F$.

It is known that every $n$-universal graph for the class of forests has $\Omega(n\log n)$ edges, as it must contain, for each integer $t$ with $1 \leq t \leq n$, at least $t$ vertices of degree at least $\lfloor n/t \rfloor-1$, in order to contain the forest formed by $t$ disjoint copies of $K_{1,\lfloor n/t \rfloor-1}$ and isolated vertices as a subgraph.
Using separator theorems, Babai, Chung, Erd\H{o}s, Graham and Spencer \cite{bcegs} and Chung \cite{c} proved that there exist $n$-universal graphs with $O(n^{3/2})$ edges for the class of planar graphs and for any proper minor-closed family, respectively.
It remains open whether there exists an $n$-universal graph with $O(n^{1+o(1)})$ edges for any fixed proper minor-closed family.
The currently most general result in this direction is due to Esperet, Joret and Morin \cite{ejm} who proved that there exists an $n$-universal graph with $O(n^{1+o(1)})$ edges for the class of $H$-minor free graphs, for any fixed apex graph $H$.
Their proof relies on a structural theorem for apex-minor free graphs.
Corollary \ref{half EP non-rooted} extends their result to more general minor-closed families.

\begin{corollary}
Let $H$ be an apex graph, and let $k$ be a positive integer.
If $\F$ is the class of graphs that do not half-integrally pack $k$ $H$-minors, then there exists an $n$-universal graph for $\F$ on $(1+o(1))n$ vertices with $n \cdot 2^{O(\sqrt{\log n \cdot \log\log n})}$ edges.
\end{corollary}

\begin{pf}
Let $H'$ be a graph obtained from $H$ by adding edges to make every isolated vertex of $H$ have degree 1 in $H'$.
Note that $H'$ is also an apex graph.
Let $\F'$ be the class of all $H'$-minor free graphs.
Let $n$ be a positive integer.
An immediate corollary of a combination of \cite[Corollary 40]{djmmuw} and \cite[Theorem 2]{ejm} implies that there exists an $n$-universal graph $G_n$ for $\F'$, where $G_n$ has $(1+o(1))n$ vertices and $n \cdot 2^{O(\sqrt{\log n \cdot \log\log n})}$ edges.
Since $H'$ has no isolated vertices, one can add isolated vertices into an $H'$-minor free graph on less than $n$ vertices to obtain an $n$-vertex $H'$-minor free graph.
So $G_n$ contains every simple graph in $\F'$ on at most $n$ vertices as a subgraph.

Let $\F_k'$ be the class of graphs that do not half-integrally pack $k$ $H'$-minors.
Recall that Corollary \ref{half EP non-rooted} implies Conjecture \ref{minor conj}.
So there exists an integer $N$ (only depending on $H'$ and $k$) such that for every graph $G$ in $\F_k'$, there exists $Z \subseteq V(G)$ with $\lvert Z \rvert \leq N$ such that $G-Z \in \F'$.
Define $G_n'$ to be the graph obtained from a disjoint union of $G_n$ and $K_N$ by adding all edges between $G_n$ and $K_N$.
Note that $\lvert V(G'_n) \rvert = \lvert V(G_n) \rvert + N = (1+o(1))n$ and $\lvert E(G'_n) \rvert \leq \lvert E(G_n) \rvert + N^2 + N\cdot \lvert V(G_n) \rvert = n \cdot 2^{O(\sqrt{\log n \cdot \log\log n})} + N^2 + N(1+o(1))n = n \cdot 2^{O(\sqrt{\log n \cdot \log\log n})}$.

Since $\F \subseteq \F_k'$, it suffices to show that $G_n'$ is an $n$-universal graph for $\F_k'$.
If $G$ is an $n$-vertex simple graph in $\F_k'$, then there exists $Z \subseteq V(G)$ with $\lvert Z \rvert \leq N$ such that $G-Z \in \F'$, so $G_n$ contains $G-Z$ as a subgraph and $K_N$ contains $G[Z]$ as a subgraph.
Since $G_n'$ contains all edges between $K_N$ and $G_n$, $G'_n$ contains $G$ as a subgraph.
\end{pf}

\section{Proof sketch of Theorem \ref{half EP} and organization of the paper}
The rest of the paper is dedicated to a proof of Theorem \ref{half EP}.
We sketch the proof in this section and explain the organization of this paper.

\subsection{Reducing to the connected case}
The proof of Theorem \ref{half EP} follows from an induction on the number of components of $H$.
The main difficulty of the proof is the base case; namely the case that $H$ is connected.
The majority of the paper devotes to a proof of Lemma \ref{connected half EP} which settles the connected case of Theorem \ref{half EP}.
The proof of the general case of Theorem \ref{half EP} uses a simple matching technique and is included in the proof of Theorem \ref{final half EP restatement}.

\subsection{Solving the connected case}
In the rest of this section, we assume that $H$ is a connected graph and sketch the proof of Lemma \ref{connected half EP}.
To simplify the sketch, we only consider the non-rooted version in this section.
That is, we assume that $R_v=V(G)$ for every $v \in V(H)$.
So every $\R$-compatible subdivision of $H$ in $G$ is simply a subdivision of $H$ in $G$.

We shall prove Lemma \ref{connected half EP} by induction on $k$.
Namely, we prove that for each $k$, there exists a number $f(k)$ such that either $G$ half-integrally packs $k$ subdivisions of $H$, or there exists $Z \subseteq V(G)$ with $\lvert Z \rvert \leq f(k)$ such that $G-Z$ does not contain any subdivision of $H$.
Now we assume that $G$ does not half-integrally pack $k$ subdivisions of $H$ and show the existence of $f(k)$.

\subsubsection{Step 1: defining a tangle}
A {\it separation} of a graph $G$ is a pair $(A,B)$ of edge-disjoint subgraphs such that $A \cup B=G$.
The {\it order} of $(A,B)$ is $\lvert V(A \cap B) \rvert$.

Fix a positive integer $\theta$.
Suppose that there exists a separation $(A,B)$ of $G$ of order less than $\theta$ such that each of $A-V(B)$ and $B-V(A)$ contains a subdivision of $H$.
Since $G$ does not half-integrally pack $k$ subdivisions of $H$, each of $B-V(A)$ and $A-V(B)$ does not half-integrally pack $k-1$ subdivisions of $H$.
By the induction hypothesis, there exist $Z_A \subseteq V(A)-V(B)$ and $Z_B \subseteq V(B)-V(A)$ with $\lvert Z_A \rvert \leq f(k-1)$ and $\lvert Z_B \rvert \leq f(k-1)$ such that $A-(V(B) \cup Z_A)$ and $B-(V(A) \cup Z_B)$ does not contain any subdivision of $H$.
Since $H$ is connected, any subdivision of $H$ in $G-(Z_A \cup Z_B)$ must intersect $V(A \cap B)$.
Hence, $G-(Z_A \cup Z_B \cup V(A \cap B))$ does not contain a subdivision of $H$.
The size of $\lvert Z_A \cup Z_B \cup V(A \cap B) \rvert$ is at most $2f(k-1)+\theta$.
Therefore, if such a separation exists, then we are done by setting $f(k) \leq 2f(k-1)+\theta$.

So we may assume that for every separation $(A,B)$ of order less than $\theta$, either $A-V(B)$ or $B-V(A)$ does not contain a subdivision of $H$.
Note that if neither $A-V(B)$ nor $B-V(A)$ contains a subdivision of $H$, then $G-V(A \cap B)$ does not contain a subdivision of $H$, and we are done by setting $f(k) \leq \theta$.

Therefore, we may assume that for every separation $(A,B)$ of order less than $\theta$, exactly one of $A-V(B)$ and $B-V(A)$ contains a subdivision of $H$.
This allows us to define a tangle $\T$ in $G$ of order $\theta$.
The notion of tangles is widely exploited in the study of graph minors.
Intuitively, a tangle of order $\theta$ indicates which one among $A$ and $B$ is the important side of each separation $(A,B)$ of order less than $\theta$.
In our application, $B$ is the important side of $(A,B)$ if $B-V(A)$ contains a subdivision of $H$.

Formally, a tangle of order $\theta$ is a collection of separations of order less than $\theta$ and will be precisely defined in Section \ref{sec: using clique minors}.
We define a collection $\T$ of separations of order less than $\theta$ such that for every separation $(A,B)$ of $G$ of order less than $\theta$, $(A,B) \in \T$ if and only if $B-V(A)$ contains a subdivision of $H$.
Then one can verify that $\T$ is a tangle in $G$ of order $\theta$.

\subsubsection{Step 2: using large clique minors}

We shall prove Lemma \ref{connected half EP} by considering two cases.
One case is that $G$ contains a large complete graph as a minor such that for every separation $(A,B) \in \T$, $V(B)$ intersects some branch set; the other case is that no such complete graph minor exists.
We deal with the first case in this subsection, and the complete proof of this case is included in Section \ref{sec: using clique minors}.
This case is relatively simpler than the other case.

By applying \cite[Theorem 6.8]{lt} by taking $H=kH$, we know that there exists $Z_0 \subseteq V(G)$ with bounded size such that for every $v \in V(G)$, exists a separation $(A_v,B_v) \in \T-Z_0$ of order less than the maximum degree of $H$ such that $v \in V(A_v)-V(B_v)$.
(Here $\T-Z_0$ is a tangle in $G-Z_0$ which is the collection of separations of $G-Z_0$ of small order ``consistent'' with $\T$.)
In particular, for every branch vertex $v$ of a subdivision of $H$ in $G-Z_0$, there exists $(A_v,B_v) \in \T-Z_0$ of small order such that $v \in V(A_v)-V(B_v)$.

Hence, for every subdivision $S$ of $H$ in $G-Z_0$, there exists a separation $(A_S,B_S) \in \T-Z_0$ of small order such that every branch vertex of $S$ is contained in $V(A_S)-V(B_S)$, where $(A_S,B_S)=(\bigcup_v A_v, \bigcap_v B_v)$ and the union and the intersection are over all branch vertices $v$ of $S$.
By using a result in \cite{rs XIII}, we show that either there exists $(A'_S,B'_S) \in \T-Z_0$ of order smaller than $(A_S,B_S)$ with $A_S \subseteq A'_S$, or there exist disjoint paths in $B_S$ connecting pairs of vertices in $V(A_S \cap B_S)$ such that these paths together with $S \cap A_S$ form a subdivision of $H$ in $G-Z_0$.
Since $G-Z_0$ does not contain $k$ disjoint subdivisions of $H$, it is impossible that there exist $k$ distinct subdivisions $S_1,S_2,...,S_k$ in $G-Z_0$ in which the later case holds for each $S_i$ such that for all distinct $i$ and $j$, $A_{S_i}$ is disjoint from $A_{S_j}$, and the disjoint paths in $B_{S_i}$ for $S_i$ are disjoint from the paths in $B_{S_j}$ for $S_j$.
By applying a result in \cite{l}, we can further delete a subset $Z_1$ of $V(G)$ with bounded size such that for every subdivision $S$ of $H$ in $G-Z_0$, either $S \cap Z_1 \neq \emptyset$, or $S$ is a subdivision of $H$ in $G-(Z_0 \cup Z_1)$ and there exists $(A_S',B_S') \in \T-(Z_0 \cup Z_1)$ of order less than $(A_S,B_S)$ such that $A_S-Z_1 \subseteq A_S'$.
In particular, for every subdivision $S$ of $H$ in $G-(Z_0 \cup Z_1)$ (and hence in $G-Z_0$), there exists a separation $(A'_S,B'_S) \in \T-(Z_0 \cup Z_1)$ of order strictly less than the order of $(A_S,B_S)$ such that $A'_S$ contains all branch vertices of $S$.

Therefore, by repeatedly applying this argument, there exists a set $Z \subseteq V(G)$ of bounded size such that either $G-Z$ has no subdivision of $H$, or for every subdivision $S$ in $G-Z$, there exists $(A_S,B_S) \in \T-Z$ of order 0 such that all branch vertices of $S$ are contained in $A_S$.
If the later case holds, then $S \subseteq A_S$ since $H$ is connected and $(A_S,B_S)$ has order 0, but it contradicts the definition of $\T$.
Hence the former case holds and we are done.

\subsubsection{Step 3: using nearly embedding}
By Step 2, we may assume that $G$ does not contain a large clique minor mentioned in the previous case. 
The rest of the paper deals with this case.
This case is much more complicated than the previous one and involves several technical definitions for new techniques and developed machinery, so we only include a very high-level sketch for this case.

By Robertson and Seymour's structure theorem for excluding such a clique minor \cite{rs XVI}, there exist $Z \subseteq V(G)$ with small size, a surface $\Sigma$ of bounded genus and disks with pairwise disjoint interior such that we can map the vertices and edges of $G-Z$ to points in the union of those disks such that no edges of $G-Z$ has ends in the interior of different disks, and all but a bounded number of disks have at most 3 vertices in the boundary of the disks.
Furthermore, for the disks having at least 4 vertices in the boundary, the subgraph of $G-Z$ contained in this disk (called a ``vortex'') has a path-decomposition with bounded adhesion.
(Vortices and path-decompositions will be formally defined in Section \ref{sec: surfaces}.)
Then for each disk not corresponding to a vortex, connect the vertices on the boundary by the curve in the boundary of the disk to obtain a drawing $G'$ in $\Sigma$.
Furthermore, one can define a metric $m_\T$ on the vertices of $G'$ with respect to $\T$.
Roughly speaking, any two vertices of $G'$ that are far apart with respect to $m_\T$ if and only if every closed curve in $\Sigma$ containing both of them in the ``inside'' of the region determined by the curve must intersect $G'$ many times.

Section \ref{sec: surfaces} is dedicated to the formal setup for the above setting and proving some results about drawing in surfaces that will be used in later sections.

Each subdivision of $H$ in $G-Z$ defines a ``drawing'' in $\Sigma$ where there are possibly some crossings inside the disks that were chosen to accommodate vertices and edges of $G-Z$.
(Such a ``drawing'' motivates the notion of ``pseudo-embeddings'' that will be formally defined in Section \ref{sec: pseudoembedding}.)
Since every disk that is not corresponding to a vortex has only at most 3 vertices in the boundary, if there exists a crossing inside a disk not corresponding to a vortex, then this disk contains a branch vertex of this subdivision.
Therefore, for each subdivision $S$ of $H$ in $G-Z$, only a bounded number of disks contain crossings, and we call such a disk a ``crossing disk'' for $S$.
Note that it implies that we can add some pairwise disjoint curves in $\Sigma$ disjoint from the interior of all crossing disks for $S$ to connect some pairs of vertices in the boundary of the crossing disks such that those curves together with the subgraph of $S$ contained in those crossing disks form a topological embedding of $H$ in $\Sigma$ with crossings only in crossing disks.
A theorem of Robertson and Seymour \cite{rs XII} (Theorem \ref{linkage on surface} in this paper) implies that if 
	\begin{itemize}
		\item the number of curves that we are supposed to add is small,
		\item there are pairwise disjoint paths in $G-Z$ from the endpoints of those curves to some vertices far away with respect to $m_\T$, and 
		\item the endpoints of those curves belonging to different crossing disks are pairwise far apart with respect to the metric $m_\T$, 
	\end{itemize}
then there exist pairwise disjoint paths in $G-Z$ connecting the pairs of the endpoints of the curves that we are supposed to add, and hence the subgraph of $S$ contained in the union of the crossing disks can be extended to another subdivision of $H$ in $G-Z$.

A similar technique shows that if there exist $k$ distinct subdivisions $S_1,...,S_k$ of $H$ in $G$ such that the subgraphs of $S_1,...,S_k$ contained in their crossing disks are pairwise disjoint, and if the curves that are supposed to be added can be realized by disjoint paths in $G-Z$, then $G-Z$ half-integrally packs $k$ subdivisions of $H$.
Since $G$ does not half-integrally pack $k$ subdivisions of $H$, one of the following obstructions appears.
	\begin{itemize}
		\item[Obstruction 1:] There exists some $S_i$ such that the number of curves that are supposed to be added for $S_i$ is too large.
		\item[Obstruction 2:] There exists some $S_i$ such that there are no disjoint paths in $G-Z$ from the endpoints of those curves to vertices far away with respect to $m_\T$.
		\item[Obstruction 3:] There exist $S_1,...,S_k$ such that some of the crossing disks for $S_1,...,S_k$ are nearby with respect to $m_\T$.
	\end{itemize}

Section \ref{sec: pseudoembedding} aims to overcome Obstruction 1.
We use a well-quasi-ordering technique and a disentanglement process to show that as long as the subgraph of $G-Z$ contained in the disks corresponding to vortices has a path-decomposition of bounded adhesion, it is sufficient to add a bounded number of curves for each $S_i$ to extend the subgraph of $S_i$ contained in its crossing disks to a subdivision of $H$.
More details and intuitions will be included in Section \ref{subsec:intuition_pseudo_embed}.

Sections \ref{sec: gauges}, \ref{sec: rerouting} and \ref{sec: vortices} are dedicated to developing tools to solve Obstructions 2 and 3.
Roughly speaking, they allow us to enlarge and merge the crossing disks to obtain a new set of crossing disks such that the union of these new disks includes more branch vertices of the subdivision if no such disjoint paths in $G-Z$ exist.
Lemmas proved in Section \ref{sec: surfaces} ensure that the subgraphs of $S_i$ contained in those new crossing disks remain under control.
The number of times for applying this process is bounded since there are only finitely many branch vertices of $S_i$.
In Section \ref{sec: gauges}, we will show that there exists a ``gauge'' that regularizes the behavior of all subdivisions of $H$ at once, which is a preparation for enlarging and merging crossing-disks.
Formal arguments about enlarging and merging those crossing disks are stated in Section \ref{sec: vortices}, with the help of the methods developed in Section \ref{sec: rerouting}.

Now it remains to deal with Obstruction 3.
Using Lemma \ref{sweeping balls into vortices}, for each $S_i$, we can repeatedly merge crossing disks for $S_i$ to obtain a new set of crossing disks such that the subgraphs of $S_i$ contained in those new crossing disks remain under control, until the new crossing disks for $S_i$ are pairwise far apart.
Note that this process will terminate since the number of crossing disks for each $S_i$ is bounded.
By further applying this technique and the method for dealing with Obstruction 2, we can show that for each $S_i$, one can find a new set of crossing disks such that the union of the interior of these new crossing disks contains at least one more branch vertex of $S_i$.
The details are included in Section \ref{sec: vortices}.
Note that this process can only be applied a bounded number of times since the number of branch vertices of each $S_i$ is bounded.

Therefore, none of Obstructions 1-3 can exist.
Hence we can construct a half-integral packing of $k$ subdivisions of $H$ in $G$ to obtain a contradiction, unless $G-Z$ has very few distinct subdivisions of $H$.
So all those subdivisions can be killed by deleting few vertices.
The proof of Theorem \ref{half EP} will be completed in Section \ref{sec: half-integral EP}.

\section{Notations}
We define some notions and notations that will be frequently used in this paper.
Graphs in this paper are allowed to have loops or parallel edges.
For every positive integer $k$, we denote the set $\{1,2,...,k\}$ by $[k]$.
For a function $f$ and a subset $X$ of its domain, $f(X)$ denotes the set $\{f(x): x \in X\}$; for a function $f$ and a sequence $s=(s_1,s_2,...,s_k)$, $f(s)$ denotes the sequence $(f(s_1),f(s_2),...,f(s_k))$.

The following alternative definition of subdivisions is useful.
Let $G,H$ be graphs.
A {\it homeomorphic embedding} from $H$ into $G$ is a function with domain $V(H) \cup E(H)$ such that the following hold.
\begin{itemize}
	\item $\pi$ maps vertices of $H$ to vertices of $G$ injectively.
	\item $\pi$ maps each loop of $H$ with end $v$ to a cycle in $G$ containing $\pi(v)$ and maps each non-loop edge with ends $u,v$ of $H$ to a path in $G$ from $\pi(u)$ to $\pi(v)$.
	\item If $e,f$ are two distinct edges of $H$, then $\pi(e) \neq \pi(f)$ and $V(\pi(e) \cap \pi(f)) = \pi(e \cap f)$.
	\item If $v$ is a vertex of $H$ not incident with an edge $e$ of $H$, then $\pi(v) \not \in V(\pi(e))$.
\end{itemize}
Clearly, the image of any homeomorphic embedding $\pi$ from $H$ into $G$ is a subdivision of $H$ in $G$ whose set of branch vertices is $\pi(V(H))$.
Let $\R=(R_v: v \in V(H))$, where $R_v$ is a subset of $V(G)$ for each $v \in V(H)$.
Then an {\it $\R$-compatible homeomorphic embedding} from $H$ into $G$ is a homeomorphic embedding $\pi$ from $H$ into $G$ such that $\pi(v) \in R_v$ for every $v \in V(H)$.
Note that the image of any $\R$-compatible homeomorphic embedding from $H$ into $G$ is an $\R$-compatible subdivision of $H$ in $G$.

\section{Using complete graph minors} \label{sec: using clique minors}
The goal of this section is to prove Lemma \ref{with clique minor} which deals with the half-integral packing in a graph that has a nice complete graph minor.

Given a graph $H$, an $H$-{\it minor} of a graph $G$ is a map $\alpha$ with domain $V(H) \cup E(H)$ such that the following hold.
\begin{itemize}
	\item For every $h \in V(H)$, $\alpha(h)$ is a nonempty connected subgraph of $G$.
	\item If $h_1$ and $h_2$ are different vertices of $H$, then $\alpha(h_1)$ and $\alpha(h_2)$ are disjoint.
	\item For each edge $e$ of $H$ with ends $h_1,h_2$, $\alpha(e)$ is an edge of $G$ with one end in $\alpha(h_1)$ and one end in $\alpha(h_2)$; furthermore, if $h_1=h_2$, then $\alpha(e) \in E(G)-E(\alpha(h_1))$ and has all ends in $\alpha(h_1)$.
	\item If $e_1, e_2$ are two different edges of $H$, then $\alpha(e_1) \neq \alpha(e_2)$.
\end{itemize}
We say that {\it $G$ contains an $H$-minor} if such a function $\alpha$ exists.

Recall that a separation of a graph $G$ is a pair $(A,B)$ of edge-disjoint subgraphs with $A \cup B=G$, and the order of $(A,B)$ is $\lvert V(A) \cap V(B) \rvert$.
A {\it tangle} $\T$ in $G$ of {\it order} $\theta$ is a set of separations of $G$, each of order less than $\theta$ such that
\begin{enumerate}
	\item[(T1)] for every separation $(A,B)$ of $G$ of order less than $\theta$, either $(A,B) \in \T$ or $(B,A) \in \T$;
	\item[(T2)] if $(A_1, B_1), (A_2,B_2), (A_3,B_3) \in \T$, then $A_1 \cup A_2 \cup A_3 \neq G$;
	\item[(T3)] if $(A,B) \in \T$, then $V(A) \neq V(G)$.
\end{enumerate}
Furthermore, for $Z \subseteq V(G)$ with $\lvert Z \rvert<\theta$, we define $\T-Z$ to be the set of all separations $(A',B')$ of $G-Z$ of order less than $\theta-\lvert Z \rvert$ such that there exists $(A,B) \in \T$ with $Z \subseteq V(A \cap B)$, $A'=A-Z$ and $B'=B-Z$.
It is proved in \cite[Theorem 8.5]{rs X} that $\T-Z$ is a tangle in $G-Z$ of order $\theta-\lvert Z \rvert$.
A tangle $\T$ in $G$ {\it controls} an $H$-minor $\alpha$ if $\alpha$ is an $H$-minor such that there do not exist $(A,B) \in \T$ of order less than $\lvert V(H) \rvert$ and $h \in V(H)$ such that $V(\alpha(h)) \subseteq V(A)$.

The following theorem is useful.

\begin{theorem}[{\cite[Theorem (5.4)]{rs XIII}}] \label{vertex linkage}
Let $G$ be a graph, and let $Z$ be a subset of $V(G)$ with $\lvert Z \rvert = \xi$.
Let $k \geq \lfloor \frac{3}{2} \xi \rfloor$, and let $\alpha$ be a $K_k$-minor in $G$.
If there is no separation $(A,B)$ of $G$ of order less than $\lvert Z \rvert$ such that $Z \subseteq V(A)$ and $A \cap \alpha(h)=\emptyset$ for some $h \in V(K_k)$, then for every partition $(Z_1,...,Z_n)$ of $Z$ into non-empty subsets, there are $n$ connected graphs $T_1, ..., T_n$ of $G$, pairwise disjoint and $V(T_i) \cap Z = Z_i$ for $1 \leq i \leq n$.
\end{theorem}

Let $G$ be a graph and $\T$ a tangle in $G$.
We say that a subset $X$ of $V(G)$ is {\it free} with respect to $\T$ if there exists no $(A,B) \in \T$ of order less than $\lvert X \rvert$ such that $X \subseteq V(A)$.

The following lemma is proved in \cite{l}.

\begin{lemma}[{\cite[Lemma 3.2]{l}}] \label{stronger mw spider}
Let $G$ be a graph and $\T$ a tangle in $G$ of order $\theta$, and let $c$ be a positive integer.
For every $1 \leq i \leq c$, let $d_i,k_i$ be positive integers, and let $\{X_{i,j} \subseteq V(G): j \in J_i\}$ be a family of subsets of $V(G)$ indexed by a set $J_i$.
Let $d,k$ be integers such that $\theta \geq (kcd)^{d+1}+d$, $d_i \leq d$ and $k_i \leq k$ for $1 \leq i \leq c$.
Let $J_i^* \subseteq J_i$ with $\lvert J_i^* \rvert \leq k_i$ for each $1 \leq i \leq c$, such that $\bigcup_{i=1}^c \bigcup_{j \in J_i^*} X_{i,j}$ is free with respect to $\T$ and $X_{i,j} \cap X_{i',j'} = \emptyset$ for distinct pairs $(i,j),(i',j')$ with $1 \leq i \leq i' \leq c$, $j \in J_i^*$ and $j' \in J_{i'}^*$.
If $\lvert X_{i,j} \rvert \leq d_i$ for every $1 \leq i \leq c$ and $j \in J_i$, then either 
	\begin{enumerate}
		\item there exist $J_1',J_2',...,J_c'$ with $J_i^* \subseteq J_i' \subseteq J_i$ and $\lvert J_i' \rvert = k_i$ for each $1 \leq i \leq c$ such that $\bigcup_{1 \leq i \leq c} \bigcup_{j \in J'_i} X_{i,j}$ is free with respect to $\T$, and $X_{i,j} \cap X_{i',j'} = \emptyset$ for all distinct pairs $(i,j),(i',j')$ with $1 \leq i \leq i' \leq c'$, $j \in J_i'$ and $j'\in J_{i'}'$, or
		\item there exist $Z \subseteq V(G)$ with $\lvert Z \rvert \leq (kcd)^{d+1}$ and integer $i$ with $1 \leq i \leq c$ and $\lvert J_i^* \rvert < k_i$ such that for every $j \in J_i$, either $X_{i,j} \cap Z \neq \emptyset$, or $X_{i,j}$ is not free with respect to $\T-Z$.
	\end{enumerate}
\end{lemma}

We need another lemma about the freeness.

\begin{lemma} \label{A side disjoint}
Let $\theta$ be a positive integer.
Let $G$ be a graph and $\T$ a tangle in $G$ of order at least $\theta$.
Let $\Q \subseteq \T$ such that $V(A) \cap V(B)$ is disjoint from $V(A') \cap V(B')$ for distinct $(A,B),(A',B') \in \Q$.
Assume that for every $(A,B) \in \Q$ and vertex $v \in V(A)-V(B)$, there exists a path in $G[A]$ from $v$ to $V(A) \cap V(B)$.
If $\bigcup_{(A,B) \in \Q} (V(A) \cap V(B))$ is free with respect to $\T$ and has size at most $\theta$, then $V(A) \cap V(A') = \emptyset$ for distinct $(A,B),(A',B') \in \Q$.
\end{lemma}

\begin{pf}
Denote $\Q=\{(A_1,B_1),...,(A_k,B_k)\}$.
For each $i$ with $1 \leq i \leq k$, define $X_i=V(A_i) \cap V(B_i)$.
We first assume that $X_1 \cap V(A_2) \neq \emptyset$.
Since $X_{1} \cap X_{2} = \emptyset$, there exists a vertex $v$ in $X_1 \cap V(A_2)-V(B_2)$.
Let $A^* = \bigcup_{i=1}^k A_i$ and $B^* = \bigcap_{i=1}^k B_i$.
Note that $v \not \in V(A^*) \cap V(B^*)$, so the order of $(A^*,B^*)$ is less than $\sum_{i=1}^k \lvert X_i \rvert$.
By (T1) and (T2), $(A^*,B^*) \in \T$.
Since $A^*$ contains $\bigcup_{i=1}^k X_i$ and $(A^*,B^*)$ has order less than $\sum_{i=1}^k \lvert X_i \rvert = \lvert \bigcup_{i=1}^k X_i \rvert$, $\bigcup_{i=1}^k X_i$ is not free with respect to $\T$, a contradiction.
Hence, $X_1 \cap V(A_2) = \emptyset$.
That is, $X_1 \subseteq V(B_2)-V(A_2)$.
Similarly, $X_i \subseteq V(B_j)-V(A_j)$ for every distinct $i,j \in [k]$.

Now suppose there exist distinct $i,j \in [k]$ such that $V(A_i) \cap V(A_j) \neq \emptyset$.
Let $w \in V(A_i) \cap V(A_j)$.
By the assumption, there exists a path $P$ in $G[A_i]$ from $w$ to $X_i$.
Since $X_i \subseteq V(B_j)-V(A_j)$, $P$ contains a vertex $u$ in $X_j$.
So $u \in X_j \cap V(A_i) \neq \emptyset$, a contradiction.
This proves the lemma.
\end{pf}

\bigskip

The following lemma is proved in \cite{lt}.

\begin{lemma}[{\cite[Theorem 3.4]{lt}}] \label{spider tangle}
Let $h$ and $d$ be positive integers.
Let $G$ be a graph, and let $S$ be a subset of vertices of degree at least $d-1$ in $G$.
Let $\T$ be a tangle in $G$ of order $\theta$.
If $\theta \geq (hd+1)^{d+1}+d$, then either 
	\begin{enumerate}
		\item there exist $h$ vertices $v_1,v_2,...,v_h \in S$ and $h$ pairwise disjoint subsets $X_1,X_2,...,X_h$ of $V(G)$, where $X_i$ consists of $v_i$ and $d-1$ neighbors of $v_i$ for each $i \in [h]$, such that $\bigcup_{i=1}^h X_i$ is free in $\T$, or
		\item there exists a set $C \subseteq V(G)$ with $\lvert C \rvert \leq (hd+1)^{d+1}$ such that for every $v \in S-C$, there exists $(A,B) \in \T-C$ of order less than $d$ such that $v \in V(A)-V(B)$.
	\end{enumerate}
\end{lemma}

Let $G$ and $H$ be graphs.
Let $\pi$ be a homeomorphic embedding from $H$ into $G$.
Let $(A,B)$ be a separation of $G$.
We define the {\it shade of $\pi$ with respect to $(A,B)$} to be the subgraph $S$ of $H$ with $V(S)=\{v \in V(H): \pi(v) \in V(A)-V(B)\}$ and $E(S)=\{e \in E(H): V(\pi(e)) \subseteq V(A)\}$.

The following lemma is the heart of the proof of the main result (Lemma \ref{with clique minor}) of this section.

\begin{lemma} \label{clique_minor_uniform} 
For any connected graph $H$ and positive integers $k,q,\xi_0$, there exist integers $\theta,t,\xi$ such that if $S$ is a (possibly empty) subgraph of $H$, $G$ is a graph, $\R=(R_v: v \in V(H))$ is a collection of subsets of $V(G)$, $\T$ is a tangle in $G$ of order at least $\theta$ controlling a $K_t$-minor, a subset $Z_0$ of $V(G)$ with $\lvert Z_0 \rvert \leq \xi_0$, a collection $\F$ of $\R$-compatible homeomorphic embeddings from $H$ into $G-Z_0$, and a collection $\Q \subseteq \T-Z_0$ of separations of order at most $q$ such that
	\begin{itemize}
		\item $G-Z_0$ does not contain $k$ disjoint $\R$-compatible subdivisions of $H$, and
		\item for every $\pi \in \F$, there exists $(A_\pi,B_\pi) \in \Q$ such that the shade of $\pi$ with respect to $(A_\pi,B_\pi)$ equals $S$,
	\end{itemize}
then there exists $Z$ with $Z_0 \subseteq Z \subseteq V(G)$ of size at most $\xi$ such that for every $\pi \in \F$, either 
	\begin{enumerate}
		\item $\pi(V(H) \cup E(H)) \cap Z \neq \emptyset$, or
		\item there exists $(A'_\pi,B'_\pi) \in \T-Z$ of order at most $q+\Delta-1$, where $\Delta$ is the maximum degree of $H$, such that either
			\begin{enumerate}
				\item the order of $(A'_\pi,B'_\pi)$ is zero and $\pi(V(H)) \subseteq V(A'_\pi)-V(B'_\pi)$, or
				\item $\pi(V(H)) \cap V(A'_\pi)-V(B'_\pi) \supset \pi(V(H)) \cap V(A_\pi)-V(B_\pi)$, or
				\item $\pi(V(H)) \cap V(A'_\pi)-V(B'_\pi) = \pi(V(H)) \cap V(A_\pi)-V(B_\pi) \neq \emptyset$ and $\lvert V(A'_\pi) \cap V(B'_\pi) \rvert < \lvert V(A_\pi) \cap V(B_\pi) \rvert$.
			\end{enumerate}
	\end{enumerate}
\end{lemma}

\begin{pf}
Let $H$ be a connected graph, and let $k,q,\xi_0$ be positive integers.
Let $\Delta$ be the maximum degree of $H$.
Let $q' = \max\{q,\Delta,2\}$.
We define the following.
	\begin{itemize}
		\item $\xi = \xi_0 + (k(\lvert V(H) \rvert+\lvert E(H) \rvert+1)q')^{q'+1} + (\Delta+1)^{\Delta+1}$.
		\item $t=\xi+\lfloor \frac{3}{2}k(\lvert V(H) \rvert+\lvert E(H) \rvert+1)q' \rfloor$.
		\item $\theta=\xi+t+1$. 
	\end{itemize}

Let $S,G,\R,\T,Z_0,\F,\Q$ be as stated in the lemma.

Let $\F'= \{\pi \in \F:$ there exists no $(A,B) \in \T-Z_0$ of order zero with $\pi(V(H)) \cap V(A) \neq \emptyset\}$.
Note that when $S$ contains at least one vertex, for every $\pi \in \F$, $\pi(V(H)) \cap V(A_\pi) \neq \emptyset$, so for every $\pi \in \F'$, $V(A_\pi) \cap V(B_\pi) \neq \emptyset$.
For each $\pi \in \F'$, let $A^0_{\pi}=A_{\pi}-\{v \in V(A_{\pi}):$ there exists no path in $G-Z_0$ from $v$ to $V(A_{\pi}) \cap V(B_{\pi})\}$, and let $B^0_\pi$ be a subgraph of $G-Z_0$ such that $(A^0_\pi,B^0_\pi)$ is a separation of $G-Z_0$ with $V(A^0_\pi) \cap V(B^0_\pi) = V(A_\pi) \cap V(B_\pi)$.
Note that $(A^0_\pi,B^0_\pi) \in \T-Z_0$ for every $\pi \in \F'$ by (T1) and (T2).

\noindent{\bf Claim 1:} For every $\pi \in \F'$, $\pi(V(H)) \cap V(A_\pi)-V(B_\pi) = \pi(V(H)) \cap V(A^0_\pi)-V(B^0_\pi)$.

\noindent{\bf Proof of Claim 1:}
Suppose to the contrary that $\pi(V(H)) \cap V(A_\pi)-V(B_\pi) \neq \pi(V(H)) \cap V(A^0_\pi)-V(B^0_\pi)$.
Since $A_\pi \supseteq A^0_\pi$, there exists $v \in \pi(V(H)) \cap V(A_\pi)-V(A^0_\pi)$.
Since $v \not \in V(A^0_\pi)$, there exists a component $C$ of $A_\pi-V(A_\pi \cap B_\pi)$ containing $v$ such that there exists no edge of $G-Z_0$ between $V(C)$ and $V(A_\pi \cap B_\pi)$.
Note that $C$ is a component of $G-Z_0$.
Since $H$ is connected and $v \in V(C) \cap \pi(V(H))$, $\pi(V(H)) \subseteq V(C)$.
Hence $(C,G-V(C))$ is a separation of $G-Z_0$ of order zero with $\pi(V(H)) \subseteq V(C)$.
Since $C \subseteq A_\pi$, $(C,G-V(C)) \in \T-Z_0$ by (T1) and (T2), contradicting that $\pi \in \F'$.
$\Box$

We define the following.
	\begin{itemize}
		\item For each $v \in V(H)-V(S)$, let $\J_v$ be the maximal collection of subsets of $V(G)-Z_0$ such that for every $Y \in \J_v$, $Y$ consists of a vertex in $R_v-Z_0$ and $\Delta-1$ of its neighbors in $G-Z_0$.
		\item For each $e \in E(H)-E(S)$, let $\J_e$ be the maximal collection of subsets of $V(G)-Z_0$ such that for every $Y \in \J_e$, $Y$ consists of a vertex in $V(G)-Z_0$ and one of its neighbors in $V(G)-Z_0$.
		\item If $V(S) \neq \emptyset$, then let $\J_0 = \{V(A^0_\pi) \cap V(B^0_\pi): \pi \in \F'\}$; otherwise, let $\J_0=\emptyset$. 
	\end{itemize}
Note that for each $v \in V(H)-V(S)$, every member of $\J_v$ has size $\Delta \leq q'$, and each member of $\J_0$ has size at most $q \leq q'$.
And every member of $\J_v \cup \J_e \cup \J_0$ is a subset of $V(G)-Z_0$.
Furthermore, if $S$ contains at least one vertex, then every member of $\J_0$ is non-empty.

\noindent{\bf Claim 2:} There exist $U_0 \subseteq V(G)$ with $\lvert U_0 \rvert \leq (k(\lvert V(H) \rvert+\lvert E(H) \rvert+1)q')^{q'+1}$ and $r \in (V(H)-V(S)) \cup (E(H)-E(S)) \cup \{0\}$ such that for every $X \in \J_r$, either $X \cap U_0 \neq \emptyset$ or $X$ is not free with respect to $\T-(U_0 \cup Z_0)$.
In addition, if $V(S)=\emptyset$, then $r \in V(H)-V(S)$. 

\noindent{\bf Proof of Claim 2:}
Suppose that this claim does not hold.
If $V(S)=\emptyset$, then let $k'=0$; otherwise, let $k'=k$.
Since the order of $\T-Z_0$ is at least $\theta-\xi_0 \geq (k(\lvert V(H) \rvert+\lvert E(H) \rvert+1)q')^{q'+1}+q'$ and this claim does not hold, by Lemma \ref{stronger mw spider}, there exists $\J_0' \subseteq \J_0$ with $\lvert \J_0' \rvert =k'$, for every $e \in E(H)-E(S)$, there exists $\J_e' \subseteq \J_e$ with $\lvert \J_e' \rvert=k'$, and for every $v \in V(H)-V(S)$, there exists $\J_v' \subseteq \J_v$ with $\lvert \J_v' \rvert=k$ such that $\J_0' \cup \bigcup_{y \in (V(H)-V(S)) \cup (E(H)-E(S))} \J_y'$ consists of pairwise disjoint members, and $\bigcup_{X \in \J_0' \cup \bigcup_{y \in (V(H)-V(S)) \cup (E(H)-E(S))}\J_y'}X$ is free with respect to $\T-Z_0$.

Let $X^*=\bigcup_{X \in \J_0' \cup \bigcup_{y \in (V(H)-V(S)) \cup (E(H)-E(S))}\J_y'}X$.
Note that $\lvert X^* \rvert \leq k(\lvert V(H) \rvert+\lvert E(H) \rvert+1)q'$.
Let $\pi_1,\pi_2,...,\pi_{k'}$ be $\R$-compatible homeomorphic embeddings from $H$ to $G-Z_0$ such that $\J_0'=\{V(A^0_{\pi_j} \cap B^0_{\pi_j}): j \in [k']\}$.
By definition, for each $i \in [k']$ and $v \in V(A^0_{\pi_i})-V(B^0_{\pi_i})$, there exists a path in $G[A^0_{\pi_i}]-Z_0$ from $v$ to $V(A^0_{\pi_i} \cap B^0_{\pi_i})=V(A_{\pi_i} \cap B_{\pi_i})$.
By Lemma \ref{A side disjoint}, $V(A^0_{\pi_1}),...,V(A^0_{\pi_{k'}})$ are pairwise disjoint.
Moreover, $X^* \subseteq \bigcap_{j=1}^{k'}V(B^0_{\pi_j})$, for otherwise $(G[X^*] \cup \bigcup_{j=1}^{k'}A^0_{\pi_j}, \bigcap_{j=1}^{k'}B^0_{\pi_j}-E(G[X^*]))$ is a separation of $G-Z_0$ in $\T-Z_0$ of order less than $\lvert X^* \rvert$ such that the first part contains $X^*$, contradicting that $X^*$ is free with respect to $\T-Z_0$.

For each $v \in V(H)-V(S)$, we denote the members of $\J_v'$ as $X_v^1,...,X_v^k$.
For each $e \in E(H)-E(S)$, we denote the members of $\J_e'$ as $X_e^1,...,X_e^{k'}$.
For each $e \in E(H)-E(S)$ and $j \in [k]-[k']$, let $X_e^j=\emptyset$.

For each $j \in [k']$, let $W_j=\{u \in V(A^0_{\pi_j} \cap B^0_{\pi_j}):$ there exist $w \in V(S)$, $e \in E(H)-E(S)$ incident with $w$ and having all ends in $V(S)$, and a path $P$ in $A^0_{\pi_j}$ from $\pi_j(w)$ to $u$ internally disjoint from $V(B^0_{\pi_j})$ such that $P \subseteq \pi_j(e)$ and $V(\pi_j(e)) \not \subseteq V(A^0_{\pi_j})\}$; let $W_j'=\{u \in V(A^0_{\pi_j} \cap B^0_{\pi_j}):$ there exists $e \in E(H)$ incident with a vertex $w \in V(S)$ and a vertex $w' \not \in V(S)$, and a path $P$ in $A^0_{\pi_j}$ from $\pi_j(w)$ to $u$ internally disjoint from $V(B^0_{\pi_j})$ such that $P \subseteq \pi_j(e)$ and $V(\pi_j(e)) \not \subseteq V(A^0_{\pi_j})\}$.
For each $j \in [k]-[k']$, let $W_j=\emptyset$ and $W'_j=\emptyset$.

So for each $j \in [k]$, there exists a partition $\P_j$ of a subset of $W_j \cup \bigcup_{y \in E(H)-E(S)}X_y^j$ such that for each part $P$ of $\P_j$, there exist $e \in E(H)-E(S)$ with all ends in $V(S)$ and a part $P'$ in $\P_j-\{P\}$ such that each of $P$ and $P'$ consists of a vertex in $W_j \cap V(\pi_j(e))$ and a vertex in $X_e^j$.
Similarly, for each $j \in [k]$, there exists a partition $\P_j'$ of a subset of $W'_j \cup \bigcup_{v \in V(H)-V(S)}X_v^j$ such that for each part of $\P_j'$, there exists $e \in E(H)-E(S)$ incident with at most one vertex in $V(S)$, such that this part either consists of one vertex in $W'_j \cap V(\pi_j(e))$ and one vertex in $X_v^j$, where $v$ is the end of $e$ not in $V(S)$, or one vertex in $X_u^j$ and one vertex in $X_v^j$, where $u,v$ are the ends of $e$.
Let $\P_0=\bigcup_{j=1}^k(\P_j \cup \P_j')$.
Let $\P = \P_0 \cup \{X^*-\bigcup_{Y \in \P_0}Y\}$.
Note that $\P$ is a partition (with possibly one empty part) of $X^*$. 
By removing the possibly empty part of $\P$, we may assume that $\P$ has no empty parts.

Denote the members of $\P$ by $L_1,L_2,...,L_{\lvert \P \rvert}$.
We say that a part of $\P$ is {\it good} if it is not equal to $X^*-\bigcup_{Y \in \P_0}Y$.

Since $\bigcup_{j'=1}^{k'}V(A^0_{\pi_{j'}} \cap B^0_{\pi_{j'}}) \subseteq X^* \subseteq \bigcap_{j'=1}^{k'}V(B^0_{\pi_{j'}})$ and each good part of $\P$ consists of two vertices, where at least one vertex is in $(\bigcap_{j'=1}^{k'}B^0_{\pi_{j'}})-(\bigcup_{j'=1}^{k'}A^0_{\pi_{j'}})$, we know that if there exist $\lvert \P \rvert$ pairwise disjoint connected subgraphs $T_1,...,T_{\lvert \P \rvert}$ of $G-Z_0$ such that $V(T_j) \cap X^* = L_j$ for each $j \in [\lvert \P \rvert]$, then by choosing those $T_j$'s to be minimal, each $T_j$ corresponding to a good part is a path in $\bigcap_{j'=1}^{k'}B^0_{\pi_{j'}}$ internally disjoint from $V(\bigcup_{j'=1}^{k'}A^0_{\pi_{j'}})$, so there exist $k$ $\R$-compatible homeomorphic embeddings from $H$ into $G-Z_0$ with pairwise disjoint images, a contradiction. 
Hence, such $\lvert \P \rvert$ pairwise disjoint connected subgraphs of $G-Z_0$ do not exist.
Note that $t \geq \lfloor \frac{3}{2}k(\lvert V(H) \rvert+\lvert E(H) \rvert+1)q' \rfloor + \xi_0$ and $\T$ controls a $K_t$-minor, so $\T-Z_0$ controls a $K_{t'}$-minor $\alpha$, for some $t' \geq \lfloor \frac{3}{2}k(\lvert V(H) \rvert+\lvert E(H) \rvert+1)q' \rfloor$.
Since $\lvert X^* \rvert \leq k(\lvert V(H) \rvert+\lvert E(H) \rvert+1)q'$ and $t' \geq \lfloor \frac{3}{2}k(\lvert V(H) \rvert+\lvert E(H) \rvert+1)q' \rfloor$, by Theorem \ref{vertex linkage}, there exists a separation $(A,B)$ of $G-Z_0$ of order less than $\lvert X^* \rvert$ such that $X^* \subseteq V(A)$ and $A \cap \alpha(w)=\emptyset$ for some  $w \in V(K_{t'})$.
Since $X^*$ is free with respect to $\T-Z_0$, $(A,B) \not \in \T-Z_0$.
Since the order of $\T-Z_0$ is greater than $\lvert X^* \rvert$, $(B,A) \in \T-Z_0$ by (T1).
But $\T-Z_0$ controls $\alpha$ and $(B,A) \in \T-Z_0$ is a separation of order less than $\lvert V(K_{t'}) \rvert$, so $V(\alpha(w)) \cap (V(A)-V(B)) \neq \emptyset$.
This implies that $\alpha(w) \cap A \neq \emptyset$, a contradiction.
$\Box$

Let $U_0$ be the set and $r$ be the element mentioned in Claim 2.

\noindent{\bf Claim 3:} $r \in (V(H)-V(S)) \cup \{0\}$.

\noindent{\bf Proof of Claim 3:}
Suppose to the contrary that $r \in E(H)-E(S)$.
By Claim 2, for each non-loop edge $e$ of $G-Z_0$, the set of its ends either intersects $U_0$ or is not free with respect to $\T-(Z_0 \cup U_0)$.
So for each non-loop edge $e$ of $G-(Z_0 \cup U_0)$, there exists a separation $(A_e,B_e) \in \T-(Z_0 \cup U_0)$ of order at most 1 such that the ends of $e$ are contained in $A_e$, and we assume that $A_e$ is maximal subject to this property.

If $G-(Z_0 \cup U_0)$ has no non-loop edge, then every component of $G-(Z_0 \cup U_0)$ has at most 1 vertex, so by letting $Z=Z_0 \cup U_0$, we know that for every $\pi \in \F$, Statement 1 or 2(a) hold, since $\lvert Z \rvert \leq \xi_0+(k(\lvert V(H) \rvert+\lvert E(H) \rvert+1)q')^{q'+1} \leq \xi$ and $H$ is connected.
So we may assume that there exists a non-loop edge $e$ of $G-(Z_0 \cup U_0)$.
By (T3), there exists $v \in V(B_e)-V(A_e)$.
We choose $v$ such that $v$ is adjacent to $V(A_e) \cap V(B_e)$ if possible.

If $v$ is not incident with a non-loop edge, then $(G[V(A_e) \cup \{v\}], B_e-\{v\})$ is a separation of $G-(Z_0 \cup U_0)$ of order at most 1 and is in $\T-(Z_0 \cup U_0)$ by (T1) and (T2), contradicting the maximality of $A_e$.
So we may assume that $v$ is incident with a non-loop edge $e'$.
We choose $e'$ such that $e'$ is incident with $V(A_e) \cap V(B_e)$ if possible.
If $e'$ is not incident with $V(A_e) \cap V(B_e)$, then by the our choices of $v$ and $e'$, either $V(A_e) \cap V(B_e)=\emptyset$, or the unique vertex in $V(A_e) \cap V(B_e)$ is not adjacent to any vertex in $V(B_e)-V(A_e)$, so $(A_e \cup A_{e'}, (B_e-V(A_e)) \cap B_{e'})$ is a separation of $G-(Z_0 \cup U_0)$ of order at most 1 and is in $\T-(Z_0 \cup U_0)$ by (T1) and (T2), contradicting the maximality of $A_e$.
So we may assume that $e'$ is incident with $V(A_e) \cap V(B_e)$.
Then $(A_e \cup A_{e'}, B_e \cap B_{e'})$ is a separation of $G-(Z_0 \cup U_0)$ of order at most 1 and is in $\T-(Z_0 \cup U_0)$ by (T1) and (T2), contradicting the maximality of $A_e$.
$\Box$

\noindent{\bf Claim 4:} We may assume that $r=0$ and $V(S) \neq \emptyset$.

\noindent{\bf Proof of Claim 4:}
Suppose to the contrary that $r \in (V(H)-V(S)) \cup (E(H)-E(S))$ or $V(S)=\emptyset$.
By Claim 3, $r \in V(H)-V(S)$ or $V(S)=\emptyset$.
Since $V(S)=\emptyset$ implies that $r \in V(H)-V(S)$ by Claim 2, $r \in V(H)-V(S)$ in either case.

Let $R_r'= \{v \in R_r-(U_0 \cup Z_0): v$ has degree at least $\Delta-1$ in $G-(U_0 \cup Z_0)\}$.
By Claim 2, there exists no subset $X$ of $V(G)-(U_0 \cup Z_0)$ consisting of a vertex in $R_r'$ and $\Delta-1$ of its neighbors in $G-(U_0 \cup Z_0)$ such that $X$ is free with respect to $\T-(U_0 \cup Z_0)$.
Since the order of $\T-(U_0 \cup Z_0)$ is at least $\theta-\xi_0-(k(\lvert V(H) \rvert+\lvert E(H) \rvert+1)q')^{q'+1} \geq (\Delta+1)^{\Delta+1}+\Delta$, by Lemma \ref{spider tangle} (taking $(h,d,G,S,\T)=(1,\Delta,G-(U_0 \cup Z_0),R'_r,\T-(U_0 \cup Z_0))$), there exists a set $C_r \subseteq V(G)-(U_0 \cup Z_0)$ with $\lvert C_r \rvert \leq (\Delta+1)^{\Delta+1}$ such that for every vertex $u \in R'_r-C_r$, there exists $(A_u,B_u) \in \T-(C_r \cup U_0 \cup Z_0)$ of order at most $\Delta-1$ such that $u \in V(A_u)-V(B_u)$.
In addition, for every vertex $u \in R_r-(U_0 \cup Z_0)$ with degree less than $\Delta-1$ in $G-(U_0 \cup Z_0)$, there exists $(A_u,B_u) \in \T-(U_0 \cup Z_0)$ of order less than $\Delta-1$ such that $u \in V(A_u)-V(B_u)$.
Hence, for every vertex $u \in R_r-(C_r \cup U_0 \cup Z_0)$, there exists $(A_u,B_u) \in \T-(C_r \cup U_0 \cup Z_0)$ of order at most $\Delta-1$ such that $u \in V(A_u)-V(B_u)$.

Let $Z = C_r \cup U_0 \cup Z_0$.
Then for every $\R$-compatible homeomorphic embedding $\pi$ from $H$ into $G-Z_0$ with $\pi(V(H) \cup E(H)) \cap Z \neq \emptyset$, $\pi(r) \in R_r-(C_r \cup U_0 \cup Z_0)$; for each such $\pi$, let $(A'_\pi,B'_\pi)=((A_\pi-Z) \cup A_{\pi(r)}, (B_\pi-Z) \cap B_{\pi(r)})$, so $(A'_\pi,B'_\pi) \in \T-Z$ has order at most $\lvert V(A_\pi) \cap V(B_\pi) \rvert + \lvert V(A_{\pi(r)}) \cap V(B_{\pi(r)}) \rvert \leq q+\Delta-1$, and $\pi(r) \in \pi(V(H)) \cap (V(A'_\pi)- V(B'_\pi)) - (V(A_\pi)-V(B_\pi))$, since $r \in V(H)-V(S)$.
Note that $V(A'_\pi)-V(B'_\pi) \supseteq V(A_\pi)-(V(B_\pi) \cup Z)$ and $\lvert Z \rvert \leq \lvert C_r \rvert + \lvert U_0 \rvert + \lvert Z_0 \rvert \leq (\Delta+1)^{\Delta+1} + (k(\lvert V(H) \rvert+\lvert E(H) \rvert+1)q')^{q'+1} + \xi_0 \leq \xi$.
So Statement 2(b) of this lemma holds.
$\Box$

Hence we may assume that $r=0$ and $V(S) \neq \emptyset$ by Claim 4.
Let $Z = U_0 \cup Z_0$.
Hence $\lvert Z \rvert \leq (k(\lvert V(H) \rvert+\lvert E(H) \rvert+1)q')^{q'+1} + \xi_0 \leq \xi$.

Suppose to the contrary that this lemma does not hold.
So there exists $\pi \in \F$ such that Statements 1 and 2 do not hold for $\pi$.
In particular, $\pi(V(H) \cup E(H)) \cap Z = \emptyset$.
So $\pi$ is a homeomorphic embedding from $H$ into $G-Z$.

\noindent{\bf Claim 5:} $\pi \in \F'$.

\noindent{\bf Proof of Claim 5:}
Suppose to the contrary that $\pi \in \F-\F'$.
Then there exists $(A,B) \in \T-Z_0$ of order zero such that $\pi(V(H)) \cap V(A) \neq \emptyset$.
So $(A-Z,B-Z) \in \T-Z$ has order zero such that $\pi(V(H)) \cap V(A-Z) \neq \emptyset$.
Since $H$ is connected, $\pi(V(H)) \subseteq V(A-Z) = V(A-Z) -V(B-Z)$.
Hence Statement 2(a) of this lemma holds, a contradiction.
$\Box$

Since $V(S) \neq \emptyset$, $\pi(V(H)) \cap V(A_\pi)-V(B_\pi) \neq \emptyset$.
If $Z \cap V(A_\pi \cap B_\pi) \neq \emptyset$, then defining $(A'_\pi,B'_\pi)=(A_\pi-Z, B_\pi-Z)$ implies that Statement 2(c) holds, a contradiction.
So $Z \cap V(A_\pi \cap B_\pi) = \emptyset$.

Let $(A'_\pi,B'_\pi)$ be the separation in $\T-Z$ such that $V(A_\pi) \cap V(B_\pi) \subseteq V(A'_\pi)$, and subject to this, the order is minimum, and subject to these, $A'_\pi-V(B'_\pi)$ is maximal.
Since $Z \cap V(A_\pi \cap B_\pi) = \emptyset$, $(A_\pi-Z,B_\pi-Z)$ is a candidate for $(A'_\pi,B'_\pi)$, so $(A'_\pi,B'_\pi)$ exists. 

By Claims 4 and 5, $V(S) \neq \emptyset$ and $\pi \in \F'$, so $V(A_\pi) \cap V(B_\pi) = V(A^0_\pi) \cap V(B^0_\pi) \in \J_0$.
Since $Z \cap V(A_\pi \cap B_\pi) = \emptyset$, by Claims 2 and 4, $V(A_\pi) \cap V(B_\pi)$ it is not free with respect to $\T-Z$, so the order of $(A'_\pi,B'_\pi)$ is less than $\lvert V(A_\pi) \cap V(B_\pi) \rvert$ by the definition of $(A'_\pi,B'_\pi)$.
Since Statements 2(b) and 2(c) do not hold, $A_\pi-(V(B_\pi) \cup Z) \not \subseteq A'_\pi-V(B'_\pi)$. 

Let $(A',B')=(A'_\pi \cup (A_\pi-Z), B'_\pi \cap (B_\pi-Z))$ and $(A'',B'')=(A'_\pi \cap (A_\pi-Z), B'_\pi \cup (B_\pi-Z))$.
Since $\theta$ is sufficiently large, both $(A',B')$ and $(A'',B'')$ belong to $\T-Z$ by (T1) and (T2).
Since $V(A_\pi \cap B_\pi) \cap Z=\emptyset$ and $V(A_\pi \cap B_\pi) \subseteq V(A'_\pi)$, we know $V(A_\pi \cap B_\pi) \subseteq V(A'' \cap B'')$.
So $\lvert V(A'' \cap B'') \rvert \geq \lvert V(A_\pi \cap B_\pi) \rvert = \lvert V((A_\pi-Z) \cap (B_\pi-Z)) \rvert$.
By the submodularity, the order of $(A',B')$ is at most the order of $(A'_\pi,B'_\pi)$.
Since $A' \supseteq A'_\pi$, by the definition of $(A'_\pi,B'_\pi)$, the order of $(A',B')$ equals the order of $(A'_\pi,B'_\pi)$.
Since $A'-V(B') \supseteq A'_\pi-V(B'_\pi)$, by the definition of $(A'_\pi,B'_\pi)$, $A'-V(B')=A'_\pi-V(B'_\pi)$. 
Hence $A_\pi-(V(B_\pi) \cup Z) \subseteq (A_\pi-Z)-(V(B_\pi)-Z) \subseteq (A_\pi-Z)-V(B'_\pi \cap (B_\pi-Z)) \subseteq A'-V(B')=A'_\pi-V(B'_\pi)$, a contradiction.
This proves the lemma.
\end{pf}

\begin{lemma} \label{clique_minor_one_step} 
For any connected graph $H$ and positive integers $k,q_0,\xi_0$, there exist integers $\theta,t,\xi,q$ such that if $G$ is a graph, $\R=(R_v: v \in V(H))$ is a collection of subsets of $V(G)$, $\T$ is a tangle in $G$ of order at least $\theta$ controlling a $K_t$-minor, a subset $Z_0$ of $V(G)$ with $\lvert Z_0 \rvert \leq \xi_0$, a collection $\F$ of $\R$-compatible homeomorphic embeddings from $H$ into $G-Z_0$, and a collection $\Q = \{(A_\pi,B_\pi): \pi \in \F\} \subseteq \T-Z_0$ of separations of order at most $q_0$ such that $G-Z_0$ does not contain $k$ disjoint $\R$-compatible subdivisions of $H$, then there exists $Z$ with $Z_0 \subseteq Z \subseteq V(G)$ of size at most $\xi$ such that for every $\pi \in \F$, either 
	\begin{enumerate}
		\item $\pi(V(H) \cup E(H)) \cap Z \neq \emptyset$, or
		\item there exists $(A'_\pi,B'_\pi) \in \T-Z$ of order at most $q$ such that either
			\begin{enumerate}
				\item the order of $(A'_\pi,B'_\pi)$ is zero and $\pi(V(H)) \subseteq V(A'_\pi)-V(B'_\pi)$, or
				\item $\pi(V(H)) \cap V(A'_\pi)-V(B'_\pi) \supset \pi(V(H)) \cap V(A_\pi)-V(B_\pi)$, or
				\item $\pi(V(H)) \cap V(A'_\pi)-V(B'_\pi) = \pi(V(H)) \cap V(A_\pi)-V(B_\pi) \neq \emptyset$ and $\lvert V(A'_\pi) \cap V(B'_\pi) \rvert < \lvert V(A_\pi) \cap V(B_\pi) \rvert$.
			\end{enumerate}
	\end{enumerate}
\end{lemma}

\begin{pf}
Let $H$ be a connected graph, and let $k,q_0,\xi_0$ be positive integers.
Let $\Delta$ be the maximum degree of $H$.
Let $\theta_1,t_1,\xi_1$ be the integers $\theta,t,\xi$ mentioned in Lemma \ref{clique_minor_uniform} by taking $(H,k,q,\xi_0)$ to be $(H,k,q_0,\xi_0)$, respectively.
Define $q=q_0+\Delta-1$, $t=t_1$, $\xi=2^{\lvert V(H) \rvert + \lvert E(H) \rvert} \cdot \xi_1$, and $\theta = \theta_1 + \xi + t$.

Let $G,\R,\T,Z_0,\F,\Q$ be as stated in the lemma.
For each (possibly empty) subgraph $S$ of $H$, let $\F_S = \{\pi \in \F:$ the shade of $\pi$ with respect to $(A_\pi,B_\pi)$ equals $S\}$, and let $\Q_S = \{(A_\pi,B_\pi): \pi \in \F_S\}$.
By Lemma \ref{clique_minor_uniform}, for each subgraph $S$ of $H$, there exists $Z_S$ with $Z_0 \subseteq Z_S \subseteq V(G)$ with $\lvert Z_S \rvert \leq \xi_1$ such that for every $\pi \in \F_S$, either 
	\begin{itemize}
		\item $\pi(V(H) \cup E(H)) \cap Z_S \neq \emptyset$, or
		\item there exists $(A'_\pi,B'_\pi) \in \T-Z_S$ of order at most $q_0+\Delta-1=q$ such that either
			\begin{itemize}
				\item the order of $(A'_\pi,B'_\pi)$ is zero and $\pi(V(H)) \subseteq V(A'_\pi)-V(B'_\pi)$, or
				\item $\pi(V(H)) \cap V(A'_\pi)-V(B'_\pi) \supset \pi(V(H)) \cap V(A_\pi)-V(B_\pi)$, or
				\item $\pi(V(H)) \cap V(A'_\pi)-V(B'_\pi) = \pi(V(H)) \cap V(A_\pi)-V(B_\pi) \neq \emptyset$ and $\lvert V(A'_\pi) \cap V(B'_\pi) \rvert < \lvert V(A_\pi) \cap V(B_\pi) \rvert$.
			\end{itemize}
	\end{itemize}

Let $Z = \bigcup_SZ_S$, where the union is over all subgraphs $S$ of $H$.
Since there are at most $2^{\lvert V(H) \rvert + \lvert E(H) \rvert}$ subgraphs of $H$, $\lvert Z \rvert \leq 2^{\lvert V(H) \rvert + \lvert E(H) \rvert} \cdot \xi_1 = \xi$.
Then $Z$ satisfies the conclusion of this lemma.
\end{pf}

\begin{lemma} \label{clique_minor_one_big_step} 
For any connected graph $H$ and positive integers $k,q_0,\xi_0$, there exist integers $\theta,t,\xi,q$ such that if $G$ is a graph, $\R=(R_v: v \in V(H))$ is a collection of subsets of $V(G)$, $\T$ is a tangle in $G$ of order at least $\theta$ controlling a $K_t$-minor, a subset $Z_0$ of $V(G)$ with $\lvert Z_0 \rvert \leq \xi_0$, a collection $\F$ of $\R$-compatible homeomorphic embeddings from $H$ into $G-Z_0$, and a collection $\Q = \{(A_\pi,B_\pi): \pi \in \F\} \subseteq \T-Z_0$ of separations of order at most $q_0$ such that $G-Z_0$ does not contain $k$ disjoint $\R$-compatible subdivisions of $H$, then there exists $Z$ with $Z_0 \subseteq Z \subseteq V(G)$ of size at most $\xi$ such that for every $\pi \in \F$, either 
	\begin{enumerate}
		\item $\pi(V(H) \cup E(H)) \cap Z \neq \emptyset$, or
		\item there exists $(A'_\pi,B'_\pi) \in \T-Z$ of order at most $q$ such that either
			\begin{enumerate}
				\item the order of $(A'_\pi,B'_\pi)$ is zero and $\pi(V(H)) \subseteq V(A'_\pi)-V(B'_\pi)$, or
				\item $\pi(V(H)) \cap V(A'_\pi)-V(B'_\pi) \supset \pi(V(H)) \cap V(A_\pi)-V(B_\pi)$. 
			\end{enumerate}
	\end{enumerate}
\end{lemma}

\begin{pf}
Let $H$ be a connected graph, and let $k,q_0,\xi_0$ be positive integers.
For every positive integer $i$, let $\theta_i,t_i,\xi_i,q_i$ be the integers $\theta,t,\xi,q$ mentioned in Lemma \ref{clique_minor_one_step} by taking $(H,k,q_0,\xi_0)$ to be $(H,k,q_0,\xi_{i-1})$, respectively.
Define $q=\sum_{i=0}^{q_0}q_i$, $t=\sum_{i=1}^{q_0}t_i$, $\xi = \xi_{q_0}$ and $\theta = \xi+t+\sum_{i=1}^{q_0}\theta_i$.

Let $G,\R,\T,Z_0,\F,\Q$ be as stated in the lemma.
Let $\F_0=\F$ and $\Q_0=\Q$.
For every $\pi \in \F_0=\F$, let $(A_\pi^0,B_\pi^0)=(A_\pi,B_\pi)$.
For every positive integer $i$, applying Lemma \ref{clique_minor_one_step} by taking $(H,k,q_0,\xi_0,G,\R,\T,Z_0,\F,\Q) = (H,k,q_0,\xi_{i-1},G,\R,\T,Z_{i-1},\F_{i-1},\Q_{i-1})$, we obtain the following.
	\begin{itemize}
		\item A set $Z_i$ with $Z_{i-1} \subseteq Z_i \subseteq V(G)$ with size at most $\xi_i$ such that for every $\pi \in \F_{i-1}$, either 
			\begin{itemize}
				\item $\pi(V(H) \cup E(H)) \cap Z_i \neq \emptyset$, or
				\item there exists $(A^i_\pi,B^i_\pi) \in \T-Z_i$ of order at most $q_i$ such that either
					\begin{itemize}
						\item the order of $(A^i_\pi,B^i_\pi)$ is zero and $\pi(V(H)) \subseteq V(A^i_\pi)-V(B^i_\pi)$, or
						\item $\pi(V(H)) \cap V(A^i_\pi)-V(B^i_\pi) \supset \pi(V(H)) \cap V(A^{i-1}_\pi)-V(B^{i-1}_\pi)$, or
						\item $\pi(V(H)) \cap V(A^i_\pi)-V(B^i_\pi) = \pi(V(H)) \cap V(A^{i-1}_\pi)-V(B^{i-1}_\pi) \neq \emptyset$ and $\lvert V(A^i_\pi) \cap V(B^i_\pi) \rvert < \lvert V(A^{i-1}_\pi) \cap V(B^{i-1}_\pi) \rvert$.
					\end{itemize}
			\end{itemize}
		\item A collection $\F_i = \{\pi \in \F_{i-1}: \pi(V(H) \cup E(H)) \cap Z_i = \emptyset, V(A^i_\pi) \cap V(B^i_\pi) \neq \emptyset, \pi(V(H)) \cap V(A^i_\pi)-V(B^i_\pi) = \pi(V(H)) \cap V(A^{i-1}_\pi)-V(B^{i-1}_\pi) \neq \emptyset\}$.
			Note that every member of $\F_i$ is an $\R$-compatible homeomorphic embedding from $H$ into $G-Z_i$.
		\item A collection $\Q_i = \{(A^i_\pi,B^i_\pi): \pi \in \F_i\} \subseteq \T-Z_i$.
			Note that for every $\pi \in \F_i$, 
			\begin{itemize}
				\item since $\pi(V(H) \cup E(H)) \cap Z_i = \emptyset$, $(A^i_\pi,B^i_\pi)$ is defined and $\pi(V(H)) \cap V(A^i_\pi)-V(B^i_\pi) \supseteq \pi(V(H)) \cap V(A^{i-1}_\pi)-V(B^{i-1}_\pi)$,
				\item since $V(A^i_\pi) \cap V(B^i_\pi) \neq \emptyset$ and $\pi(V(H)) \cap V(A^i_\pi)-V(B^i_\pi) = \pi(V(H)) \cap V(A^{i-1}_\pi)-V(B^{i-1}_\pi) \neq \emptyset$, the order of $(A^i_\pi,B^i_\pi)$ is less than the order of $(A^{i-1}_\pi,B^{i-1}_\pi)$.
			\end{itemize}
			So every member of $\Q_i$ has order at most $q_0-i$, and for every $\pi \in \F_i$, $\pi(V(H)) \cap V(A^i_\pi)-V(B^i_\pi) \supseteq \pi(V(H)) \cap V(A^0_\pi)-V(B^0_\pi) = \pi(V(H)) \cap V(A_\pi)-V(B_\pi)$.
	\end{itemize}
Note that for every $\pi \in \F_{q_0}$, $(A^{q_0}_\pi,B^{q_0}_\pi)$ has order zero and satisfies $V(A^{q_0}_\pi) \cap V(B^{q_0}_\pi) \neq \emptyset$. 
Hence $\F_{q_0}=\emptyset$.

Let $Z = \bigcup_{i=0}^{q_0}Z_0$.
Note that $Z=Z_{q_0}$.
So $\lvert Z \rvert \leq \xi$.
For every $\pi \in \F$, let $i_\pi$ be the least positive integer such that $\pi \not \in \F_{i_\pi}$, so $\pi \in \F_{i_\pi-1}$ and $i_\pi \leq q_0$, and hence either 
	\begin{itemize}
		\item[(i)] $\pi(V(H) \cup E(H)) \cap Z_{i_\pi} \neq \emptyset$, or
		\item[(ii)] the order of $(A^{i_\pi}_\pi,B^{i_\pi}_\pi)$ is zero and $\pi(V(H)) \subseteq V(A^{i_\pi}_\pi)-V(B^{i_\pi}_\pi)$, or 
		\item[(iii)] $\pi(V(H)) \cap V(A^{i_\pi}_\pi)-V(B^{i_\pi}_\pi) \supset \pi(V(H)) \cap V(A^{i_\pi-1}_\pi)-V(B^{i_\pi-1}_\pi)$, or
		\item[(iv)] $\pi(V(H)) \cap V(A^{i_\pi}_\pi)-V(B^{i_\pi}_\pi) = \pi(V(H)) \cap V(A^{i_\pi-1}_\pi)-V(B^{i_\pi-1}_\pi) \neq \emptyset$, $\lvert V(A^{i_\pi}_\pi) \cap V(B^{i_\pi}_\pi) \rvert < \lvert V(A^{i_\pi-1}_\pi) \cap V(B^{i_\pi-1}_\pi) \rvert$, and $V(A^{i_\pi}) \cap V(B^{i_\pi})=\emptyset$.
	\end{itemize}
Note that (iv) implies (ii) since $H$ is connected.
Therefore, for every $\pi \in \F$, either 
	\begin{itemize}
		\item $\pi(V(H) \cup E(H)) \cap Z \neq \emptyset$, or 
		\item there exists $(A'_\pi,B'_\pi) = (A^{i_\pi}_\pi-Z,B^{i_\pi}_\pi-Z) \in \T-Z$ such that $(A'_\pi,B'_\pi)$ has order at most $q_{i_\pi} \leq q$, and either the order of $(A'_\pi,B'_\pi)$ is zero and $\pi(V(H)) \subseteq V(A'_\pi)-V(B'_\pi)$, or $\pi(V(H)) \cap V(A'_\pi)-V(B'_\pi) \supset \pi(V(H)) \cap V(A^{i_\pi-1}_\pi)-V(B^{i_\pi-1}_\pi) \supseteq \pi(V(H)) \cap V(A_\pi)-V(B_\pi)$.
	\end{itemize}
This proves the lemma.
\end{pf}

\bigskip

Now we are ready to state and prove the main result of this section.

\begin{lemma} \label{with clique minor}
For any connected graph $H$ and positive integer $k$, there exist integers $\theta,t,\xi$ such that if $G$ is a graph, $\R=(R_v: v \in V(H))$ is a collection of subsets of $V(G)$, and $\T$ is a tangle in $G$ of order at least $\theta$ controlling a $K_t$-minor, then one of the following holds.
	\begin{enumerate}
		\item $G$ contains $k$ disjoint $\R$-compatible subdivisions of $H$.
		\item There exists $Z \subseteq V(G)$ with $\lvert Z \rvert \leq \xi$ such that $\pi(V(H) \cup E(H)) \cap Z \neq \emptyset$ for every $\R$-compatible homeomorphic embedding $\pi$ from $H$ into $G$.
		\item There exists a separation $(A,B) \in \T$ such that $A-V(B)$ contains an $\R$-compatible subdivision of $H$.
	\end{enumerate}
\end{lemma}

\begin{pf}
Let $H$ be a connected graph, and let $k$ be a positive integer.
Let $\xi_0=1$ and $q_0=1$.
For every positive integer $i$, let $\theta_i$, $t_i$, $\xi_i$, $q_i$ be the integers $\theta,t,\xi,q$ mentioned in Lemma \ref{clique_minor_one_big_step} by taking $(H,k,q_0,\xi_0)=(H,k,q_{i-1},\xi_{i-1})$.
Define $\theta = \sum_{i=1}^{\lvert V(H) \rvert+1}\theta_i$, $t = \sum_{i=1}^{\lvert V(H) \rvert+1}t_i$ and $\xi = \xi_{\lvert V(H) \rvert+1}$.

Let $G,\R,\T$ be as stated in the lemma.
We may assume that $G$ does not contain $k$ disjoint $\R$-compatible subdivisions of $H$, for otherwise we are done.

Let $Z_0=\emptyset$.
Note that $\lvert Z_0 \rvert \leq \xi_0$.
Let $\F_0$ be the collection of all $\R$-compatible homeomorphic embeddings from $H$ into $G-Z_0=G$.
For every $\pi \in \F_0$, let $(A^0_\pi,B^0_\pi) = (\emptyset, G)$.
Let $\Q_0 = \{(A^0_\pi,B^0_\pi): \pi \in \F_0\}$.
Note that $\Q_0$ is a subset of $\T$, and every member of $\Q_0$ has order at most $0 \leq q_0$.

For every positive integer $i$, we apply Lemma \ref{clique_minor_one_big_step} by taking $(H,k,q_0,\xi_0,G,\R,\T,Z_0,\F,\Q) \allowbreak =(H,k,q_{i-1},\xi_{i-1},G,\R,\T,Z_{i-1},\F_{i-1},\Q_{i-1})$, we obtain the following.
	\begin{itemize}
		\item A set $Z_i$ with $Z_{i-1} \subseteq Z_i \subseteq V(G)$ with $\lvert Z_i \rvert \leq \xi_i$ such that for every $\pi \in \F_{i-1}$, either
			\begin{itemize}
				\item $\pi(V(H) \cup E(H)) \cap Z_i \neq \emptyset$, or
				\item there exists $(A^i_\pi,B^i_\pi) \in \T-Z_i$ of order at most $q_i$ such that either
					\begin{itemize}
						\item the order of $(A^i_\pi,B^i_\pi)$ is zero and $\pi(V(H)) \subseteq V(A^i_\pi)-V(B^i_\pi)$, or
						\item $\pi(V(H)) \cap V(A^i_\pi)-V(B^i_\pi) \supset \pi(V(H)) \cap V(A^{i-1}_\pi)-V(B^{i-1}_\pi)$. 
					\end{itemize}
			\end{itemize}
		\item A collection $\F_i = \{\pi \in \F_{i-1}: \pi(V(H) \cup E(H)) \cap Z_i = \emptyset$ and there exists no separation $(A,B) \in \T-Z_i$ of order zero with $\pi(V(H)) \subseteq V(A)-V(B)\}$.
			Note that each member of $\F_i$ is an $\R$-compatible homeomorphic embedding from $H$ into $G-Z_i$.
		\item A collection $\Q_i = \{(A^i_\pi,B^i_\pi): \pi \in \F_i\}$.
			Note that every member $(A^i_\pi,B^i_\pi)$ has order at most $q_i$ such that $\lvert \pi(V(H)) \cap V(A^i_\pi)-V(B^i_\pi) \rvert \geq i$.
	\end{itemize}
Note that $\Q_{\lvert V(H) \rvert+1}=\emptyset$, so $\F_{\lvert V(H) \rvert+1} = \emptyset$.

Define $Z=Z_{\lvert V(H) \rvert+1}$.
So $\lvert Z \rvert \leq \xi$.

Suppose that this lemma does not hold.
Then there exists an $\R$-compatible homeomorphic embedding $\pi$ from $H$ to $G-Z$ such that there exists no $(A,B) \in \T$ such that $A-V(B)$ contains $\pi(V(H) \cup E(H))$.
Note that $\pi \in \F_0$.
So there exists $i_\pi$ such that $i_\pi$ is the least positive integer $i$ such that $\pi \not \in \F_i$ and $\pi \in \F_{i-1}$.
Since $Z \supseteq Z_{i_\pi}$ and $\pi \in \F_{i_\pi-1}-\F_{i_\pi}$, there exists a separation $(A,B) \in \T-Z_{i_\pi}$ of order zero with $\pi(V(H)) \subseteq V(A)-V(B)$.
Since $H$ is connected and $(A,B)$ has order zero, $\pi(V(H) \cup E(H)) \subseteq A=A-V(B)$.
Hence there exists $(A',B') \in \T$ with $V(A')=V(A) \cup Z$ and $V(B')=V(B) \cup Z$ such that $A'-V(B')=A-V(B)$ contains $\pi(V(H) \cup E(H))$, a contradiction.
This proves the lemma.
\end{pf}

\section{Surfaces and vortices} \label{sec: surfaces}

In this section, we introduce and extend some machinery about drawings in surfaces, which will be extensively used in later sections.

\subsection{Preliminaries}
In this subsection, we define some notions and state some known results about surfaces and drawings.

\subsubsection{Drawings}

A {\it surface} is a nonnull compact $2$-manifold without boundary.
A {\it line} in a surface $\Sigma$ is a subset of $\Sigma$ homeomorphic to $[0,1]$.
An {\it O-arc} in a surface $\Sigma$ is a subset of $\Sigma$ homeomorphic to a circle.
For every subset $\Delta$ of a surface $\Sigma$, we denote the closure of $\Delta$ by $\overline{\Delta}$, and the boundary of $\Delta$ by $\partial\Delta$.

A {\it drawing} $\Gamma$ in a surface $\Sigma$ is a pair $(U,V)$, where $V \subseteq U \subseteq \Sigma$, $U$ is closed, $V$ is finite, $U-V$ has only finitely many arc-wise connected components, called {\it edges}, and for every edge $e$, either the closure $\bar{e}$ of $e$ is a line whose set of ends is $\bar{e} \cap V$, or $\bar{e}$ is an O-arc and $\lvert \bar{e} \cap V \rvert =1$.
Every component of $\Sigma-U$ is called a {\it region}.
Every member of $V$ is called a {\it vertex}.
For a drawing $\Gamma=(U,V)$, we write $U=U(\Gamma), V= V(\Gamma)$, and define the set of edges to be $E(\Gamma)$.
If $v$ is a vertex of a drawing $\Gamma$ and $e$ is an edge or a region of $\Gamma$, we say that $e$ is {\it incident with} $v$ if $v \in \overline{e}$.
Note that the incidence relation between $V(\Gamma)$ and $E(\Gamma)$ defines a graph, and we say that $\Gamma$ is a {\it drawing of $G$} in $\Sigma$ if $G$ is defined by this incidence relation.
In this case, we say that $G$ is {\it embeddable} in $\Sigma$, or $G$ can be {\it drawn} in $\Sigma$.

Let $\Sigma$ be a surface and $\Gamma$ a drawing in $\Sigma$.
The sets $\{v\}$, for all $v \in V(\Gamma)$, the edges and regions of $\Gamma$ are called the {\it atoms} of $\Gamma$.
We say that a drawing $\Gamma'$ is a {\it subdrawing} of $\Gamma$ if $V(\Gamma') \subseteq V(\Gamma)$ and $E(\Gamma') \subseteq E(\Gamma)$.
We write $\Gamma' \subseteq \Gamma$ if $\Gamma'$ is a subdrawing of $\Gamma$.
If $\Delta \subseteq \Sigma$ is a closed set such that either $\bar{e} \subseteq \Delta$ or $e \cap \Delta=\emptyset$ for each $e \in E(\Gamma)$, then we define $\Gamma \cap \Delta$ to be the drawing $(U(\Gamma) \cap \Delta, V(\Gamma) \cap \Delta)$.
A subset $Z$ of $\Sigma$ is {\it $\Gamma$-normal} if $Z \cap U(\Gamma) \subseteq V(\Gamma)$.
If $\Sigma$ is connected and not a sphere, we say that $\Gamma$ is {\it $\theta$-representative} if $\lvert F \cap V(\Gamma) \rvert \geq \theta$ for every non-null-homotopic $\Gamma$-normal O-arc $F$ in $\Sigma$.
We say $\Gamma$ is {\it $2$-cell} if $\Sigma$ is connected and every region of $\Gamma$ is an open disk.

\subsubsection{Metrics in drawings}

Let $\Sigma$ be a surface, and let $\Gamma$ be a drawing of a graph $G$ in $\Sigma$.
A {\it tangle} in $\Gamma$ and a {\it separation} of $\Gamma$ are a tangle in $G$ and a separation of $G$, respectively.
A tangle $\T$ in $\Gamma$ of order $\theta$ is said to be {\it respectful (towards $\Sigma$)} if $\Sigma$ is connected and for every $\Gamma$-normal O-arc  $F$ in $\Sigma$ with $\lvert F \cap V(\Gamma) \rvert < \theta$, there is a closed disk $\Delta \subseteq \Sigma$ with $\partial\Delta=F$ such that $(\Gamma \cap \Delta, \Gamma \cap \overline{\Sigma-\Delta}) \in \T$.
It is clear that $\Delta$ has to be unique, and we write $\Delta = \ins(F)$; the function {\it ins} is called the {\it inside function} of $\T$.

Let $\Gamma$ be a $2$-cell drawing in a surface $\Sigma$.
We say that a drawing $K$ in $\Sigma$ is a {\it radial drawing} of $\Gamma$ if it satisfies the following.
	\begin{itemize}
		\item $U(\Gamma) \cap U(K) = V(\Gamma) \subseteq V(K)$.
		\item Each region $r$ of $\Gamma$ contains a unique vertex of $K$.
		\item $K$ is a drawing of a bipartite graph, and $(V(\Gamma), V(K)-V(\Gamma))$ is a bipartition of it.
		\item For every $v \in V(\Gamma)$, the edges of $K \cup \Gamma$ incident with $v$ belong alternately to $\Gamma$ and to $K$ (in their cyclic order around $v$).
	\end{itemize}

Let $\Gamma$ be a $2$-cell drawing in a surface $\Sigma$, and let $K$ be a radial drawing of $\Gamma$.
Assume there exists a respectful tangle $\T$ in $\Gamma$ of order $\theta$.
So $\T$ defines the inside function $\ins(\cdot)$ of $\T$.
If $W$ is a closed walk of $K$, we define $K|W$ to be the subdrawing of $K$ formed by the vertices and the edges in $W$.
If the length of $W$ is less than $2\theta$, then we define $\ins(W)$ to be the union of $U(K|W)$ and $\ins(C)$, taken over all cycles $C$ of $K|W$.
For any two atoms $a,b$ of $K$, define a function $m_\T(a,b)$ as follows:
	\begin{itemize}
		\item if $a=b$, then $m_\T(a,b)=0$;
		\item if $a \neq b$ and $a,b \subseteq \ins(W)$ for some closed walk $W$ of $K$ of length less than $2 \theta$, then $m_\T(a,b) = \min\frac{1}{2} \lvert E(W) \rvert$, taking over all such closed walks $W$;
		\item otherwise, $m_\T(a,b) = \theta$.
	\end{itemize}
Note that $K$ is bipartite, so $m_\T$ is integral.
In addition, for every atom $c$ of $\Gamma$, we define $a(c)$ to be an atom of $K$ such that
	\begin{itemize}
		\item $a(c)=c$ if $c \subseteq V(\Gamma)$; 
		\item $a(c)$ is the region of $K$ including $c$ if $c$ is an edge of $\Gamma$; 
		\item $a(c)=\{v\}$, where $v$ is the vertex of $K$ in $c$, if $c$ is a region of $\Gamma$.
	\end{itemize}
For atoms $b,c$ of $\Gamma$, we define $m_\T(b,c) = m_\T(a(b),a(c))$.
Note that \cite[Theorem 9.1]{rs XI} implies that if $\T$ is a respectful tangle in $\Gamma$, then $m_\T$ is a metric on the atoms of $\Gamma$.

Let $\T$ be a respectful tangle in a 2-cell drawing $G$ in a surface $\Sigma$.
If $X,Y$ are sets of atoms of $G$, then we define $m_\T(X,Y)=\min\{m_\T(x,y): x \in X, y \in Y\}$.
When one of $X$ and $Y$, say $Y$, has size one, we denote $m_\T(X,Y)$ by $m_\T(X,y)$, where $y$ is the unique element of $Y$.

The following is a restatement of a result in \cite{lt}.

\begin{theorem}[{\cite[Theorem 5.3]{lt}}] \label{A distance} 
Let $\Sigma$ be a surface, and let $\Gamma$ be a $2$-cell drawing of a graph in $\Sigma$ with $E(\Gamma) \neq \emptyset$.
Let $\T$ be a respectful tangle of order $\theta$ in $\Gamma$.
Let $x \in V(\Gamma)$.
If $(A,B) \in \T$ is a separation of $\Gamma$ such that $x \in V(A)-V(B)$ and there exists a path $P$ in $A$ from $x$ to a vertex $y \in V(A)$ internally disjoint from $V(B)$, then $m_\T(x,y) \leq \lvert V(A) \cap V(B) \rvert$.
\end{theorem}

We prove a variation of Theorem \ref{A distance} as follows.

\begin{lemma} \label{A distance set}
Let $\Sigma$ be a surface, and let $\Gamma$ be a $2$-cell drawing of a graph in $\Sigma$ with $E(\Gamma) \neq \emptyset$.
Let $\T$ be a respectful tangle of order $\theta$ in $\Gamma$.
Let $X \subseteq V(\Gamma)$.
If $(A,B) \in \T$ is a separation of $\Gamma$ of order less than $\lvert X \rvert$ such that $X \subseteq V(A)$, and subject to this, $A$ is minimal, then $m_\T(X,y) \leq \lvert V(A) \cap V(B) \rvert$ for every $y \in V(A)$.
\end{lemma}

\begin{pf}
Since $\lvert V(A \cap B) \rvert < \lvert X \rvert$ and $X \subseteq V(A)$, $X-V(B) \neq \emptyset$.
By the minimality of $A$, every component of $A-V(A \cap B)$ intersects $X-V(B)$, and each vertex in $V(A \cap B)$ is either in $X$ or adjacent to some vertex in $V(A)-V(B)$.
So for every $y \in V(A)$, either $y \in X$, or there exists a path in $A$ from a vertex $x_y \in X-V(B)$ to $y$ internally disjoint from $V(B)$.
If $y \in X$, then $m_\T(X,y)=0$; otherwise, $m_\T(x_y,y) \leq \lvert V(A \cap B) \rvert$ by Theorem \ref{A distance}.
Hence $m_\T(X,y) \leq \lvert V(A \cap B) \rvert$ for every $y \in V(A)$.
\end{pf}

\subsubsection{Conformal tangles and zones}

Let $G$ and $H$ be graphs.
Let $\T'$ be a tangle in $H$ of order $\theta \geq 2$.
Let $\alpha$ be an $H$-minor in $G$.
Let $\T$ be the set of separations $(A,B)$ of $G$ of order less than $\theta$ such that there exists $(A',B') \in \T'$ with $\alpha(E(A')) = E(A) \cap \alpha(E(H))$.
It was proved in \cite[Theorem (6.1)]{rs X} that $\T$ is a tangle in $G$ of order $\theta$, and we say that $\T$ is the {\it tangle induced by $\T'$}.
We say that $\T'$ is {\it conformal} with a tangle $\T''$ in $G$ if $\T \subseteq \T''$.

Let $\Gamma$ be a $2$-cell drawing in a surface $\Sigma$, and let $\T$ be a respectful tangle of order $\theta$ in $\Gamma$.
Let $x$ be an atom of $\Gamma$.
A {\it $\lambda$-zone around $x$} is an open disk $\Delta$ in $\Sigma$ with $x \subseteq \Delta$ such that $\partial\Delta$ is an O-arc, $\partial\Delta \subseteq U(\Gamma)$, $m_\T(x,y) \leq \lambda$ for every atom $y$ of $\Gamma$ with $y \subseteq \overline{\Delta}$, and if $x \in E(\Gamma)$, then $\lambda \geq 2$.
A {\it $\lambda$-zone} is a $\lambda$-zone around some atom.
Let $\Delta$ be a $\lambda$-zone.
Note that $U(\Gamma) \cap \partial\Delta$ is a cycle, and the drawing $\Gamma' = \Gamma \cap (\Sigma - \Delta)$ is $2$-cell in $\Sigma$.
We say that $\Gamma'$ is the {\it drawing obtained from $\Gamma$ by clearing $\Delta$}.
We say that $\T'$ is a {\it tangle of order $\theta-4\lambda-2$ obtained by clearing $\Delta$} if $\T'$ is a tangle in $\Gamma'$ of order $\theta-4\lambda-2$ satisfying the following.
\begin{itemize}
	\item $\T'$ is respectful with a metric $m_{\T'}$.
	\item $\T'$ is conformal with $\T$.
	\item If $x,y$ are atoms of $\Gamma$ and $x',y'$ are atoms of $\Gamma'$ with $x \subseteq x'$ and $y \subseteq y'$, then $m_\T(x,y) \geq m_{\T'}(x',y') \geq m_\T(x,y)-4\lambda-2$.
\end{itemize}

The following are some useful results.

\begin{theorem}[{\cite[Theorem~(7.10)]{rs XII}}] \label{clean a zone}
Let $\Delta$ be a $\lambda$-zone.
If $\theta \geq 4\lambda+3$, then there exists a unique respectful tangle of order $\theta-4\lambda-2$ obtained by clearing $\Delta$.
\end{theorem}

\begin{theorem}[{\cite[Theorem~(9.2)]{rs XIV}}] \label{big zone contains ball}
Let $\Gamma$ be a $2$-cell drawing in a surface $\Sigma$, and let $\T$ be a respectful tangle in $\Gamma$ of order $\theta$.
Let $x$ be an atom of $\Gamma$, and $\lambda$ an integer with $2 \leq \lambda \leq \theta-4$.
Then there exists a $(\lambda+3)$-zone $\Delta$ around $x$ such that $x' \subseteq \Delta$ for every atom $x'$ of $\Gamma$ with $m_\T(x,x') \leq \lambda$.
\end{theorem}

\begin{lemma}[{\cite[Lemma 5.7]{lt}}] \label{disjoint boundary of zone}
Let $\Gamma$ be a $2$-cell drawing in a surface, $z$ an atom, and $\T$ a respectful tangle in $\Gamma$ of order $\theta$.
Let $\lambda$ be a nonnegative integer, and let $C$ be the cycle of the boundary of a $\lambda$-zone around $z$.
If $\theta \geq \lambda+8$, then there exists a $(\lambda+7)$-zone $\Lambda$ around $z$ such that the cycle bounding $\Lambda$ is disjoint from $C$, and $\Lambda$ contains the $\lambda$-zone bounded by $C$.
\end{lemma}

\subsubsection{Vortices and segregations}

A {\it tree-decomposition} of a graph $G$ is a pair $(T,\X)$, where $T$ is a tree and $\X=\{X_t: t \in V(T)\}$ such that the following hold.
	\begin{itemize}
		\item $\bigcup_{t \in V(T)} X_t = V(G)$.
		\item For every edge $e$ of $G$, there exists $t \in V(T)$ such that $X_t$ contains all ends of $e$.
		\item For every vertex $v$ of $G$, the subgraph of $T$ induced by the set $\{t \in V(T): v \in X_t\}$ is connected.
	\end{itemize}
The {\it width} of $(T,\X)$ is $\max\{\lvert X_t \rvert-1: t \in V(T)\}$.
The {\it tree-width} of $G$ is the minimum width of a tree-decomposition of $G$.
The {\it adhesion} of $(T,\X)$ is $\max\{\lvert X_t \cap X_{t'} \rvert: tt' \in E(T)\}$.
We say that a tree-decomposition $(T,\X)$ is a {\it path-decomposition} if $T$ is path.

A {\it society} is a pair $(S,\Omega)$, where $S$ is a graph and $\Omega$ is a cyclic permutation of a subset $\overline{\Omega}$ of $V(S)$.
Let $\rho$ be a nonnegative integer.
A society $(S,\Omega)$ is a {\it $\rho$-vortex} if for all distinct $u,v \in \overline{\Omega}$, there do not exist $\rho+1$ mutually disjoint paths of $S$ between $I \cup \{u\}$ and $J \cup \{v\}$, where $I$ is the set of vertices in $\overline{\Omega}$ after $u$ and before $v$ in $\Omega$, and $J$ is the set of vertices in $\overline{\Omega}$ after $v$ and before $u$ in $\Omega$.

Let $(S,\Omega)$ be a society with $\overline{\Omega} = \{v_1,v_2,...,v_{\lvert\overline{\Omega}\rvert}\}$ in order.
A {\it vortical decomposition} of $(S,\Omega)$ is a path-decomposition $(t_1t_2...t_{\lvert \overline{\Omega} \rvert}, \{X_{t_i}:1 \leq i \leq \lvert \overline{\Omega} \rvert\})$ such that $v_i \in X_{t_i}$ for $1 \leq i \leq  \lvert \overline{\Omega} \rvert$.
The following theorem ensures the existence of a vortical decomposition.

\begin{theorem}[{\cite[Theorem (8.1)]{rs IX}}] \label{path decomp of vortex}
For every positive integer $\rho$, every $\rho$-vortex has a vortical decomposition with adhesion at most $\rho$.
\end{theorem}

A {\it segregation} of a graph $G$ is a set $\Se$ of societies such that the following conditions hold. 
\begin{itemize}
	\item $S$ is a subgraph of $G$ for every $(S, \Omega) \in \Se$, and $\bigcup_{(S,\Omega) \in \Se}S=G$.
	\item For any distinct $(S,\Omega)$ and $(S', \Omega') \in \Se$, $V(S \cap S') \subseteq \overline{\Omega} \cap \overline{\Omega'}$ and $E(S \cap S') = \emptyset$.
\end{itemize}
We write $V(\Se) = \bigcup_{(S,\Omega) \in \Se}\overline{\Omega}$.
For a tangle $\T$ in $G$, we say that a segregation $\Se$ of $G$ is {\it $\T$-central} if for every $(S,\Omega) \in \Se$, there is no $(A,B) \in \T$ of order at most half of the order of $\T$ with $B \subseteq S$.

Let $\Se$ be a segregation of a graph $G$.
Let $\Se_1,\Se_2$ be subsets of $\Se$ with $\Se_1 \cap \Se_2=\emptyset$ and $\Se_1 \cup \Se_2 = \Se$ such that $\lvert \overline{\Omega} \rvert \leq 3$ for every $(S,\Omega) \in \Se_1$.
We say $\Se$ is {\it effective with respect to $(\Se_1,\Se_2)$} if for every $(S,\Omega) \in \Se_1$ and every $v \in \overline{\Omega}$, there exist $\lvert \overline{\Omega} \rvert-1$ paths in $S$ from $v$ to $\overline{\Omega}-\{v\}$ such that the intersection of them is $\{v\}$.

Let $\Sigma$ be a surface and $\Se = \{(S_1, \Omega_1), ..., (S_k, \Omega_k)\}$ a segregation of $G$.
An {\it arrangement} of $\Se$ in $\Sigma$ is a function $\alpha$ with domain $\Se \cup V(\Se)$, such that the following conditions hold.
\begin{itemize}
	\item For $1 \leq i \leq k$, $\alpha(S_i, \Omega_i)$ is a closed disk $\Delta_i \subseteq \Sigma$, and $\alpha(x) \in \partial\Delta_i$ for each $x \in \overline{\Omega_i}$.
	\item For $1 \leq i<j \leq k$, if $x \in \Delta_i \cap \Delta_j$, then $x=\alpha(v)$ for some $v \in \overline{\Omega_i} \cap \overline{\Omega_j}$.
	\item For all distinct $x,y \in V(\Se)$, $\alpha(x) \neq \alpha(y)$.
	\item For $1 \leq i \leq k$, $\Omega_i$ is mapped by $\alpha$ to a natural order of $\alpha(\overline{\Omega_i})$ determined by $\partial\Delta_i$.
\end{itemize}
An arrangement is {\it proper} if $\Delta_i \cap \Delta_j = \emptyset$ for all $1 \leq i < j \leq k$ such that $\lvert \overline{\Omega_i} \rvert, \lvert \overline{\Omega_j} \rvert >3$.

Let $\Se$ be a segregation of a graph $G$.
Let $\Se_1,\Se_2$ be subsets of $\Se$ with $\Se_1 \cap \Se_2 = \emptyset$ and $\Se_1 \cup \Se_2=\Se$ such that $\lvert \overline{\Omega} \rvert \leq 3$ for every $(S,\Omega) \in \Se_1$.
Let $\alpha$ be a proper arrangement of $\Se$ in a surface $\Sigma$.
The {\it skeleton of $\alpha$ with respect to $(\Se_1,\Se_2)$} is the drawing $\Gamma = (U,V)$ in $\Sigma$ with $V(\Gamma)=\{\alpha(v): v \in V(\Se)\}$ such that $U(\Gamma)$ consists of the boundary of $\alpha(S,\Omega)$ for each $(S,\Omega) \in \Se_1$ with $\lvert \overline{\Omega} \rvert = 3$, and a line in the boundary of $\alpha(S',\Omega')$ with ends $\overline{\Omega'}$ for each $(S',\Omega') \in \Se_1$ with $\lvert \overline{\Omega'} \rvert=2$.
Note that we do not add any edges into the skeleton for $(S,\Omega) \in \Se_1$ with $\lvert \overline{\Omega} \rvert \leq 1$ and for $(S,\Omega) \in \Se_2$.
Note that if $\Se$ is an effective segregation with respect to $(\Se_1,\Se_2)$, then the skeleton of $\alpha$ with respect to $(\Se_1,\Se_2)$ is a minor of $G$: contracting some edges in $S$ for $(S,\Omega) \in \Se_1$ and deleting redundant vertices and edges; we call this minor a {\it natural minor}.

We shall prove a lemma that allows us to modify a segregation but still keep the $\T$-central property.
We will need the following lemma proved in \cite{lt}.

\begin{lemma}[{\cite[Lemma 6.4]{lt}}] \label{make segregation central}
Let $\rho$ be an integer, $G$ a graph, $\T$ a tangle in $G$ of order at least $2\rho+2$, and $\Se$ a segregation of $G$.
If $(S,\Omega) \in \Se$ is a $\rho$-vortex and there exists no $(A,B) \in \T$ of order at most $2\rho+1$ such that $B \subseteq S$, then there exists no $(A',B') \in \T$ of order at most the half of the order of $\T$ such that $B' \subseteq S$.
\end{lemma}

\begin{lemma} \label{modifying making central}
Let $\rho$ be a positive integer and $\Sigma$ a surface.
Let $G$ be a graph and $\T$ a tangle in $G$ of order at least $2\rho+2$.
Let $\Se$ be a $\T$-central segregation of $G$, and let $\Se_1$ be a subset of $\Se$ such that $\lvert\overline{\Omega} \rvert \leq 3$ for every $(S,\Omega) \in \Se_1$.
Let $\Se_1'$ be a subset of $\Se_1$, and let $\Se'$ be a segregation of $G$ such that $\Se_1' \subseteq \Se'$ and $(S,\Omega)$ is a $\rho$-vortex for every $(S,\Omega) \in \Se'-\Se_1'$.
Let $\alpha$ be a proper arrangement of $\Se'$ in $\Sigma$, and let $G'$ be the skeleton of $\alpha$ with respect to $(\Se'_1,\Se'-\Se_1')$.
Assume that there exists a natural $G'$-minor $\beta$ in $G$ and a tangle $\T'$ in $G'$ conformal with $\T$.
If $\T'$ has order at least $2\rho+2$, then $\Se'$ is $\T$-central.
\end{lemma}

\begin{pf}
Since $\Se'_1 \subseteq \Se_1$ and $\Se$ is $\T$-central, by Lemma \ref{make segregation central}, it suffices to prove that there exists no $(A,B) \in \T$ with order at most $2\rho+1$ such that $B \subseteq S$ for some $(S,\Omega) \in \Se'-\Se'_1$.
Suppose to the contrary that there exist $(S,\Omega) \in \Se'-\Se_1'$ and $(A,B) \in \T$ of order at most $2\rho+1$ with $B \subseteq S$.
Define $(A',B')$ to be a separation of $G'$ such that $\beta(E(A'))=E(A) \cap \beta(E(G'))$, $\beta(E(B'))=E(B) \cap \beta(E(G'))$, and for every $C \in \{A,B\}$, a vertex $v \in V(G')$ is in $V(C')$ if and only if $\beta(v) \cap V(C) \neq \emptyset$. 
Note that $\lvert V(A' \cap B') \rvert \leq \lvert V(A \cap B) \rvert \leq 2\rho+1$.
Since $\T'$ is conformal with $\T$ and $(A,B) \in \T$, $(A',B') \in \T'$.
Since $V(S) \cap V(G') \subseteq \overline{\Omega}$, $V(B') \subseteq \overline{\Omega}$.
So every vertex in $V(B')-V(A')$ has degree zero in $B'$.
By (T1), we may repeatedly remove vertices in $V(B')-V(A')$ from $B'$ and add them into $A'$ to obtain a separation $(A'',B'') \in \T'$ with $V(A'' \cap B'')=V(A' \cap B')$ such that $V(B'') \subseteq V(A'' \cap B'')$.
Hence, $\lvert V(B'') \rvert \leq 2\rho+1$, so $(B'',A'') \in \T'$ by (T1) and (T3), a contradiction.
This proves that $\Se'$ is $\T$-central.
\end{pf}

\subsubsection{Enlarging a vortex}

Let $\Gamma$ be a 2-cell drawing in a surface having a respectful tangle $\T$.
Let $\Lambda$ be a $\lambda$-zone (with respect to $m_\T$) around some atom of $\Gamma$ for some nonnegative integer $\lambda$.
For every $v \in V(\Gamma) \cap \partial\Lambda$, a {\it loose component with respect to $(\Lambda,v,\T)$} is a component $L$ of $\Gamma-v$ such that some vertex of $L$ is adjacent to $v$, and there exists no separation $(A,B) \in \T$ with $V(A \cap B)=\{v\}$ and $V(B)=V(L) \cup \{v\}$; we call $v$ the {\it attachment} of $L$.
{\it A loose component with respect to $(\Lambda,\T)$} is a loose component with respect to $(\Lambda,v,\T)$ for some $v \in V(\Gamma) \cap \partial\Lambda$.

The following result proved in \cite{lt} will be useful.

\begin{lemma}[{\cite[Lemma 6.1]{lt}}] \label{extend a vortex}
Let $t,\rho,\theta$ be nonnegative integers.
Let $G$ be a graph.
Let $\Se=\Se_1 \cup \Se_2$ with $\Se_1 \cap \Se_2=\emptyset$ be a segregation of $G$ such that $\lvert \overline{\Omega} \rvert \leq 3$ for every $(S,\Omega) \in \Se_1$.
Let $\alpha$ be a proper arrangement of $\Se$ with respect to $(\Se_1,\Se_2)$ of $G$ in a surface $\Sigma$.
Let $(S,\Omega) \in \Se$ be a $\rho$-vortex.
Let $G'$ be the skeleton of $\alpha$.
Let $\T'$ be a respectful tangle in $G'$ of order $\theta$.
If $G'$ is $2$-cell and $\theta \geq 4t+59$, then there exists a cycle $C$ such that the following hold.
	\begin{enumerate}
		\item $C$ bounds a $(t+14)$-zone $\Lambda$ in $G'$ around some vertex in $\bar{\Omega}$.
		\item $\Lambda$ contains every atom $x$ of $G'$ with $m_{\T'}(x,y) \leq t$ for some $y \in \bar{\Omega}$.
		\item The closure of $\Lambda$ contains $\alpha(S,\Omega)$.
		\item\label{def_C_vortex} Let $S'$ be the union of $S''$ over all societies $(S'',\Omega'') \in \Se$ with 
			\begin{itemize}
				\item either $\alpha(S'',\Omega'') \subseteq \overline{\Lambda}$, or 
				\item $\lvert \overline{\Omega''} \rvert=2$ and $\alpha(S'',\Omega'') \cap E(C) \neq \emptyset$, or
				\item $\overline{\Omega''}$ is contained in the union of some loose component with respect to $(\Lambda,\T')$ and its attachment, or 
				\item $\lvert \overline{\Omega''} \rvert=1$ and $\overline{\Omega''} \subseteq V(C)$.
			\end{itemize}
			Let $\overline{\Omega'}=V(C)-\{x \in V(C):$ every edge of $G'$ incident with $x$ is either contained in $\overline{\Lambda}$ or incident in $G'$ with a vertex in a loose component with respect to $(\Lambda,\T')\}$, and let $\Omega'$ be a cyclic ordering consistent with the cyclic ordering of $C$.
If every $(S'',\Omega'') \in \Se_2$ with $\alpha(S'',\Omega'') \subseteq \overline{\Lambda}$ is a $\rho_{S''}$-vortex for some nonnegative integer $\rho_{S''}$, then $(S',\Omega')$ is a $(\rho+t+8+ \sum_{S''} \rho_{S''})$-vortex, where the sum is over all societies $(S'',\Omega'') \in \Se_2-\{(S,\Omega)\}$ with $\alpha(S'',\Omega'') \subseteq \overline{\Lambda}$.
		\item Let $\Se_1^*=\Se_1-\{(S'',\Omega'') \in \Se_1: S'' \subseteq S'\}$ and $\Se^*_2=(\Se_2-\{(S'',\Omega'') \in \Se_2: S'' \subseteq S'\}) \cup \{(S',\Omega')\}$.
			If $m_{\T'}(x,y) \geq 3$ for every atom $x \subseteq \partial\Lambda$ and $y \in V(\overline{\Omega''})$ with $(S'',\Omega'') \in \Se_2^*-\{(S',\Omega')\}$, then $\Se^*$ is a segregation, and there exists a proper arrangement $\alpha^*$ of $\Se^*_1 \cup \Se^*_2$ with respect to $(\Se_1^*,\Se_2^*)$ such that the skeleton $G^*$ of $\alpha^*$ 
			\begin{itemize}
				\item is obtained from $G'$ by clearing $\Lambda$ and deleting some edges in $E(C)$ and all loose components with respect to $(\Lambda,\T')$ and deleting all resulting isolated vertices, 
				\item is 2-cell, and 
				\item has a respectful tangle $\T^*$ conformal with $\T'$ of order at least $\theta-4t-58$ such that $m_{\T'}(x',y') \geq m_{\T^*}(x,y) \geq m_{\T'}(x',y')-4t-58$ for all atoms $x,y$ of $G^*$, where $x',y'$ are atoms of $G'$ with $x' \subseteq x$ and $y' \subseteq y$.
			\end{itemize}
	\end{enumerate}
\end{lemma}

Let $G$ be a graph.
Let $\Se=\Se_1 \cup \Se_2$ with $\Se_1 \cap \Se_2=\emptyset$ be a segregation of $G$ such that $\lvert \overline{\Omega} \rvert \leq 3$ for every $(S,\Omega) \in \Se_1$.
Let $\alpha$ be a proper arrangement of $\Se$ with respect to $(\Se_1,\Se_2)$ of $G$ in a surface $\Sigma$.
Let $G'$ be the skeleton of $\alpha$.
Let $\T'$ be a respectful tangle in $G'$ of order $\theta$.
For a cycle $C$ bounds a $\lambda$-zone $\Lambda$ in $G'$ for some $\lambda \geq 0$, we define the {\it $C$-vortex} to be the society $(S',\Omega')$ mentioned in Statement \ref{def_C_vortex} in Lemma \ref{extend a vortex}.

\subsection{Better arrangements}

Let $\kappa,\rho$ be nonnegative integers.
For subsets $\Se_1,\Se_2$ of a segregation $\Se$, we say that $(\Se_1,\Se_2)$ is a {\it $(\kappa,\rho)$-witness} of $\Se$ if $\Se_1 \cup \Se_2=\Se$, $\Se_1 \cap \Se_2=\emptyset$, $\lvert \overline{\Omega} \rvert \leq 3$ for every $(S,\Omega) \in \Se_1$, $\lvert \Se_2 \rvert \leq \kappa$ and every member of $\Se_2$ is a $\rho$-vortex.
(Note that we do not require that $\lvert \overline{\Omega} \rvert >3$ for every member $(S,\Omega) \in \Se_2$.)
We say that $\Se$ is a {\it $(\kappa,\rho)$-segregation} if it has a $(\kappa,\rho)$-witness.

Let $\Sigma$ be a surface, $\theta$ an integer and $\phi$ a nondecreasing function with domain $\Z$.
We say that an arrangement $\alpha$ of a segregation $\Se$ of a graph $G$ with a $(\kappa,\rho)$-witness $(\Se_1,\Se_2)$ in $\Sigma$ is a {\it $(\Sigma, \theta,\phi)$-arrangement with respect to $(\Se_1,\Se_2)$} if the following conditions hold.
	\begin{itemize}
		\item $\alpha$ is a proper arrangement of $\Se$ in $\Sigma$.
		\item The skeleton of $\alpha$ with respect to $(\Se_1,\Se_2)$ is a 2-cell drawing in $\Sigma$.
		\item There exists a respectful tangle $\T'$ of order at least $\theta$ in the skeleton of $\alpha$ with respect to $(\Se_1,\Se_2)$.
		\item Let $\rho'$ be the smallest nonnegative integer such that every member of $\Se_2$ is a $\rho'$-vortex. 
			Then $m_{\T'}(\overline{\Omega}, \overline{\Omega'}) \geq \phi(\rho')$ for all distinct members $(S,\Omega),(S',\Omega')$ of $\Se_2$.
	\end{itemize}
We further say this $(\Sigma, \theta,\phi)$-arrangement $\alpha$ is a {\it $(\Sigma, \theta,\phi,\T)$-arrangement with respect to $(\Se_1,\Se_2)$} if the following conditions hold.
	\begin{itemize}
		\item $\T$ is a tangle in $G$ such that $\Se$ is a $\T$-central segregation, 
		\item The skeleton of $\alpha$ with respect to $(\Se_1,\Se_2)$ is a natural minor of $G$.
		\item The respectful tangle $\T'$ of the skeleton of $\alpha$ with respect to $(\Se_1,\Se_2)$ mentioned above is conformal with $\T$.
	\end{itemize}

\begin{lemma} \label{make vortices far apart}
For any positive integers $\kappa, \rho$ and nondecreasing function $\phi$ with domain $\Z$, there exists an integer $\rho^*=\rho^*(\kappa,\rho,\phi)$ such that for every positive integer $\theta^*$, there exists an integer $\theta=\theta(\kappa,\rho,\phi,\theta^*)$ such that if $\T$ is a tangle in a graph $G$, and $\Se$ is a $\T$-central segregation of $G$ with a $(\kappa,\rho)$-witness $(\Se_1,\Se_2)$ with a proper arrangement $\alpha$ in a surface $\Sigma$ such that the skeleton of $\alpha$ with respect to $(\Se_1,\Se_2)$ is a natural minor of $G$ and is a 2-cell drawing in a surface $\Sigma$ with a respectful tangle with order at least $\theta$ conformal with $\T$, then there exists a $\T$-central segregation $\Se^*$ of $G$ with a $(\kappa,\rho^*)$-witness $(\Se_1^*,\Se_2^*)$ such that the following statements hold.
	\begin{enumerate}
		\item $\Se_1^* \subseteq \Se_1$ and $\bigcup_{(S,\Omega) \in \Se_2}S \subseteq \bigcup_{(S,\Omega) \in \Se_2^*}S$.
		\item $\Se^*$ has a $(\Sigma,\theta^*,\phi,\T)$-arrangement with respect to $(\Se_1^*,\Se_2^*)$ in $\Sigma$.
	\end{enumerate}
\end{lemma}

\begin{pf}
Let $\kappa,\rho$ be positive integers, and let $\phi$ be a nondecreasing function with domain $\Z$.
Let $(\tau_m: n \geq 0)$, $(a^*_{m,n}: m, n \geq 0)$ and $(b^*_{m,n}: m,n \geq 0)$ be sequences defined as follows.
	\begin{itemize}
		\item $\tau_0=\rho$.
		\item For each $m \geq 0$, define $a^*_{m,0}=\phi(\tau_m)$.
		\item For each $m \geq 0, n \geq 0$, define $b^*_{m,n}=2 (a^*_{m,n}+31)$. 
		\item For each $m \geq 0, n \geq 1$, define $a^*_{m,n}=\kappa b^*_{m,n-1}$.
		\item For each $m \geq 1$, define $\tau_m= (\kappa+1)\tau_{m-1}+a^*_{m-1,\kappa}+8$.
	\end{itemize}
Define $\rho^*=\tau_{\kappa-1}$.
Let $\theta^*$ be a positive integer.
Define $\theta=\theta^*+(4b^*_{\kappa,\kappa}+58\kappa)(\kappa-1)+2\rho^*+2$.

Let $G$ be a graph, and let $\T$ be a tangle in $G$.
Let $\Sigma$ be a surface.

We say that $(\Se,\Se_1,\Se_2)$ is {\it $i$-good} if $i$ is a nonnegative integer such that $\Se$ is a $\T$-central segregation of $G$ with a $(\kappa-i, \tau_i)$-witness $(\Se_1,\Se_2)$ and has a proper arrangement of $\Se$ in $\Sigma$ whose skeleton with respect to $(\Se_1,\Se_2)$ is a natural minor of $G$ and is a 2-cell drawing in $\Sigma$ with a respectful tangle with order at least $\theta-(4b^*_{\kappa,\kappa}+58\kappa)i$ conformal with $\T$.

Let $\Se$ be a $\T$-central segregation of $G$ with a $(\kappa, \rho)$-witness $(\Se_1,\Se_2)$.
Let $\alpha$ be a proper arrangement of $\Se$ in $\Sigma$ such that the skeleton $G'$ of $\alpha$ with respect to $(\Se_1,\Se_2)$ is a natural minor of $G$ and is a 2-cell drawing in $\Sigma$ with a respectful tangle $\T'$ with order at least $\theta$ conformal with $\T$.
Note that $(\Se,\Se_1,\Se_2)$ is 0-good.

Let $i$ be the maximal integer with $0 \leq i \leq \kappa$ such that $(\Se,\Se_1,\Se_2)$ is $i$-good.
Suppose that this lemma does not hold.
We further assume that $\Se,\Se_1,\Se_2$ are chosen so that $i$ is maximum among all counterexamples.
In other words, there exists no segregation $\Se^*$ satisfying the conclusions of this lemma with respect to $(\Se,\Se_1,\Se_2)$, but for every triple $(\Se',\Se_1',\Se_2')$ that is $j$-good for some $j>i$, there exists a segregation $\Se^*$ satisfying the conclusions of this lemma with respect to $(\Se',\Se_1',\Se_2')$.

Note that this lemma holds if $\kappa-i \leq 1$, since we can take $(\Se^*,\Se^*_1,\Se^*_2)=(\Se,\Se_1,\Se_2)$.
So $i \leq \kappa-2$.
By the maximality of $i$, there does not exist a $\T$-central segregation $\Se'$ of $G$ and subsets $\Se_1',\Se_2'$ of $\Se'$ with $\Se_1' \cup \Se_2'=\Se'$, $\lvert \Se_2' \rvert < \lvert \Se_2 \rvert$, $\Se_1' \subseteq \Se_1$ and $\bigcup_{(S,\Omega) \in \Se_2}S \subseteq \bigcup_{(S,\Omega) \in \Se_2'}S$ such that $(\Se',\Se_1',\Se_2')$ is $(i+1)$-good.

Let $\rho_0$ be the smallest integer such that every member of $\Se_2$ is a $\rho_0$-vortex.
Note that $\rho_0 \leq \tau_i$.
Define two sequences $(a_n: n \geq 0)$ and $(b_n: n \geq 0)$ such that $a_0=\phi(\rho_0), b_0=2(a_0+31)$ and for $n \geq 1$, $a_n=\kappa b_{n-1}$ and $b_n=2(a_n+31)$.
Since $\rho_0 \leq \tau_i$, it is easy to show that for every $n \geq 0$, $a_n \leq a^*_{i,n}$ and $b_n \leq b^*_{i,n}$ by induction on $n$.

\noindent{\bf Claim 1:} There exists a collection $\C$ of cycles in $G'$ with the following properties.
	\begin{itemize}
		\item Each member of $\C$ bounds a disk in $\Sigma$.
		\item For each member $C$ of $\C$, there exists $(S_C,\Omega_C) \in \Se_2$ such that $C$ satisfies the conclusions of Lemma \ref{extend a vortex} for choosing $(S,\Omega)=(S_C,\Omega_C)$ and $t=a_{\lvert \Se_2 \rvert-\lvert \C \rvert}$.
		\item For each $(S,\Omega) \in \Se_2$, there exists $C \in \C$ such that $\alpha(S,\Omega)$ is contained in the closed disk bounded by $C$.
		\item For any two distinct members $C_1,C_2$ of $\C$, the closure of the disks bounded by $C_1$ and $C_2$ are disjoint, and $m_{\T'}(V(C_1),V(C_2)) \geq 3$. 
	\end{itemize}

\noindent{\bf Proof of Claim 1:}
For each member $(S,\Omega) \in \Se_2$, let $C_S$ be the cycle $C$ mentioned in Lemma \ref{extend a vortex} by choosing $t=\phi(\rho_0)$ and $\rho=\rho_0$.
Then the collection $\{C_S: (S,\Omega) \in \Se_2\}$ satisfies the first three bullets of the conditions for $\C$ in this claim.
Let $\C'$ be a collection of cycles in $G'$ satisfying the first three bullets of this claim, and subject to this, $\lvert \C' \rvert$ is minimum.
Since $\{C_S: (S,\Omega) \in \Se_2\}$ is a candidate for $\C'$, $\lvert \C' \rvert \leq \lvert \Se_2 \rvert \leq \kappa$.

We shall prove that $\C'$ satisfies all bullets of the conditions for $\C$ in this claim.

For each $D \in \C'$, let $\Lambda_D$ be the open disk in $\Sigma$ bounded by $D$.
Since $\C'$ satisfies this claim, for every $D \in \C'$, there exists $(S_D,\Omega_D) \in \Se_2$ such that $D$ satisfies the conclusions of Lemma \ref{extend a vortex} for choosing $(S,\Omega)=(S_D,\Omega_D)$ and $t=a_{\lvert \Se_2 \rvert-\lvert \C' \rvert}$.
For each $D \in \C'$, let $\hat{D}$ be the cycle $C$ mentioned in Lemma \ref{extend a vortex} for choosing $(S,\Omega)=(S_D,\Omega_D)$ and $t=a_{\lvert \Se_2 \rvert-\lvert \C' \rvert}+17$, and let $\Lambda_{\hat{D}}$ be the open disk bounded by $\hat{D}$.

Let $\Delta_1,\Delta_2,...,\Delta_\ell$ (for some $\ell$) be the connected components of $\bigcup_{C \in \C'} \overline{\Lambda_{\hat{C}}}$ in $\Sigma$, where $\overline{\Lambda_{\hat{C}}}$ is the closure of $\Lambda_{\hat{C}}$.
We may assume that $\overline{\Lambda_{\hat{C}}} \cap \overline{\Lambda_{\hat{C'}}} \neq \emptyset$ for some distinct members $C,C' \in \C'$, for otherwise $\C'$ satisfies the fourth bullet of this claim and we are done.
Hence $\ell < \lvert \C' \rvert$ and $\lvert \C' \rvert \geq 2$.
For each $j \in [\ell]$, let $k_j$ be the number of members $C \in \C'$ such that $\Lambda_{\hat{C}} \subseteq \Delta_j$, and let $(S_j,\Omega_j)$ be a member of $\Se_2$ such that it equals $(S_C,\Omega_C)$ for some $C \in \C'$ with $\Lambda_{\hat{C}} \subseteq \Delta_j$.
Note that for each $j \in [\ell]$ and each atom $x$ of $G' \cap \Delta_j$, $m_{\T'}(x,\overline{\Omega_j}) \leq k_j \cdot 2(a_{\lvert \Se_2 \rvert - \lvert \C' \rvert}+17+14)=k_jb_{\lvert \Se_2 \rvert - \lvert \C' \rvert} \leq \kappa b_{\lvert \Se_2 \rvert - \lvert \C' \rvert}=a_{\lvert \Se_2 \rvert - \lvert \C' \rvert +1} \leq a_{\lvert \Se_2 \rvert - \ell}$.
For each $j \in [\ell]$, define $D_j$ to be the cycle in $G'$ bounding an $(a_{\lvert \Se_2 \rvert - \ell}+14)$-zone $\Lambda_j$ mentioned in Lemma \ref{extend a vortex} by choosing $(S,\Omega)=(S_j,\Omega_j)$ and $t=a_{\lvert \Se_2 \rvert - \ell}$.
Note that for each $j \in [\ell]$, $\overline{\Lambda_j} \supseteq \Delta_j$.
So $\bigcup_{j=1}^\ell\overline{\Lambda_j} \supseteq \bigcup_{C \in \C'} \overline{\Lambda_{\hat{C}}} \supseteq \bigcup_{(S,\Omega) \in \Se_2}\alpha(S,\Omega)$.
Hence the collection $\{D_j: 1 \leq j \leq \ell\}$ satisfies the first three bullets of the conditions for $\C$ but has size less than $\lvert \C' \rvert$, a contradiction.
This proves the claim.
$\Box$

Let $\C$ be a collection satisfying Claim 1 with $\lvert \C \rvert$ maximum.
If $\lvert \C \rvert = \lvert \Se_2 \rvert$, then $m_{\T'}(\overline{\Omega},\overline{\Omega'}) \geq a_0=\phi(\rho_0)$ for any distinct members $(S,\Omega),(S',\Omega')$ of $\Se_2$, so $\alpha$ is a $(\Sigma,\theta^*,\phi,\T)$-arrangement of $\Se$, a contradiction.
So $\lvert \C \rvert < \lvert \Se_2 \rvert$.
For each member $C \in \C$, let $(S'_C,\Omega'_C)$ be the society $(S',\Omega')$ mentioned in Lemma \ref{extend a vortex} by choosing $(S,\Omega)=(S_C,\Omega_C)$ and $t=a_{\lvert \Se_2 \rvert - \lvert \C \rvert}$, where $(S_C,\Omega_C)$ is the society mentioned in the second bullet in Claim 1.
By Lemma \ref{extend a vortex}, each $(S'_C,\Omega'_C)$ is a $(\rho_0+a_{\lvert \Se_2 \rvert-\lvert \C \rvert}+8+\lvert \Se_2 \rvert\rho_0)$-vortex.
Note that $\rho_0+a_{\lvert \Se_2 \rvert-\lvert \C \rvert}+8+\lvert \Se_2 \rvert\rho_0 \leq (\kappa+1)\rho_0+a_{\kappa}+8 \leq (\kappa+1)\tau_i+a^*_{i,\kappa}+8 \leq \tau_{i+1}$.

If $V(S'_C) \cap V(S'_{C'}) \neq \emptyset$ for distinct members $C,C'$ of $\C$, then since the closed disks bounded by $C$ and $C'$ are disjoint, some loose component for $C$ intersects some loose component for $C'$, but it is impossible. 
Hence $V(S'_C) \cap V(S'_{C'}) = \emptyset$ for distinct members $C,C'$ of $\C$.

Define $\Se_2'=\{(S'_C,\Omega'_C): C \in \C\}$ and $\Se_1'=\{(S_1,\Omega_1) \in \Se_1: S_1 \not\subseteq \bigcup_{(S'_C,\Omega'_C) \in \Se_2'}S'_C\}$.
Recall that $\lvert \Se_2' \rvert = \lvert \C \rvert < \lvert \Se_2 \rvert \leq \kappa-i$.
Let $\Se'=\Se_1' \cup \Se_2'$.
Since $m_{\T'}(V(C),V(C')) \geq 3$ for all distinct members $C,C'$ of $\C$, Lemma \ref{extend a vortex} implies that $\Se'$ is a segregation of $G$ with a $(\kappa-i-1, \tau_{i+1})$-witness $(\Se_1',\Se_2')$, and there exists a proper arrangement $\alpha'$ in $\Sigma$ such that the skeleton $G''$ of $\alpha'$ with respect to $(\Se_1',\Se_2')$ is a natural minor of $G$, is 2-cell, and has a respectful tangle $\T''$ of order at least $\theta-(4b^*_{\kappa,\kappa}+58\kappa)i-\kappa(4a_{\kappa-1}+58) \geq \theta-(4b^*_{\kappa,\kappa}+58\kappa)i-(4b_\kappa+58\kappa) \geq \theta-(4b^*_{\kappa,\kappa}+58\kappa)(i+1)$ in $G''$ conformal with $\T$.
Furthermore, $\Se_1' \subseteq \Se_1$, $\bigcup_{(S,\Omega) \in \Se_2} S \subseteq \bigcup_{(S,\Omega) \in \Se_2'}S$ and $\lvert \Se'_2 \rvert=\lvert \C \rvert < \lvert \Se_2 \rvert$.
In addition, Lemma \ref{modifying making central} implies that $\Se'$ is $\T$-central.
This proves that $(\Se',\Se'_1,\Se'_2)$ is $(i+1)$-good, a contradiction. 
\end{pf}

\bigskip

The following lemma provides a tool to sweep atoms in a small zone of the skeleton of an arrangement of a segregation into vortices.

\begin{lemma} \label{sweeping balls into vortices}
For any positive integers $\kappa,k$, there exists $\kappa^*=\kappa^*(\kappa,k)$ such that for all positive integers $\rho,\lambda$ and nondecreasing function $\phi$ with domain $\Z$, there exists $\rho^*=\rho^*(\kappa,k,\rho,\lambda,\phi)$ such that for every integer $\theta^*$, there exists $\theta=\theta(\kappa,k,\rho,\lambda,\phi,\theta^*)$ such that the following hold.
Let $\T$ be a tangle in a graph $G$, and let $\Se$ be a $\T$-central segregation of $G$ with a $(\kappa,\rho)$-witness $(\Se_1,\Se_2)$ with a proper arrangement $\alpha$ in a surface $\Sigma$ such that the skeleton $G'$ of $\alpha$ with respect to $(\Se_1,\Se_2)$ is a natural minor of $G$ and is a 2-cell drawing in a surface $\Sigma$ with a respectful tangle with order at least $\theta$ that is conformal with $\T$.
If $\Lambda_1,\Lambda_2,...,\Lambda_k$ are $\lambda$-zones around some atoms of $G'$, then there exists a $\T$-central segregation $\Se^*$ of $G$ with a $(\kappa^*,\rho^*)$-witness $(\Se_1^*,\Se_2^*)$ such that the following hold.
	\begin{enumerate}
		\item $\Se_1^* \subseteq \Se_1$ and $(\bigcup_{(S,\Omega) \in \Se_2}S) \cup (\bigcup_{(S,\Omega) \in \Se_1, \alpha(S,\Omega) \subseteq \bigcup_{i=1}^k \overline{\Lambda_i}}S) \subseteq \bigcup_{(S,\Omega) \in \Se_2^*}S$.
		\item $\Se^*$ has a $(\Sigma,\theta^*,\phi,\T)$-arrangement with respect to $(\Se_1^*,\Se_2^*)$ in $\Sigma$.
	\end{enumerate}
\end{lemma}

\begin{pf}
Let $\kappa,k$ be positive integers.
Define $\kappa^*=\kappa+k$.

Let $\rho,\lambda$ be positive integers, and let $\phi$ be a nondecreasing function with domain $\Z$.
Let $\lambda_0=\lambda+5$ and $\lambda_i = 3\lambda_{i-1}+34$ for every $i \geq 1$.
Let $\rho'=\lambda_{\kappa^*}+22+(\kappa+1) \rho$.
Define $\rho^*$ to be the number $\rho^*$ mentioned in Lemma \ref{make vortices far apart} by taking $\kappa=\kappa^*$, $\rho=\rho'$ and $\phi=\phi$.

Let $\theta^*$ be an integer.
Let $\theta'$ be the number $\theta$ mentioned in Lemma \ref{make vortices far apart} by taking $\kappa=\kappa^*$, $\rho=\rho'$, $\phi=\phi$ and $\theta^*=\theta^*$.
Define $\theta=\theta' + 4\rho'\kappa^*+2\rho'+2$.

Let $\T$ be a tangle in a graph $G$, and let $\Se$ be a $\T$-central segregation of $G$ with a $(\kappa,\rho)$-witness $(\Se_1,\Se_2)$ with a proper arrangement $\alpha$ in a surface $\Sigma$ such that the skeleton $G'$ of $\alpha$ with respect to $(\Se_1,\Se_2)$ is a natural minor of $G$ and is a 2-cell drawing in a surface $\Sigma$ with a respectful tangle $\T'$ with order at least $\theta$ that is conformal with $\T$.
Let $\Lambda_1,\Lambda_2,...,\Lambda_k$ be $\lambda$-zones around some atoms of $G'$.

For each $i \in [k]$, let $v_i$ be an atom of $G'$ such that $\Lambda_i$ is a $\lambda$-zone around $v_i$.
For each $(S,\Omega) \in \Se_2$, let $v_S$ be a vertex in $\overline{\Omega}$, and let $\Lambda_S$ be a 5-zone around $v_S$ containing $\overline{\Omega}$.
Note that the existence of $\Lambda_S$ follows from Lemma \ref{big zone contains ball}.
Let $\W_0 = \{\Lambda_i,\Lambda_S: i \in [k], (S,\Omega) \in \Se_2\}$.

\noindent{\bf Claim 1:} There exist $t \in [k+\lvert \Se_2 \rvert] \cup \{0\}$, a set $R_t \subseteq \{v_i,v_S: i \in [k], (S,\Omega) \in \Se_2\}$ with $\lvert R_t \rvert \leq k+\lvert \Se_2 \rvert-t$ and a collection $\W_t$ of $\lambda_t$-zones around some elements in $R_t$ such that $\bigcup_{W \in \W_t}W \supseteq (\bigcup_{i=1}^k\Lambda_i) \cup \bigcup_{(S,\Omega) \in \Se_2}\alpha(S,\Omega)$ and $m_{\T'}(r_1,r_2) \geq 2\lambda_t+32$ for every distinct elements $r_1,r_2$ of $R_t$.

\noindent{\bf Proof of Claim 1:}
There exist $j \in [k+\lvert \Se_2 \rvert] \cup \{0\}$, a set $R_j \subseteq \{v_i,v_S: i \in [k], (S,\Omega) \in \Se_2\}$ with $\lvert R_j \rvert \leq k+ \lvert \Se_2 \rvert-j$ and a collection $\W_j$ of $\lambda_j$-zones around some elements in $R_j$ such that $\bigcup_{W \in \W_j}W \supseteq (\bigcup_{i=1}^k\Lambda_i) \cup \bigcup_{(S,\Omega) \in \Se_2}\alpha(S,\Omega)$, as we can take $j=0$ and $R_j=\{v_i,v_S:i\in [k], (S,\Omega) \in \Se_2\}$ and $\W_j=\W_0$.
We assume that $j$ is as large as possible.
We may assume that there exist distinct $a,b \in R_j$ such that $m_{\T'}(a,b) \leq 2\lambda_j+31$, for otherwise we are done.
In particular, $R_j \neq \emptyset$, so $j \leq k+\lvert \Se_2 \rvert-1$.
By Lemma \ref{big zone contains ball}, there exists a $(3\lambda_j+34)$-zone $W$ around $a$ containing every atom $x$ of $G'$ with $m_{\T'}(x,a)\leq 3\lambda_j+31$.
For each $v \in \{a,b\}$, if there exists a member of $\W_j$ around $v$, then let $W_v$ be this member; otherwise, let $W_v=\emptyset$.
Since $m_{\T'}(a,b) \leq 2\lambda_j+31$, $W$ contains $W_a \cup W_b$.
Define $\W_{j+1}=(\W_j-\{W_a,W_b\}) \cup \{W\}$ and $R_{j+1}=R_j-\{b\}$.
Since $\lambda_{j+1}=3\lambda_j+34$, the existence of the sets $\W_{j+1}$ and $R_{j+1}$ contradicts the maximality of $j$.
This proves the claim.
$\Box$

Note that for every vertex $v \in V(G')$, there exists $(S_v,\Omega_v) \in \Se$ with $v \in \overline{\Omega_v}$, and $(S_v,\Omega_v)$ is a $\max\{\rho,1\}$-vortex.
By Lemma \ref{extend a vortex}, for each $x \in R_t$, there exists a $(\lambda_t+14)$-zone $\Lambda_x$ around $x$ in $G'$ such that the following hold.
	\begin{itemize}
		\item $\Lambda_x$ contains every atom $y$ in $G'$ with $m_{\T'}(x,y) \leq \lambda_t$.
		\item $\bigcup_{x \in R_t} \Lambda_x \supseteq (\bigcup_{i=1}^k \Lambda_i) \cup (\bigcup_{(S,\Omega) \in \Se_2}\alpha(S,\Omega))$.
		\item Each $\Lambda_x$ defines a $(\lambda_t+22+(\kappa+1) \rho)$-vortex $(S_x',\Omega_x')$ as the vortex $(S',\Omega')$ mentioned in Lemma \ref{extend a vortex}.
	\end{itemize}
Since $\lvert \Se_2 \rvert \leq \kappa$, $t \leq k+\kappa \leq \kappa^*$ and $\lambda_t+22+(\kappa+1) \rho \leq \rho'$.

Define $\Se_2'=\{(S_x',\Omega_x'): x \in R_t\}$.
Define $\Se_1'=\{(S,\Omega) \in \Se_1: S \not \subseteq \bigcup_{(S',\Omega') \in \Se_2'}S'\}$, and $\Se'=\Se_1' \cup \Se'_2$.
Since $m_{\T'}(r_1,r_2) \geq 2\lambda_t+32 = 2(\lambda_t+14)+4$ for any distinct $r_1,r_2 \in R_t$, Lemma \ref{extend a vortex} implies that $\Se'$ is a segregation of $G$ with a $(\kappa^*,\rho')$-witness $(\Se_1',\Se_2')$, and there exists a proper arrangement of $\Se'$ such that its skeleton with respect to $(\Se_1',\Se_2')$ is a natural minor of $G$, is 2-cell and has a respectful tangle of order at least $\theta-(4\lambda_t+58)\kappa^* \geq \theta'+2\rho'+2$ that is conformal with $\T$. 
By Lemma \ref{modifying making central}, $\Se'$ is $\T$-central.
Then a desired $\T$-central segregation $\Se^*$ of $G$ with a $(\kappa^*,\rho^*)$-witness $(\Se_1^*,\Se_2^*)$ and a $(\Sigma, \theta^*, \phi,\T)$-arrangement with respect to $(\Se_1^*,\Se_2^*)$ in $\Sigma$ can be obtained by applying Lemma \ref{make vortices far apart} to $\Se'$.
This completes the proof.
\end{pf}

\section{Pseudo-embeddings} \label{sec: pseudoembedding}

Let $H$ be a graph and $\Sigma$ a surface.
A {\it pseudo-embedding} of $H$ in $\Sigma$ is a mapping $\pi$ with domain $V(H) \cup E(H)$ such that the following hold. 
\begin{itemize}
	\item $\pi|_{V(H)}$ is an injection from $V(H)$ to $\Sigma$.
	\item For each edge $e$ with ends $u,v$ (not necessarily distinct), $\pi(e)$ is a simple curve $c:[0,1] \rightarrow \Sigma$ such that $c(0)=\pi(u)$, $c(1)=\pi(v)$, $c((0,1)) \cap \pi(V(H))=\emptyset$, and $c([0,1]) \neq \{c(0),c(1)\}$.
	\item If $e,e'$ are distinct edges of $H$, then $\pi(e) \cap \pi(e')$ consists of finitely many points, and for every $x \in \pi(e) \cap \pi(e')$, there exists an open set $B \subseteq \Sigma$ such that $B \cap \pi(e) \cap \pi(e')=\{x\}$.
\end{itemize}
We say that a point $x \in \Sigma$ is a {\it crossing-point} if there exist distinct edges $e,e'$ of $H$ such that $x \in \pi(e) \cap \pi(e')- \pi(V(e) \cup V(e'))$.
And for any curve $\gamma: [0,1] \rightarrow \Sigma$, we also denote the image of $\gamma$ as $\gamma$ when there is no danger for creating confusion.

We say that a disk $\Delta$ in a surface $\Sigma$ is {\it oriented} if $\Delta$ is equipped with a linear ordering $\Omega(p)$ of $\partial\Delta$ for some $p \in \Delta$ such that $p$ is the first element in $\Omega(p)$, and $\Omega(p)$ is consistent with a linear ordering of $\partial \Delta$ obtained by tracing $\partial \Delta$.

A {\it march} in a graph $G$ is a sequence of distinct vertices of $G$.
A {\it rooted graph} is a pair $(G,\sigma)$, where $G$ is graph and $\sigma$ is a march in $G$.
We say two rooted graphs $(G_1,\sigma_1)$ and $(G_2,\sigma_2)$ are {\it isomorphic} if there exists an isomorphism $\iota$ from $G_1$ to $G_2$ such that $\sigma_2=\iota(\sigma_1)$.
(Note that it implies that the number of entries of $\sigma_1$ equals the number of entries of $\sigma_2$.)
If some vertices of $(G_1,\sigma_1)$ and $(G_2,\sigma_2)$ are labelled, then we say that they are {\it isomorphic} if there exists an isomorphism $\iota$ from $G_1$ to $G_2$ with $\sigma_2=\iota(\sigma_1)$ such that 
	\begin{itemize}
		\item for every $v \in V(G_1)$, $v$ has a label if and only if $\iota(v)$ has a label, and
		\item for every $v \in V(G_1)$ that has a label, the labels of $v$ and $\iota(v)$ are the same.
	\end{itemize}

\subsection{General intuition} \label{subsec:intuition_pseudo_embed}

In this subsection, we give a rough intuition for the notions and lemmas proved in this section.
More precise descriptions will be included in other subsections when they are ready to be stated.

Let $G$ be a graph that has a segregation $\Se$ with a $(\kappa,\rho)$-witness $(\Se_1,\Se_2)$ with a proper arrangement $\alpha$ in a surface $\Sigma$.
So for each $(S,\Omega) \in \Se$, we can map each vertex of $S$ to a point in $\alpha(S,\Omega)$ and each edge of $S$ to a line in $\alpha(S,\Omega)$ such that each vertex $v$ in $\overline{\Omega}$ is mapped to $\alpha(v)$.
Hence we get a pseudo-embedding of $G$ in $\Sigma$. 
So every homeomorphic embedding from $H$ into $G$ also leads to a pseudo-embedding of $H$ in $\Sigma$.
Our goal is to show that if there exists no small set in $G$ hitting all subdivisions of $H$, then we can construct a half-integral packing of many subdivisions of $H$.
Note that the intersection of each society in $\Se$ and each homeomorphic embedding from $H$ into $G$ gives a ``small piece'' of a subdivision of $H$.
We will try to combine those small pieces with some extra paths to construct a desired half-integral packing of subdivisions of $H$.

Let $\pi$ be a homeomorphic embedding from $H$ into $G$.
Note that every crossing-point in the pseudo-embedding of $H$ in $\Sigma$ given by $\pi$ is contained in $\alpha(S,\Omega)$ for some $(S,\Omega) \in \Se$.
For each $(S,\Omega) \in \Se_1$, $\lvert \overline{\Omega} \rvert \leq 3$, so it is relatively easy to analyze the intersection of $\alpha(S,\Omega)$ and the image of $\pi$.
The main difficulty lies on the intersection of the image of $\pi$ and $\alpha(S,\Omega)$ for $(S,\Omega) \in \Se_2$.
So we will work on $\Se_1$ after we work on $\Se_2$, and it is the motivation of the notions about addenda defined in Section \ref{subsec:addenda_template}.

Note that $\lvert \Se_2 \rvert \leq \kappa$.
So we restrict ourselves to the intersection of the image of $\pi$ and at most $\kappa$ disks with disjoint closure.
To simplify the situation, we further assume that for each edge $e$ of $H$, if we trace $\pi(e)$ from an endpoint to the other endpoint, then as long as it touches the boundary of a disk, then it must pass through it.
It is the intuition behind the definition of legality introduced in Section \ref{subsec:dive_depth}.
Note that this extra assumption does not make us lose of generality, by using machineries developed in Sections \ref{sec: gauges}-\ref{sec: vortices}.

Fix a disk $\Delta$ among those $\kappa$ disks.
Note that for each edge $e$ of $H$, $\pi(e)$ can pass though the boundary of $\Delta$ arbitrarily many (but still finitely many) times.
So there can be arbitrarily many connected components in the intersection of $\pi(E(H))$ and $\Delta$.
And it causes complications when we try to assemble them in order to construct a half-integral packing.

The key strategy is to find a finite set of ``small pieces'' such that for any homeomorphic embedding $\pi'$ from $H$ to $G$ and any disk $\Delta'$, the intersection of $\pi'(E(H))$ and $\Delta'$ ``simulates'' one of those small pieces, and such that we can find pairwise disjoint curves in $\Sigma$ internally disjoint from those $\kappa$ disks to connect those small pieces to obtain a homeomorphic embedding from $H$ to $G$.
In Section \ref{sec: vortices}, we will prove that the existence of this finite set implies the existence of a desired half-integral packing.
The members of this finite set are essentially the ``templates'' defined in Section \ref{subsec:addenda_template}.

Hence the difficulty is to find such a finite set of small pieces.
The objective of this section is to overcome this difficulty.

Note that even though there are arbitrarily many connected components of the intersection of $\pi(E(H))$ and $\Delta$, there are only at most $\lvert V(H) \rvert$ of them intersecting $\pi(V(H))$.
So there are only a bounded number of pieces coming from such a connected component, and we can collect all of them into our set of small pieces without losing the finiteness.
Hence it suffices to focus on the remaining connected components. 

That is, now we get many curves in $\overline{\Delta}$, where each of them has both endpoints in $\partial \Delta$, and the intersection of any two different curves is contained in the interior of $\Delta$.
We can treat the union of them as a 1-regular graph $D$ whose vertices are the endpoints of those curves and whose edges are given by the curves, and the vertices of this graph are ordered by the natural order of $\partial \Delta$, so we get a rooted graph whose march is given by this natural order.
Such a 1-regular graph $D$ is the $\Delta$-dive defined in Section \ref{subsec:dive_depth}; a $\Delta$-dive together with those finitely many connected components intersecting $\pi(V(H))$ that we ignored form the $\Delta$-miniature defined in Section \ref{subsec:dive_depth}.

Recall that our goal is to find a finite set of rooted graphs such that any such 1-regular rooted graph (dive) $D$ contains some member of this finite set (template) as a sub-rooted graph, and such that we can add pairwise disjoint curves in $\Sigma$ internally disjoint from $\overline{\Delta}$ to connect vertices in this sub-rooted graph to construct a homeomorphic embedding from $H$ to $G$.

For adding pairwise disjoint curves to connect vertices in the sub-rooted graph mentioned in the previous paragraph, it can be shown (see Sections \ref{subsec:linking_disentanglement} and \ref{subsec:addenda_template}) that it is essentially equivalent to finding at most $\rho$ curves $\gamma_1,\gamma_2,...,\gamma_\rho$ in $\Sigma$ such that the union of those curves contains $D$, and $\bigcup_{i=1}^\rho \gamma_i - D$ is a union of pairwise disjoint curves in $\Sigma$ internally disjoint from $\overline{\Delta}$.
Recall that $\rho$ is given by the depth of the vortices in $\Se_2$, and there are no $\rho+1$ pairwise intersecting curves in $D$.

We construct those curves in the following procedure.
(See Figure \ref{fig_disentanglement} for an illustration.)
Let $\Delta'$ be a disk containing $\Delta$ such that $\partial \Delta' \cap \partial \Delta = \emptyset$.
Let $p_1,p_2,...,p_\rho$ be distinct points in $\partial \Delta'$ appearing in the order listed if we trace $\partial \Delta'$ in the natural order consistent with the natural order of $\partial \Delta$.
For each $i \in [\rho]$, let $\gamma_i$ be the curve that consists of the single point $p_i$.
In the process, we will grow each $\gamma_i$ by adding a curve in $\Delta'$ internally disjoint from $\partial \Delta' \cup \bigcup_{j} \gamma_j$ to connect a vertex of $D$ and the endpoint $q_i$ of $\gamma_i$, where $q_i \neq p_i$ unless $\gamma_i$ consists of a single point.
Now we trace $\partial \Delta$ starting from a point in $\partial \Delta$ according to the natural order.
As long as we touch the first end $v$ of an edge of $D$, we first choose $i \in [\rho]$ such that either $p_i=q_i$ or $q_i$ is the second end of an edge of $D$, and subject to this, $q_i$ is as close to $v$ as possible, and then we grow $\gamma_i$ by adding a curve in $\Delta'$ to connect $q_i$ and $v$ internally disjoint from $\Delta \cup \partial \Delta' \cup \bigcup_{j}\gamma_j$; as long as we touch the second end $v$ of an edge of $D$, say the first end of this edge is in $\gamma_i$, we grow $\gamma_i$ by adding this edge of $D$ (note that every edge of $D$ is a curve).
This is the intuition of Phase 1 stated in Section \ref{subsubsec:def_disentanglement}.

This will construct desired curves due to the choice of $i$ in the procedure. 
However, it is difficult to analyze those curves to get the desired finite set of small pieces because those curves are ``entangled'' in the sense that any connected component of $\bigcup_j\gamma_j \cap \Delta$ can consist of arbitrarily many curves.
We want to simplify them by ``disentangling'' them by doing the following procedure.
During the above process, we have some ``transition periods''.
That is, each transition period consists of the points in $\partial\Delta$ visited during the time period between the time we just finished seeing a maximal sequence of consecutive second ends of edges, and the time that we were about to see a new first end of some edge.
And for any transition period, each curve $\gamma_i$ has a connected component in either $\gamma_i \cap \overline{\Delta}$ or $\gamma_i-\Delta$ such that the endpoints of this connected component define an interval in $\partial\Delta$ containing the entire transition period, and we say $\gamma_i$ is visiting $\Delta$ if this connected component is in $\overline{\Delta}$. 
For each transition period, we ``lift up'' each curve $\gamma_i$ visiting $\Delta$ by replacing the connected component $C$ of $\gamma_i \cap \overline{\Delta}$ showing that $\gamma_i$ is visiting $\Delta$ by the union of three curves: a curve in $\Delta$ from one endpoint of $C$ to a point $z$ in the transition period, a curve in $\Delta'-\Delta$ from $z$ to another point $z'$ in the transition period, and a curve in $\Delta$ from $z'$ to the other endpoint of $C$; and we do those ``lifting'' in a way that we do not create any crossing outside $\Delta$.
This is the intuition of the disentanglement defined in Section \ref{subsec:disentanglement} and essentially the Phase 2 stated in Section \ref{subsubsec:def_disentanglement}.
The disentanglement has the property that we can easily partition the intersection of the curves and $\Delta$ into parts (called ``bundles'' in Section \ref{subsec:disentanglement}), where each part only involves at most $\rho$ connected components, and the minimal interval in $\partial \Delta$ containing the endpoints of the curves in one part is disjoint from the minimal interval in $\partial \Delta$ containing the endpoints of the curves in any other part.
This makes the analysis simpler, and a well-quasi-ordering technique ensures the existence of the finite set of small pieces.

Finally, we will show how to resemble ``essential'' parts obtained by disentanglements to obtain a sub-rooted graph of the original 1-regular graph which is sufficient to construct a desired half-integral packing in Sections \ref{subsec:linking_disentanglement}.

\subsection{Dives, depth, and miniatures} \label{subsec:dive_depth}

Let $\kappa$ be a nonnegative integer, and $\Delta_1,\Delta_2,...,\Delta_\kappa$ be oriented open disks in $\Sigma$ with pairwise disjoint closure. 
For each $i \in [\kappa]$, let $p_i$ be the point in $\partial\Delta_i$ and $\Omega_i(p_i)$ the linear ordering witnessing that $\Delta_i$ is an oriented disk.
We say that a pseudo-embedding $\pi$ of a graph $H$ in $\Sigma$ is {\it legal with respect to $\{\Delta_1,...,\Delta_\kappa\}$} if the following conditions hold.
	\begin{itemize}
		\item $\pi(V(H)) \cap (\bigcup_{i=1}^\kappa \partial \Delta_i) = \emptyset$.
		\item $\pi(E(H)) \cap \{p_i: i \in [\kappa]\} = \emptyset$.
		\item For every $i \in [\kappa]$ and $x \in \partial \Delta_i \cap \pi(E(H))$, there exists an open set $B \subseteq \Sigma$ with $x \in B$ such that 
			\begin{itemize}
				\item $B \cap (\Delta_i \cup \{x\}) \cap \pi(E(H))$ has exactly one connected component, and this component contains a point other than $x$, and
				\item $B \cap \Delta_i \cap \pi(E(H))$ has exactly one connected component. 
			\end{itemize} 
		\item $\pi(E(H)) \cap \partial \Delta_i$ is finite.
	\end{itemize}
Note that the fourth bullet implies that for each $i \in [\kappa]$ and $x \in \partial \Delta_i \cap \pi(E(H))$, there exists an open set $B \subseteq \Sigma$ with $x \in B$ such that $B \cap \partial\Delta_i \cap \pi(E(H)) = \{x\}$.
For each $i \in [\kappa]$, let $\Gamma_i$ be the union of the connected components of $\overline{\Delta_i} \cap \pi(E(H))$ disjoint from $\overline{\Delta_i} \cap \pi(V(H))$; define $S_i$ to be the graph such that $V(S_i)=\Gamma_i \cap \partial \Delta_i$ and for each pair of distinct vertices $x,y$ of $S_i$, they are adjacent if and only if there exists $e \in E(H)$ such that some connected component of $\overline{\Delta_i} \cap \pi(e)$ links $x,y$ in $\Gamma_i$.
Note that $S_i$ is a simple 1-regular graph.
As vertices of $S_i$ are points in $\partial\Delta_i$, we may assume that $\Omega_i(p_i)$ is an ordering on $V(S_i)$.
We define the following.
	\begin{itemize}
		\item Define the {\it $\Delta_i$-dive} of $\pi$ to be the rooted graph $(S_i,\mu_i)$, where $\mu_i$ is the march with $\lvert V(S_i) \rvert$ entries such that for every $j \in [\lvert V(S_i) \rvert]$, the $j$-th entry is the $j$-th vertex of $V(S_i)$ in $\Omega_i(p_i)$.
		\item Define the {\it $\Delta_i$-depth} of $\pi$ to be the minimum $\rho$ such that there exists no partition $\{I',I''\}$ of $\Omega_i(p_i)$ into two cyclic intervals such that there exist $\rho+1$ disjoint paths in $S_i$ from $I'$ to $I''$.
		\item Define $R_i$ to be the graph such that $V(R_i)=\{v \in V(H): \pi(v) \in \Delta_i\} \cup (\partial \Delta_i \cap \pi(E(H))-V(S_i))$, and $E(R_i)=\{e \in E(H): \pi(e) \subseteq \Delta_i\} \cup \{xv: v \in V(H), \pi(v) \in \Delta_i, x \in \partial\Delta_i$ and there exists $e \in E(H)$ such that some connected component of $\overline{\Delta_i} \cap \pi(e)$ contains $\{x,\pi(v)\}\}$.
		\item Define the {\it $\Delta_i$-miniature} of $\pi$ in $\Sigma$ to be the rooted graph $(S_i \cup R_i, \nu_i)$, where $\nu_i$ is the march with $\lvert V(S_i) \rvert + \lvert \partial \Delta_i \cap \pi(E(H))-V(S_i) \rvert = \lvert \partial \Delta_i \cap \pi(E(H)) \rvert$ entries such that for each $j \in [\lvert \partial \Delta_i \cap \pi(E(H)) \rvert]$, the $j$-th entry in $\nu_i$ is the $j$-th vertex in $V(S_i) \cup (\partial \Delta_i \cap \pi(E(H))-V(S_i)) = \partial \Delta_i \cap \pi(E(H))$ in $\Omega_i(p_i)$.
	\end{itemize}
Note that $\lvert V(R_i) \rvert \leq \lvert V(H) \rvert+2\lvert E(H) \rvert$, $\lvert E(R_i) \rvert \leq 2\lvert E(H) \rvert$, and $V(R_i) \cap \partial\Delta_i$ naturally partitions $\partial\Delta_i$ into $\lvert V(R_i) \cap \partial\Delta_i \rvert+1$ (non-cyclic) intervals.
Moreover, for each edge of $S_i$, we call an end of this edge the {\it first end} if it appears earlier than the other end; otherwise, we call this end the {\it second end}.

\subsection{Disentanglements, foundations, highlights and signatures} \label{subsec:disentanglement}

In this subsection, we assume the following.
	\begin{itemize}
		\item $\kappa$ is a nonnegative integer, and $\Delta_1,\Delta_2,...,\Delta_\kappa$ are oriented open disks in $\Sigma$ with pairwise disjoint closure.
		\item For each $i \in [\kappa]$, $p_i$ is the point in $\partial\Delta_i$ such that $\Omega_i(p_i)$ is the linear ordering witnessing that $\Delta_i$ is an oriented disk.
		\item $\pi$ is a pseudo-embedding of a graph $H$ in $\Sigma$ legal with respect to $\{\Delta_1,...,\Delta_\kappa\}$.
		\item For each $i \in [\kappa]$, $(S_i,\mu_i)$ is the $\Delta_i$-dive of $\pi$, and $R_i$ is the graph disjoint from $S_i$ such that $(S_i \cup R_i,\nu_i)$ is the $\Delta_i$-miniature of $\pi$ in $\Sigma$.
	\end{itemize}
When there is no danger for creating confusion, we do not distinguish an edge $e$ of $S_i$ and a curve contained in $\overline{\Delta_i}$ internally disjoint from $\partial\Delta_i$ whose endpoints are the ends of $e$.

\subsubsection{Definition of disentanglements} \label{subsubsec:def_disentanglement}

Let $i$ be a fixed element in $[\kappa]$, and let $\rho$ be the $\Delta_i$-depth.
The {\it $\Delta_i$-disentanglement} of $\pi$ is a rooted graph $(D_i,\sigma_i)$ together with a labelling function defined on the vertices, where $D_i$ is a 1-regular simple graph whose every vertex is a point in $\partial \Delta_i$ and labelled with a number in $[\rho]$ and a word which is either ``real'' or ``fake'' constructed by the following procedure, and $\sigma_i$ is the march that orders all vertices of $D_i$ consistent with $\Omega_i(p_i)$.

The procedure will be divided into 2 phases.
At the beginning of the procedure, $D_i$ is empty.
Let $\Delta_i'$ be an open disk in $\Sigma-\bigcup_{j \in [\kappa]-\{i\}}\overline{\Delta_j}$ containing $\Delta_i$ with $\partial\Delta_i \cap \partial\Delta_i' = \emptyset$.

\medskip

\noindent{\bf Phase 1: Adding ``real'' vertices and defining curves $\zeta_j$'s.}
During this phase, we will add some vertices into $D_i$, label them with ``real'' and some numbers, and define some disjoint curves $\zeta_1.\zeta_2,\cdots$.
At the beginning of the phase, none of the curves $\zeta_j$'s is defined.
For each $j$, once $\zeta_j$ is defined, it might be repeatedly updated by concatenating with other curves during this phase, and the following properties are maintained during the entire phase:
	\begin{itemize}
		\item $\zeta_j$ has one endpoint in $\Delta_i'-\overline{\Delta_i}$ and one endpoint in $\partial\Delta_i$.
		\item Every point in $\zeta_j \cap \partial\Delta_i$ is a vertex of $S_i$ and is labelled with number $j$ when it is added into $D_i$.
		\item Every component of $\zeta_j \cap \overline{\Delta_i}$ is an edge of $S_i$.
		\item If $\zeta_j$ contains a vertex $v$ of $S_i$, then $\zeta_j$ contains the edge of $S_i$ incident with $v$.
		\item If we trace $\zeta_j$ from the endpoint in $\Delta_i'-\overline{\Delta_i}$, then the first point in $\zeta_j \cap \partial\Delta_i$ is the first end of some edge of $S_i$, and every component of $\zeta_j-\Delta_i$ is a curve from the second end $v$ of an edge of $S_i$ to the first end of another edge of $S_i$ appearing later than $v$ in $\Omega_i(p_i)$.
	\end{itemize}
In addition, during this phase, if $\zeta_x$ is undefined for some positive integer $x$, then $\zeta_y$ is undefined for all integers $y \geq x$.

Now we start this phase.
We trace $\partial \Delta$ by starting from $p_i$ according to $\Omega_i(p_i)$, and whenever we encounter a vertex $v$ in $S_i$, we do the following operations: (see Figure \ref{fig_disentanglement} for an example) 
	\begin{itemize}
		\item If $v$ is the first end of an edge of $S_i$, then do the following operations. 
			\begin{itemize}
				\item If there exists no positive integer $x$ such that $\zeta_x$ is defined and the endpoint of $\zeta_x$ in $\partial\Delta_i$ is the second end of an edge of $S_i$, then 
					\begin{itemize}
						\item let $x$ be the smallest positive integer such that $\zeta_x$ is undefined, and 
						\item define $\zeta_x$ to be a curve in $\Delta_i'-\Delta_i$ from a point in $\Delta_i'-(\overline{\Delta_i} \cup \bigcup_{j=1}^{x-1}\zeta_j)$ to $v$ internally disjoint from $\overline{\Delta_i}$;
					\end{itemize}
				\item otherwise, 
					\begin{itemize}
						\item let $x$ be the integer such that $\zeta_x$ is defined, the endpoint of $\zeta_x$ in $\partial\Delta_i$ is the second end of an edge of $S_i$, and subject to these, the endpoint of $\zeta_x$ in $\partial\Delta_i$ is as close to $v$ in $\Omega_i(p_i)$ as possible, and
						\item concatenate $\zeta_x$ with a curve $\zeta$ in $\Delta_i'$ from the endpoint of $\zeta_x$ in $\partial\Delta_i$ to $v$ internally disjoint from $\overline{\Delta_i}$ and disjoint from other defined curves $\zeta_j$'s.

							(Note that the existence of $\zeta$ follows from the choice of $x$ and the properties of $\zeta_j$'s maintained during this phase.)
					\end{itemize}
				\item Add $v$ into $V(D_i)$ and label $v$ with $x$ and ``real''. 
			\end{itemize}
		\item If $v$ is the second end of an edge of $S_i$, then do the following operations. 
			\begin{itemize}
				\item Let $x$ be the number labelled on the neighbor of $v$ in $S_i$.
				\item Add $v$ into $V(D_i)$, and label $v$ with $x$ and ``real''.
				\item Add the edge of $S_i$ incident with $v$ into $E(D_i)$.
				\item Concatenate $\zeta_x$ with the edge of $S_i$ incident with $v$. 
			\end{itemize}
	\end{itemize}
This is the end of Phase 1.
It is clear that the procedure maintains the properties of those $\zeta_j$'s mentioned above.
The following lemma shows that the numbers used in the labels of the vertices of $D_i$ are in $[\rho]$.

\begin{lemma} \label{num_zeta}
$\zeta_{j}$ is undefined for every $j \geq \rho+1$. 
\end{lemma}

\begin{pf}
By the properties maintained in Phase 1, it suffices to show that $\zeta_{\rho+1}$ is undefined.
Suppose to the contrary that $\zeta_{\rho+1}$ is defined.
Then when $\zeta_{\rho+1}$ is about to be defined (say, when we visit vertex $v$), for every $j \in [\rho]$, $\zeta_j$ is defined and the endpoint of $\zeta_j$ in $\partial\Delta_i$ is the first end of some edge $e_j$ of $S_i$ such that the first end of $e_j$ appears before $v$ in $\Omega_i(p_i)$ and the other end of $e_j$ appears after $v$ in $\Omega_i(p_i)$.
In addition, $v$ is the first end of some edge of $S_i$.
Hence there exist $\rho+1$ disjoint paths in $S_i$ between $I'$ and $I''$, where $I'$ is the interval of $\Omega_i(p_i)$ from the first point in $\Omega_i(p_i)$ to $v$, and $I''$ is the interval of $\Omega_i(p_i)$ such that $\{I',I''\}$ is a partition of $\Omega_i(p_i)$.
So the $\Delta_i$-depth is greater than $\rho$, a contradiction.
\end{pf}

\bigskip

\noindent{\bf Phase 2: Adding ``fake'' vertices and defining curves $\gamma_j$'s.}
During this phase, we will add vertices into $D_i$, label them ``fake'' and some numbers, and define curves $\gamma_j$'s by modifying the curves $\zeta_j$'s.
At beginning of this phase, for each $j$ for which $\zeta_j$ is defined, define $\gamma_j$ to be $\zeta_j$.

Before we start this phase of the procedure, we define some terminology.
We say that a vertex $v$ of $S_i$ is a {\it closer} if it is not the last entry in $\mu_i$, and $v$ is the second end of some edge of $S_i$, and the vertex of $S_i$ appeared right after $v$ in $\Omega_i(p_i)$ is the first end of an edge of $S_i$.

Now we start Phase 2.
We trace $\partial\Delta_i$ by starting from $p_i$ according to $\Omega_i(p_i)$, and whenever we encounter a closer $v$, we do the following operations: (see Figure \ref{fig_disentanglement} for an example) 
	\begin{itemize}
		\item Let $k$ be the number of indices $j$ such that $\gamma_j \cap \overline{\Delta_i}$ has a connected component whose one endpoint (say $u_j$) is strictly before $v$ and one endpoint (say $v_j$) is strictly after $v$ in $\Omega_i(p_i)$.
			And let $j_1,j_2,...,j_k$ be those indices such that $u_{j_1},u_{j_2},...,u_{j_k}$ appear in $\Omega_i(p_i)$ in the order listed.
		\item Pick $2k$ points $z_1,z_2,...,z_k,z_k',z_{k-1}',...,z_1'$ in $\partial\Delta_i$ strictly between $v$ and the vertex of $S_i$ right after $v$ in $\Omega_i(p_i)$, where $z_1,z_2,...,z_k,z_k',z_{k-1}',...,z_1'$ appear in $\Omega_i(p_i)$ in the order listed.
		\item For each $\ell \in [k]$, modify $\gamma_{j_\ell}$ by 
			\begin{itemize}
				\item first deleting the connected component of $\gamma_{j_\ell} \cap \overline{\Delta_i}$ from $u_{j_\ell}$ to $v_{j_\ell}$,
				\item then concatenating a curve in $\overline{\Delta_i}$ from $u_{j_\ell}$ to $z_{k-\ell+1}$ internally disjoint from $\partial\Delta_i$,
				\item then concatenating a curve in $\Delta_i'-(\Delta_i \cup \bigcup_{\alpha \neq j_\ell}\gamma_\alpha)$ from $z_{k-\ell+1}$ to $z'_{k-\ell+1}$ internally disjoint from $\overline{\Delta_i}$, and
				\item then concatenating a curve in $\overline{\Delta_i}$ from $z'_{k-\ell+1}$ to $v_{j_\ell}$ internally disjoint from $\partial\Delta_i$.
			\end{itemize}
			Note that we can choose those curves in the concatenation process such that the sub-curves of the modified $\gamma_{j_\ell}$'s from $u_{j_\ell}$ to $z'_{k-\ell+1}$ are pairwise disjoint.
		\item For each $\ell \in [k]$, add $z_\ell$ and $z'_\ell$ into $V(D_i)$ and label both of them with ``fake'' and $j_{k-\ell+1}$, and add $u_{j_{k-\ell+1}}z_\ell$ and $z'_\ell v_{j_{k-\ell+1}}$ into $E(D_i)$.

			(Note that we can treat each of $u_{j_{k-\ell+1}}z_\ell$ and $z'_\ell v_{j_{k-\ell+1}}$ a sub-curve of $\gamma_{j_{k-\ell+1}}$ contained in $\overline{\Delta_i}$.)
		\item If $k>0$, then call $z_k$ a {\it fake closer} and $z'_k$ a {\it fake opener}; if $k=0$, then call $v$ a {\it fake closer} and the vertex of $S_i$ right after $v$ in $\Omega_i(p_i)$ a {\it fake opener}.
	\end{itemize}
This is the end of Phase 2 and completes the definition of $D_i$.

\medskip

Define $\sigma_i$ to be the march consisting of the vertices of $D_i$ ordered according to $\Omega_i(p_i)$, starting at $p_i$.
This completes the definition of the $\Delta_i$-disentanglement $(D_i,\sigma_i)$ (with a labelling function). 
See Figure \ref{fig_disentanglement} for an example.

We also call the first entry of $\sigma_i$ a fake opener and the last entry of $\sigma_i$ a fake closer.
A {\it life period} is a set of consecutive entries of $\sigma_i$ starting from a fake opener, ending at a fake closer, and containing exactly one fake opener and one fake closer.

A {\it bundle} in the $\Delta_i$-disentanglement is a labelled rooted graph whose underlying graph is a subgraph of $D_i$ induced by a life period, and whose march contains all its vertices and is a subsequence of $\sigma_i$, and the label on each vertex is the same as its label in $D_i$.
See Figure \ref{fig_disentanglement} for an example.
Note that the collection of the vertex-sets of all bundles is a partition of $V(D_i)$.
The {\it type} of a bundle $(B,\sigma_B)$ is the sequence with $\lvert V(B) \rvert$ entries such that for each $j \in [\lvert V(B) \rvert]$, the $j$-th entry is the label of the $j$-th entry of $\sigma_B$.

\begin{figure} 
	\begin{picture}(100,360) (-35,-180)
		\thinlines
		\multiput(5,200)(5,0){79}{\line(1,0){3}}
		\multiput(5,197)(0,-5){12}{\line(0,1){3}}
		\multiput(400,197)(0,-5){12}{\line(0,1){3}}
		\multiput(5,141)(5,0){79}{\line(1,0){3}}
		\put(10,205){{\scriptsize$p_i$}}
		\put(-17,160){{$\partial\Delta_i$}}
		\put(190,150){{$\Delta_i$}}

		\thicklines
		\put(30,200){\circle*{5}} \put(28,205){$v_1$}
		\put(50,200){\circle*{5}} \put(48,205){$v_2$}
		\put(70,200){\circle*{5}} \put(68,205){$v_3$}
		\put(90,200){\circle*{5}} \put(88,205){$v_4$}
		\qbezier(50,200)(70,140)(90,200)
		\put(190,200){\circle*{5}} \put(188,205){$v_5$}
		\put(210,200){\circle*{5}} \put(208,205){$v_6$}
		\put(230,200){\circle*{5}} \put(228,205){$v_7$}
		\qbezier(30,200)(130,140)(230,200)
		\put(250,200){\circle*{5}} \put(248,205){$v_8$}
		\qbezier(70,200)(160,140)(250,200)
		\put(270,200){\circle*{5}} \put(269,205){$v_9$}
		\qbezier(210,200)(240,140)(270,200)
		\put(330,200){\circle*{5}} \put(328,205){$v_{10}$}
		\put(350,200){\circle*{5}} \put(348,205){$v_{11}$}
		\qbezier(190,200)(270,140)(350,200)
		\put(370,200){\circle*{5}} \put(368,205){$v_{12}$}
		\qbezier(330,200)(350,140)(370,200)

		\thinlines
		\multiput(5,80)(5,0){79}{\line(1,0){3}}
		\multiput(5,77)(0,-5){15}{\line(0,1){3}}
		\multiput(400,77)(0,-5){15}{\line(0,1){3}}
		\multiput(5,6)(5,0){79}{\line(1,0){3}}
		\put(10,85){{\scriptsize$p_i$}}
		\put(-17,40){{$\partial\Delta_i$}}
		\put(190,15){{$\Delta_i$}}

		\thicklines
		\put(30,80){\circle*{5}} \put(28,85){$v_1$}\put(28,93){R}\put(28,103){$1$}
		\put(50,80){\circle*{5}} \put(48,85){$v_2$}\put(48,93){R}\put(48,103){$2$}
		\put(70,80){\circle*{5}} \put(68,85){$v_3$}\put(68,93){R}\put(68,103){$3$}
		\put(90,80){\circle*{5}} \put(88,85){$v_4$}\put(88,93){R}\put(88,103){$2$}
		\qbezier(50,80)(70,20)(90,80)
		\put(110,80){\circle{5}} \put(108,93){F}\put(108,103){$3$}
		\put(130,80){\circle{5}} \put(128,93){F}\put(128,103){$1$}
		\put(150,80){\circle{5}} \put(148,93){F}\put(148,103){$1$}
		\put(170,80){\circle{5}} \put(168,93){F}\put(168,103){$3$}
		\qbezier(70,80)(85,50)(110,80)
		\qbezier(30,80)(71,40)(130,80)
		\put(190,80){\circle*{5}} \put(188,85){$v_5$}\put(188,93){R}\put(188,103){$2$}
		\put(210,80){\circle*{5}} \put(208,85){$v_6$}\put(208,93){R}\put(208,103){$4$}
		\put(230,80){\circle*{5}} \put(228,85){$v_7$}\put(228,93){R}\put(228,103){$1$}
		\qbezier(150,80)(190,40)(230,80)
		\put(250,80){\circle*{5}} \put(248,85){$v_8$}\put(248,93){R}\put(248,103){$3$}
		\qbezier(170,80)(197,40)(250,80)
		\put(270,80){\circle*{5}} \put(268,85){$v_9$}\put(268,93){R}\put(268,103){$4$}
		\qbezier(210,80)(240,20)(270,80)
		\put(290,80){\circle{5}} \put(288,93){F}\put(288,103){$2$}
		\put(310,80){\circle{5}} \put(308,93){F}\put(308,103){$2$}
		\qbezier(190,80)(240,60)(290,80)
		\put(330,80){\circle*{5}} \put(328,85){$v_{10}$}\put(328,93){R}\put(328,103){$4$}
		\put(350,80){\circle*{5}} \put(348,85){$v_{11}$}\put(348,93){R}\put(348,103){$2$}
		\qbezier(310,80)(330,40)(350,80)
		\put(370,80){\circle*{5}} \put(368,85){$v_{12}$}\put(368,93){R}\put(368,103){$4$}
		\qbezier(330,80)(350,20)(370,80)

		\thinlines
		\put(23,30){\line(1,0){115}}
		\put(23,28){\line(0,1){4}}
		\put(138,28){\line(0,1){4}}
		\put(60,35){{\small{bundle 1}}}

		\put(143,30){\line(1,0){155}}
		\put(143,28){\line(0,1){4}}
		\put(298,28){\line(0,1){4}}
		\put(195,35){{\small{bundle 2}}}

		\put(303,30){\line(1,0){75}}
		\put(303,28){\line(0,1){4}}
		\put(378,28){\line(0,1){4}}
		\put(320,35){{\small{bundle 3}}}

		\thinlines
		\multiput(5,-90)(5,0){79}{\line(1,0){3}}
		\multiput(5,-93)(0,-5){15}{\line(0,1){3}}
		\multiput(400,-93)(0,-5){15}{\line(0,1){3}}
		\multiput(5,-164)(5,0){79}{\line(1,0){3}}
		\put(10,-85){{\scriptsize$p_i$}}
		\put(-17,-130){{$\partial\Delta_i$}}
		\put(190,-155){{$\Delta_i$}}

		\multiput(-20,-40)(5,0){86}{\line(1,0){3}}
		\multiput(-20,-43)(0,-5){26}{\line(0,1){3}}
		\multiput(410,-43)(0,-5){26}{\line(0,1){3}}
		\multiput(-20,-170)(5,0){86}{\line(1,0){3}}	
		\put(370,-60){{$\Delta'_i$}}

		\thicklines
		\put(30,-90){\circle*{5}} 
		\put(50,-90){\circle*{5}} 
		\put(70,-90){\circle*{5}} 
		\put(90,-90){\circle*{5}} 
		\qbezier(50,-90)(70,-150)(90,-90)
		\put(110,-90){\circle{5}} 
		\put(130,-90){\circle{5}} 
		\put(150,-90){\circle{5}} 
		\put(170,-90){\circle{5}} 
		\qbezier(70,-90)(85,-120)(110,-90)
		\qbezier(30,-90)(71,-130)(130,-90)
		\put(190,-90){\circle*{5}} 
		\put(210,-90){\circle*{5}} 
		\put(230,-90){\circle*{5}} 
		\qbezier(150,-90)(190,-130)(230,-90)
		\put(250,-90){\circle*{5}} 
		\qbezier(170,-90)(197,-130)(250,-90)
		\put(270,-90){\circle*{5}} 
		\qbezier(210,-90)(240,-150)(270,-90)
		\put(290,-90){\circle{5}} 
		\put(310,-90){\circle{5}} 
		\qbezier(190,-90)(240,-110)(290,-90)
		\put(330,-90){\circle*{5}} 
		\put(350,-90){\circle*{5}} 
		\qbezier(310,-90)(330,-130)(350,-90)
		\put(370,-90){\circle*{5}} 
		\qbezier(330,-90)(350,-150)(370,-90)

		\qbezier(25,-80)(27,-85)(30,-90) \put(22,-75){$\gamma_1$}
		\qbezier(45,-80)(47,-85)(50,-90) \put(42,-75){$\gamma_2$}
		\qbezier(65,-80)(67,-85)(70,-90) \put(62,-75){$\gamma_3$}
		\qbezier(90,-90)(140,-10)(190,-90)
		\qbezier(110,-90)(140,-30)(170,-90)
		\qbezier(130,-90)(140,-50)(150,-90)
		\qbezier(205,-80)(207,-85)(210,-90) \put(202,-75){$\gamma_4$}
		\qbezier(290,-90)(300,-50)(310,-90)
		\qbezier(270,-90)(300,-30)(330,-90)

		\thinlines
		\put(23,-140){\line(1,0){115}}
		\put(23,-142){\line(0,1){4}}
		\put(138,-142){\line(0,1){4}}
		\put(60,-135){{\small{bundle 1}}}

		\put(143,-140){\line(1,0){155}}
		\put(143,-142){\line(0,1){4}}
		\put(298,-142){\line(0,1){4}}
		\put(195,-135){{\small{bundle 2}}}

		\put(303,-140){\line(1,0){75}}
		\put(303,-142){\line(0,1){4}}
		\put(378,-142){\line(0,1){4}}
		\put(320,-135){{\small{bundle 3}}}
	\end{picture}
	\caption{An example of the $\Delta_i$-dive $(S_i,\mu_i)$, the $\Delta_i$-disentanglement $(D_i,\sigma_i)$, and a $\Delta_i$-foundation. (1) The top side of the picture includes the graph $S_i$ and the disk $\Delta_i$, where $V(S_i)=\{v_1,v_2,...,v_{12}\}$. (2) The middle part of the picture includes the graph $D_i$ at the end of the procedure, together with the labels on the vertices of $D_i$, and the disk $\Delta_i$. The vertices in $V(D_i)-V(S_i)$ are indicated by empty circles. Each vertex in $D_i$ is labelled with a number and a word which is either ``real'' or ``fake''. In the picture, ``R'' denotes ``real'', and ``F'' denotes ``fake''. (3) The bottom side of the picture includes the graph $D_i$, a $\Delta_i$-foundation $\{\gamma_1,\gamma_2,\gamma_3,\gamma_4\}$, and disks $\Delta_i$ and $\Delta_i'$.} \label{fig_disentanglement}
\end{figure}

\begin{lemma} \label{disentanglement_basic_3}
Every bundle has at most $\rho$ edges and contains at least one vertex labelled with ``real''.
\end{lemma}

\begin{pf}
Note that every edge of a bundle is contained in $\gamma_j$ for some $j$, and no $\gamma_j$ contains at least two edges of a bundle.
So every bundle contains at most $\rho$ edges by Lemma \ref{num_zeta}.
Also, if there exists no closer, then Phase 2 is void, so $V(D_i)=V(S_i)$ and hence every bundle contains a vertex labelled with ``real''; otherwise, every bundle contains a closer (i.e. the vertex $v$ in Phase 2), and this closer is a vertex of $S_i$, so it is labelled with ``real''.
\end{pf}

\bigskip

Let $b$ be the largest integer such that $\gamma_b$ is defined.
So $\{\gamma_j: j \in [b]\}$ is a set of curves satisfying the following:
	\begin{itemize}
		\item For distinct $j,j' \in [b]$, $\gamma_j \cap \gamma_{j'} \subseteq \Delta_i$.
		\item For each $j \in [b]$, every connected component of $\gamma_j \cap \overline{\Delta_i}$ is an edge of $D_i$ whose both ends are labelled with $j$.
		\item For each $j \in [b]$, every vertex of $D_i$ contained in $\gamma_j$ is labelled with $j$.
		\item For each $j \in [b]$, $\gamma_j$ passes through all vertices of $V(D_i)$ with label $j$ in the order as the appearance of those vertices in $\Omega_i(p_i)$.
	\end{itemize}
A {\it $\Delta_i$-foundation} is a set consisting of $b$ curves satisfying the above conditions.
See Figure \ref{fig_disentanglement} for an example.

Note that $\Omega_i(p_i)$ is a linear ordering on $V(D_i) \cup (V(R_i) \cap \partial\Delta_i)$.
We define the {\it $\Delta_i$-highlight} of $\pi$ to be the rooted graph whose underlying graph is $D_i \cup R_i$ and whose march is formed by the vertices in $V(D_i) \cup (V(R_i) \cap \partial\Delta_i)$ ordered by $\Omega_i(p_i)$.

Recall that $V(R_i) \cap \partial \Delta_i$ defines $\lvert V(R_i) \cap \partial \Delta_i \rvert+1$ intervals of $\Omega_i(p_i)$.
Define the {\it $\Delta_i$-signature} of $\pi$ to be the sequence $(Q,s)$, where
	\begin{itemize}
		\item $Q$ is the rooted graph with some vertices labelled such that 
			\begin{itemize}
				\item the underlying graph of $Q$ is the union of $R_i$ and the underlying graphs of the bundles intersecting at least two intervals of $\Omega_i(p_i)$ determined by $V(R_i) \cap \partial\Delta_i$, 
				\item the march of $Q$ is formed by the vertices in $V(Q) \cap \partial \Delta_i$ ordered by $\Omega_i(p_i)$, and
				\item the vertices of $Q$ that are labelled are the vertices that are vertices of some bundles and the vertices in $V(R_i)-\partial \Delta_i$, where
					\begin{itemize}
						\item the labels of the vertices of bundles are the same as their labels in the bundles, and
						\item each vertex in $V(R_i)-\partial\Delta_i$ is a vertex $v$ of $H$ and its label is $v$.
					\end{itemize}
			\end{itemize}
		\item $s=(s_1,s_2,...,s_{\lvert V(R_i) \cap \partial \Delta_i \rvert+1})$, and for each $j \in [\lvert V(R_i) \cap \partial \Delta_i \rvert+1]$, $s_j$ is the finite sequence such that  
			\begin{itemize}
				\item each entry of $s_j$ is a bundle with vertex-set contained in the $j$-th interval determined by $V(R_i) \cap \partial \Delta_i$, where the first interval is the one that contains $p_i$, and 
				\item for each $r$, the $r$-th entry of $s_j$ is the $r$-th bundle in that interval (ordered by $\Omega_i(p_i)$).
			\end{itemize}
	\end{itemize}
Note that the number intervals of $\Omega_i(p_i)$ determined by $V(R_i) \cap \partial\Delta_i$ is upper bounded by a function of $\lvert V(H) \rvert + \lvert E(H) \rvert$, and for any two distinct intervals of $\Omega_i(p_i)$ determined by $V(R_i) \cap \partial\Delta_i$, at most one bundle can intersect both of them.
Recall that each bundle contains at most $\rho$ edges by Lemma \ref{disentanglement_basic_3}.
So the number of non-isomorphic $Q$'s is bounded by a function of $\lvert V(H) \rvert+\lvert E(H) \rvert$ and the $\Delta_i$-depth of $\pi$.
And there exists a finite set only depending on $\lvert V(H) \rvert+\lvert E(H) \rvert$ and the $\Delta_i$-depth such that each entry of any $s_j$ is isomorphic to a rooted graph in this set.
Note that for each bundle, the numbers in its type determine the graph structure and the march of this bundle.
Moreover, $s$ records the information of the types of the bundles, so the $\Delta_i$-depth can be told by simply seeing the $\Delta_i$-signature.

We define a binary relation $\preceq'$ on the set of all $\Delta_i$-signatures as follows: for $\Delta_i$-signatures $(Q,s)$ and $(Q',s')$, where $s=(s_1,s_2,...,s_{\lvert s \rvert})$ and $s'=(s'_1,s'_2,...,s'_{\lvert s' \rvert})$, define $(Q,s) \preceq' (Q',s')$ if and only if 
	\begin{itemize}
		\item $Q$ and $Q'$ are isomorphic (so $\lvert s \rvert = \lvert s' \rvert$), 
		\item for each $j \in [\lvert s \rvert]$, there exist positive integers $t^j_1<t^j_2<...<t^j_{\lvert s_j \rvert}$ such that for each $r \in [\lvert s_j \rvert]$, the type of the $r$-th entry in $s_j$ equals the type of the $t^j_r$-th entry in $s'_j$, and
		\item the $\Delta_i$-depth corresponding to $(Q,s)$ equals the $\Delta_i$-depth corresponding to $(Q',s')$.
	\end{itemize}
Clearly, $\preceq'$ is a quasi-ordering.
(Recall that a {\it quasi-ordering} is a reflexive and transitive binary relation.)

A quasi-ordering $\sqsubseteq$ on a set $X$ is a {\it well-quasi-ordering} if for every infinite sequence $a_1,a_2,...$ in $X$, there exist $1 \leq i <j$ such that $a_i \sqsubseteq a_j$.
We will need the following Higman's Lemma for well-quasi-ordering.

\begin{theorem}[\cite{h}] \label{Higman lemma}
If $\sqsubseteq$ is a well-quasi-ordering on a set $X$, then $\sqsubseteq'$ is a well-quasi-ordering on the set of finite sequences on $X$, where $\sqsubseteq'$ is the relation such that two finite sequences $a=(a_1,a_2,...,a_{\lvert a \rvert})$ and $b=(b_1,b_2,...,b_{\lvert b \rvert})$ over $X$ satisfy $a \sqsubseteq' b$ if and only if there exist $i_1<i_2<...<i_{\lvert a \rvert}$ such that $a_j \sqsubseteq b_{i_j}$ for every $j \in \lvert a \rvert$.
\end{theorem}

\begin{lemma} \label{wqo Delta-signature}
Let $\Sigma$ be a surface, and let $W$ be a set of oriented open disks in $\Sigma$ with pairwise disjoint closure.
Let $H$ be a connected graph and $\Delta \in W$.
For every positive integer $\rho$, the set of $\Delta$-signatures of pseudo-embeddings of $H$ in $\Sigma$ legal with respect to $W$ with $\Delta$-depth at most $\rho$ is well-quasi-ordered by $\preceq'$.
\end{lemma}

\begin{pf}
Since the $\Delta$-depth is bounded and $H$ is a fixed graph, there exist finite sets $A, B$ such that for any $\Delta$-signature $(Q,(s_1,s_2,...,s_k))$ of a pseudo-embedding of $H$ in $\Sigma$ with $\Delta$-depth at most $\rho$,  $Q$ is chosen from $A$, and each entry of each $s_j$ is chosen from $B$.
Hence $\preceq'$ is a well-quasi-ordering by Theorem \ref{Higman lemma}.
\end{pf}

\subsection{Linking disentanglements} \label{subsec:linking_disentanglement}

Let $\rho$ be a positive integer.
Let $\Sigma$ be a surface and $\Delta$ an oriented open disk in $\Sigma$.
Let $H$ be a connected graph, and let $\pi$ be a pseudo-embedding of $H$ in $\Sigma$ legal with respect to $\{\Delta\}$.
Let $\P=\{\gamma_1,...,\gamma_\rho\}$ be a $\Delta$-foundation.
Let $I$ be a subset of $[\rho]$.
We call the induced subgraph of a bundle an {\it $I$-subbundle} if it is the intersection of this bundle and $\bigcup_{i \in I}\gamma_i$.

\begin{lemma} \label{bundle_index_frist}
Let $\Sigma$ be a surface and $\Delta$ an oriented open disk in $\Sigma$.
Let $H$ be a connected graph, and let $\pi$ be a pseudo-embedding of $H$ in $\Sigma$ legal with respect to $\{\Delta\}$.
Let $\Delta'$ be an open disk in $\Sigma$ with $\Delta' \supset \Delta$ and $\partial\Delta \cap \partial \Delta'=\emptyset$.
Let $\P=\{\gamma_1,...,\gamma_\rho\}$ be a $\Delta$-foundation with $\bigcup_{i=1}^\rho\gamma_i \subseteq \Delta'$.
Let $I$ be a subset of $[\rho]$.
If $B_1,B_2$ are two distinct $I$-subbundles appearing in $\Delta$ in the order listed such that no $I$-subbundle appearing in $\Delta$ strictly between $B_1$ and $B_2$ has more than $\lvert E(B_1) \rvert$ edges, then there exists $I' \subseteq I$ with $\lvert I' \rvert = \lvert E(B_1) \rvert$ such that 
	\begin{enumerate}
		\item $B_1$ is an $I'$-subbundle,
		\item every $I$-subbundle strictly between $B_1$ and $B_2$ is an $I'$-subbundle, 
		\item if $\lvert E(B_1) \rvert \geq \lvert E(B_2) \rvert$, then $B_2$ is an $I'$-subbundle, and
		\item if $\lvert E(B_1) \rvert < \lvert E(B_2) \rvert$, then the $I'$-subbundle of $B_2$ has $\lvert I' \rvert$ edges.
	\end{enumerate}
\end{lemma}

\begin{pf}
Let $D_1,D_2$ be the bundles such that $B_1,B_2$ are the $I$-subbundles of $D_1,D_2$, respectively.
Let $I'$ be the minimal subset of $[\rho]$ such that $B_1$ is an $I'$-subbundle of $D_1$.
Note that $I' \subseteq I$ and $\lvert I' \rvert = \lvert E(B_1) \rvert$.
So Statement 1 holds.

Note that by the process of defining the $\Delta$-disentanglement, for every $j \in I-I'$, no connected component of $\gamma_j \cap \overline{\Delta}$ with one endpoint before $D_1$ and one endpoint after $D_1$. 

Suppose that there exists a bundle $D$ that appears strictly between $D_1$ and $D_2$ or equals $D_2$ such that the $I$-subbundle of $D$ is not an $I'$-subbundle.
We further assume that $D$ is as close to $D_1$ as possible.
So there exists $j \in I-I'$ such that $D$ contains a connected component of $\gamma_j \cap \overline{\Delta}$, and for every $I$-bundle $D'$ that is strictly between $D_1$ and $D$ or equal to $D_1$, $D'$ does not contain a connected component of $\gamma_j \cap \overline{\Delta}$.
By the process of defining the $\Delta$-disentanglement, $D$ must contain a connected component of $\gamma_{j'} \cap \overline{\Delta}$ for every $j' \in I'$.
Hence the $I$-subbundle of $D$ has more edges than $B_1$.
This shows that Statements 2 and 3 hold, and $\lvert E(B_2) \rvert > \lvert E(B_1) \rvert$.
It also implies that the $I'$-subbundle of $B_2$ contains $\lvert I' \rvert$ edges, so Statement 4 holds.
\end{pf}

\begin{lemma} \label{subbundle}
For every positive integer $\rho$, there exists an integer $c=c(\rho)$ such that the following hold. 
Let $\Sigma$ be a surface and $\Delta$ an oriented open disk in $\Sigma$.
Let $H$ be a connected graph, and let $\pi$ be a pseudo-embedding of $H$ in $\Sigma$ legal with respect to $\{\Delta\}$.
Let $\Delta'$ be an open disk in $\Sigma$ with $\Delta' \supset \Delta$ and $\partial\Delta \cap \partial \Delta'=\emptyset$.
Let $\P=\{\gamma_1,...,\gamma_\rho\}$ be a $\Delta$-foundation with $\bigcup_{i=1}^\rho\gamma_i \subseteq \Delta'$.
Let $I$ be a subset of $[\rho]$.
Let $D$ be the subgraph of the $\Delta$-highlight consisting of all $I$-subbundles.
If $B_1,B_2$ are two distinct $I$-subbundles with $\lvert E(B_1) \rvert = \lvert E(B_2) \rvert$ appearing in $\Delta$ in the order listed such that no $I$-subbundle appearing in $\Delta$ strictly between $B_1,B_2$ has more than $\lvert E(B_1) \rvert$ edges, then there exists a set $\Gamma$ of pairwise disjoint curves such that the following statements hold.
	\begin{enumerate}
		\item Each curve in $\Gamma$ is contained in $\Delta'$ and internally disjoint from $\overline{\Delta}$.
		\item Each curve in $\Gamma$ connects two vertices in $W$ with the same labels, where $W$ is the set consisting of the last $\lvert E(B_1) \rvert$ vertices in $B_1$, the first $\lvert E(B_2) \rvert$ vertices in $B_2$, and all vertices of $D$ appearing in $\Delta$ between them. 
		\item Every vertex that is one of the last $\lvert E(B_1) \rvert$ vertices of $B_1$ or one of the first $\lvert E(B_2) \rvert$ vertices of $B_2$ is an endpoint of a curve in $\gamma$.
		\item If $\gamma$ is a member of $\Gamma$ such that both its endpoints are labelled with ``fake'', then $\gamma \subseteq \bigcup_{i \in I}\gamma_i$.
		\item For each integer $i \in I$, if $v_i$ is the vertex in the last $\lvert E(B_1) \rvert$ vertices of $B_1$ contained in $\gamma_i$, then there exists a connected component of $E(D) \cup \bigcup_{\gamma \in \Gamma}\gamma$ containing $v_i$ and the vertex in the first $\lvert E(B_2) \rvert$ vertices of $B_2$ contained in $\gamma_i$. 
		\item The number of members of $\Gamma$ with both endpoints labelled with ``real'' is at most $c$.
	\end{enumerate}
\end{lemma}

\begin{pf}
Let $a_1=b_1=c_1=1$, and for $i \in [\rho-1]$, let 
	\begin{itemize}
		\item $a_{i+1}= c_i+i+1$. 
		\item $b_{i+1,0}=a_{i+1}$. 
		\item for every $j \in [2^{4\rho} \cdot (\rho !)^4]$, $b_{i+1,j}= b_{i+1,j-1}+a_{i+1}$, and 
		\item $c_{i+1}=b_{i+1,2^{4\rho} \cdot (\rho !)^4}$.
	\end{itemize}
Define $c=c_\rho$. 

Let Statement 6' be the statement obtained by replacing the number $c$ in Statement 6 in this lemma by $c_{\lvert E(B_1) \rvert}$.
We shall prove that Statements 1-5 and 6' hold by induction on $\lvert E(B_1) \rvert$.
Note that it implies the lemma since $c=c_\rho \geq c_{\lvert E(B_1) \rvert}$.

Note that by Lemma \ref{bundle_index_frist}, there exists $I' \subseteq I$ with $\lvert I' \rvert = \lvert E(B_1) \rvert$ such that every $I$-subbundle between $B_1,B_2$ (including $B_1,B_2$) is an $I'$-subbundle.
For any two vertices $x,y$ in $V(D) \cap \partial \Delta$, where $x$ appears earlier than $y$ in the linear order given by $\Delta$, we denote the set of all vertices of $D$ appearing in $\partial \Delta$ strictly between $x,y$ by $(x,y)$, and define $[x,y]=(x,y) \cup \{x,y\}$.

\noindent{\bf Claim 1:} Statements 1-5 and 6' hold when $\lvert E(B_1) \rvert=1$. 

\noindent{\bf Proof of Claim 1:}
Assume that $\lvert E(B_1) \rvert=1$.
Let $x$ be the last vertex of $B_1$ and $y$ be the first vertex of $B_2$.
If all vertices in $[x,y]$ are labelled with ``fake'', then it is straightforward to see that the set consisting of the connected components of $(\bigcup_{i \in I'}\gamma_i)-\Delta$ having both endpoints labelled with ``fake'' is desired.
So we may assume that $[x,y]$ contains at least one vertex labelled with ``real''.
Let $x'$ be the vertex in $[x,y]$ labelled with ``real'' closest to $x$, and let $y'$ be the vertex in $[x,y]$ labelled with ``real'' closest to $y$.
Since $\lvert I' \rvert = \lvert E(B_1) \rvert=1$, $x'$ and $y'$ have the same labels.
By the definition of the $\Delta$-disentanglement, $x'$ is the last vertex of some $I'$-subbundle and $y'$ is the first vertex of some $I'$-subbundle.
Then it is straightforward to see that the set consisting of the connected components of $(\bigcup_{i \in I'}\gamma_i)-\Delta$ disjoint from $(x',y')$ having both endpoints labelled with ``fake'' and a curve in $\Delta'$ connecting $x',y'$ internally disjoint from $\overline{\Delta}$ satisfies the conclusions.
$\Box$

Claim 1 proves the base of the induction.
Now we assume that $\lvert E(B_1) \rvert \geq 2$ and this lemma holds if $\lvert E(B_1) \rvert$ is smaller.
In addition, we may assume that there exists an $I$-subbundle strictly between $B_1$ and $B_2$, for otherwise the set consisting of some connected component of $\gamma_j$ for each $j \in I'$ is a desired set $\Gamma$.

\noindent{\bf Claim 2:} If no $I$-subbundle strictly between $B_1,B_2$ has $\lvert E(B_1) \rvert$ edges, then there exists $\Gamma$ satisfying Statements 1-5, and there are at most $a_{\lvert E(B_1) \rvert}$ members of $\Gamma$ with both endpoints labelled with ``real''.

\noindent{\bf Proof of Claim 2:}
Assume that no $I$-subbundle strictly between $B_1,B_2$ has $\lvert E(B_1) \rvert$ edges.
Let $M$ be the smallest number such that every $I$-subbundle strictly between $B_1,B_2$ has at most $M$ edges.
So $M \leq \lvert E(B_1) \rvert-1$.
Let $A_1,A_2,...,A_k$ (for some positive integer $k$) be the $I$-subbundles strictly between $B_1,B_2$ with $M$ edges.
Let $I_0$ be the minimal subset of $I$ such that $A_1$ is an $I_0$-subbundle.
By Lemma \ref{bundle_index_frist}, $I_0 \subseteq I'$, and $A_j$ is an $I_0$-subbundle for every $j \in [k]$.

Let $q_1$ be the number of vertices in the last $\lvert E(B_1) \rvert$ vertices of $B_1$ labelled with ``fake''.
Let $q_2$ be the number of vertices in the first $\lvert E(B_2) \rvert$ vertices of $B_2$ labelled with ``fake''.
Note that $M \geq q_1$ and $M \geq q_2$ by the construction of the $\Delta$-disentanglement.
Let $s_1$ be the sequence of formed by the first $\lvert I' \rvert - \lvert I_0 \rvert$ vertices in the last $\lvert E(B_1) \rvert$ vertices of $B_1$.
Let $s_2$ be the sequence of formed by the last $\lvert I' \rvert - \lvert I_0 \rvert$ vertices in the first $\lvert E(B_1) \rvert$ vertices of $B_2$.
So $s_1,s_2$ only contain vertices with labelled ``real''.
Let $s_2'$ be the sequence that lists the entries of $s_2$ in the reverse order.
By the construction of the $\Delta$-disentanglement, $s_1=s_2'$.
So there exist $\lvert I' \rvert - \lvert I_0 \rvert$ pairwise disjoint curves $r_1,r_2,...,r_{\lvert I' \rvert-\lvert I_0 \rvert}$ in $\Delta'$ internally disjoint from $\overline{\Delta}$ such that for every $j \in [\lvert I' \rvert - \lvert I_0 \rvert]$, $r_j$ connects the $j$-th vertex in the last $\lvert E(B_1) \rvert$ vertices of $B_1$ and the $j$-th last vertex in the first $\lvert E(B_2) \rvert$ vertices of $B_2$.
Note that for each $j$, the endpoints of $r_j$ have the same labels.

Let $B_1',B_2'$ be the $I_0$-subbundle of $B_1,B_2$, respectively.
Since $I_0 \subseteq I'$, $\lvert E(B_1') \rvert = \lvert E(B_2') \rvert = \lvert I_0 \rvert = M$.
By the construction of the $\Delta$-disentanglement, the set of the last $\lvert E(B_1') \rvert$ vertices of $B_1'$ equals the set of the last $\lvert E(B_1) \rvert$ vertices of $B_1$, and the set of the first $\lvert E(B_2') \rvert$ vertices of $B_2'$ equals the set of the first $\lvert E(B_2) \rvert$ vertices of $B_2$.
Since $M<\lvert E(B_1) \rvert$, by the induction hypothesis, there exists a set $\Gamma'$ of curves satisfying Statements 1-5 and 6' of this lemma, where $B_1$ and $B_2$ are replaced by $B_1'$ and $B_2'$.
Hence the set $\Gamma' \cup \{r_j: j \in [\lvert I' \rvert - \lvert I_0 \rvert]\}$ satisfies Statements 1-5 and contains at most $c_M+\lvert I' \rvert - \lvert I_0 \rvert \leq c_{\lvert E(B_1) \rvert-1}+ \lvert E(B_1) \rvert \leq a_{\lvert E(B_1) \rvert}$ members with both endpoints labelled with ``real''.
$\Box$

Let $A_1,A_2,...,A_k$ (for some nonnegative integer $k$) be the $I$-subbundles strictly between $B_1,B_2$ with $\lvert E(B_1) \rvert$ edges.
By Claim 2, we may assume that $k \geq 1$.
By Lemma \ref{bundle_index_frist}, $A_j$ is an $I'$-subbundle for every $j \in [k]$.

We say that an ordered pair $(B,B')$ of $I$-subbundles is {\it tight} if $B$ appears prior to $B'$, $\{B,B'\} \subseteq \{B_1,B_2,A_j: j \in [k]\}$, and every $I$-subbundle strictly between $B,B'$ has less than $\lvert E(B) \rvert$ edges.
The {\it type} of a tight pair $(B,B')$ is the sequence $(t_B,t_{B'})$, where $t_B$ is the type of $B$ and $t_{B'}$ is the type of $B'$.
Note that there are at most $2^{2\lvert I' \rvert} \cdot (\lvert I' \rvert !)^2$ different types of $I'$-subbundles.
So there are at most $2^{4\lvert I' \rvert} \cdot (\lvert I' \rvert !)^4$ different types of tight pairs.

\noindent{\bf Claim 3:} If the type of $(B_1,A_1)$ and $(A_k,B_2)$ are the same, and some vertex in the last $\lvert E(B_1) \rvert$ vertices of $B_1$ is labelled with ``real'', then there exists $\Gamma$ satisfying Statements 1-5, and there are at most $c_{\lvert E(B_1) \rvert-1}+1$ members of $\Gamma$ having both endpoints labelled with ``real''.

\noindent{\bf Proof of Claim 3:}
Let $v$ be the first vertex in the last $\lvert E(B_1) \rvert$ vertices of $B_1$.
Let $u$ be the last vertex in the first $\lvert E(B_2) \rvert$ vertices of $B_2$.
Since every $I$-subbundle strictly between $B_1$ and $A_1$ has less than $\lvert E(B_1) \rvert$ edges, $v$ and the last vertex in the first $\lvert E(A_1) \rvert$ vertices of $A_1$ have the same label.
Since $A_1$ and $B_2$ have the same type, $u$ and $v$ have the same label.

Since some vertex in the last $\lvert E(B_1) \rvert$ vertices of $B_1$ is labelled with ``real'', $v$ is labelled with ``real''.
Let $j$ be the number labelled on $v$.
Let $I_1=I'-\{j\}$.
By the induction hypothesis, there exists a set $\Gamma_1$ of curves satisfying Statements 1-5 (where $B_1,B_2$ are replaced by the $I_1$-subbundles of $B_1,B_2$, respectively), such that there are at most $c_{\lvert I_1 \rvert}$ curves in $\Gamma_1$ with both endpoints labelled with ``real''.
Then adding a curve connecting $u$ and $v$ into $\Gamma_1$ gives a set of curves satisfying Statements 1-5, and there are at most $c_{\lvert I_1 \rvert}+1 \leq c_{\lvert E(B_1) \rvert-1}+1$ curves having both endpoints labelled with ``real''.
$\Box$

Let $M$ be the number of types of the tight pairs whose both entries are between $B_1,B_2$ (including $B_1,B_2$). 
So $M \leq 2^{4\lvert I' \rvert} \cdot (\lvert I' \rvert !)^4$.

\noindent{\bf Claim 4:} There exists $\Gamma$ satisfying Statements 1-5, and there are at most $b_{\lvert E(B_1) \rvert,M}$ members of $\Gamma$ with both endpoints labelled with ``real''.

\noindent{\bf Proof of Claim 4:}
We shall prove this claim by induction on $M$.

If for every $B \in \{B_1,A_j: j \in [k]\}$, every vertex in the last $\lvert E(B_1) \rvert$ vertices of $B$ is labelled with ``fake'', then all vertices in the closed interval between the first vertex in the last $\lvert E(B_1) \rvert$ vertices of $B_1$ and the last vertex in the first $\lvert E(B_2) \rvert$ vertices of $B_2$ are labelled with ``fake'', so collecting connected components of $\bigcup_{i \in I'}\gamma_i-\Delta$ gives a desired collection of curves.
So we may assume that there exists $B^* \in \{B_1,A_j: j \in [k]\}$ such that some vertex in the last $\lvert E(B_1) \rvert$ vertices of $B^*$ is labelled with ``real'', and we choose such $B^*$ to be as close to $B_1$ as possible.

Let $A_0=B_1$.
Hence there exists $j_0 \in [k] \cup \{0\}$ such that $B^*=A_{j_0}$.
Let $x_0$ be the first vertex in the last $\lvert E(B_1) \rvert$ vertices of $B_1$.
Let $y_0$ be the last vertex in the first $\lvert E(A_{j_0}) \rvert$ vertices of $A_{j_0}$.
If $j_0=0$, then let $\Gamma_0=\emptyset$; otherwise, let $\Gamma_0$ be the collection of the connected components of $\bigcup_{i \in I'}\gamma_i-\Delta$ incident with a vertex in $[x_0,y_0]$.
Note that no curve in $\Gamma_0$ has both endpoints labelled with ``real'' by the choice of $j_0$.

Let $A_{k+1}=B_2$.
Let $j_1$ be the maximum in $[k]-[j_0-1]$ such that the type of $(A_{j_0},A_{j_0+1})$ equals the type of $(A_{j_1},A_{j_1+1})$.
Note that $j_1$ exists since $j_0$ is a candidate.
If $j_1 \neq j_0$, then by Claim 3, there exists $\Gamma_2$ satisfying Statements 1-5 (where $B_1,B_2$ are replaced by $A_{j_0},A_{j_1+1}$, respectively), such that there are at most $c_{\lvert E(B_1) \rvert-1}+1 \leq a_{\lvert E(B_1) \rvert}$ members of $\Gamma_2$ having both endpoints labelled with ``real''; if $j_1=j_0$, then by Claim 2, there exists $\Gamma_2$ satisfying Statements 1-5 (where $B_1,B_2$ are replaced by $A_{j_0},A_{j_1+1}=A_{j_0+1}$, respectively), such that there are at most $a_{\lvert E(B_1) \rvert}$ members of $\Gamma_2$ having both endpoints labelled with ``real''.

When $M=1$, $j_1=k$, so $\Gamma_0 \cup \Gamma_2$ satisfies Statements 1-5 and contains at most $a_{\lvert E(B_1) \rvert} \leq b_{\lvert E(B_1) \rvert,0}$ curves with both endpoints labelled with ``real''.
So we may assume that $M \geq 2$ and this claim holds if $M$ is smaller.

By the maximality of $j_1$, there are at most $M-1$ types of the tight pairs which both entries are between $A_{j_1+1},B_2$ (including $A_{j_1+1}$ and $B_2$).
By the induction hypothesis, there exists $\Gamma_3$ satisfying Statements 1-5 (where $B_1,B_2$ are replaced by $A_{j_1+1},B_2$), such that there are at most $b_{\lvert E(B_1) \rvert,M-1}$ members of $\Gamma_3$ having both endpoints labelled with ``real''.
Then $\Gamma_0 \cup \Gamma_2 \cup \Gamma_3$ satisfies Statements 1-5 and contains at most $a_{\lvert E(B_1) \rvert}+b_{\lvert E(B_1) \rvert,M-1} \leq b_{\lvert E(B_1) \rvert,M}$ curves with both endpoints labelled with ``real''.
This proves the claim.
$\Box$

Recall that $M \leq 2^{4\lvert I' \rvert} \cdot (\lvert I' \rvert !)^4 \leq 2^{4\rho} \cdot (\rho !)^4$.
So $b_{\lvert E(B_1) \rvert,M} \leq c_{\lvert E(B_1) \rvert}$.
Hence this lemma follows from Claim 4.
\end{pf}

\begin{lemma} \label{simulating miniature 0}
Let $\Sigma$ be a surface and $\Delta$ an oriented open disk in $\Sigma$.
Let $\Delta'$ be an open disk in $\Sigma$ with $\Delta' \supset \Delta$ and $\partial \Delta \cap \partial \Delta' = \emptyset$.
Let $H$ be a connected graph.
Then for every $\Delta$-signature $(Q,s)$ of some pseudo-embedding $\pi$ of $H$ in $\Sigma$ legal with respect to $\{\Delta\}$, there exists an integer $\tau=\tau(Q,s)$ such that the following statement holds:
if $\pi'$ is a pseudo-embedding of $H$ in $\Sigma$ legal with respect to $\{\Delta\}$ with $(Q,s) \preceq' (Q',s')$, where $(Q',s')$ is the $\Delta$-signature of $\pi'$, then there exist 
	\begin{itemize}
		\item a set $\Gamma$ of at most $\tau$ pairwise disjoint curves in $\Delta'-\Delta$, each connecting a pair of points in $\partial\Delta$ and internally disjoint from $\overline{\Delta}$, and 
		\item a homeomorphic embedding $\eta$ from $M$ to $M'+E_\Gamma$, where $(M,\sigma)$ and $(M',\sigma')$ are the $\Delta$-miniatures of $\pi$ and $\pi'$, respectively, and $E_\Gamma=\{xy: x,y \in V(M'),$ some curve in $\Gamma$ has endpoints $x,y\}$
	\end{itemize}
such that 
	\begin{enumerate}
		\item $\eta(\sigma)$ is a subsequence of $\sigma'$, and
		\item $\eta(V(M)-\partial \Delta) = V(M')-\partial \Delta$.
	\end{enumerate}
\end{lemma}

\begin{pf}
Let $L$ be the $\Delta$-highlight of $\pi$ and let $L'$ be the $\Delta$-highlight of $\pi'$.
By the definition of $\preceq'$, there exists an isomorphism $\iota$ from $L$ to an induced subgraph of $L'$ preserving the labels and marches.
Let $\rho$ be the maximum number of edges of a bundle of $L$.
Let $c$ be the number $c(\rho)$ mentioned in Lemma \ref{subbundle}.
Let $A_1,A_2,...,A_k$ be the bundles in $L$ (for some nonnegative integer $k$).
Define $\tau= (k-1)c$.
Note that $\tau$ only depends on $(Q,s)$.

Since $\iota$ is an isomorphism from $L$ to an induced subgraph of $L'$, there exist bundles $B_1,B_2,...,B_k$ in $L'$ such that $B_j=\iota(A_j)$ for each $j \in [k]$.
Let $\{\gamma'_i: i \in [\rho']\}$ be a $\Delta$-foundation for the $\Delta$-disentanglement of $\pi'$ with $\gamma_i' \subseteq \Delta'$ for each $i \in [\rho']$, where $\rho'$ is the maximum number of edges of a bundle of $L'$.
Note that $\rho \leq \rho'$.

\noindent{\bf Claim 1:} For every $i \in [k-1]$, there exists a set $\Y_i$ of pairwise disjoint curves contained in $\Delta'-\Delta$ and internally disjoint from $\overline{\Delta}$ such that the following statements hold.
	\begin{itemize}
		\item[(i)] Each curve in $\Y_i$ connects two vertices in $W$ with the same labels, where $W$ is the set of vertices consisting of the last $\lvert E(B_i) \rvert$ vertices of $B_i$, the first $\lvert E(B_{i+1}) \rvert$ vertices of $B_{i+1}$, and all vertices of $L'$ appearing in $\partial \Delta$ strictly between them.
		\item[(ii)] If $\gamma \in \Y_i$ connects two vertices labelled with ``fake'', then $\gamma \subseteq \bigcup_{i \in [\rho']}\gamma_i'$.
		\item[(iii)] For each $j \in [\rho']$, if the vertex $v_j$ in the last $\lvert E(B_i) \rvert$ vertices of $B_i$ contained in $\gamma_j$ is labelled with ``fake'', then there exists a connected component of $E(L') \cup \bigcup_{\gamma \in \Y_i}\gamma$ containing $v_j$ and the vertex in the first $\lvert E(B_{i+1}) \rvert$ vertices of $B_{i+1}$ contained in $\gamma_j$.
		\item[(iv)] There are at most $c$ members of $\Y_i$ with both endpoints labelled with ``real''.
	\end{itemize}

\noindent{\bf Proof of Claim 1:}
Let $i \in [k-1]$.
Let $Z$ be the set consisting of the last $\lvert E(B_i) \rvert$ vertices of $B_i$ labelled with ``fake''.
Let $I$ be the minimal subset of $[\rho']$ such that the $I$-subbundle of $B_i$ contains $Z$.
Let $B'_i$ and $B'_{i+1}$ be the $I$-subbundles of $B_i$ and $B_{i+1}$, respectively.
Note that $\lvert E(B'_i) \rvert = \lvert E(B'_{i+1}) \rvert = \lvert I \rvert$, and there exists no $I$-subbundle appearing in $\Delta$ strictly between $B_i'$ and $B_{i+1}'$ having more than $\lvert I \rvert = \lvert E(B_i') \rvert$ edges.
So there exists a set $\Gamma$ of curves satisfying Lemma \ref{subbundle} (where $B_1,B_2$ are replaced by $B_i'$ and $B_{i+1}'$).
Then $\Gamma$ is a desired set $\Y_i$ for this claim.
$\Box$

For each $i \in [k-1]$, let $\Y_i$ be the set mentioned in Claim 1.
Let $\Gamma = \{\gamma \in \bigcup_{i \in [k-1]}\Y_i:$ both endpoints of $\gamma$ are labelled with ``real''$\}$.
So $\lvert \Gamma \rvert \leq (k-1)c = \tau$ by (iv).

Let $(M,\sigma)$ and $(M',\sigma')$ be the $\Delta$-miniatures of $\pi$ and $\pi'$, respectively.
Note that $V(M) \subseteq V(L)$ and $V(M') \subseteq V(L')$.
For every $v \in V(M)$, define $\eta(v) = \iota(v)$.
So $\eta(\sigma) = \iota(\sigma)$ is a subsequence of $\sigma'$, and $\eta(V(M)-\partial \Delta) = \iota(V(M)-\partial \Delta) = V(M')-\partial \Delta$.

For every edge $e$ of $M$ incident with at least one vertex in $V(M)-\partial \Delta$, it is an edge of $L$, and we define $\eta(e)=\iota(e)$, so $\eta(e)$ is an edge of $L'$ between $\eta(x)$ and $\eta(y)$, where $x,y$ are the ends of $e$.
Let $(S,\sigma_S)$ and $(S',\sigma_{S'})$ be the $\Delta$-dives of $\pi$ and $\pi'$, respectively.
Let $\{\gamma_i: i \in [\rho]\}$ be a $\Delta$-foundation for the $\Delta$-disentanglement of $\pi$ with $\gamma_i \subseteq \Delta'$ for each $i \in [\rho]$.

We say that an edge $e \in E(M)$ is {\it tricky} if $e$ is not incident with $V(M)-\partial \Delta$.
Note that for every tricky $e \in E(M)$, $e$ is an edge of $S$, the ends of $e$, say $x_e$ and $y_e$, are vertices in $L$ labelled with ``real'', and there exist $j_e \in [\rho]$ and a curve $\gamma_e$ contained in $\gamma_{j_e}$ connecting $x_e$ and $y_e$ such that every internal point of $\gamma_e$ intersecting $\partial \Delta$ is labelled with ``fake''.
So for every tricky $e \in E(M)$, by (iii), some connected component $X_e$ of $E(L') \cup \bigcup_{\gamma \in \bigcup_{i=1}^{k-1}\Y_i}\gamma$ contains $\iota(x_e)$ and $\iota(y_e)$.
Moreover, for every tricky $e \in E(M)$, by statements (i) and (ii) and seeing the labels of the endpoints of the curves in $\bigcup_{\gamma \in \bigcup_{i=1}^{k-1}\Y_i}\gamma$, $X_e$ is a connected component of $E(L') \cup \bigcup_{\gamma \in \Gamma}\gamma \cup \bigcup_{\gamma''}\gamma''$, where the last union is over all curves $\gamma''$ contained in $\bigcup_{i=1}^{\rho'}\gamma_i'$ connecting two vertices labelled with ``fake'' internally disjoint from $\Delta$.

For every tricky $e \in E(M)$, let $X_e'$ be the curve contained in $X_e$ between $\iota(x_e)$ and $\iota(y_e)$, and let $X'_{e,1},X'_{e,2},...,X'_{e,\ell_e}$ (for some positive integer $\ell_e$) be the internally disjoint curves contained in $X_e'$ such that $X_e' = \bigcup_{j=1}^{\ell_e}X'_{e,j}$, and for each $j \in [\ell_e]$, the endpoints of $X'_{e,j}$ are vertices of $L'$ labelled with ``real'', and every internal point of $X'_{e,j}$ intersecting $\partial \Delta$ is labelled with ``fake''.
Note that for each tricky $e \in E(M)$ and $j \in [\ell_e]$, either $X'_{e,j}$ is a member $\gamma_{e,j}$ of $\Gamma$, or $X'_{e,j}$ is not a member of $\Gamma$ and there exists an edge $f_{e,j}$ of $S'$ whose ends are the endpoints of $X'_{e,j}$.

Let $E_\Gamma = \{xy: x,y \in V(M')$, some curve in $\Gamma$ has endpoints $x,y\}$.
For any tricky $e$ and $j \in [\ell_e]$, when $\gamma_{e,j}$ is defined, we define $f_{e,j}$ to be the element in $E_\Gamma$ corresponding to $\gamma_{e,j}$.
For each tricky $e \in E(M)$, define $\eta(e) = \bigcup_{j \in [\ell_e]}f_{e,j}$.
Note that for each tricky $e$, $\eta(e)$ is a path in $S'+E_\Gamma$ between $\iota(x_e)$ and $\iota(y_e)$.
In addition, if $e_1,e_2$ are distinct tricky edges, then points in the intersection of $X_{e_1}$ and $X_{e_2}$ are crossing points in $\Delta$, so $\eta(e_1)$ and $\eta(e_2)$ are disjoint.
Therefore, $\eta$ is a desired homeomorphic embedding from $M$ to $M'+\Gamma_E$.
This proves the lemma.
\end{pf}

\subsection{Addenda, templates and port numbers} \label{subsec:addenda_template}

An {\it addendum} of a pseudo-embedding $\pi$ of a connected graph $H$ in a surface $\Sigma$ legal with respect to a set $\{\Delta_1,...,\Delta_{\kappa}\}$ of oriented open disks with pairwise disjoint closure is a set $\A$ of pairwise disjoint open disks in $\Sigma-\bigcup_{i=1}^\kappa \overline{\Delta_i}$ with pairwise disjoint closure such that the following hold.
	\begin{itemize}
		\item $\bigcup_{\Delta \in \A}\overline{\Delta} \cap \bigcup_{i=1}^\kappa \overline{\Delta_i} = \emptyset$.
		\item $\pi$ has no crossing-point in $\Sigma-(\bigcup_{i=1}^\kappa \Delta_i \cup \bigcup_{\Delta \in \A}\Delta)$.
		\item For every $\Delta \in \A$, $\partial \Delta \cap \pi(V(H)) = \emptyset$, $\Delta \cap \pi(V(H)) \neq \emptyset$, and if $x$ is a point in $\partial \Delta \cap \pi(e)$ for some $e \in E(H)$ with $\pi(e) \not \subseteq \overline{\Delta}$, then there exists an open set $B$ containing $x$ such that $B \cap \overline{\Delta} \cap \pi(e)$ has only one connected component.
	\end{itemize}
Since $\Delta \cap \pi(V(H)) \neq \emptyset$ for every $\Delta \in \A$, we have $\lvert \A \rvert \leq \lvert V(H) \rvert$.
For each $\Delta \in \A$, $v \in V(H)$ with $\pi(v) \in \Delta$, and $e \in E(H)$ incident with $v$ with $\pi(e) \not \subseteq \overline{\Delta}$, we define the following:
	\begin{itemize}
		\item If $e$ is a non-loop edge, then we define $e_v$ to be the point $x$ in $\partial\Delta \cap \pi(e)$ such that some connected component of $\pi(e) \cap (\Delta \cup \{x\})$ contains $\{\pi(v),x\}$.
		\item If $e$ is a loop, then we define $e_v',e_v''$ to be the points $x',x''$, respectively, in $\partial\Delta \cap \pi(e)$ such that some connected component of $\pi(e) \cap (\Delta \cup \{x',x''\})$ contains $\{\pi(v),x',x''\}$.
	\end{itemize}
For each $\Delta \in \A$, define the {\it $\Delta$-caption} to be the labelled rooted graph $(H_\Delta,\Omega_\Delta)$ as follows.
	\begin{itemize}
		\item $V(H_\Delta)=\{v \in V(H): \pi(v) \in \Delta\} \cup \{e_v: v \in V(H), \pi(v) \in \Delta, e$ is a non-loop edge of $H$ incident with $v, \pi(e) \not \subseteq \overline{\Delta}\} \cup \{e_v',e_v'': v \in V(H), \pi(v) \in \Delta, e$ is a loop in $H$ incident with $v, \pi(e) \not \subseteq \overline{\Delta}\}$.
		\item For each vertex $x \in V(H_\Delta) \cap \Delta$, $x$ is labelled with $v$, where $v$ is the vertex of $H$ such that $v=x$.
		\item $E(H_\Delta)=\{e \in E(H): \pi(e) \subseteq \overline{\Delta}\} \cup \{ve_v: v \in V(H), e_v \in V(H_\Delta)\} \cup \{ve_v',ve_v'': v \in V(H), e_v',e_v'' \in V(H_\Delta)\}$.
		\item $\Omega_\Delta$ is a march that orders the vertices $V(H_\Delta)-\{v \in V(H): \pi(v) \in \Delta\}$ in a way that is consistent with a natural ordering of $\partial \Delta$. 
	\end{itemize}
Note that there are at most $2^{\lvert V(H) \rvert + \lvert E(H) \rvert} \cdot (2\lvert E(H) \rvert)!$ non-isomorphic $\Delta$-captions.
The {\it $(\Sigma,(\Delta_1, \allowbreak ...,\Delta_\kappa),\A)$-signature} of $\pi$ is the tuple $(Q_1,Q_2,...,Q_\kappa, \{(H_\Delta,\Omega_\Delta): \Delta \in \A\})$, where $Q_i$ is the $\Delta_i$-signature of $\pi$ for each $i \in [\kappa]$. 
Note that for distinct $\Delta,\Delta' \in \A$, $V(H_\Delta) \cap V(H_{\Delta'}) = \emptyset$.
So $\{(H_\Delta,\Omega_\Delta): \Delta \in \A\}$ is a set with size $\lvert \A \rvert$.

Let $\Sigma$ be a surface and $\Delta_1,...,\Delta_\kappa$ be oriented open disks in $\Sigma$ with pairwise disjoint closure.
Let $H$ be a connected graph.
For any two legal pseudo-embeddings $\pi_1,\pi_2$ of $H$ in $\Sigma$ with respect to $\{\Delta_1,...,\Delta_\kappa\}$ and for every addendum $\A_1$ of $\pi_1$ and every addendum $\A_2$ of $\pi_2$, we define $(Q_1^1,Q_2^1,...,Q_\kappa^1,\A_1) \preceq (Q_1^2,Q_2^2,...,Q_\kappa^2,\A_2)$, where the former is the $(\Sigma,(\Delta_i: i \in [\kappa]),\A_1)$-signature of $\pi_1$ and the latter is the $(\Sigma,(\Delta_i: i \in [\kappa]),\A_2)$-signature of $\pi_2$, if and only if the following conditions hold.
	\begin{itemize}
		\item For each $i \in [\kappa]$, $Q_i^1 \preceq' Q_i^2$.
		\item $\{(H_\Delta,\Omega_\Delta): \Delta \in \A_1\}$ and $\{(H_\Delta,\Omega_\Delta): \Delta \in \A_2\}$ are equal (up to isomorphism for labelled rooted graphs).
	\end{itemize}
Note that $\preceq$ is a quasi-order.
A {\it $(\Sigma,(\Delta_i: i \in [\kappa]))$-template} is a minimal $(\Sigma,(\Delta_i: i \in [\kappa]),\A)$-signature of some pseudo-embedding of $H$ that is legal with respect to $\{\Delta_i: i \in [\kappa]\}$, where $\A$ is an addendum, with respect to $\preceq$.
For $i \in [\kappa]$, we say that $(Q,s)$ is a {\it $\Delta_i$-template} if it is the $i$-th entry of some $(\Sigma,(\Delta_1,...,\Delta_\kappa))$-template.
For any pseudo-embedding $\pi$ of $H$ in $\Sigma$ legal respect to $\{\Delta_i: i \in [\kappa]\}$ and for $i \in [\kappa]$, we say that $(Q,s)$ is a {\it $\Delta_i$-template for $\pi$} if there exists an addendum $\A$ of $\pi$ such that $(Q,s)$ is the $i$-th entry of some $(\Sigma,(\Delta_1,...,\Delta_\kappa))$-template that is less or equal to the $(\Sigma,(\Delta_i: i \in [\kappa],\A))$-signature of $\pi$ with respect to $\preceq$.
Recall that the definition of $\preceq'$ implies that the $(\Sigma, (\Delta_1,...,\Delta_\kappa),\A)$-signatures of two pseudo-embeddings with different $\Delta_i$-depth for some $i \in [\kappa]$ are incomparable with respect to $\preceq$.
Hence for every $\Delta_i$-template $(Q,s)$, there exists a unique positive integer $\rho$ such that $(Q,s)$ can only be an entry of some $(\Sigma,(\Delta_1,...,\Delta_\kappa),\A)$-template of a pseudo-embedding with $\Delta_i$-depth $\rho$. 
We call this number $\rho$ the ${\it depth}$ of a $\Delta_i$-template $(Q,s)$.

\begin{lemma} \label{finitely many templates}
Let $\kappa$ be a positive integer.
Let $\Sigma$ be a surface and $\Delta_1,...,\Delta_\kappa$ be oriented open disks in $\Sigma$ with pairwise disjoint closure.
Let $H$ be a connected graph.
For every positive integer $\rho$, there exists a positive integer $M$ such that there are at most $M$ different $(\Sigma,(\Delta_1,...,\Delta_\kappa))$-templates among all pseudo-embeddings of $H$ in $\Sigma$ that are legal with respect to $\{\Delta_1,...,\Delta_\kappa\}$ with $\Delta_i$-depth at most $\rho$ for all $i \in [\kappa]$.
\end{lemma}

\begin{pf}
Note that $\lvert \A \rvert \leq \lvert V(H) \rvert$ and there are at most $2^{\lvert V(H) \rvert+\lvert E(H) \rvert}(2\lvert E(H) \rvert)!$ non-isomorphic $\Delta$-captions for any addendum $\A$ and $\Delta\in \A$.
By Lemma \ref{wqo Delta-signature}, $\preceq'$ is a well-quasi-ordering since the $\Delta_i$-depth is at most $\rho$ for all $i \in [\kappa]$.
Hence $\preceq$ is a well-quasi-ordering by Theorem \ref{Higman lemma}.
So the number of minimal elements with respect to $\preceq$ is finite.
\end{pf}

\begin{lemma} \label{simulating miniature}
Let $\kappa$ be a positive integer.
Let $\Sigma$ be a surface and $\Delta_1,...,\Delta_\kappa$ be oriented open disks in $\Sigma$ with pairwise disjoint closure.
Let $H$ be a connected graph.
Then for every positive integer $\rho$, there exists an integer $\tau$ such that the following holds.

Let $\pi,\pi'$ be pseudo-embeddings of $H$ in $\Sigma$ legal with respect to $\{\Delta_i: i \in [\kappa]\}$.
Let $i\in[\kappa]$, and let $(Q,s),(Q',s')$ be the $\Delta_i$-signatures of $\pi,\pi'$, respectively.
If $(Q,s) \preceq' (Q',s')$, and $(Q,s)$ is a $\Delta_i$-template with depth at most $\rho$, then for every open disk $\Delta'_i$ in $\Sigma$ with $\Delta'_i \supset \Delta_i$ and $\partial\Delta_i \cap \partial\Delta_i'=\overline{\Delta'_i} \cap (\bigcup_{j \in [\kappa]-\{i\}}\overline{\Delta_j})=\emptyset$, there exists a set $\Gamma$ of at most $\tau$ pairwise disjoint curves in $\Delta'_i-\Delta_i$, each connecting a pair of points in $\partial\Delta_i$ and internally disjoint from $\overline{\Delta_i}$, such that there exists a homeomorphic embedding $\eta$ from $M$ to $M'+E_\Gamma$, where $(M,\sigma)$ and $(M',\sigma')$ are the $\Delta_i$-miniatures of $\pi$ and $\pi'$, respectively, and $E_\Gamma=\{xy: x,y \in V(M'),$ some curve in $\Gamma$ has endpoints $x,y\}$ such that 
	\begin{enumerate}
		\item $\eta(\sigma)$ is a subsequence of $\sigma'$, and
		\item $\eta(V(M)-\partial \Delta_i) = V(M')-\partial \Delta_i$.
	\end{enumerate}
\end{lemma}

\begin{pf}
Define $\tau$ to be the maximum of $\tau_{\ref{simulating miniature 0}}(Q^*,s^*)$ among all entries $(Q^*,s^*)$ of $(\Sigma,(\Delta_1,...,\Delta_\kappa))$-templates of pseudo-embeddings legal with respect to $\{\Delta_j: j \in [\kappa]\}$ with $\Delta_j$-depth at most $\rho$ for all $j \in [\kappa]$, where $\tau_{\ref{simulating miniature 0}}$ is the number $\tau$ mentioned in Lemma \ref{simulating miniature 0}.
Note that the maximum is over finitely many terms by Lemma \ref{finitely many templates}, so $\tau$ exists.
Then this lemma immediately follows from Lemma \ref{simulating miniature 0}.
\end{pf}

\bigskip

Note that the number $\tau$ mentioned in Lemma \ref{simulating miniature} depends on $H, \Sigma, \kappa, \rho$ but does not depend on the set $\{\Delta_1,...,\Delta_\kappa\}$ as we may homotopically move the disks without violating the conclusion of Lemma \ref{simulating miniature}.

Let $H$ be a connected graph, $\Sigma$ be a surface, $\kappa,\rho$ be positive integers, and $\Delta_1,...,\Delta_\kappa$ be oriented open disks in $\Sigma$ with pairwise disjoint closure.
Let $p$ be the minimum number $\tau$ mentioned in Lemma \ref{simulating miniature}. 
Let $i \in [\kappa]$, and let $\pi$ be a pseudo-embedding of $H$ in $\Sigma$ legal with respect to $\{\Delta_j: j \in [\kappa]\}$ such that the $\Delta_i$-depth is at most $\rho$.
Let $(Q,s)$ be a $\Delta_i$-template for $\pi$, and let $\pi_0$ be a pseudo-embedding of $H$ in $\Sigma$ legal with respect to $\{\Delta_i: i \in [\kappa]\}$ such that the $\Delta_i$-signature of $\pi_0$ is $(Q,s)$.
By Lemma \ref{simulating miniature}, there exists a rooted graph $(H_i,\sigma_{H_i})$ such that 
	\begin{itemize}
		\item $H_i$ is a subgraph of $M$ and $\sigma_{H_i}$ is a subsequence of $\sigma$, where $(M,\sigma)$ is the $\Delta_i$-miniature of $\pi$,
		\item for every open disk $\Delta'$ with $\Delta' \supset \Delta_i$ and $\partial \Delta_i \cap \partial\Delta' = \overline{\Delta'} \cap \bigcup_{j \in [\kappa]-\{i\}}\overline{\Delta_j} = \emptyset$, there exists a set $\Gamma$ of pairwise disjoint curves contained in $\Delta'-\Delta_i$ such that $\lvert \Gamma \rvert \leq p$, and each member of $\Gamma$ connects two vertices in $V(H_i) \cap \partial\Delta_i$ and is internally disjoint from $\overline{\Delta_i}$, and 
		\item there exists a homeomorphic embedding $\eta$ from $M_0$ into $H_i+E_\Gamma$, where $(M_0,\sigma_0)$ is the $\Delta_i$-miniature of $\pi_0$ and $E_\Gamma = \{xy: x,y \in V(H_i) \cap \partial \Delta_i,$ some member $\gamma \in \Gamma$ has endpoints $x,y\}$, such that 
			\begin{itemize}
				\item $\eta(\sigma_0)$ is a subsequence of $\sigma_{H_i}$, 
				\item $\eta(V(M_0)-\partial \Delta_i) = V(H_i)-\partial \Delta_i$,
				\item every entry in $\sigma_{H_i}$ is either an entry in $\eta(\sigma_0)$ or incident with an edge in $E_\Gamma$, 
				\item $\sigma_{H_i}$ orders its entries in a way that is consistent with the ordering of $\Delta_i$, and
				\item every vertex of $H_i$ is either in $V(M)-\partial \Delta_i$ or an entry of $\sigma_{H_i}$.
			\end{itemize}
	\end{itemize}
We call such a rooted graph $(H_i,\sigma_{H_i})$ a {\it $(Q,s)$-witness} for $\pi$.
Note that $\lvert V(H_i) \rvert = \lvert V(H_i)-\partial \Delta_i \rvert + \lvert \sigma_{H_i} \rvert$ only depends on $p$, $(Q,s)$ and $H$.

Let $H$ be a connected graph, $\Sigma$ a surface and $\kappa,\rho$ positive integers. 
Let $\Delta_1,...,\Delta_\kappa$ be oriented open disks in $\Sigma$ with pairwise disjoint closure.
For each $i \in [\kappa]$, let $q_i = \max\lvert V(H_i)\rvert$, where the maximum is over all $(Q,s)$-witnesses $(H_i,\sigma_{H_i})$ and all $\Delta_i$-templates $(Q,s)$ with $\Delta_i$-depth at most $\rho$. 
Note that $q_i$ exists since there are only finitely many $\Delta_i$-templates with $\Delta_i$-depth at most $\rho$. 
Let $q=\sum_{i=1}^\kappa q_i$.
Note that $q$ depends on $H, \Sigma, \kappa, \rho$ but does not depend on the set $\{\Delta_1,...,\Delta_\kappa\}$ as we may homotopically move the disks.
We define the {\it $(H,\Sigma,\kappa,\rho)$-port number} to be $q$.

\section{Lifting and gauges} \label{sec: gauges}

In this section, we formally deal with the relationship between homeomorphic embeddings from $H$ into $G$ and pseudo-embeddings of $H$ in $\Sigma$, when $G$ has a segregation with a proper arrangement in $\Sigma$.
We will first show that the legality assumption defined in Section \ref{sec: pseudoembedding} does not affect the generality by using the notion lifting defined below.
Then we will show that there is a ``gauge'' in the skeleton of the arrangement in the surface regulating the behavior of all homeomorphic embeddings from $H$ into $G$ at once.

\subsection{Lifting}

Let $H$ be a connected graph and let $\Sigma$ be a surface.
Let $\kappa$ be a positive integer.
Let $\Delta_1,\Delta_2,...,\Delta_\kappa$ be open disks in $\Sigma$ with disjoint closure.
Let $\pi$ be a pseudo-embedding from $H$ into $\Sigma$ with $\pi(V(H)) \cap \bigcup_{i \in [\kappa]}\partial \Delta_i = \emptyset$ such that there are only finitely many connected components of $\pi(E(H)) \cap \bigcup_{i=1}^\kappa\partial \Delta_i$.
We define a {\it lifting of $\pi$ with respect to $\{\Delta_i: i \in [\kappa]\}$} to be the pseudo-embedding from $H$ into $\Sigma$ obtained from $\pi$ by repeatedly applying the following operations until none of them is applicable:
	\begin{itemize}
		\item[(LIFT1)] If $\gamma$ is a curve with $\gamma \subseteq \pi(e) \cap \bigcup_{i \in [\kappa]}\overline{\Delta_i}$ for some $e \in E(H)$ such that 
			\begin{itemize}
				\item the endpoints of $\gamma$ are two points in $\partial\Delta_i$ for some $i \in [\kappa]$, 
				\item $\gamma$ does not contain any crossing point, and 
				\item there exists a curve $\gamma'$ connecting the endpoints of $\gamma$ internally disjoint from $\bigcup_{i \in [\kappa]}\overline{\Delta_i} \cup \pi(E(H))$,
			\end{itemize}
			then replace $\gamma$ by $\gamma'$.
		\item[(LIFT2)] If $z \in \pi(e) \cap \partial\Delta_i$ for some $e \in E(H)$ and $i \in [\kappa]$ is a point such that there exists an open set $B_z$ containing $z$ such that $B_z \cap \pi(e) \cap \Delta_i=\emptyset=B_z \cap \pi(V(H))$ and $B_z \cap \pi(e)$ has exactly one connected component, and this component is not contained in $\partial \Delta_i$, then replace the curve in $B_z \cap \pi(e)$ by a curve with the same pair of endpoints contained in $B_z$ internally disjoint from $\overline{\Delta_i}$.
		\item[(LIFT3)] If $z \in \pi(e) \cap \partial\Delta_i$ for some $e \in E(H)$ and $i \in [\kappa]$ is a point such that there exists an open set $B_z$ containing $z$ such that $B_z \cap \pi(e) -\overline{\Delta_i}=\emptyset=B_z \cap \pi(V(H))$ and $B_z \cap \pi(e)$ has exactly one connected component, then replace the curve in $B_z \cap \pi(e)$ by a curve with the same pair of endpoints contained in $B_z$ internally disjoint from $\Sigma-\Delta_i$.
	\end{itemize}
Note that for every lifting $\pi'$ of $\pi$ with respect to $\{\Delta_i: i \in [\kappa]\}$, $\pi'(E(H)) \cap \bigcup_{i \in [\kappa]}\partial\Delta_i \subseteq \pi(E(H)) \cap \bigcup_{i \in [\kappa]}\partial\Delta_i$.
Moreover, if for each $i \in [\kappa]$, there exists a point $p_i \in \partial \Delta_i-\pi(E(H))$ witnessing that $\Delta_i$ is an oriented disk, then every lifting of $\pi$ with respect to $\{\Delta_i: i \in [\kappa]\}$ is legal with respect to $\{\Delta_i: i \in [\kappa]\}$.

Let $H$ be a connected graph and $\kappa,\rho$ be positive integers.
Let $\alpha$ be a proper arrangement of a segregation $\Se$ with a $(\kappa,\rho)$-witness $(\Se_1,\Se_2)$ of a graph $G$ in a surface $\Sigma$.
Let $G'$ be the skeleton of $\alpha$.
For each $(S,\Omega) \in \Se_2$, let $C_S$ be a cycle of $G'$ such that the open disk $\Delta_S$ bounded by $C_S$ contains $\alpha(S,\Omega)$.
Assume that $\overline{\Delta_S} \cap \overline{\Delta_{S'}}=\emptyset$ for distinct $(S,\Omega),(S',\Omega') \in \Se_2$.
Let $\eta$ be a homeomorphic embedding from $H$ into $G$ with $\eta(V(H)) \cap \bigcup_{(S,\Omega) \in \Se_2}V(C_S) = \emptyset$.
We define a {\it lifting of $\eta$ with respect to $\{C_S: (S,\Omega) \in \Se_2\}$} to be a lifting of a pseudo-embedding $\pi$ from $H$ into $\Sigma$ with respect to $\{\Delta_S: (S,\Omega) \in \Se_2\}$ such that there exists a pseudo-embedding $\pi_G$ from $G$ into $\Sigma$ such that 
	\begin{itemize}
		\item $\pi_G(v)=\alpha(v)$ for every $v \in \bigcup_{(S,\Omega) \in \Se}\overline{\Omega}$,
		\item for every $v \in V(G)-\bigcup_{(S,\Omega) \in \Se}\overline{\Omega}$, $\pi_G(v)$ is in the interior of $\alpha(S_v,\Omega_v)$, where $(S_v,\Omega_v)$ is the unique member in $\Se$ with $v \in V(S_v)$,
		\item for every $e \in E(G)$, $\pi_G(e)$ is contained in $\alpha(S_e,\Omega_e)$ and internally disjoint from $\partial\alpha(S_e,\Omega_e)$, where $(S_e,\Omega_e)$ is the unique member in $\Se$ with $e \in E(S_e)$, 
		\item $\pi(v)=\pi_G(\eta(v))$ for every $v \in V(H)$, and
		\item $\pi(e)=\bigcup_{e' \in E(\eta(e))}\pi_G(e')$.
	\end{itemize}

\begin{lemma} \label{lifting keeps depth}
Let $H$ be a connected graph and $\kappa,\rho$ be positive integers.
Let $G'$ be the skeleton of a proper arrangement $\alpha$ of a segregation $\Se$ with a $(\kappa,\rho)$-witness $(\Se_1,\Se_2)$ of a graph $G$ in a surface $\Sigma$.
For each $(S,\Omega) \in \Se_2$, let $C_S$ be a cycle of $G'$ such that the open disk $\Delta_S$ bounded by $C_S$ contains $\alpha(S,\Omega)$, and assume that $\Delta_S$ is an oriented disk given by a point $p_S \in \partial \Delta_S-V(G')$.
Assume that $\overline{\Delta_S} \cap \overline{\Delta_{S'}}=\emptyset$ for distinct $(S,\Omega),(S',\Omega') \in \Se_2$.
If $\eta$ is a homeomorphic embedding from $H$ into $G$ with $\eta(V(H)) \cap \bigcup_{(S,\Omega) \in \Se_2}V(C_S)=\emptyset$, then for every lifting $\pi'$ of $\eta$ with respect to $\{C_{S'}:(S',\Omega') \in \Se_2\}$, $\pi'$ is legal with respect to $\{\Delta_S: (S,\Omega) \in \Se_2\}$, and for every $(S,\Omega) \in \Se_2$, 
	\begin{enumerate}
		\item if $(M',\sigma')$ is the $\Delta_S$-miniature of $\pi'$, then there exists a homeomorphic embedding $\eta^*$ from $M'$ to $\eta(E(H)) \cap \bigcup_{(S'',\Omega'') \in \Se, \alpha(S'',\Omega'') \subseteq \overline{\Delta_S}}S''$ such that
			\begin{enumerate}
				\item for every $v \in V(H)$ with $\pi'(v) \in V(M')-\partial\Delta_S$, $\eta^*(\pi'(v))=\eta(v)$,
				\item $\eta^*(\sigma')$ is a subsequence of the march formed by the elements of $V(\eta(E(H))) \cap V(C_S)$ ordered by the order of $\Delta_S$,
				\item for every $x \in V(M') \cap \partial\Delta_S$, if $e \in E(H)$ is the edge of $H$ with $x \in \pi'(e)$, then $\eta^*(x) \in V(\eta(e))$, and 
				\item for every $x \in E(M')$, if $e \in E(H)$ is the edge of $H$ with $x \subseteq \pi'(e)$, then $\eta^*(x) \subseteq \eta(e)$, and
			\end{enumerate}
		\item if $\eta(V(H)) \cap \bigcup_{(S'',\Omega'') \in \Se, \alpha(S'',\Omega'') \subseteq \overline{\Delta_S}}S'' \subseteq V(S)$, then the $\Delta_S$-depth of $\pi'$ is at most $\rho+2$. 
	\end{enumerate}
\end{lemma}

\begin{pf}
Let $\pi'$ be a lifting of $\eta$ with respect to $\{C_{S'}:(S',\Omega') \in \Se_2\}$.
By the definition of a lifting, $\pi'$ is legal with respect to $\{\Delta_S: (S,\Omega) \in \Se_2\}$.

Fix $(S,\Omega)$ to be any member of $\Se_2$.
Statement 1 follows from the definition of a lifting.
So it suffices to prove Statement 2.

Let $\rho'$ be the $\Delta_S$-depth of $\pi'$.
We may assume that $\rho' \geq 3$ for otherwise $\rho' \leq \rho+2$ and we are done.
Hence there exists a partition of $\partial\Delta_S$ into two intervals $I,J$ such that there are $\rho'$ pairwise disjoint paths $P_1,P_2,...,P_{\rho'}$ in the subset of the image of $\pi'$ contained in $\overline{\Delta_S}$ from $I$ to $J$.
By taking short-cuts, we may assume that $P_1,P_2,...,P_{\rho'}$ are internally disjoint from $\partial\Delta_S$.
For each $i \in [\rho']$, let $x_i$ be the end of $P_i$ in $I$.
By changing the indices, we may assume that $x_1,x_2,...,x_{\rho'}$ appear in $I$ in the order listed.
For each $i \in [\rho']$, let $W_i = \eta^*(E(P_i))$, so $W_i$ is a path in $\eta(E(H)) \cap \bigcup_{(S'',\Omega'') \in \Se, \alpha(S'',\Omega'') \subseteq \overline{\Delta_S}}S''$ by Statement 1.

Since $\Delta_S$ is a disk and $\eta(V(H)) \cap \bigcup_{(S'',\Omega'') \in \Se, \alpha(S'',\Omega'') \subseteq \overline{\Delta_S}}S'' \subseteq V(S)$, by (LIFT1) and (LIFT2) mentioned in the definition of a lifting, $W_i$ intersects $E(S)$ for each $i \in [\rho']$.
For each $i \in [\rho']$, let $Q_i$ be the component of $W_i-E(S)$ containing $x_i$, and let $q_i$ be the end of $Q_i$ other than $x_i$.
Since $\Delta_S$ is a disk, there exists the minimal interval $I_S$ of $\Omega$ containing all ends of $Q_i$ in $\overline{\Omega}$ such that $q_1,q_2,...,q_{\rho'}$ appear in $\Omega$ in the order listed.
Note that $I_S$ is unique since $\rho' \geq 3$.
By the minimality of $I_S$, $q_1$ and $q_{\rho'}$ are the ends of $I_S$.

For each $i \in [\rho']$, let $y_i$ be the end of $W_i$ in $J$, let $R_i$ be the component of $W_i-E(S)$ containing $y_i$, and let $r_i$ be the end of $R_i$ other than $y_i$.
Since $\rho' \geq 3$, the planarity of $\Delta_S$ implies that none of $r_i$ is contained in $I_S$.
Let $J_S$ be the interval of $\Omega$ disjoint from $I_S$ such that $I_S \cup J_S=\overline{\Omega}$.
So $\{r_i: i \in [\rho']\} \subseteq J_S$.

Note that $\bigcup_{i=2}^{\rho'-1}W_i$ contains $\rho'-2$ disjoint paths in $S$ from $I_S$ to $J_S$ since $q_1,q_{\rho'}$ are the ends of $I_S$ and $W_i \cap E(S) \neq \emptyset$ for every $i \in [\rho']$.
So $\rho \geq \rho'-2$ since $(S,\Omega)$ is a $\rho$-vortex.
Hence $\rho' \leq \rho+2$.
This proves the lemma.
\end{pf}

\subsection{Gauges}

Let $G$ be a graph.
Let $\Se$ be a segregation of $G$, and let $\Se_1,\Se_2$ be subsets of $\Se$ with $\Se_1 \cap \Se_2=\emptyset$ and $\Se_1 \cup \Se_2=\Se$ such that $\lvert \overline{\Omega} \rvert \leq 3$ for every $(S,\Omega) \in \Se_1$.
Let $\alpha$ be a proper arrangement of $\Se$ in a surface $\Sigma$, and let $G'$ be the skeleton of $\alpha$ with respect to $(\Se_1,\Se_2)$.
Let $k$ be a positive integer.
For $(S, \Omega) \in \Se$, a {\it buffer system around $(S,\Omega)$ with $k$ layers} is a collection $\C_S$ that consists of $k$ disjoint cycles $C_1,C_2,...,C_k$ in $G'$ such that $C_1$ is a cycle bounding an open disk containing the disk $\alpha(S,\Omega)$, and for each $i \in [k]-\{1\}$, $C_i$ is a cycle bounding an open disk containing the disk bounding by $C_{i-1}$ in $\Sigma$. 
A {\it buffer system with $k$ layers of $\alpha$ with respect to $(\Se_1,\Se_2)$} is a collection $\{\C_{S}: (S,\Omega) \in \Se_2\}$ such that for each $(S,\Omega) \in \Se_2$, $\C_S$ is a buffer system around $(S,\Omega)$ with $k$ layers.
For any vertex $v$ of $G'$, a {\it buffer system around $v$ with $k$ layers} is a collection of $k$ pairwise disjoint cycles $C_1,C_2,...,C_k$ in $G'$ such that $v \in \Delta_i$ for each $i \in [k]$, and $\Delta_1 \subset \Delta_2 \subset ... \subset \Delta_k$, where $\Delta_i$ is the open disk in $\Sigma$ bounded by $C_i$ for each $i \in [k]$.

Recall the definition of the $C$-vortex defined below Lemma \ref{extend a vortex}.
Figure \ref{fig_gague_example} might be helpful for understanding Lemma \ref{gauge}.

\begin{figure} 
	\begin{picture}(200,300) (-50,-150)
		\put(200,0){\oval(50,50)} \put(182,0){$\alpha(S,\Omega)$} 
		\put(200,0){\oval(70,70)}
		\put(200,0){\oval(90,90)}
		\put(200,0){\oval(110,110)}
		\thicklines
		\put(200,0){\oval(130,130)}
		\thinlines
		\put(200,0){\oval(150,150)}
		\put(200,0){\oval(170,170)}
		\put(200,0){\oval(190,190)}
		\put(200,0){\oval(210,210)}
		\thicklines
		\put(200,0){\oval(230,230)}
		\thinlines
		\put(200,0){\oval(250,250)}
		\put(200,0){\oval(270,270)}

		\put(150,0){\circle*{5}} 
		\put(255,30){\circle*{5}} 
		\put(190,50){\circle*{5}} 
		\put(200,-10){\circle*{5}}
		\linethickness{0.3mm}
		\qbezier(150,0)(200,-10)(190,50)
		\qbezier(150,0)(150,-50)(200,-10)
		\qbezier(200,-10)(250,-10)(190,50)
		\qbezier(200,-10)(270,-30)(255,30)

		\put(80,-50){\circle*{5}} 
		\put(80,90){\circle*{5}} 
		\put(230,130){\circle*{5}}
		\put(350,-20){\circle*{5}} 
		\qbezier(150,0)(100,-180)(80,-50)
		\qbezier(150,0)(200,50)(80,90)
		\qbezier(80,-50)(130,30)(80,90)
		\qbezier(80,-50)(55,-50)(30,-50)
		\qbezier(80,-50)(60,-10)(30,0)
		\qbezier(190,50)(200,80)(230,130)
		\qbezier(200,-10)(50,-130)(350,-20)

		\qbezier(80,90)(55,105)(30,120)
		\qbezier(80,90)(55,90)(30,90)
		\qbezier(80,90)(55,70)(30,50)

		\qbezier(230,130)(200,-10)(350,50)
		\qbezier(350,50)(450,50)(280,0)
		\qbezier(280,0)(260,-10)(350,-20)

		\qbezier(350,-20)(400,-10)(380,-50)

	\end{picture}
	\caption{An example for Statement \ref{vortex_picture} in Lemma \ref{gauge}. $(S,\Omega)$ is a member of $\Se_2$. The cycles drawn by thin lines are the members of $\C_S$. The outer cycle drawn by a thick line is $C_S$. The inner cycle drawn by a think line is $C_S'$. Other curves and vertices denote $\pi$.} \label{fig_gague_example}
\end{figure}

\begin{lemma} \label{gauge}
Let $H$ be a connected graph, $\kappa,\rho$ be positive integers, $\theta_0,\mu$ be non-decreasing and nonnegative functions with domain ${\mathbb Z}$, and $\Sigma$ be a surface.
Then there exist positive integers $\theta,\tau,\rho^*,\lambda^*$ such that the following hold.
Let $G$ be a graph and $\Se$ a segregation of $G$ with a $(\kappa,\rho)$-witness $(\Se_1,\Se_2)$.
Let $\alpha$ be a proper arrangement of $\Se$ in $\Sigma$ such that the skeleton $G'$ of $\alpha$ has a respectful tangle $\T'$ of order at least $\theta$ such that $m_{\T'}(\overline{\Omega},\overline{\Omega'}) \geq \tau$ for all distinct $(S,\Omega),(S',\Omega') \in \Se_2$.
Then there exists a buffer system $\{\C_S: (S,\Omega) \in \Se_2\}$ of $\alpha$ with respect to $(\Se_1,\Se_2)$ such that for any $(S,\Omega) \in \Se_2$, $\lvert \C_S \rvert \leq \lambda^*$ and each member of $\C_S$ is a cycle bounding a $\lambda^*$-zone in $G'$, and if we fix a linear ordering of each of the boundaries of those $\lambda^*$-zones, where the first point is not in $V(G')$, then for every homeomorphic embedding $\pi$ from $H$ into $G$, there exist $\{C_S \in \C_S: (S,\Omega) \in \Se_2\}$ and $\{C'_S \in \C_S: (S,\Omega) \in \Se_2\}$ such that the following hold.
	\begin{enumerate}
		\item $V(C_S) \cap V(C_{S'}) =\emptyset$ for all $(S,\Omega),(S',\Omega') \in \Se_2$.
		\item $V(C_S) \cap \pi(V(H))=\emptyset$ for all $(S,\Omega) \in \Se_2$.
		\item For each $(S,\Omega) \in \Se_2$, $C_S$ bounds a $\lambda_S$-zone $\Lambda_S$ in $G'$ for some $\lambda_S \leq \lambda^*$, and the $C_S$-vortex $(S_{C_S},\Omega_{C_S})$ is a $\rho_S$-vortex for some $\rho_S \leq \rho^*$. 
		\item If $e \in E(H)$ and $(S,\Omega) \in \Se_2$ with $V(C_S) \cap V(\pi(e)) \neq \emptyset$, then $\pi(e) \not \subseteq S_{C_S}$. 
		\item\label{vortex_picture} For each $(S,\Omega) \in \Se_2$, 
				\begin{enumerate}
					\item there exist at least $\mu(\sum_{(S',\Omega') \in \Se_2}\rho_{S'})$ disjoint cycles in $G'$ contained in $\overline{\Lambda_S}$, where the open disks bounded by those cycles are nested and each of them contains the closed disk bounded by $C_S'$, 
					\item $C_S'$ bounds a $\lambda_S$-zone in $G'$ such that $(\pi(V(H)) \cap V(S_{C_S})) \cup \bigcup_{e \in E(H), \pi(e) \subseteq S_{C_S}}\pi(e)$ is contained in $S_{C_S'}-V(C_S')$, where $(S_{C_S'},\Omega_{C_S'})$ is the $C_S'$-vortex, and 
					\item $m_{\T'}(V(C_S),V(C_S')) \geq \mu(\sum_{(S',\Omega') \in \Se_2}\rho_{S'})$.
				\end{enumerate}
		\item\label{gague_X} There exists a set $X_\pi \subseteq V(G')$ such that for every $x \in X_\pi$, there exists a buffer system $\C_x$ around $x$ and cycles $C_x,C'_x \in \C_x$ such that the following hold.
			\begin{enumerate}
				\item For each $x \in X_\pi$, $C_x$ bounds a $\lambda_x$-zone $\Lambda_x$ around $x$ in $G'$, where $\lambda_x \leq \lambda^*$, such that $\pi(V(H))\cap \bigcup_{(S,\Omega) \in \Se_1, \alpha(S,\Omega) \subseteq \overline{\Lambda_x}}V(S)-\bigcup_{(S,\Omega) \in \Se_2}V(S) \neq \emptyset$. 
				\item $\overline{\Lambda_x} \cap \overline{\Lambda_{x'}} = \emptyset$ for distinct $x, x' \in X_\pi$.
				\item $V(C_x) \cap \pi(V(H))=\emptyset$ for every $x \in X_\pi$.
				\item If $x \in X_\pi$ and $e \in E(H)$ with $V(C_x) \cap V(\pi(e)) \neq \emptyset$, then $\pi(e) \not \subseteq S_{C_x}$, where $(S_{C_x},\Omega_{C_x})$ is the $C_x$-vortex. 
				\item There exist at least $\mu(\sum_{(S',\Omega') \in \Se_2}\rho_{S'})$ disjoint cycles contained in $\overline{\Lambda_x}$, where the open disks bounded by those cycles are nested and each of them contains the closed disk bounded by $C_x'$.
				\item $C_x'$ bounds a $\lambda_x$-zone in $G'$ such that $(\pi(V(H)) \cap V(S_{C_x})) \cup \bigcup_{e \in E(H),\pi(e) \subseteq S_{C_x}}\pi(e)$ is contained in $S_{C_x'}-V(C_x')$, where $(S_{C_x'},\Omega_{C_x'})$ is the $C_x'$-vortex. 
				\item $m_{\T'}(V(C_x), V(C_x')) \geq \mu(\sum_{(S',\Omega') \in \Se_2}\rho_{S'})$.
				\item\label{gague_X_Nx} For each $x \in X_\pi$, let $N_x$ be the union of the connected components of $\pi(E(H)) \cap S_{C_x}$ intersecting $\pi(V(H))$, and let $Y_x= \{y \in V(N_x) \cap V(C_x):$ there exists a path in $\pi(e)$ between $y$ and $\pi(V(H))$ internally disjoint from $V(C_x)$, where $e \in E(H)$ is the edge with $y \in V(\pi(e))\}$. 
		Then for each $x \in X_\pi$, $\lvert Y_x \rvert \leq 2\lvert E(H) \rvert$. 
			\end{enumerate}
		\item $\pi(V(H)) \subseteq \bigcup_{x \in X_\pi} (V(S_{C_x})-V(C_x)) \cup \bigcup_{(S,\Omega) \in \Se_2} (V(S_{C_S})-V(C_S))$. 
		\item For distinct $y_1,y_2 \in \Se_2 \cup X_\pi$, $S_{C_{y_1}}$ and $S_{C_{y_2}}$ are disjoint.
		\item If $\pi^*$ is a lifting of $\pi$ with respect to $\{C_S: (S,\Omega) \in \Se_2\}$, then the following hold.
			\begin{enumerate}
				\item $\pi^*$ is legal with respect to $\{\Lambda_S: (S,\Omega) \in \Se_2\}$.
				\item For every $(S,\Omega) \in \Se_2$, the $\Lambda_S$-depth of $\pi^*$ is at most $\rho_S$.
				\item\label{Y_distance_gague} Let $\T''$ be the tangle obtained from $\T'$ by clearing $\bigcup_{(S',\Omega') \in \Se_2}\Lambda_{S'} \cup \bigcup_{x' \in X_\pi}\Lambda_{x'}$.
					For each $(S,\Omega) \in \Se_2$ and for each $\Lambda_S$-template $(Q_S,s_S)$ such that $(Q_S,s_S)$ is at most the $\Lambda_S$-signature of $\pi^*$ with respect to $\preceq'$, let $(H_S,\sigma_{H_S})$ be a $(Q_S,s_S)$-witness for $\pi^*$ and $Y_S=V(C_S) \cap V(H_S)$.
				Then the following hold.
					\begin{enumerate}
						\item The order of $\T''$ is at least $\theta_0(\sum_{(S',\Omega')\in \Se_2}\rho_{S'})+ \sum_{(S',\Omega') \in \Se_2}\lvert Y_{S'} \rvert+\sum_{x \in X_\pi}\lvert Y_x \rvert$.
						\item $m_{\T''}(Y_x,Y_{x'}) \geq \theta_0(\sum_{(S',\Omega')\in \Se_2}\rho_{S'})$ for distinct $x,x' \in X_\pi$.
						\item $m_{\T''}(Y_S,Y_x) \geq \theta_0(\sum_{(S',\Omega')\in \Se_2}\rho_{S'})$ for all $(S,\Omega) \in \Se_2$ and $x \in X_\pi$.
						\item $m_{\T''}(Y_S,Y_{S'}) \geq \theta_0(\sum_{(S',\Omega')\in \Se_2}\rho_{S'})$ for distinct $(S,\Omega),(S',\Omega') \in \Se_2$.
					\end{enumerate}
			\end{enumerate}
	\end{enumerate}
\end{lemma}

\begin{pf}
Let $H$ be a connected graph, $\kappa,\rho$ be positive integers, $\theta_0,\mu$ be non-decreasing and nonnegative functions with domain ${\mathbb Z}$, and $\Sigma$ be a surface.
Let $\lambda_0=5$ and $\lambda_i=\lambda_{i-1}+21$ for $i \geq 1$. 
For every nonnegative integer $i$, let $\rho_i=4(\rho+\lambda_i+15)$. 
Let $f_1$ be the function such that $f_1(x)=x+\theta_0(\kappa\rho_x)+\mu(\kappa  \rho_x)+(\kappa+\lvert V(H) \rvert)(4\lambda_x+2)+\lambda_x$ for every nonnegative integer $x$, and for $i \geq 2$, let $f_i$ be the function such that $f_i(x)=f_1(f_{i-1}(x))$. 
Let $h=(2\lvert V(H) \rvert+ \lvert E(H) \rvert+1)(\kappa+\lvert V(H) \rvert+1)$.
Let $\tau_0= 3f_h(\kappa+\lvert V(H) \rvert)$.
Define $\rho^*=\rho_{\tau_0}+2\lambda_{\tau_0}$, $\lambda^*=\lambda_{\tau_0}$ and $\tau=\theta_0(\kappa\rho^*)+2\lambda_{\tau_0}+(4\lambda_{\tau_0}+2)\tau_0$.
Let $r$ be the $(H,\Sigma,\kappa,\rho^*)$-port number.
Define $\theta=\theta_0(\kappa\rho^*)+ r+ 2\lvert V(H) \rvert \lvert E(H) \rvert+(4\lambda^*+2)(\kappa + \lvert V(H) \rvert) + \tau$.

Let $G$ be a graph and $\Se$ a segregation of $G$ with a $(\kappa,\rho)$-witness $(\Se_1,\Se_2)$.
Let $\alpha$ be a proper arrangement of $\Se$ in $\Sigma$ such that the skeleton $G'$ of $\alpha$ has a respectful tangle $\T'$ of order at least $\theta$ such that $m_{\T'}(\overline{\Omega},\overline{\Omega'}) \geq \tau$ for all distinct $(S,\Omega),(S',\Omega') \in \Se_2$.

For each $(S,\Omega) \in \Se_2$, define $\C_S$ to be $\{C_{S,i}: 1 \leq i \leq \tau_0\}$ as follows.
Let $x_S$ be a vertex in $\overline{\Omega}$.
Let $C_{S,0}$ be the cycle of a 5-zone $\Lambda_{S,0}$ around $x_S$ containing $\alpha(S,\Omega)$.
The existence of $C_{S,0}$ follows from Lemma \ref{big zone contains ball}. 
For $i \geq 1$, define $C_{S,i-1}'$ to be the cycle bounding a $(\lambda_{i-1}+7)$-zone around $x_S$ such that $C_{S,i-1}'$ is disjoint from $C_{S,i-1}$ and the closed disk bounded by $C_{S,i-1}'$ contains $\Lambda_{S,i-1}$, and define $C_{S,i}$ to be the cycle $C$ mentioned in Lemma \ref{extend a vortex} by taking $(S,\Omega)=(S,\Omega)$ and $t=\lambda_{i-1}+7$, such that the closed disk bounded by $C_{S,i}$ contains $C_{S,i-1}'$.
Note that the existence of $C_{S,i-1}'$ follows from Lemma \ref{disjoint boundary of zone}.
By Lemma \ref{extend a vortex}, the cycle $C_{S,i}$ bounds a $(\lambda_{i-1}+7+14)$-zone $\Lambda_{S,i}$, and since $m_{\T'}(\overline{\Omega_1},\overline{\Omega_2}) \geq \tau \geq 2\lambda_{\tau_0}$ for all distinct $(S_1,\Omega_1),(S_2,\Omega_2) \in \Se_2$, if $i \leq \tau_0$, then the $C_{S,i}$-vortex, denoted by $(S_{S,i},\Omega_{S,i})$, is a $\rho_i$-vortex.
Note that $\lambda_i=\lambda_{i-1}+21$, so $\Lambda_{S,i}$ is a $\lambda_i$-zone and contains all atoms $y$ with $m_{\T'}(x_S,y) \leq \lambda_{i-1}+7$.

Note that for any $(S,\Omega) \in \Se_2$, $\lvert \C_S \rvert = \tau_0 \leq \lambda_{\tau_0} = \lambda^*$ and each member of $\C_S$ bounds a $\lambda_{\tau_0}$-zone in $G'$.
Since $\tau > 2\lambda_{\tau_0}$, no cycle belonging to $\C_S$ intersects a cycle belonging to $\C_{S'}$, for distinct $(S,\Omega),(S',\Omega') \in \Se_2$.
We shall prove that the buffer system $\{\C_S: (S,\Omega) \in \Se_2\}$ satisfies the conclusions of this lemma.

For each cycle in $\{\C_S: (S,\Omega) \in \Se_2\}$, fix a linear ordering on the boundary of the disk bounded by it, where the first point is not in $V(G')$.
Let $\pi$ be a homeomorphic embedding from $H$ into $G$.

For each $v \in V(H)$ with $\pi(v) \in V(S_1)$ for some $(S_1,\Omega_1) \in \Se_1$, we fix such an $(S_1,\Omega_1)$ and let $x_v$ be a vertex in $\overline{\Omega_1}$ and $\C_{x_v}$ be a buffer system $\{C_{x_v,i}: i \in [\tau_0]\}$ such that for every $i \in [\tau_0]$, the cycle $C_{x_v,i}$ bounds a $\lambda_i$-zone $\Lambda_{x_v,i}$ around $x_v$ with $\alpha(S_1,\Omega_1) \subseteq \overline{\Lambda_{x_v,i}}$ such that the $C_{x_v,i}$-vortex, denoted by $(S_{x_v,i},\Omega_{x_v,i})$, is a $\rho_i$-vortex. 
Note that $\C_{x_v}$ exists by a similar argument as the existence of $\C_S$ for $(S,\Omega) \in \Se_2$.
Let $X=\{x_v: v \in V(H)\}$.

For every $y \in \Se_2 \cup X$, let $b_y = 1$ and $\rho_y = \rho_1$.
For simplicity, for each $(S,\Omega) \in \Se_2$ and $i \in [\tau_0]$, we also write $C_{S,i}$, $\Lambda_{S,i}$ and $(S_{S,i},\Omega_{S,i})$ as $C_{(S,\Omega),i}$, $\Lambda_{(S,\Omega),i}$ and $(S_{(S,\Omega),i},\Omega_{(S,\Omega),i})$, respectively, and we write $b_{(S,\Omega)}$ as $b_S$.
We modify $b_y$ and $\rho_y$ for all $y \in \Se_2 \cup X$ and modify $X$ by repeatedly applying the following operations until none of them is applicable.
	\begin{itemize}
		\item[(OP1)] If there exists $y \in \Se_2 \cup X$ such that $V(C_{y,b_y}) \cap \pi(V(H)) \neq \emptyset$, then replace $b_y$ by $b_y+1$, and replace $\rho_y$ by $\rho_{b_y+1}$.
		\item[(OP2)] If there exist $y \in \Se_2 \cup X$ and $e \in E(H)$ such that $\pi(e) \subseteq S_{y,b_y}$ and $\pi(e) \cap V(C_{y,b_y}) \neq \emptyset$, then replace $b_y$ by $b_y+1$, and replace $\rho_y$ by $\rho_{b_y+1}$. 
		\item[(OP3)] If there exist $y \in \Se_2 \cup X$ and $x_v \in X-\{y\}$ such that either $\pi(v) \in V(S_{y,b_y})$ or $m_{\T'}(V(C_{y,b_y}),V(C_{x_v,b_{x_v}})) \leq j$, where $j=\theta_0(\sum_{(S',\Omega') \in \Se_2}\rho_{b_{S'}}) + \mu(\sum_{(S',\Omega') \in \Se_2}\rho_{b_{S'}})+\sum_{y \in \Se_2 \cup X}(4\lambda_{b_y}+2)$, then replace $b_y$ by $b_y + j+\lambda_{b_{x_v}}$, replace $\rho_y$ by $\rho_{b_y+j+\lambda_{b_{x_v}}}$, and remove $x_v$ from $X$. 
		\item[(OP4)] If none of (OP1)-(OP3) is applicable, and there exists $y \in \Se_2 \cup X$ such that either 
			\begin{itemize}
				\item $b_y \leq j$, or 
				\item $\pi(V(H)) \cap V(S_{y,b_y})-(V(S_{y,b_y-j})-V(C_{y,b_y-j})) \neq \emptyset$, or 
				\item there exists $e \in E(H)$ with $\pi(e) \subseteq S_{y,b_y}$ but $\pi(e) \not \subseteq S_{y,b_y-j}-V(C_{y,b_y-j})$, 
			\end{itemize}
		where $j=\mu(\sum_{(S',\Omega') \in \Se_2}\rho_{b_{S'}})$, then 
			\begin{itemize}
				\item[(i)] if either $\pi(V(H)) \cap V(S_{y,b_y+j})-(V(S_{y,b_y})-V(C_{y,b_y})) \neq \emptyset$, or there exists $e \in E(H)$ such that $\pi(e) \subseteq S_{y,b_y+j}$ but $\pi(e) \not \subseteq S_{y,b_y}-V(C_{y,b_y})$, then replace $b_y$ by $b_y+j$, and replace $\rho_y$ by $\rho_{b_y+j}$,
				\item[(ii)] otherwise, replace $b_y$ by $b_y+j$.
			\end{itemize}
	\end{itemize}

\noindent{\bf Claim 1:} (OP1)-(OP4) are applied at most $(2\lvert V(H) \rvert + \lvert E(H) \rvert+1)(\kappa+\lvert V(H) \rvert+1)=h$ times. 

\noindent{\bf Proof of Claim 1:}
Let $W = \bigcup_{y \in \Se_2 \cup X}S_{y,b_y}$. 
Let $D = \lvert \pi(V(H)) \cap V(W) - \bigcup_{y \in \Se_2 \cup X}V(C_{y,b_y}) \rvert  + \lvert \{e \in E(H): \pi(e) \subseteq W, V(\pi(e)) \cap \bigcup_{y \in \Se_2 \cup X}V(C_{y,b_y}) = \emptyset\}\rvert + \lvert V(H) \rvert-\lvert X \rvert$.
It is clear that $D$ strictly increases whenever any of (OP1)-(OP3) is applied, and $D$ does not decrease when (OP4) is applied.
So (OP1)-(OP3) can be applied at most $2 \lvert V(H) \rvert + \lvert E(H) \rvert$ times.

Now we prove that (OP4) cannot be applied consecutively more than $\kappa+\lvert V(H) \rvert$ times without increasing $D$.
Assume that (OP4) is just applied on $y \in \Se_2 \cup X$, and $D$ does not increase.
If (i) happens, then since $D$ does not increase, (OP4) cannot be applied again until (OP1) or (OP2) is applied on this $y$.
So we may assume that (ii) happens, and hence $\mu(\sum_{(S',\Omega') \in \Se_2}\rho_{b_{S'}})$ remains the same.
Since (ii) happens, (OP4) cannot be applied on the same $y$ until either any of (OP1)-(OP3) is applied or (OP4) is applied on another member $y' \in \Se_2 \cup X$.
If (i) happens when applying (OP4) on $y'$, then either $D$ strictly increase, or no more (OP4) cannot be applied until (OP1) or (OP2) is applied on $y'$; if (ii) happens when applying (OP4) on $y'$, then $\mu(\sum_{(S',\Omega') \in \Se_2}\rho_{b_{S'}})$ remains the same and (OP4) cannot be applied to $y$ and $y'$.
By repeating this argument, we know that (OP4) can be consecutively applied at most $\lvert \Se_2 \rvert+\lvert X \rvert \leq \kappa+\lvert V(H) \rvert$ times without increasing $D$.

Therefore, the total number of times that (OP1)-(OP4) are applied is at most $(2\lvert V(H) \rvert + \lvert E(H) \rvert+1)(\kappa+\lvert V(H) \rvert+1)$.
$\Box$

By Claim 1, the process must terminate.
Then it is easy to prove by induction that whenever any of (OP1)-(OP4) is applied, $\sum_{y \in \Se_2 \cup X}b_y$ is replaced by a number at most $f_1(\sum_{y \in \Se_2 \cup X}b_y)$. 
Hence by Claim 1, $\sum_{y \in \Se_2 \cup X}b_y \leq f_{h}(\kappa+\lvert V(H) \rvert) \leq \tau_0/3$. 

For every $(S,\Omega) \in \Se_2$, define $\rho_S = \rho_{b_S}$.
Let $\beta = \mu(\sum_{(S',\Omega') \in \Se_2}\rho_{b_{S'}})$.

For each $(S,\Omega) \in \Se_2$, define $C_S = C_{S,b_S}$ and $C'_S = C_{S,b_S-\beta}$.
Clearly, Statement 1 of this lemma holds.
Note that for each $(S,\Omega) \in \Se_2$, the $C_S$-vortex, denoted by $(S_{C_S},\Omega_{C_S})$, is a $\rho_{b_S}$-vortex. 
Since $\rho^* \geq \rho_{\tau_0} \geq \rho_{b_S}=\rho_S$ for each $(S,\Omega)$, Statement 3 holds.
Since (OP1) and (OP2) are not applicable, Statements 2 and 4 hold.

Since (OP4) is not applicable, Statement 5(a) holds.
Since (OP1), (OP2) and (OP4) are not applicable, Statement 5(b) holds.
Suppose that $m_{\T'}(V(C_S),V(C_S')) \leq \beta-1$ for some $(S,\Omega) \in \Se_2$.
Then $\beta \geq 2$.
So $m_{\T'}(x_S,V(C_S)) \leq m_{\T'}(x_S,V(C_S'))+m_{\T'}(V(C_S),V(C_S')) \leq m_{\T'}(x_S,V(C_S'))+\beta-1 \leq \lambda_{b_S-\beta}+\beta-1 \leq \lambda_{b_S-2}+7$.
Since $\Lambda_{S,b_S-1}$ contains all atoms $x$ of $G'$ with $m_{\T'}(x,x_S) \leq \lambda_{b_S-2}+7$ by the definition of $C_{S,b_S-1}$, so $V(C_S)=V(C_{S,b_S})$ is not disjoint from $V(C_{S,b_S-1})$, a contradiction.
Hence Statement 5(c) holds.

Define $X_\pi=X$.
For each $x \in X$, define $C_x=C_{x,b_x}$ and $C'_x=C_{x,b_x-\beta}$.
Then Statement 6(a) holds.
Since (OP1)-(OP3) are not applicable, Statements 6(b)-6(d) hold.
Statements 6(e)-6(g) follow from a similar argument for showing Statements 5(a)-5(c).
For each $x \in X_\pi$, let $N_x$ be the union of the connected components of $\pi(E(H)) \cap S_{C_x}$ intersecting $\pi(V(H))$, and let $Y_x= \{y \in V(N_x) \cap V(C_x):$ there exists a path in $\pi(e)$ between $y$ and $\pi(V(H))$ internally disjoint from $V(C_x)$, where $e \in E(H)$ is the edge with $y \in V(\pi(e))\}$.
Since each edge $e \in E(H)$ has at most two ends, $\lvert V(\pi(e)) \cap Y_x \rvert \leq 2$.
So Statement 6(h) holds.

Statement 7 in this lemma clearly holds.
Since (OP3) is not applicable and $m_{\T'}(\overline{\Omega},\overline{\Omega'}) \geq \tau > 2\lambda^*+2$, Statement 8 holds.
In particular, the disks bounded by $C_S$ for distinct $(S,\Omega) \in \Se_2$ are disjoint.

Let $\pi^*$ be a lifting of $\pi$ with respect to $\{C_S: (S,\Omega) \in \Se_2\}$.
Since Statements 7 and 8 in this lemma hold, we can apply Lemma \ref{lifting keeps depth} for each $(S,\Omega) \in \Se_2$ (by taking $(S,\Omega)$ in Lemma \ref{lifting keeps depth} to be $(S_{C'_S},\Omega_{C'_S})$), so Statement 9(a) holds.
Since Statements 5(b) and 6(f) hold, and since $(S_{C'_S},\Omega_{C'_S})$ is a $\rho_{b_S-\beta}$-vortex and $\rho_S=\rho_{b_S} \geq \rho_{b_S-\beta}+2$, we know that Statement 9(b) holds. 

Let $\T''$ be the tangle obtained from $\T'$ by clearing $\bigcup_{(S',\Omega') \in \Se_2}\Lambda_{S'} \cup \bigcup_{x' \in X_\pi}\Lambda_{x'}$.
For each $(S,\Omega) \in \Se_2$, let $(Q_S,s_S)$ be a $\Lambda_S$-template such that $(Q_S,s_S)$ is at most the $\Lambda_S$-signature of $\pi^*$ with respect to $\preceq'$, and we let $(H_S,\sigma_{H_S})$ be a $(Q_S,s_S)$-witness for $\pi^*$ and $Y_S=V(C_S) \cap V(H_S)$.
Note that $(Q_S,s_S)$ exists since Statement 9(b) holds.
Since $\rho_S \leq \rho^*$ for every $(S,\Omega) \in \Se_2$, $\sum_{(S,\Omega) \in \Se_2}\lvert Y_S \rvert \leq \sum_{(S,\Omega) \in \Se_2}\lvert V(H_S) \rvert$ is at most the $(H,\Sigma,\kappa,\rho^*)$-port number, which is denoted by $r$.
Since $\Lambda_y$ is a $\lambda^*$-zone for each $y \in \Se_2 \cup X$, the order of $\T''$ is at least $\theta-(4\lambda^*+2)(\lvert \Se_2 \rvert + \lvert X \rvert) \geq \theta-(4\lambda^*+2)(\kappa+\lvert V(H) \rvert) \geq \theta_0(\kappa\rho^*)+r+ 2\lvert V(H) \rvert \lvert E(H) \rvert+ \tau \geq \theta_0(\sum_{(S',\Omega')\in \Se_2}\rho_{S'})+ \sum_{(S',\Omega') \in \Se_2}\lvert Y_{S'} \rvert+\sum_{x \in X_\pi}\lvert Y_x \rvert$.
So Statement 9(c)(i) holds.

For each $(S,\Omega) \in \Se_2$, we write $Y_S$ as $Y_{(S,\Omega)}$ and write $\Lambda_{S}$ as $\Lambda_{(S,\Omega)}$.
Note that for each $y \in \Se_2 \cup X$, $\Lambda_y$ is a $\lambda_{b_y}$-zone.
For distinct $y_1,y_2 \in \Se_2 \cup X$ with $\{y_1,y_2\} \not \subseteq \Se_2$, since (OP3) is not applicable, $m_{\T''}(Y_{y_1},Y_{y_2}) \geq m_{\T''}(C_{y_1},C_{y_2}) \geq m_{\T'}(C_{y_1},C_{y_2})-\sum_{y \in \Se_2 \cup X}(4\lambda_{b_y}+2) \geq \theta_0(\sum_{(S',\Omega') \in \Se_2}\rho_{S'}) + \mu(\sum_{(S',\Omega') \in \Se_2}\rho_{S'})$.
So Statements 9(c)(ii) and 9(c)(iii) hold.
For distinct $(S,\Omega),(S',\Omega') \in \Se_2$, $m_{\T''}(Y_S,Y_{S'}) \geq m_{\T'}(Y_S,Y_{S'})-(4\lambda^*+2)(\kappa+\lvert V(H) \rvert) \geq (m_{\T'}(\overline{\Omega},\overline{\Omega'})-2\lambda^*)-(4\lambda^*+2)(\kappa+\lvert V(H) \rvert) \geq (\tau-2\lambda^*)-(4\lambda^*+2)(\kappa+\lvert V(H) \rvert) \geq \theta_0(\kappa\rho^*)$.
So Statement 9(c)(iv) holds.
This proves the lemma.
\end{pf}

\bigskip

We call the buffer system $\{\C_S: (S,\Omega) \in \Se_2\}$ mentioned in Lemma \ref{gauge} an {\it $(H,\theta_0,\mu)$-gauge} with respect to $(\Se_1,\Se_2,\alpha)$.
Note that the numbers $\rho^*,\lambda^*$ mentioned in Lemma \ref{gauge} only depend on $H,\Sigma,\theta_0,\mu,\kappa,\rho$, and are independent of $\Se$.
For every homeomorphic embedding $\pi$ from $H$ into $G$, let $X_\pi$ be the set mentioned in Statement \ref{gague_X} in Lemma \ref{gauge}, and for every $x \in X_\pi$, let $N_x$ be the graph and $Y_x$ be the set mentioned in Statement \ref{gague_X_Nx} in Lemma \ref{gauge}, and we define $(H_x,\Omega_{H_x})$ to be the labelled rooted graph as follows.
	\begin{itemize}
		\item $V(H_x)=\{v \in V(H): \pi(v) \in V(N_x)\} \cup Y_x$,
		\item $E(H_x) = \{e \in E(H): \pi(e) \subseteq N_x\} \cup \{yv: y \in Y_x, v \in V(H_x)-Y_x,$ there exists an edge $e'$ of $H$ incident with $v$ such that some path in $N_x$ between $y$ and $\pi(v)$ is a subset of $\pi(e')\}$, 
		\item the vertices that are labelled are the vertices in $V(H_x)-Y_x$, and for each $v \in V(H_x)-Y_x$, $v$ is a vertex of $H$ and we label $v$ with $v$, and
		\item $\Omega_{H_x}$ is the march that orders vertices in $Y_x$ in a way that is consistent with a natural ordering of $V(C_x)$, where $C_x$ is the cycle mentioned in Statement \ref{gague_X} in Lemma \ref{gauge}.
	\end{itemize}
Note that elements in $\{(H_x,\Omega_{H_x}): x \in X_\pi\}$ are pairwise non-isomorphic because of the labelling.
We treat the set $\{(H_x,\Omega_{H_x}): x \in X_\pi\}$ as an isomorphic class; that is, for distinct $\pi',\pi''$, we say $\{(H_x,\Omega_{H_x}): x \in X_{\pi'}\}$ and $\{(H_x,\Omega_{H_x}): x \in X_{\pi''}\}$ are equal if there exists a bijection $\phi$ between these two sets such that $Z$ is isomorphic to $\phi(Z)$ for each $Z$.
For each $(S,\Omega) \in \Se_2$, let $C_S,C_S'$ be the cycles mentioned in Lemma \ref{gauge}.
Let $\pi^*$ be a lifting of $\pi$ with respect to $\{C_S: (S,\Omega) \in \Se_2\}$. 
Note that the set $\{\Lambda_x: x \in X_{\pi}\}$, denoted by $\A_{\pi^*}$, is an addendum of $\pi^*$ by Lemma \ref{gauge}, where $\Lambda_x$ is defined in Statement \ref{gague_X} in Lemma \ref{gauge}.
And for each $x \in X_{\pi}$, $(H_x,\Omega_{H_x})$ is isomorphic to the $\Lambda_x$-caption of $\pi^*$.
For each $(S,\Omega) \in \Se_2$, let $\Lambda_S$ be the disk mentioned in Lemma \ref{gauge}.
Fix an ordering $\sigma_{\Se_2}$ of the members of $\Se_2$.
Then for each $(S,\Omega) \in \Se_2$, there exists a $\Lambda_S$-template $(Q_S,s_S)$ such that $((Q_{S'},s_{S'}): (S',\Omega') \in \Se_2)$, where the sequence is ordered according to $\sigma_{\Se_2}$, is a subsequence of a $(\Sigma,(\Lambda_{S'}: (S',\Omega') \in \Se_2))$-template which is a minimal $(\Sigma,(\Lambda_{S'}: (S',\Omega') \in \Se_2),\A_{\pi^*})$-signature.
We define the {\it $(H,\theta_0,\mu)$-snapshot} of $\pi$ with respect to this $(H,\theta_0,\mu)$-gauge and $\sigma_{\Se_2}$ to be the sequence \linebreak $\big(((Q_S,s_S): (S,\Omega) \in \Se_2), (C_S: (S,\Omega) \in \Se_2),(C_S': (S,\Omega) \in \Se_2), \{(H_x,\Omega_{H_x}): x \in X_\pi\} \big)$.

\section{Rerouting} \label{sec: rerouting}

We will prove two lemmas about rerouting paths or curves in this section.
They will be used in later sections when we connect small pieces of homeomorphic embeddings to construct a half-integral packing.

We say that $(S,\Omega, \Omega_0)$ is a {\it neighborhood} if $S$ is a graph and $\Omega, \Omega_0$ are cyclic permutations with $\overline{\Omega}, \overline{\Omega_0} \subseteq V(S)$.
A neighborhood $(S,\Omega, \Omega_0)$ is {\it rural} if $S$ has a drawing $\Gamma$ in the plane and there are disks $\Delta_0 \subseteq \Delta$ such that 
	\begin{itemize}
		\item $\Gamma$ uses no point outside $\Delta$ and none in the interior of $\Delta_0$,  
		\item $\overline{\Omega}$ are the vertices in $\Gamma \cap \partial\Delta$, and $\overline{\Omega_0}$ are the vertices in $\Gamma \cap \Delta_0$, and 
		\item the cyclic permutations of $\overline{\Omega}$ and $\overline{\Omega_0}$ coincide with the natural cyclic orders on $\Delta$ and $\Delta_0$, respectively.
	\end{itemize}
In this case, we say that $(\Gamma, \Delta, \Delta_0)$ is a {\it presentation} of $(S,\Omega, \Omega_0)$.
For a fixed presentation $(\Gamma, \Delta, \Delta_0)$ of a neighborhood $(S,\Omega, \Omega_0)$ and an integer $s \geq 0$, an {\it $s$-nest} for $(\Gamma, \Delta, \Delta_0)$ is a sequence $(C_1, C_2, ..., C_s)$ of pairwise disjoint cycles of $S$ such that $\Delta_0 \subseteq \Delta_1 \subseteq ... \subseteq \Delta_s \subseteq \Delta$, where $\Delta_i$ is the closed disk in the plane bounded by $C_i$, for $1 \leq i \leq s$, in the drawing $\Gamma$, and $\overline{\Delta_0} \cap \partial\Delta_1=\emptyset$.

If $(S,\Omega, \Omega_0)$ is a neighborhood and $(S_0,\Omega_0)$ is a society, then $(S \cup S_0, \Omega)$ is a society and we call this society the {\it composition} of the society $(S_0,\Omega_0)$ with the neighborhood $(S,\Omega, \Omega_0)$.
A society $(S,\Omega)$ is {\it $s$-nested} if it is the composition of a society with a rural neighborhood that has an $s$-nest for some presentation of it.

A subgraph $F \subseteq S$ for a rural neighborhood $(S,\Omega,\Omega_0)$ with presentation $(\Gamma,\Delta,\Delta_0)$ is {\it perpendicular} to an $s$-nest $(C_1,C_2,...,C_s)$ for $(\Gamma,\Delta,\Delta_0)$ if for every component $P$ of $F$
	\begin{itemize}
		\item $P$ is a path with one end in $\overline{\Omega}$ and one end in $\overline{\Omega_0}$, and 
		\item $P \cap C_i$ is a path for each $i=1,2,...,s$.
	\end{itemize}

A {\it linkage} $L$ in a graph $G$ is a subgraph of $G$ such that every component of $L$ is a path.
We say that two linkages $L,L'$ in $G$ are {\it equivalent} if for every two vertices $u,v$ of $G$, $u,v$ are the ends of a component of $L$ if and only if $u,v$ are the ends of a component of $L'$.
We say that a linkage $L$ in $G$ is {\it vital} in $G$ if $V(L)=V(G)$, and there exists no linkage $L'$ in $G$ with $L' \neq L$ equivalent to $L$.

\begin{theorem}[{\cite[(1.1)]{rs XXI}}] \label{unique linkage}
For every positive integer $p$, there exists an integer $w$ such that every graph that has a vital linkage with $p$ components has tree-width less than $w$.
\end{theorem}

We shall prove the following lemma, where its proof included in this paper is inspired by a proof of a theorem in \cite{kntw}.

\begin{lemma} \label{perpendicular}
For positive integers $k,s$, there exists an integer $s'$ with $s' \geq s$ such that the following holds.
If $(S,\Omega)$ is a composition of a society $(S_0,\Omega_0)$ with an $s'$-nested rural neighborhood $(S',\Omega,\Omega_0)$ that has an $s'$-nest $(C_1,C_2,...,C_{s'})$ for some presentation $(\Gamma,\Delta,\Delta_0)$, and $F$ is a linkage in $S$ with at most $k$ components such that every component of $F$ is either
	\begin{itemize}
		\item a path $P$ with both ends in $\overline{\Omega}$ such that there exists another component $P'$ of $F$ such that $P'$ is a path with both ends in $\overline{\Omega}$ but belonging to different intervals of $\Omega$ determined by the ends of $P$, or
		\item a path with one end in $\overline{\Omega}$ and one end in $V(S_0)$,
	\end{itemize}
then there exist a linkage $F'$ in $S$ equivalent to $F$ and a society $(S_0',\Omega_0')$ with $S_0 \subseteq S_0'$ such that 
	\begin{enumerate}
		\item $(S,\Omega)$ is a composition of $(S_0',\Omega_0')$ with an $s$-nested rural neighborhood $(S'',\Omega,\Omega_0')$ that has an $s$-nest $(C_1',C_2',...,C_s')$ for some its presentation such that $F' \cap S''$ is perpendicular to $(C_1',C_2',...,C_s')$, and
		\item $S'' \cap S_0'$ is a cycle with vertex-set $\overline{\Omega_0'}$ passing through $\overline{\Omega_0'}$ in the order $\Omega_0'$.
	\end{enumerate}
\end{lemma}

\begin{pf}
Let $w$ be the number $w$ mentioned in Theorem \ref{unique linkage} by taking $p=k$.
Define $s'=2w+4+s$.

Suppose that $S$ and $F$ form a counterexample of this lemma with $\lvert V(S) \rvert + \lvert E(S) \rvert$ minimum.
By deleting vertices and the minimality of $S$, we know that $V(S)=V(F \cup \bigcup_{i=1}^{s'}C_i)$.
Furthermore, by contracting edges in $E(F) \cap E(\bigcup_{i=1}^{s'}C_i)$ and the minimality of $S$, we know that $E(F) \cap E(\bigcup_{i=1}^{s'}C_i)=\emptyset$.
If $V(S)-V(F) \neq \emptyset$, say $v \in V(S)-V(F)$, then $v \in (\bigcup_{i=1}^{s'}V(C_i))-V(F)$, and contracting an edge in $\bigcup_{i=1}^{s'}C_i$ incident with $v$ contradicts the minimality of $S$.
So $V(S)=V(F)$.

\noindent{\bf Claim 1:} The tree-width of $S$ is less than $w$.

\noindent{\bf Proof of Claim 1:}
Suppose that the tree-width of $S$ is at least $w$.
Then by Theorem \ref{unique linkage}, $F$ is not a vital linkage in $S$.
So there exists a linkage $L$ in $S$ with $L \neq F$ equivalent to $F$.
Hence there exists an edge $e \in E(F)-E(L)$.
Since $E(F) \cap E(\bigcup_{i=1}^{s'}C_i)=\emptyset$, $e \not \in E(\bigcup_{i=1}^{s'}C_i)$.
So $(C_1,C_2,...,C_{s'})$ is an $s'$-nest for $S-e$, and $L \subseteq S-e$.
Since $L$ is equivalent to $F$ in $S$, the minimality of $S$ implies that conclusion of this lemma holds for $S-e$ and hence in $S$, a contradiction. 
So the tree-width of $S$ is less than $w$.
$\Box$

Let $W$ be a component of $F \cap S'$.
By definition, $W$ is a path.
If there exists a path $R$ in $S$ with ends $x,y$ such that $\{x,y\} \subseteq V(W)$ and $R \cap V(F)=\{x,y\}$, then we say the graph $F'$ is obtained from $F$ by {\it rerouting through $R$} if $F'$ is obtained from $F$ by deleting the edges and internal vertices of the subpath in $W$ with ends $x,y$ and adding $R$.
Note that $F'$ is a linkage in $S$ equivalent to $F$.
If such $R$ exists, we say that $F$ can be {\it rerouted through $R$}.

If there exists a path $R$ in $C_i$ for some $i \in [s']$ such that $F'$ can be obtained from $F$ by rerouting through $R$, then one can delete the edges in $F-F'$ from $S$ and apply induction, contradicting the minimality of $S$.
So $F$ cannot be rerouted through any path in $C_i$ for any $i \in [s']$.

\noindent{\bf Claim 2:} If $Z$ is a component of $F \cap S'$ such that $Z$ intersects $C_i$ for some $i$ with $2w+3 \leq i \leq s'$, then $Z$ has exactly one end in $\overline{\Omega}$. 

\noindent{\bf Proof of Claim 2:}
We say that any path in $F \cap S'$ from $\overline{\Omega_0}$ to $\overline{\Omega_0}$ with at least one edge is a {\it bounce}.
The {\it height} of a bounce $P$ is the maximum $h$ such that $V(P) \cap V(C_h) \neq \emptyset$.
Let $\Delta_0$ be the disk in which $S_0$ is located.
Note that $P \cup \Delta_0$ bounds a disk containing $\Delta_0$, and we denote this disk by $\Delta_P$.

If a bounce $P$ has height $h$ for some $h \geq 2$, then $\lvert V(P) \cap V(C_{h-1}) \rvert \geq 2$ by the planarity of the rural neighborhood, and there exists a bounce $Q$ disjoint from $P$ such that $\Delta_Q \subseteq \Delta_P$ and $V(Q) \cap V(C_{h-1}) \neq \emptyset$; otherwise there exists a path $R$ in $C_{h-1}$ with ends in $V(P) \cap V(C_{h-1})$ such that rerouting $F$ through $R$ contradicts the minimality of $S$.
In other words, if there exists a bounce $P_h$ with height $h$ for some $h \geq 2$, then there exists a bounce $P_{h-1}$ disjoint from $P_h$ with height at least $h-1$ such that $\Delta_{P_{h-1}} \subseteq \Delta_{P_h}$.
Therefore, if there exists a bounce with height at least $2w+3$, then there exist bounces $P_i$ for every $i \in [w+1]$ such that the following hold.
	\begin{itemize}
		\item $P_i$ are pairwise disjoint for $i \in [w+1]$.
		\item $\Delta_{P_j} \subseteq \Delta_{P_{j+1}}$ for every $j \in [w]$.
		\item $\lvert V(P_i) \cap V(C_{2i}) \rvert \geq 2$ for every $i \in [w+1]$.
	\end{itemize}
Since $E(F) \cap E(\bigcup_{i=1}^{s'}C_i)=\emptyset$, if there exists a bounce with height at least $2w+3$, then $\{V(P_i) \cup V(C_{2i}): i \in [w+1]\}$ is a set of $w+1$ edge-disjoint but pairwise intersecting subgraphs of $S$, so $S$ has tree-width at least $w$, contradicting Claim 1.
Hence every bounce has height at most $2w+2$.
That is, if $P$ is a component of $F \cap S'$ such that $P$ intersects $C_i$ for some $i$ with $2w+3 \leq i \leq s'$, then $P$ has at least one end in $\overline{\Omega}$.
Observe that no component of $F \cap S'$ contains both ends in $\overline{\Omega}$ by the definition of $F$ and the planarity of a rural neighborhood.
So if $P$ is a component of $F \cap S'$ such that $P$ intersects $C_i$ for some $i$ with $2w+3 \leq i \leq s'$, then $P$ has exactly one end in $\overline{\Omega}$. 
$\Box$

For any linkage $L$ in $S$ equivalent to $F$, we say that an $(s+1)$-nest $(C_0',C_1',C_2',...,C_s')$ for $(\Gamma,\Delta,\Delta_0)$ is {\it $L$-compatible} if any component of $L \cap S'$ that intersects $\bigcup_{i=0}^s C_i'$ is a path with one end in $\overline{\Omega}$ and one end in $V(S_0)$. 
Note that $(C_{2w+3},C_{2w+4},C_{2w+5},..., \allowbreak C_{2w+3+s})$ is $F$-compatible by Claim 2.

We choose a linkage $F'$ in $S$ equivalent to $F$ and an $F'$-compatible $(s+1)$-nest $(C_0',C_1',C_2', \allowbreak ..., \allowbreak C_s')$ such that the number of components of $F' \cap (\bigcup_{i=0}^s C_s')$ is minimum.

\noindent{\bf Claim 3:} For every component $Q$ of $F' \cap S'$ and $i \in [s] \cup \{0\}$, $Q \cap C_i'$ is a path.

\noindent{\bf Proof of Claim 3:}
Suppose to the contrary that there exist a component $Q$ of $F' \cap S'$ and $i \in [s] \cup \{0\}$ such that $Q \cap C_i'$ is not a path.
So there exist $x,y \in V(Q) \cap V(C_i')$ belonging to different components of $Q \cap C_i'$, the subpath $R$ of $Q$ with ends $x,y$, and a path $P \subseteq C_i'$ with ends $x,y$ such that $R \cup P$ bounds an open disk disjoint from $\Delta_0$.
We further choose the pair $Q$ and $i$ such that $i$ is maximum, and subject to this, $R$ is internally disjoint from the closure of the disk $\Delta_i'$ bounded by $C_i'$ if possible.

Since $Q$ is a path with one end in $\overline{\Omega}$ and one end in $V(S_0)$, the maximality of $i$ implies that $R$ is disjoint from $C'_j$ for every $j \geq i+1$.
If $R$ is internally disjoint from the closure of $\Delta_i'$, then define $C_i''$ to be the cycle obtained by $C_i'$ by deleting all edges and internal vertices of $P$ and adding $R$, and define $C_j''=C_j'$ for $j \in ([s] \cup \{0\})-\{i\}$.
Then $(C_0'',C_1'',...,C_s'')$ is $F'$-compatible, but the number of components of $F' \cap (\bigcup_{i=0}^s C_i'')$ is smaller than the number of components of $F' \cap (\bigcup_{i=0}^s C_i')$, a contradiction.

Hence $R$ is not internally disjoint from the closure of $\Delta_i'$.
So $R$ is contained in the closure of $\Delta_i'$.
If some internal vertex of $P$ belongs to $V(F')$, then there exists a component $Q'$ of $F' \cap S'$ (possibly equal to $Q$) such that $Q' \cap C_i'$ is not a path and the corresponding subpath $R'$ of $Q'$ is different from $R$ and is internally disjoint from the closure of $\Delta_i'$ by the planarity of $(\Gamma,\Delta,\Delta_0)$, contradicting the choice of $Q,i$.
So no internal vertex of $P$ belongs to $V(F')$.
Define $F''$ to be the linkage in $S$ obtained by rerouting $F'$ through $P$.
Then $(C_0',C_1',...,C_s')$ is $F''$-compatible, but the number of components of $F'' \cap (\bigcup_{i=0}^s C_i')$ is less than the number of components of $F' \cap (\bigcup_{i=0}^s C_i')$, a contradiction.
This proves the claim.
$\Box$

Let $\Delta_0'$ be the open disk bounded by $C_0'$.
Let $\Omega_0'$ be a cyclic permutation of $V(C_0')$ consistent with the natural cyclic order of $\partial \Delta_0'$.
Let $(S'',\Omega,\Omega_0')$ be the rural neighborhood that has presentation $(\Gamma-\Delta_0', \Delta, \Delta_0')$.
So $(C_1',C_2',...,C_s')$ is an $s$-nest for some presentation of $(S'',\Omega,\Omega_0')$.
Let $S_0'$ be the minimal subgraph of $S$ such that $S_0' \supseteq C_0'$ and $(S,\Omega)$ is a composition of $(S_0',\Omega_0')$ with $(S'',\Omega,\Omega_0')$.
Since $S'' \subseteq S'$, $S_0' \supseteq S_0$.
Note that $S'' \cap S_0'$ is a cycle with vertex-set $\overline{\Omega_0'}$ passing through $\overline{\Omega_0'}$ in the order $\Omega_0'$.

It suffices to show that $F' \cap S''$ is perpendicular to $(C_1',C_2',...,C_s')$.

\noindent {\bf Claim 4:} Every component of $F' \cap S''$ intersects $\bigcup_{i=0}^sC_i'$.

\noindent{\bf Proof of Claim 4:}
Suppose to the contrary that there exists a component $Q$ of $F' \cap S''$ disjoint from $\bigcup_{i=0}^sC_i'$.
In particular, $Q$ is disjoint from $V(C_0')=\overline{\Omega_0'}$.
Since $F'$ is equivalent to $F$ and $\Delta_0' \supseteq \Delta_0$, $Q$ is a path with at least one end in $\overline{\Omega}$.
Since $Q$ is disjoint from $C'_s$, $Q$ is a path with both ends in $\overline{\Omega}$ disjoint from $\bigcup_{i=0}^sC_i'$.
Since $F'$ is equivalent to $F$, there exists another component $Q''$ of $F'$ such that $Q''$ is a path with both ends in $\overline{\Omega}$ belonging to different intervals of $\Omega$ determined by the ends of $Q$, contradicting the planarity.
$\Box$

Let $Q$ be a component of $F' \cap S''$.
Since $S'' \subseteq S'$, $Q$ is contained in some component $Q'$ of $F' \cap S'$.
By Claim 4, $Q' \supseteq Q$ intersects $\bigcup_{i=0}^sC_i'$.
Since $(C_0',C_1',...,C'_s)$ is $F'$-compatible, $Q'$ is a path with one end in $\overline{\Omega}$ and one end in $V(S_0)$, and by Claim 3, $Q' \cap C_i'$ is a path for every $i \in [s] \cup \{0\}$, so $Q$ is a path with one end in $\overline{\Omega}$ and one end in $\overline{\Omega_0'}$ such that $Q \cap C_i'$ is a path for every $i \in [s]$.
Therefore, $F' \cap S''$ is perpendicular to $(C_1',C_2',...,C_s')$.
\end{pf}

\bigskip

Let $\Sigma$ be a connected surface, and let $\Delta_1, ..., \Delta_t$ be pairwise disjoint closed disks in $\Sigma$.
Let $\Gamma$ be a drawing in $\Sigma$ such that $U(\Gamma) \cap \Delta_i = V(\Gamma) \cap \partial\Delta_i$ for $1 \leq i \leq t$.
Let $Z = \bigcup_{i=1}^t V(\Gamma) \cap \partial\Delta_i$.
We say that a partition $(Z_1, Z_2, ..., Z_p)$ of $Z$ satisfies the {\it topological feasibility condition} if there exist pairwise disjoint disks $D_1, D_2, ..., D_p$ in $\Sigma$ such that $D_j \cap (\bigcup_{i=1}^t \Delta_i) = Z_j$ for $j \in [p]$.

\begin{theorem}[{\cite[Theorem~(3.2)]{rs XII}}] \label{linkage on surface}
	For every connected surface $\Sigma$ and all integers $t \geq 0$ and $z \geq 0$, there exists a positive integer $\theta$ such that the following is true.
	Let $\Delta_1, ..., \Delta_t$ be pairwise disjoint closed disks in $\Sigma$, and let $\Gamma$ be a $2$-cell drawing in $\Sigma$ such that $U(\Gamma) \cap \Delta_i = V(\Gamma) \cap \partial\Delta_i$ for $1 \leq i \leq t$.
	Let $\lvert Z \rvert \leq z$, where $Z = \bigcup_{i=1}^t (V(\Gamma) \cap \partial\Delta_i)$, and let $(Z_1, Z_2,...,Z_p)$ be a partition of $Z$ satisfying the topological feasibility condition.
	Let $\T$ be a respectful tangle of order at least $\theta$ in $\Gamma$ with metric $m_\T$ such that $m_\T(r_i,r_j) \geq \theta$ for $1 \leq i< j \leq t$, where $r_i$ is the region of $\Gamma$ meeting $\Delta_i$ for $1 \leq i \leq t$, and $V(\Gamma) \cap \partial\Delta_i$ is free for $1 \leq i \leq t$.
	Then there are mutually disjoint connected subdrawings $\Gamma_1, \Gamma_2, ..., \Gamma_p$ of $\Gamma$ with $V(\Gamma_j) \cap Z = Z_j$ for $1 \leq j \leq p$.
\end{theorem}

\begin{lemma} \label{bounding crossings}
For every positive integer $k$ and for every surface $\Sigma$, there exists an integer $\tau$ such that if $\Delta_1,...,\Delta_k$ are open disks in $\Sigma$ with pairwise disjoint closure and $\Gamma_1,\Gamma_2,...,\Gamma_k$ are $k$ sets of curves in $\Sigma-\bigcup_{i=1}^k\Delta_i$ such that for each $i \in [k]$, $\lvert \Gamma_i \rvert \leq k$ and members of $\Gamma_i$ are pairwise disjoint curves connecting pairs of distinct points in $\bigcup_{j=1}^k\partial\Delta_j$, then there exist $k$ sets $\Gamma_1',\Gamma_2',...,\Gamma'_k$ of simple curves in $\Sigma-\bigcup_{i=1}^k\Delta_i$ such that the following hold.
	\begin{enumerate}
		\item For each $i \in [k]$, $\lvert \Gamma_i' \rvert = \lvert \Gamma_i \rvert$ and curves in $\Gamma_i'$ are pairwise disjoint.
		\item For each $i \in [k]$ and $\gamma \in \Gamma_i$, there uniquely exists a curve $\gamma' \in \Gamma_i'$ such that the endpoints of $\gamma$ are the same as the ends of $\gamma'$.
		\item For distinct members $\gamma,\gamma'$ of $\bigcup_{i=1}^k\Gamma_i'$ and for each point $x \in \gamma \cap \gamma'$, there exists an open set $B_x \subseteq \Sigma$ with $x \in B_x$ such that $B_x-\{x\}$ does not contain any point that belongs to at least two distinct members of $\bigcup_{i=1}^k\Gamma_i'$.
		\item Every point in $\Sigma$ belongs to at most two members of $\bigcup_{i=1}^k\Gamma'_i$.
		\item There exist at most $\tau$ points in $\Sigma$ belonging to at least two members of $\bigcup_{i=1}^k\Gamma_i'$.
	\end{enumerate}
\end{lemma}

\begin{pf}
Let $k$ be a positive integer, and let $\Sigma$ be a surface.
We may assume that $\Sigma$ is connected since we can prove this lemma componentwisely.
Let $\theta$ be the number $\theta$ mentioned in Theorem \ref{linkage on surface} by taking $\Sigma=\Sigma$, $t=k$ and $z=2k^2$.
Let $\Delta_1,...,\Delta_k$ be $k$ open disks with pairwise disjoint closure.
It is not hard to see that there exists a number $g$ (only depending on $\Sigma$ and $k$) such that one can construct a 2-cell drawing $G$ in $\Sigma$ with $\lvert V(G) \rvert =g$ such that $U(G) \cap \Delta_i=V(G) \cap \partial\Delta_i$ and $\lvert V(G) \cap \partial\Delta_i \rvert \geq k^2$ for every $i \in [k]$, and there exists a respectful tangle $\T$ in $G$ of order at least $\theta$ with metric $m_\T$ such that $m_\T(r_i,r_j) \geq \theta$ for $1 \leq i<j \leq k$, where $r_i$ is the region of $G$ meeting $\Delta_i$ for $1 \leq i \leq k$, and every subset of $V(G) \cap \partial\Delta_i$ with size at most $2k^2$ is free with respect $\T$ for every $1 \leq i \leq k$.
Define $\tau=k^4g$.

Let $\Gamma_1,...,\Gamma_k$ be $k$ sets of curves mentioned in the statement of this lemma.
By moving the curves in $\bigcup_{i=1}^k\Gamma_i$ homotopically, it suffices to prove the case that the endpoints of any curve in $\bigcup_{i=1}^k\Gamma_i$ are in $V(G)$.
For each $i \in [k]$, let $Z_i$ be the set of the endpoints of the members of $\Gamma_i$, and let $(Z^i_\gamma: \gamma \in \Gamma_i)$ be the partition of $Z_i$ such that $Z^i_\gamma$ consists of the endpoints of $\gamma$ for every $\gamma \in \Gamma_i$.
Since for each $i \in [k]$, members of $\Gamma_i$ are pairwise disjoint, the partition $(Z^i_\gamma: \gamma \in \Gamma_i)$ satisfies the topological feasibility condition.
So by Theorem \ref{linkage on surface}, for each $i \in [k]$, there exist pairwise disjoint subdrawings $\Gamma^i_\gamma$ in $\Sigma$ such that $V(\Gamma^i_\gamma) \cap Z_i = Z^i_\gamma$.
Since each $Z^i_\gamma$ has size 2, we may assume that each $\Gamma^i_\gamma$ is a path between the two vertices in $Z^i_\gamma$.
Let ${\Gamma^i_\gamma}'$ be the curve in $\Sigma$ with ${\Gamma^i_\gamma}'=U(\Gamma^i_\gamma)$. 
For each $i \in [k]$, let $\Gamma_i' = \{{\Gamma^i_\gamma}': \gamma \in \Gamma_i\}$.
Hence every point $x \in \Sigma$ belonging to at least two members of $\bigcup_{i=1}^k\Gamma_i'$ belongs to $U(G)$. 
By slightly perturbing members of $\bigcup_{i=1}^k\Gamma_i'$, we may assume that every point $x \in \Sigma$ belonging to at least two members of $\bigcup_{i=1}^k\Gamma_i'$ belongs to $V(G)$. 
So $\Gamma_1',...,\Gamma_k'$ satisfy Conclusions 1-3, and there are at most $\lvert V(G) \rvert$ points in $\Sigma$ belonging to at least two members of $\bigcup_{i=1}^k\Gamma_i'$.
By perturbing the members of $\bigcup_{i=1}^k\Gamma_i'$, we may assume that $\Gamma_1',...,\Gamma_k'$ satisfy Conclusions 1-4, and there are at most $(\sum_{i=1}^k\lvert \Gamma_i \rvert)^2 \lvert V(G) \rvert \leq \tau$ points in $\Sigma$ belonging to at least two distinct members of $\bigcup_{i=1}^k\Gamma_i'$.
This proves the lemma.
\end{pf}

\section{Reconnecting} \label{sec: vortices}

The goal of this section is to prove Lemma \ref{arrangement EP} which deals with the half-integral Erd\H{o}s-P\'{o}sa property when $H$ is connected and $G$ has a segregation with a proper arrangement in a surface.

Let $\rho$ be a positive integer.
Let $(S,\Omega)$ be a vortex with depth at most $\rho$, and let $(P,\X)$ be a vortical decomposition of $(S,\Omega)$ with adhesion at most $\rho$.
We denote $V(P)$ by $[\lvert V(P) \rvert]$ and any edge of $P$ is between two consecutive integers.
Let $W$ be a subgraph of $S$.
Define $\I_W$ to be a collection of closed intervals in the real line as follows.
\begin{itemize}
	\item Each member of $\I_W$, denoted by $I_C$, corresponds to a component $C$ of $W$.
	\item For each component $C$ of $W$, $I_C=[a_C,b_C]$, where $a_C$ is the minimum element $i \in [\lvert V(P) \rvert]$ such that $V(C) \cap X_i \neq \emptyset$, and $b_C$ is the maximum $i \in [\lvert V(P) \rvert]$ such that $V(C) \cap X_i \neq \emptyset$. 
\end{itemize}
We call $\I_W$ the {\it canonical interval labeling of $W$ with respect to $(P,\X)$}.

For a positive integer $t$, a {\it $t$-interval} is a union of $t$ (not necessarily disjoint) closed intervals in ${\mathbb R}$.

\begin{theorem}[{\cite[Theorem 1]{a}}] \label{t-interval stable}
Let $t,k$ be positive integers.
If $\Se$ is a finite family of $t$-intervals with no $k+1$ pairwise disjoint members, then there exists a set of at most $2t^2k$ points in ${\mathbb R}$ intersecting all members of $\Se$.
\end{theorem}

Recall that a segregation $\Se$ of a graph $G$ with a $(\kappa,\rho)$-witness is effective with respect to $(\Se_1,\Se_2)$ if for every $(S,\Omega) \in \Se_1$ and every $v \in \overline{\Omega}$, there exist $\lvert \overline{\Omega} \rvert-1$ paths in $S$ from $v$ to $\overline{\Omega}-\{v\}$ such that the intersection of them is $\{v\}$.
The following lemma is the key step toward proving the main result of this section (Lemma \ref{arrangement EP}).

\begin{lemma} \label{same imprint}
For every connected graph $H$, positive integers $k,\kappa,\rho$, and every surface $\Sigma$, there exist nonnegative nondecreasing functions $\theta_0^*,\mu^*$ with domain ${\mathbb Z}$ and integers $\theta^*,\xi^*,\lambda^*,\eta^*$ such that the following hold.
Let $G$ be a graph and $\T$ a tangle in $G$.
Let $\Se$ be an effective $\T$-central segregation of $G$ with respect to a $(\kappa,\rho)$-witness $(\Se_1,\Se_2)$.
Let $\alpha$ be a proper arrangement of $\Se$ in $\Sigma$ such that there exists a tangle $\T'$ of order at least $\theta^*$ in the skeleton $G'$ of $\alpha$ with respect to $(\Se_1,\Se_2)$ conformal with $\T$, and there exists an $(H,\theta_0^*,\mu^*)$-gauge $\C$ with respect to $(\Se_1,\Se_2,\alpha)$.
Let $\pi_1,\pi_2,...,\pi_\ell$ be distinct homeomorphic embeddings from $H$ into $G$.
Let $\sigma_{\Se_2}$ be an ordering of the members of $\Se_2$.
If the $(H,\theta_0^*,\mu^*)$-snapshots of $\pi_i$ with respect to $\C$ and $\sigma_{\Se_2}$ are identical for all $i \in [\ell]$, then one of the following statement holds.
	\begin{enumerate}
		\item $G$ half-integrally packs $k$ subdivisions of $H$, where for each $v \in V(H)$, their branch vertices corresponding to $v$ are contained in $\{\pi_i(v): i \in [\ell]\}$.
		\item There exists $Z \subseteq V(G)$ with $\lvert Z \rvert \leq \xi^*$ such that $Z \cap \pi_i(V(H) \cup E(H)) \neq \emptyset$ for each $i \in [\ell]$.
		\item There exist $W^* \subseteq V(G')$ with $\lvert W^* \rvert \leq \eta^*$ and a set $\{\Lambda_w^*: w \in W^*\}$, where each $\Lambda_w^*$ is a $\lambda^*$-zone around $w \in W^*$ in $G'$, such that for each $i \in [\ell]$, $\bigcup_{(S,\Omega)\in\Se_1, \alpha(S,\Omega) \subseteq \bigcup_{w \in W^*}\overline{\Lambda_w^*}}V(S)$ contains $\pi_i(v_i)$ for some $v_i \in V(H)$ with $\pi_i(v_i) \not \in V(\bigcup_{(S,\Omega) \in \Se_2}S)$.
		\item There exist $i \in [\ell]$ and $(A,B) \in \T$ such that $\pi_i(E(H)) \subseteq A$.
	\end{enumerate}
\end{lemma}

\begin{pf}
Let $H$ be a connected graph, $k,\kappa,\rho$ be positive integers, and $\Sigma$ be a surface.
Let $q_0=k$, and for $i \geq 1$, let $q_i=k+\sum_{j=0}^{i-1}q_j$.
For every nonnegative integer $x$, we define the following. 
	\begin{itemize}
		\item Let $r(x)$ be the $(H,\Sigma,\kappa,x)$-port number.
		\item Let $\tau_{\ref{bounding crossings}}(x)$ be the number $\tau$ mentioned in Lemma \ref{bounding crossings} by taking $k=(\kappa \cdot r(x) + 2\lvert E(H) \rvert \cdot \lvert V(H) \rvert q_{\lvert V(H) \rvert}+\kappa) \cdot k$ and $\Sigma=\Sigma$.
		\item Let $k_0(x)= 5\tau_{\ref{bounding crossings}}(x)$.
		\item Let $\phi(x)$ be the number $\theta$ mentioned in Theorem \ref{linkage on surface} by taking $\Sigma=\Sigma$, $t=\kappa+q_{\lvert V(H) \rvert}\lvert V(H) \rvert+k_0(x)$ and $z=2k\kappa r(x) + 2\lvert E(H) \rvert \lvert V(H) \rvert q_{\lvert V(H) \rvert} + 2k_0(x)+2$. 
			Note that we may assume that $\phi(x) > 2k\kappa r(x) + 2\lvert E(H) \rvert \lvert V(H) \rvert q_{\lvert V(H) \rvert} + 2k_0(x)+2$. 
		\item Let $\theta'(x)=2k_0(x)+\phi(\kappa x)$, 
		\item Let $\xi(x) = 2(r(x)+\lvert V(H) \rvert)^2(k-1)(2x+1)$. 
		\item Let $s_{\ref{perpendicular}}(x)$ be the number $s'$ mentioned in Lemma \ref{perpendicular} by taking $k=k \cdot r(x)$ and $s=x+1$.
			Note that $s_{\ref{perpendicular}}(x) \geq x+1$ by Lemma \ref{perpendicular}.
	\end{itemize}
We define the following.
	\begin{itemize}
		\item Define $\mu^*(\cdot)$ to be the function $s_{\ref{perpendicular}}(2\phi(\cdot))+r(\cdot)+2\lvert E(H) \rvert$.
		\item Define $\theta_0^*(\cdot)$ to be the function $\theta'(\cdot)+\mu^*(\cdot)+4$. 
		\item Let $\lambda_1^*, \rho^*$ be the numbers $\lambda^*,\rho^*$, respectively, mentioned in Lemma \ref{gauge} by taking $H=H, \kappa=\kappa,\rho=\rho, \theta_0=\theta_0^*,\mu=\mu^*,\Sigma=\Sigma$.
		\item Let $k_0^*=k_0(\rho^*)$.
		\item Define $\eta^*=q_{\lvert V(H) \rvert} \lvert V(H) \rvert+\kappa+k_0^*$.
		\item Let $\phi^*=2\phi(\kappa\rho^*)+(4\lambda_1^*+2)(\kappa+q_{\lvert V(H) \rvert}\lvert V(H) \rvert)+2k_0^*$. 
		\item Define $\xi^*=\xi(\rho^*)+q_{\lvert V(H) \rvert}$.
		\item Define $\lambda^*=3\lambda_1^*+\phi^*+(4\lambda_1^*+2)(\kappa+q_{\lvert V(H) \rvert}\lvert V(H) \rvert)$.
		\item Define $\theta^*= \phi^*+\xi^*+(4(\phi^*+\lambda_1^* + (4\lambda_1^*+2)(\kappa+q_{\lvert V(H) \rvert}\lvert V(H) \rvert))+2)\eta^*+2k_0^*+2\rho^*$.
	\end{itemize}

Let $G$, $\T$, $\Se$, $(\Se_1,\Se_2)$, $\alpha$, $\C$, $\pi_1,\pi_2,...,\pi_\ell,\sigma_{\Se_2}$ be the ones mentioned in the statement of this lemma.
Suppose that this lemma does not hold.

For each $i \in [\ell]$ and $(S,\Omega) \in \Se_2$, let $C_{S,i},C_{S,i}',\allowbreak \Lambda_{S,i},X_{\pi_i},\rho_{S,i}$ be the cycles $C_S,C_S'$, sets $\Lambda_S,X_\pi$ and number $\rho_S$, respectively, if we further take $\pi=\pi_i$ in Lemma \ref{gauge}.  
Since the $(H,\theta_0^*,\mu^*)$-snapshots of $\pi_i$ are identical for all $i \in [\ell]$, we know that $C_{S,i}=C_{S,j}$ and $C_{S,i}'=C'_{S,j}$ for all $i,j \in [\ell]$ and $(S,\Omega) \in \Se_2$.
So for every $(S,\Omega) \in \Se_2$, we denote $C_{S,i},C_{S,i}',\Lambda_{S,i}$ and $\rho_{S,i}$ for all $i \in [\ell]$ by $C_S,C_S',\Lambda_S$ and $\rho_S$, respectively.

For each $x \in \bigcup_{i=1}^\ell X_{\pi_i}$, let $Y_x,\Lambda_x,\lambda_x$ be the sets and number as mentioned in the statement of Lemma \ref{gauge} when taking $\pi=\pi_i$.

As the $(H,\theta_0^*,\mu^*)$-snapshots for all $\pi_i$ with respect to $\C$ and $\sigma_{\Se_2}$ are identical, we denote the terms $((Q_S,s_S): (S,\Omega) \in \Se_2)$ and $\{(H_x,\Omega_{H_x}): x \in X_\pi\}$ mentioned in the definition of the $(H,\theta_0^*,\mu^*)$-snapshot of $\pi_i$ for each $i \in [\ell]$ as $((Q_S,s_S): (S,\Omega) \in \Se_2)$ and $\{(H_x,\Omega_{H_x}): x \in X_{\pi_i}\}$, respectively.
Note that there exists a pseudo-embedding $\pi_0$ from $H$ into $\Sigma$ legal with respect to $\{\Lambda_S: (S,\Omega) \in \Se_2\}$ such that $\pi_0$ has an addendum $\A_0$ such that the $(\Sigma, (\Lambda_S: (S,\Omega) \in \Se_2),\A_0)$-signature of $\pi_0$ is a $(\Sigma,(\Lambda_S: (S,\Omega) \in \Se_2))$-template such that for each $(S,\Omega) \in \Se_2$, $(Q_S,s_S)$ is the entry of the $(\Sigma, (\Lambda_S: (S,\Omega) \in \Se_2),\A_0)$-signature of $\pi_0$ corresponding to $(S,\Omega)$, and the set $\{(H_\Delta,\Omega_\Delta): \Delta \in \A_0\}$ of $\Delta$-captions is isomorphic to $\{(H_x,\Omega_{H_x}): x \in X_{\pi_1}\}$. 

For each $i \in [\ell]$ and $(S,\Omega) \in \Se_2$, let $Y_{S,i}$ be the set $Y_S$ if we further take $\pi=\pi_i$, $\pi^*$ to be a lifting of $\pi_i$ with respect to $\{C_{S'}: (S',\Omega') \in \Se_2\}$ defining the $(H_0,\theta_0^*,\mu_0^*)$-snapshot of $\pi_i$, and $(Q_S,s_S)=(Q_S,s_S)$ in Lemma \ref{gauge}.

Let $\F_1=\{(H_x,\Omega_{H_x}): x \in X_{\pi_1}\}$.
For each $F \in \F_1$, define $\J_F$ to be the set $\{Y_x: x\in \bigcup_{i=1}^\ell X_{\pi_i}, (H_x,\Omega_{H_x}) \cong F\}$. 

\noindent{\bf Claim 1:} For each $F \in \F_1$, there exists $\J_F^* \subseteq \J_F$ with $\lvert \J_F^* \rvert = q_{\lvert \F_1 \rvert}$ such that $m_{\T'}(Y,Y') \geq \phi^*$ for every distinct $Y,Y' \in \bigcup_{F' \in \F_1}\J^*_{F'}$.

\noindent{\bf Proof of Claim 1:}
We may assume that $\F_1 \neq \emptyset$, for otherwise this claim holds.
For each $F \in \F_1$, choose $\J_F' \subseteq \J_F$ with $\lvert \J_F' \rvert \leq q_{\lvert \F_1 \rvert}$ such that $m_{\T'}(Y,Y') \geq \phi^*$ for every distinct $Y,Y' \in \bigcup_{F' \in \F_1}\J'_{F'}$.
We further choose those $\J'_F$'s such that $\sum_{F \in \F_1}\lvert \J'_F \rvert$ is as large as possible.
We are done if $\lvert \J'_F \rvert=q_{\lvert \F_1 \rvert}$ for every $F \in \F_1$.
So we may assume that there exists $F^* \in \F_1$ such that $\lvert \J_{F^*}' \rvert <q_{\lvert \F_1 \rvert}$.
By the maximality of $\sum_{F \in \F_1}\lvert \J'_F \rvert$, for every $Y \in \J_{F^*}-\bigcup_{F \in \F_1}\J'_F$, there exists $Y' \in \bigcup_{F \in \F_1}\J'_F$ such that $m_{\T'}(Y,Y') \leq \phi^*-1$.

By the definition of a gauge, for each $Y \in \bigcup_{F \in \F_1}\J_{F}$, there exist $x_Y \in V(G')$ with $m_{\T'}(x_Y,Y) \leq \lambda_1^*$ and a $\lambda_1^*$-zone $\Lambda_{x_Y}$ around $x_Y$ such that $\bigcup_{(S,\Omega) \in \Se_1, \alpha(S,\Omega) \subseteq \overline{\Lambda_{x_Y}}}V(S)$ intersects $\bigcup_{i=1}^\ell \pi_i(V(H))$. 
Let $W^*=\{x_Y: Y \in \bigcup_{F \in \F_1}\J_{F}'\}$.
So $\lvert W^* \rvert \leq q_{\lvert \F_1 \rvert}\lvert \F_1 \rvert \leq q_{\lvert V(H) \rvert}\lvert V(H) \rvert$.
Note that for every $Y \in \J_{F^*}$, there exists $w \in W^*$ such that $m_{\T'}(w,Y) \leq \lambda^*_1+\phi^*-1$.
So for each $Y \in \J_{F^*}$, $m_{\T'}(W^*,x_Y) \leq 2\lambda^*_1+\phi^*-1$.
By Lemma \ref{big zone contains ball}, for each $w \in W^*$, there exists a $(3\lambda^*_1+\phi^*+2)$-zone $\Lambda^*_w$ containing every atom $x$ with $m_{\T'}(w,x) \leq 3\lambda^*_1+\phi^*-1$.
Therefore, there exist at most $q_{\lvert V(H) \rvert}\lvert V(H) \rvert$ $(3\lambda^*_1+\phi^*+2)$-zones $\Lambda^*_w$ in $G'$ around vertices $w$ in $W^*$ such that $\bigcup_{w \in W^*}\Lambda^*_w \supseteq \bigcup_{Y \in \J_{F^*}}\Lambda_{x_Y}$.
Since for each $i \in [\ell]$, there exists $u_i \in V(H)$ such that $\pi_i(u_i) \in \bigcup_{Y \in \J_{F^*}}\bigcup_{(S,\Omega) \in \Se_1, \alpha(S,\Omega) \subseteq \overline{\Lambda_{x_Y}}}V(S)-\bigcup_{(S,\Omega) \in \Se_2}V(S)$ by the definition of $\F_1$ and a gauge, Statement 3 of this lemma holds, a contradiction.
$\Box$

For every $Y \in \bigcup_{F \in \F_1}\J_{F}^*$, let $x_Y \in \bigcup_{i=1}^\ell X_{\pi_i}$ be a vertex of $G'$ such that $Y$ is contained in the boundary cycle of $\Lambda_{x_Y}$.
Let $X=\{x_Y: Y \in \bigcup_{F \in \F_1} \J_F^*\}$.
Note that $\lvert X \rvert \leq \lvert \F_1 \rvert q_{\lvert \F_1 \rvert} \leq \lvert V(H) \rvert q_{\lvert V(H) \rvert}$.

Define $G''$ and $\T''$ to be the drawing and tangle obtained from $G'$ and $\T'$, respectively, by clearing $\bigcup_{(S,\Omega) \in \Se_2}\Lambda_S \cup \bigcup_{x \in X}\Lambda_x$.
Note that each $\Lambda_S$ and $\Lambda_x$ is a $\lambda_1^*$-zone, so the order of $\T''$ is at least $\theta^*-(4\lambda_1^*+2)(\kappa+q_{\lvert V(H) \rvert}\lvert V(H) \rvert)$.
Define $\M=\{e \in E(G''): m_{\T''}(e, \bigcup_{Y \in \bigcup_{F \in \F_1}\J_F^*}Y \cup \bigcup_{(S,\Omega) \in \Se_2}\bigcup_{i=1}^{\ell}Y_{S,i}) \geq \phi^*\}$.

\noindent{\bf Claim 2:} There exists $\M^* \subseteq \M$ with $\lvert \M^* \rvert = k_0^*$ such that every member of $\M^*$ is free with respect to $\T''$, and for distinct $Y,Y' \in \M^*$, $m_{\T''}(Y,Y') \geq \phi^*$.

\noindent{\bf Proof of Claim 2:}
Let $\M'$ be a maximal subset of $\M$ such that every member of $\M'$ is free with respect to $\T''$ and for distinct $Y,Y' \in \M'$, $m_{\T''}(Y,Y') \geq \phi^*$.
We are done if $\lvert \M' \rvert \geq k_0^*$.
So we may assume that $\lvert \M' \rvert < k_0^*$.
For each $(S,\Omega) \in \Se_2$, let $x_S$ be a vertex in $\overline{\Omega}$.
Let $W^*=\{x_Y: Y \in \bigcup_{F \in \F_1}\J_{F}^*\} \cup \{x_S: (S,\Omega) \in \Se_2\} \cup \M'$.
Hence, there exist $\lvert W^* \rvert$ $(\lambda_1^*+\phi^*+4)$-zones $\Lambda_w^*$ in $G''$ where each is around an element $w$ of $W^*$ such that $V(e) \subseteq \bigcup_{w \in W^*} \Lambda_w^*$ for every $e \in E(G'')$ with $V(e)$ free with respect to $\T''$.
Since each $\Lambda_w^{*}$ is a $(\phi^*+\lambda_1^*+4)$-zone in $G''$, it is contained in a $(\phi^*+\lambda_1^*+4+ (4\lambda_1^*+2)(\kappa+q_{\lvert V(H) \rvert}\lvert V(H) \rvert))$-zone in $G'$.
Since $\lvert W^* \rvert \leq q_{\lvert V(H) \rvert}\lvert V(H) \rvert+\kappa+k_0^* \leq \eta^*$, every edge $e$ of $G''$ with $V(e)$ free with respect to $\T''$ is contained in a union of at most $\eta^*$ $(\phi^*+\lambda_1^* + (4\lambda_1^*+2)(\kappa+q_{\lvert V(H) \rvert}\lvert V(H) \rvert))$-zones in $G'$.
Let $G_1$ and $\T_1$ be the drawing and tangle obtained from $G'$ by clearing these at most $\eta^*$ $(\phi^*+\lambda_1^* + (4\lambda_1^*+2)(\kappa+q_{\lvert V(H) \rvert}\lvert V(H) \rvert))$-zones in $G'$, respectively.
So the tangle $\T_1$ of order at least $\theta^*-(4(\phi^*+\lambda_1^* + (4\lambda_1^*+2)(\kappa+q_{\lvert V(H) \rvert}\lvert V(H) \rvert))+2)\eta^* \geq 2$.
Since every edge $e$ of $G_1$ with $V(e)$ free with respect to $\T_1$ is free with respect to $\T''$, there exists no edge $e$ of $G_1$ with $V(e)$ free with respect to $\T_1$.
But by the proof of \cite[Theorem (2.6)]{rs X}, since $\T_1$ has order at least two, there exists an edge $e^*$ of $G_1$ such that $e^* \in B$ for every $(A,B) \in \T_1$ with order at most one, a contradiction.
$\Box$

For each $(S,\Omega) \in \Se_2$, define $(G^S,\Omega^S)$ to be the $C_S$-vortex. 
Note that $G^S$ is a subgraph of $G$ and $\bigcup_{(S',\Omega') \in \Se, \alpha(S',\Omega') \subseteq \overline{\Lambda_S}} S' \subseteq G^S$.
By Lemma \ref{gauge}, $(G^S,\Omega^S)$ is a $\rho_S$-vortex for each $(S,\Omega) \in \Se_2$, so there exists a vortical decomposition $(P^S,\X^S)$ of $(G^S,\Omega^S)$ with adhesion at most $\rho_S$.
For each $(S,\Omega) \in \Se_2$ and $i \in [\ell]$, define $H^{S,i}$ to be the subgraph of $G^S \cap \pi_i(E(H))$ consisting of the components of $G^S \cap \pi_i(E(H))$ intersecting $Y_{S,i} \cup (\pi_i(V(H)) \cap V(G^S))$, and define $\I_{S,i}$ to be the canonical interval labelling of $H^{S,i}$ with respect to $(P^S,\X^S)$.
Since $H$ is connected and $\lvert Y_{S,i} \rvert \leq r(\rho_S) \leq r(\rho^*)$, each $\bigcup_{I \in \I_{S,i}}I$ is a $(r(\rho^*)+\lvert V(H) \rvert)$-interval.

\noindent{\bf Claim 3:} For each $(S,\Omega) \in \Se_2$ with $\pi_1(V(H) \cup E(H)) \cap S \neq \emptyset$, there exists $I_S \subseteq [\ell]$ with size $k$ such that $\bigcup_{I \in \I_{S,i}}I \cap \bigcup_{I \in \I_{S,j}}I=\emptyset$ for distinct $i,j \in I_S$.

\noindent{\bf Proof of Claim 3:}
Let $(S,\Omega) \in \Se_2$ with $\pi_1(V(H) \cup E(H)) \cap S \neq \emptyset$.
Since $\pi_1(V(H) \cup E(H)) \cap S \neq \emptyset$, $\pi_i(V(H) \cup E(H)) \cap S \neq \emptyset$ for every $i \in [\ell]$.
So $\bigcup_{I \in \I_{S,i}}I \neq \emptyset$ for each $i \in [\ell]$.

Suppose to the contrary that such $I_S$ does not exist.
Since $\bigcup_{I \in \I_{S,i}}I$ is a $(r(\rho^*)+\lvert V(H) \rvert)$-interval for each $i \in [\ell]$, by Theorem \ref{t-interval stable}, there exists a set $D \subseteq {\mathbb R}$ with $\lvert D \rvert \leq 2(r(\rho^*)+\lvert V(H) \rvert)^2(k-1)$ such that $D$ intersects $\bigcup_{I \in \I_{S,i}}I$ for every $i \in [\ell]$.
Since every endpoint of an interval in $\I_{S,i}$ is a positive integer, we may assume that every point in $D$ is a positive integer.
Denote the bags of $\X^S$ as $X_1,X_2,...,X_{\lvert V(P^S) \rvert}$ such that $X_i$ contains the $i$-th vertex of $\Omega$ for each $i \in [\lvert \overline{\Omega} \rvert]$.
Let $X_0$ and $X_{\lvert V(P^S) \rvert+1}$ be the empty set.
For each $i \in [\lvert V(P^S) \rvert]$, let $Z_i=X_i \cap (X_{i-1} \cup X_{i+1} \cup \overline{\Omega})$.
Since $\Se$ is $\T$-central and $H$ is connected, there do not exist $i \in [\ell]$ and $j \in [\lvert V(P^S) \rvert]$ such that $V(\pi_i(V(H) \cup E(H))) \cap X_j \neq \emptyset$ and $V(\pi_i(V(H) \cup E(H))) \cap Z_j=\emptyset$, for otherwise there exists a separation $(A,B) \in \T$ of order at most $2\rho^*+1$ with $V(A)=X_j$ and $V(A \cap B)=Z_j$ such that $\pi_i(E(H)) \subseteq A$ and Statement 4 holds.
Define $Z=\bigcup_{d \in D}Z_d$.
Then $Z \cap V(\pi_i(V(H) \cup E(H)) \cap S) \neq \emptyset$ for every $i \in [\ell]$.
But $\lvert Z \rvert \leq \lvert D \rvert \cdot (2\rho^*+1) \leq \xi^*$.
So Statement 2 holds, a contradiction.
$\Box$

For each $(S,\Omega) \in \Se_2$ with $\pi_1(V(H) \cup E(H)) \cap S \neq \emptyset$, let $I_S$ be the set mentioned in Claim 3, and we denote the elements of $I_S$ by $i_{S,1},i_{S,2},...,i_{S,k}$. 

(Figure \ref{fig_claim4_vortex} is a schematic for Claim 4.)

\begin{figure} 
	\begin{picture}(200,260) (-50,-130)
		
		\put(190,0){$\alpha(S,\Omega)$}  
		\put(200,0){\oval(110,110)}
		\thicklines
		\put(200,0){\oval(130,130)}
		\thinlines
		\put(200,0){\oval(150,150)}
		\put(200,0){\oval(170,170)}
		\put(200,0){\oval(190,190)}
		\put(200,0){\oval(210,210)}
		\thicklines
		\put(200,0){\oval(230,230)}

		\put(150,0){\circle*{5}} 
		\put(255,30){\circle*{5}} 
		\put(190,50){\circle*{5}} 
		\put(200,-10){\circle*{5}}
		\linethickness{0.3mm}
		\qbezier(150,0)(200,-10)(190,50)
		\qbezier(150,0)(150,-50)(200,-10)
		\qbezier(200,-10)(250,-10)(190,50)
		\qbezier(200,-10)(270,-30)(255,30)

		\put(85,90){\circle{5}} 
		\qbezier(85,90)(200,50)(150,0)
		\put(85,-50){\circle{5}} 
		\qbezier(85,-50)(100,-90)(150,0)
		\put(230,115){\circle{5}} 
		\qbezier(230,115)(230,110)(230,105)
		\qbezier(230,105)(220,105)(210,105)
		\qbezier(210,105)(210,100)(210,95)
		\qbezier(210,95)(200,90)(190,85)
		\qbezier(190,85)(195,85)(230,85)
		\qbezier(230,85)(225,75)(190,50)
	
		\put(200,-115){\circle{5}} 
		\put(315,-30){\circle{5}} 
		\qbezier(200,-115)(230,0)(315,-30)
		\put(150,-115){\circle{5}} 
		\put(250,-115){\circle{5}} 
		\qbezier(150,-115)(220,70)(250,-115)
	\end{picture}
	\caption{The inner thick cycle is $C_S^*$. The outer thick cycle is $C_S$. The thick curves form $R_1$. The empty circles are the vertices in $Y_{S,i_{S,1}}$.} \label{fig_claim4_vortex}
\end{figure}

\noindent{\bf Claim 4:} For each $(S,\Omega) \in \Se_2$ with $\pi_1(V(H) \cup E(H)) \cap S \neq \emptyset$, there exists a cycle $C_S^*$ in $G'$ bounding an open disk $\Lambda_S^*$ with $\alpha(S,\Omega) \subset \Lambda_S^* \subset \Lambda_S$ such that
	\begin{itemize}
		\item there exist at least $2\phi(\sum_{(S,\Omega) \in \Se_2}\rho_S)$ pairwise disjoint cycles in $G'$, each bounding a disk strictly contained in $\Lambda_S$ and strictly containing $\Lambda_S^*$, 
		\item there exist $k$ pairwise disjoint subgraphs $R_1,R_2,...,R_k$ of $G^S$ such that for each $j \in [k]$, there exists a homeomorphic embedding $\eta_j$ from the $(Q_S,\Omega_S)$-witness for the lifting of $\pi_{i_{S,j}}$ defining the snapshot into $R_j$ preserving the labelling such that $\eta_j|_{Y_{S,i_{S,j}}}$ is an identity map, and
		\item for each $j \in [k]$, every component of $R_j \cap \bigcup_{(S',\Omega') \in \Se, \alpha(S',\Omega') \subseteq \overline{\Lambda_S}-\Lambda_S^*}S'$ is a path in $G$ from $V(C_S^*)$ to $\bigcup_{i \in I_S}Y_{S,i}$ internally disjoint from $V(C_S)$. 
	\end{itemize}

\noindent{\bf Proof of Claim 4:}
Let $(S,\Omega) \in \Se_2$ with $\pi_1(V(H) \cup E(H)) \cap S \neq \emptyset$.
Let $U=\bigcup_{i \in I_S}(\pi_i(V(H)) \cap V(G^S)) \allowbreak \cup \allowbreak \bigcup_{i \in I_S} \bigcup_{e \in E(H), \pi_i(e) \subseteq G^S} \pi_i(e)$.
For each $j \in I_S$, let $W_j = \pi_j(E(H_{S,j}))$, where $H_{S,j}$ is the underlying graph of the $(Q_S,\Omega_S)$-witness for the lifting of $\pi_j$ defining the snapshot.
Note that $Y_{S,j} \subseteq V(H_{S,j})$, and $W_j \cap W_{j'}=\emptyset$ for distinct $j,j' \in I_S$. 
Let $W = \bigcup_{j \in I_S}W_j$.
Note that $U \subseteq W$.
Let $L=W-V(U)$.

Let $\Lambda'$ be the open disk bounded by $C_S'$. 
Let $(D,\Omega^D)$ be the $C'_S$-vortex.
By the definition of a lifting and the definition of an $(H,\theta_0^*,\mu^*)$-gauge, $V(L) \cap V(C_S) \subseteq \bigcup_{i \in I_S}Y_{S,i}$, and $L$ is a linkage in $G^S-(V(U) \cup (V(C_S)-V(L)))$ with at most $k\lvert Y_{S,1} \rvert \leq kr(\rho_S)$ components such that every component is either a path with one end in $\overline{\Omega^S}$ and one end in $V(D)$, or a path with both ends in $V(C_S)$ and there exists another component of $L$ intersecting the both intervals of $\Omega^S$ determined by the ends of this path.

Let $\widehat{G^S}$ be the graph obtained from $G^S$ by replacing $G^S-(V(D)-V(C_S'))$ by $G' \cap \overline{\Lambda_S}-\Lambda'$.
By the definition of an $(H,\theta_0^*,\mu^*)$-gauge, $(\widehat{G^S},\Omega^S)$ is a composition of $(D,\Omega^D)$ with a $\mu^*(\sum_{(S,\Omega) \in \Se_2}\rho_S)$-nested rural neighborhood, and $U \subseteq D-V(C_S') \subseteq \widehat{G^S}$.
So $(\widehat{G^S}-(V(U) \cup (V(C_S)-V(L))),\Omega^S|_{\overline{\Omega^S}-(V(C_S)-V(L))})$ is a composition of $(D-V(U),\Omega^D)$ with a $(\mu^*(\sum_{(S,\Omega) \in \Se_2}\rho_S)-1)$-nested rural neighborhood.

Since $\Se$ is effective with respect to $(\Se_1,\Se_2)$, we can replace each maximal subpath of a path in $L-(V(D)-V(C_S'))$ contained in some member of $\Se_1$ by an edge in $\widehat{G^S}-E(G^S)$ to obtain a linkage $\widehat{L}$ in $\widehat{G^S}-(V(U) \cup (V(C_S)-V(L)))$ with at most $k\lvert Y_{S,1} \rvert \leq kr(\rho_S)$ components such that every component is either a path with one end in $\overline{\Omega^S}$ and one end in $V(D)$, or a path with both ends in $V(C_S)$ and there exists another component of $\widehat{L}$ intersecting the both intervals of $\Omega^S$ determined by the ends of this path.
Note that $\widehat{L}$ connects the same pairs of vertices as $L$.

Since $\mu^*(\sum_{(S,\Omega) \in \Se_2}\rho_S)-1 \geq s_{\ref{perpendicular}}(2\phi(\sum_{(S,\Omega) \in \Se_2}\rho_S))$, by Lemma \ref{perpendicular}, there exist a linkage $\widehat{L'}$ in $\widehat{G^S}-(V(U) \cup (V(C_S)-V(L)))$ equivalent to $\widehat{L}$ and a society $(D',\Omega^{D'})$ with $D-V(U) \subseteq D'$ such that $(\widehat{G^S}-(V(U) \cup (V(C_S)-V(L))),\Omega^S|_{\overline{\Omega^S}-(V(C_S)-V(L))})$ is a composition of $(D',\Omega^{D'})$ with a rural neighborhood $(S'',\Omega,\Omega^{D'})$ that has a $2\phi(\sum_{(S,\Omega) \in \Se_2}\rho_S)$-nest $(C_1',C_2',...,C_{2\phi(\sum_{(S,\Omega) \in \Se_2}\rho_S)}')$ for some presentation of it such that $\widehat{L'} \cap S''$ is perpendicular to $(C_1',C_2',...,C_{2\phi(\sum_{(S,\Omega) \in \Se_2}\rho_S)}')$, and $S'' \cap D'$ is a cycle in $G'$ with vertex-set $\overline{\Omega^{D'}}$ passing through $\overline{\Omega^{D'}}$ in the order $\Omega^{D'}$.
Since $\Se$ is effective with respect to $(\Se_1,\Se_2)$, we can replace each edge of $\widehat{L'}$ contained in $\widehat{G^S}-E(G^S)$ by a path in a member of $\Se_1$ to obtain a linkage $L'$ in $G^S$ connecting the same pairs of vertices as $\widehat{L'}$.
Note that $L'$ is equivalent to $L$.

Define $C_S^*=S'' \cap D'$ and define $\Lambda_S^*$ to be the open disk bounded by $C_S^*$.
Since $L'$ is equivalent to $L$, $L' \cup U$ can be partitioned into subgraphs $R_1,R_2,...,R_k$ such that for each $j \in [k]$, there exists a homeomorphic embedding $\eta_j$ from $H_{S,i_{S,j}}$ into $R_j$ preserving the labelling such that $\eta_j|_{Y_{S,i_{S,j}}}$ is an identity map. 
Since $\widehat{L'} \cap S''$ is perpendicular to $(C_1',C_2',...,C_{2\phi(\sum_{(S,\Omega) \in \Se_2}\rho_S)}')$ and $V(C_S) \cap V(L')=V(C_S) \cap V(L) \subseteq \bigcup_{i \in I_S}Y_{S,i}$, we have that for each $j \in [k]$, every component of $R_j \cap \bigcup_{(S',\Omega') \in \Se, \alpha(S',\Omega') \subseteq \overline{\Lambda_S}-\Lambda_S^*}S'$ is a path from $V(C^*_S)$ to $\bigcup_{i \in I_S}Y_{S,i}$ internally disjoint from $V(C_S)$. 
$\Box$

For each $(S,\Omega) \in \Se_2$ with $\pi_1(V(H) \cup E(H)) \cap S \neq \emptyset$, define $Y_S^*$ to be the set of ends of the components of $\bigcup_{j \in [k]}(R_j \cap \bigcup_{(S',\Omega') \in \Se, \alpha(S',\Omega') \subseteq \overline{\Lambda_S}-\Lambda_S^*}S')$ in $V(C^*_S)$, where $R_j$ and $C_S^*$ are the graph and the cycle mentioned in Claim 4, respectively.
Note that the third conclusion of Claim 4 implies that $\lvert Y_S^* \rvert \leq \lvert \bigcup_{i \in I_S}Y_{S,i} \rvert \leq k r(\rho_S)$ for each $(S,\Omega) \in \Se_2$ with $\pi_1(V(H) \cup E(H)) \cap S \neq \emptyset$.

For each $x \in X$, define $(G^x,\Omega^x)$ to be the $C_x$-vortex. 
The following claim can be proved in the same way as Claim 4, so we omit the details.

\noindent{\bf Claim 5:} For each $x \in X$ with $x=x_Y$ for some $Y \in \bigcup_{F \in \F_1}\J^*_F$, there exists a cycle $C_x^*$ in $G'$ bounding an open disk $\Lambda_x^*$ with $\alpha(x) \subset \Lambda_x^* \subset \Lambda_x$ such that
	\begin{itemize}
		\item there exist at least $2\phi(\sum_{(S,\Omega) \in \Se_2}\rho_S)$ pairwise disjoint cycles in $G'$, each bounding a disk strictly contained in $\Lambda_x$ and strictly containing $\Lambda_x^*$, 
		\item there exist a subgraph $R$ of $G^x$ and a homeomorphic embedding $\eta$ from $H_x$ into $R$ preserving the labelling such that $\eta|_{Y_{x}}$ is an identity map, and
		\item every component of $R \cap \bigcup_{(S',\Omega') \in \Se, \alpha(S',\Omega') \subseteq \overline{\Lambda_x}-\Lambda_x^*}S'$ is a path in $G$ from $V(C_x^*)$ to $Y_x$ internally disjoint from $V(C_x)$.
	\end{itemize}

For each $x \in X$, define $Y_x^*$ to be the set of ends of the components of $\linebreak R \cap \bigcup_{(S',\Omega') \in \Se, \alpha(S',\Omega') \subseteq \overline{\Lambda_x}-\Lambda_x^*}S'$ in $V(C_x^*)$, where $R$ and $C_x^*$ are the graph and the cycle mentioned in Claim 5, respectively.
Note that $\lvert Y^*_x \rvert \leq \lvert Y_x \rvert \leq 2\lvert E(H) \rvert$ for every $x \in X$.

Define $G''',\T'''$ to be the graph and tangle obtained from $G',\T'$ by clearing $\bigcup_{(S,\Omega) \in \Se_2}\Lambda_S^* \cup \bigcup_{x \in X} \Lambda_x^*$, respectively.
Since each $\Lambda_S^*$ is contained in the $\lambda_1^*$-zone $\Lambda_S$, $\Lambda_S^*$ is a $\lambda_1^*$-zone.
Similarly, each $\Lambda_x^*$ is a $\lambda_1^*$-zone contained in $\Lambda_x$.
Hence the order of $\T'''$ is at least $\theta^*-(4\lambda^*_1+2)(\kappa+q_{\lvert \F_1 \rvert}\lvert \F_1 \rvert)$, and for any two atoms $a,b$ of $G'''$, $m_{\T'''}(a,b) \geq m_{\T'}(a,b)-(4\lambda^*_1+2)(\kappa+q_{\lvert \F_1 \rvert}\lvert \F_1 \rvert)$.

Let $\Y=\{Y_S^*: (S,\Omega) \in \Se_2$ with $\pi_1(V(H) \cup E(H)) \cap S \neq \emptyset\} \cup \{Y^*_x: x \in X\}$.
Let $\M^*$ be the set mentioned in Claim 2.
For any member $e \in \M^*$, define $\Delta_e$ to be a disk in $\Sigma$ such that $\Delta_e \cap G'=V(e)$.

\noindent{\bf Claim 6:} The following statements hold.
	\begin{itemize}
		\item $m_{\T'''}(Y,Y') \geq 2\phi(\sum_{(S,\Omega)\in \Se_2}\rho_S)$ for any two distinct members $Y,Y'$ of $\Y$. 
		\item $m_{\T'''}(Y,V(e)) \geq 2\phi(\kappa\rho^*)$ for every $Y \in \Y$ and $e \in \M^*$.
		\item Every member $Y$ of $\Y$ is free with respect to $\T'''$.
	\end{itemize}

\noindent{\bf Proof of Claim 6:}
Let $Y,Y'$ be distinct members of $\Y$.

We first assume that $Y=Y_S^*$ and $Y'=Y_{S'}^*$ for some distinct $(S,\Omega),(S',\Omega') \in \Se_2$ with $\pi_1(V(H) \cup E(H)) \cap S \neq \emptyset \neq \pi_1(V(H) \cup E(H)) \cap S'$.
Since every closed walk in the radial drawing of $G'''$ whose inside contains some atom of $Y$ and some atom of $Y'$ also contains some atom of $V(C_S)$ and some atom of $V(C_{S'})$, we have $m_{\T'''}(Y,Y') \geq m_{\T'''}(V(C_S),V(C_{S'}))$.
Since $G''$ is a subdrawing of $G'''$, $m_{\T'''}(V(C_S),V(C_{S'})) \geq m_{\T''}(V(C_S),V(C_{S'}))$.
By the definition of an $(H,\theta_0^*,\mu^*)$-gauge (or Statement \ref{Y_distance_gague} in Lemma \ref{gauge}), $m_{\T''}(V(C_S),V(C_{S'})) \geq \theta^*_0(\sum_{(S'',\Omega'')\in \Se_2}\rho_{S''})-2$.
Note that $\theta^*_0(\sum_{(S'',\Omega'')\in \Se_2}\rho_{S''})-2 \geq \mu^*(\sum_{(S'',\Omega'')\in \Se_2}\rho_{S''}) \geq \linebreak \allowbreak s_{\ref{perpendicular}}(2\phi(\sum_{(S'',\Omega'')\in \Se_2}\rho_{S''})) \geq 2\phi(\sum_{(S'',\Omega'')\in \Se_2}\rho_{S''})$.
So $m_{\T'''}(Y,Y')\geq 2\phi(\sum_{(S'',\Omega'')\in \Se_2}\rho_{S''})$.

Similarly, if $Y=Y_S^*$ and $Y'=Y^*_x$ for some $(S,\Omega) \in \Se_2$ with $\pi_1(V(H) \cup E(H)) \cap S \neq \emptyset$ and for some $x \in X$, then $m_{\T'''}(Y,Y') \geq 2\phi(\sum_{(S',\Omega')\in \Se_2}\rho_{S'})$. 

If $Y=Y^*_x$ and $Y'=Y^*_{x'}$ for some distinct $x,x' \in X$, then by Claim 1, $m_{\T'''}(Y,Y') \geq m_{\T'}(Y,Y')-(4\lambda^*_1+2)(\kappa+q_{\lvert \F_1 \rvert}\lvert \F_1 \rvert) \geq \phi^*-(4\lambda^*_1+2)(\kappa+q_{\lvert \F_1 \rvert}\lvert \F_1 \rvert) \geq 2\phi(\kappa\rho^*) \geq 2\phi(\sum_{(S,\Omega)\in \Se_2}\rho_S)$.

Furthermore, for every $e \in \M^*$, since $G''$ is a subdrawing of $G'''$, by the definition of $\M^*$, $m_{\T'''}(V(e), Y) \geq m_{\T'''}(V(e), Y'')-2\lambda_1^* \geq m_{\T''}(V(e), Y'')-2\lambda_1^* \geq \phi^*-2\lambda_1^* \geq 2\phi(\kappa\rho^*)$, where $Y''$ is the member in $\bigcup_{F \in \F_1}\J^*_F \cup \{Y_{S''}: (S'',\Omega'') \in \Se_2\}$ such that $Y'' \cup Y$ is contained in a member of $\{\overline{\Lambda_{S''}},\overline{\Lambda_x}: (S'',\Omega'') \in \Se_2, x \in X\}$.

Hence, to prove this claim, it suffices to prove that $Y$ is free with respect to $\T'''$.
Suppose that $Y$ is not free with respect to $\T'''$.
By Lemma \ref{A distance set}, there exists $(A,B) \in \T'''$ with order less than $\lvert Y \rvert$ such that $Y \subseteq V(A)$ and $m_{\T'''}(Y,z) \leq \lvert V(A \cap B) \rvert<\lvert Y \rvert \leq \max\{kr(\rho_S),2\lvert E(H) \rvert: (S,\Omega) \in \Se_2\} < 2\phi(\sum_{(S',\Omega') \in \Se_2}\rho_{S'})$ for every $z \in V(A)$.
We assume that $Y=Y^*_S$ for some $(S,\Omega) \in \Se_2$ with $\pi_1(V(H) \cup E(H)) \cap S \neq \emptyset$.
(The case that $Y=Y^*_x$ for some $x \in X$ can be proved analogously, so we omit the proof.)

Since $\lvert Y \rvert < 2\phi(\sum_{(S',\Omega') \in \Se_2}\rho_{S'})$, Claim 4 implies that there exist disjoint cycles $D_1,D_2,..., \allowbreak D_{\lvert Y \rvert}$ in $G'$ such that each $D_i$ bounds an open disk $\Delta_i$, and $\Lambda^*_S \subset \Delta_1 \subset \Delta_2 \subset ... \subset \Delta_{\lvert Y \rvert} \subset \Lambda_S$.
Since there are $\lvert Y \rvert$ disjoint paths in $G'''$ from $Y$ to $V(C_S)$, $V(C_S) \not \subseteq V(B)$.
So there exists a closed walk $W$ in the radial drawing of $G'''$ with length less than $2\lvert Y \rvert$ such that the inside of $W$ contains some vertex in $Y$ and some vertex in $V(C_S)$.

If $W$ is a walk in the radial drawing of $G'$, then the existence of the $\lvert Y \rvert$ disjoint cycles $D_1,D_2,...,D_{\lvert Y \rvert}$ and the length condition of $W$ imply that the inside of $W$ contains a vertex in $V(C_S)$ and a vertex of $V(C_S')$, so $m_{\T'}(V(C_S),V(C_S'))<2\phi(\sum_{(S',\Omega') \in \Se_2}\rho_{S'})$, contradicting the definition of an $(H,\theta_0^*,\mu^*)$-gauge.
So $W$ is not a radial drawing of $G'$.
Since $2\phi(\kappa\rho^*) \geq 2\phi(\sum_{(S',\Omega') \in \Se_2}\rho_{S'})>\lvert Y \rvert$, $W$ is disjoint from $\bigcup_{(S',\Omega') \in \Se_2-\{(S,\Omega)\}}\Lambda^*_{S'} \cup \bigcup_{x \in X}\Lambda^*_x \cup \bigcup_{e \in \M^*}\Delta_e$ by the first two statements of this claim.
Since $W$ is not a radial drawing of $G'$, $W$ intersects $\Lambda^*_S$.
Since the inside of $W$ contains a vertex in $V(C_S)$, $W$ intersects $V(D_i)$ for every $i \in [\lvert Y \rvert]$.
So the length of $W$ is at least $2 \lvert Y \rvert$, a contradiction.
This proves that $Y$ is free with respect to $\T'''$.
$\Box$

Let $\Y^* = \{Y_S^*: (S,\Omega) \in \Se_2$ with $\pi_1(V(H) \cup E(H)) \cap S \neq \emptyset\} \cup \{Y^*_x, V(e): x \in X, e \in \M^*\}$.
Note that $\Y^*=\Y \cup \{V(e): e \in \M^*\}$.
Since $G''$ is a subgraph of $G'''$, each member of $\M^*$ is free with respect to $\T'''$ by Claim 2.
So every member of $\Y^*$ is free with respect to $\T'''$ by Claim 6. 
Define $Y^*=\bigcup_{Y \in \Y^*}Y$.
Note that $\lvert Y^* \rvert \leq 2k\sum_{(S,\Omega) \in \Se_2}r(\rho_S) + 2\lvert E(H) \rvert \cdot \lvert \F_1 \rvert q_{\lvert \F_1 \rvert} + 2k_0^* \leq 2k\sum_{(S,\Omega) \in \Se_2}r(\rho_S)+2\lvert E(H) \rvert \cdot \lvert V(H) \rvert q_{\lvert V(H) \rvert}+2k_0^*$.

By the definition of a $(\Sigma, (\Lambda_S: (S,\Omega) \in \Se_2))$-template, there exist $k$ collections $\P_1,\P_2,...,\allowbreak \P_k$ and $k$ sets $\Gamma_1,...,\Gamma_k$ of pairwise disjoint simple curves such that the following statements hold.
	\begin{itemize}
		\item For each $i \in [k]$, $\P_i$ consists of a member of $\J_F^*$ for each $F \in \F_1$ and a subset of $Y_{S}$ for each $(S,\Omega) \in \Se_2$ with $\pi_1(V(H) \cup E(H)) \cap S \neq \emptyset$, such that $\bigcup_{Y \in \P_i}Y$ and $\bigcup_{Y \in \P_j}Y$ are pairwise disjoint for $i \neq j$.
		\item For each $j \in [k]$ and for each curve $\gamma \in \Gamma_j$, the endpoints of $\gamma$ are two vertices in $\bigcup_{Y \in \P_j}Y$.
		\item For each $j \in [k]$, the union of $\bigcup_{\gamma \in \Gamma_j}\gamma$ and the subgraphs $R_{S,i}$ and $R_x$ for all $(S,\Omega) \in \Se_2$ with $\pi_1(V(H) \cup E(H)) \cap S \neq \emptyset$ and for all $i \in [k]$ and $x \in X$, where $R_{S,i}$ is the subgraph $R_i$ mentioned in Claim 4, and $R_x$ is the subgraph $R$ mentioned in Claim 5, gives a pseudo-embedding of $H$ in $\Sigma$ legal with respect to $\{\Lambda^*_{S'},\Lambda^*_x: (S',\Omega') \in \Se_2,x \in X\}$ with no crossing-points in $\Sigma-(\bigcup_{(S',\Omega') \in \Se_2}\Lambda^*_{S'} \cup \bigcup_{x \in X}\Lambda^*_x)$.
	\end{itemize}
Note that each $\Gamma_i$ contains at most $\sum_{(S,\Omega) \in \Se_2}r(\rho_S)+2\lvert E(H) \rvert\lvert V(H) \rvert$ members.
By Lemma \ref{bounding crossings}, there exist $k$ sets $\Gamma_1',...,\Gamma_k'$ of pairwise disjoint simple curves such that the following statements hold.
	\begin{itemize}
		\item For each $j \in [k]$ and for each curve $\gamma \in \Gamma'_j$, the endpoints of $\gamma$ are two vertices in $\bigcup_{Y \in \P_j}Y$.
		\item For each $j \in [k]$, the union of $\bigcup_{\gamma \in \Gamma'_j}\gamma$ and the subgraphs $R_{S,i}$ and $R_x$ for all $(S,\Omega) \in \Se_2$ with $\pi_1(V(H) \cup E(H)) \cap S \neq \emptyset$ and for all $i \in [k]$ and $x \in X$, gives a pseudo-embedding of $H$ in $\Sigma$ legal with respect to $\{\Lambda^*_{S'},\Lambda^*_x: (S',\Omega') \in \Se_2,x \in X\}$ with no crossing-points in $\Sigma-(\bigcup_{(S',\Omega') \in \Se_2}\Lambda^*_{S'} \cup \bigcup_{x \in X}\Lambda^*_x)$. 
		\item For each point $x \in \gamma \cap \gamma'$ for some distinct members $\gamma,\gamma'$ of $\bigcup_{i=1}^k\Gamma_i'$, there exists an open set $B_x \subseteq \Sigma$ with $x \in B_x$ such that $B_x-\{x\}$ does not contain any point that belongs to at least two distinct members of $\bigcup_{i=1}^k\Gamma_i'$.
		\item Every point in $\Sigma$ belongs to at most two members of $\bigcup_{i=1}^k\Gamma_i'$.
		\item There exists at most $\tau_{\ref{bounding crossings}}(\rho^*)$ points belonging to at least two members of $\bigcup_{i=1}^k\Gamma_i'$.
	\end{itemize}

Note that each endpoint $v$ of a member of $\bigcup_{i=1}^k\Gamma'_i$ is a vertex in $R^v$, where $R^v$ is either equal to $R_{S,i}$ for some $(S,\Omega) \in \Se_2$ and $i \in [k]$ with $\pi_1(V(H) \cup E(H)) \cap S \neq \emptyset$, or equal to $R_x$ for some $x \in X$.
By Claims 4 and 5, every component of $R^v \cap (\bigcup_{(S,\Omega) \in \Se_2}\bigcup_{(S',\Omega') \in \Se,\alpha(S',\Omega') \subseteq \overline{\Lambda_{S}}-\Lambda^*_{S}}S' \cup \bigcup_{x \in X}\bigcup_{(S',\Omega') \in \Se,\alpha(S',\Omega') \subseteq \overline{\Lambda_{x}}-\Lambda^*_{x}}S')$ is a path from $v$ to some vertex in $Y^*$.
We define $\beta(v)$ to be the aforementioned vertex in $Y^*$.

We say a point $x$ in $\Sigma$ is a {\it $\bigcup_{i=1}^k \Gamma'_i$-crossing} if $x \in \gamma \cap \gamma'$ for some distinct $\gamma,\gamma' \in \bigcup_{i=1}^k \Gamma'_i$.
For each $\bigcup_{i=1}^k \Gamma'_i$-crossing $x$, define $\Delta_x$ to be a closed disk contained in $B_x$ and containing $x$.
Since there are only finitely many $\bigcup_{i=1}^k \Gamma'_i$-crossings, we may choose $\Delta_x$ for all $\bigcup_{i=1}^k \Gamma'_i$-crossings $x$ such that they are pairwise disjoint.

We say two $\bigcup_{i=1}^k \Gamma'_i$-crossings $x,y$ are {\it adjacent by $\gamma$} if $\gamma \in \bigcup_{i=1}^k \Gamma'_i$ and there exists $\gamma' \subseteq \gamma$ such that $x,y \in \gamma'$ and $\gamma'$ does not contain any $\bigcup_{i=1}^k \Gamma'_i$-crossing other than $x,y$.
For each unordered pair of adjacent $\bigcup_{i=1}^k \Gamma'_i$-crossings $\{x,y\}$ by some $\gamma \in \bigcup_{i=1}^k \Gamma'_i$, we define $\Delta_{x,y,\gamma}$ to be a closed disk contained in an open disk around a point in $\gamma$ between $x,y$ disjoint from $\bigcup_x \Delta_x$ (where the union is over all $\bigcup_{i=1}^k \Gamma_i'$-crossings $x$) such that $\Delta_{x,y,\gamma} \cap \gamma$ has exactly one component, $\Delta_{x,y,\gamma}$ is disjoint from all members of $(\bigcup_{i=1}^k \Gamma'_i)-\{\gamma\}$, and $\Delta_{x,y,\gamma}$ is disjoint from $\Delta_{x',y',\gamma'}$ for any unordered pair of adjacent $\bigcup_{i=1}^k \Gamma'_i$-crossings $\{x',y'\}$ different from $\{x,y\}$.
(See Figure \ref{fig_crossing_disk_linking}.)

Note that each member of $\bigcup_{i=1}^k\Gamma'_i$ is a mapping from $[0,1]$ to $\Sigma$.
Define $\Xi=\{\Delta_x: x$ is a $\bigcup_{i=1}^k\Gamma_i'$-crossing$\} \cup \{\Delta_{x,y,\gamma}: \{x,y\}$ is an unordered pair of adjacent $\bigcup_{i=1}^k\Gamma'_i$-crossings by some $\gamma \in \bigcup_{i=1}^k\Gamma_i'\}$.
Since there are at most $\tau_{\ref{bounding crossings}}(\rho^*)$ $\bigcup_{i=1}^k\Gamma_i'$-crossings, $\lvert \Xi \rvert \leq 5\tau_{\ref{bounding crossings}}(\rho^*)=\lvert \M^* \rvert$.
By deleting members of $\M^*$, we may assume that $\lvert \M^* \rvert = \lvert \Xi \rvert$.
Define $\eta$ to be a bijection from $\Xi$ to $\M^*$, and for each $\Delta \in \Xi$, we write the ends of $\eta(\Delta)$ as $\eta(\Delta)_0$ and $\eta(\Delta)_1$.

\begin{figure} 
	\begin{picture}(200,300) (-50,-200)
	
		\put(0,60){\circle*{5}}
		\put(-15,40){$\gamma(0)$} 
		\put(390,60){\circle*{5}} 
		\put(375,40){$\gamma(1)$} 
		\put(0,60){\line(1,0){390}}
		\put(120,50){$x$} 
		\put(250,50){$y$} 

		\put(130,120){\circle*{5}}
		\put(100,120){$\gamma'(0)$} 
		\put(130,0){\circle*{5}} 
		\put(100,0){$\gamma'(1)$} 
		\put(130,0){\line(0,1){120}}

		\put(260,120){\circle*{5}}
		\put(230,120){$\gamma''(0)$} 
		\put(260,0){\circle*{5}} 
		\put(230,0){$\gamma''(1)$} 
		\put(260,0){\line(0,1){120}}

		\put(130,60){\circle{40}} 
		\put(105,80){$\Delta_x$} 
		\put(260,60){\circle{40}} 
		\put(235,80){$\Delta_y$} 
		\put(195,60){\oval(30,30)}
		\put(185,80){$\Delta_{x,y}$}

		\multiput(-40,-20)(5,0){90}{\line(1,0){3}}

		\put(0,-110){\circle*{5}}
		\put(-15,-130){$\beta(\gamma(0))$} 
		\put(-15,-145){$=v_{x,\gamma,1}$} 
		\put(390,-110){\circle*{5}} 
		\put(375,-130){$\beta(\gamma(1))$}
		\put(375,-145){$=v_{y,\gamma,0}$} 
		
		\put(130,-50){\circle*{5}}
		\put(80,-50){$\beta(\gamma'(0))$}
		\put(80,-65){$=v_{x,\gamma',1}$} 
		\put(130,-170){\circle*{5}} 
		\put(80,-170){$\beta(\gamma'(1))$} 
		\put(80,-185){$=v_{x,\gamma',0}$} 
	
		\put(260,-50){\circle*{5}}
		\put(210,-50){$\beta(\gamma''(0))$} 
		\put(210,-65){$=v_{y,\gamma'',1}$} 
		\put(260,-170){\circle*{5}} 
		\put(210,-170){$\beta(\gamma''(1))$} 
		\put(210,-185){$=v_{y,\gamma',0}$} 

		\put(110,-110){\circle*{5}} 
		\put(150,-110){\circle*{5}} 
		\thicklines
		\put(110,-110){\line(1,0){40}}

		\put(260,-90){\circle*{5}} 
		\put(260,-130){\circle*{5}}
		\put(260,-130){\line(0,1){40}}

		\put(180,-110){\circle*{5}} 
		\put(210,-110){\circle*{5}} 
		\put(167,-125){$v_{x,\gamma,0}$}
		\put(207,-125){$v_{y,\gamma,1}$} 
		\put(180,-110){\line(1,0){30}}

		\put(85,-100){$\eta(\Delta_x)_0$} 
		\put(135,-100){$\eta(\Delta_x)_1$}

		\put(268,-135){$\eta(\Delta_y)_1$} 
		\put(268,-92){$\eta(\Delta_y)_0$} 

		\thinlines
		\put(-22,-35){\line(1,0){160}}
		\put(-22,-35){\line(0,-1){80}}
		\put(-22,-115){\line(1,0){140}}
		\put(118,-115){\line(1,4){20}}

		\put(127,-90){\line(1,0){65}}
		\put(127,-90){\line(0,-1){95}}
		\put(127,-185){\line(1,0){65}}
		\put(192,-90){\line(0,-1){95}}

		\put(200,-35){\line(1,0){65}}
		\put(200,-35){\line(0,-1){82}}
		\put(200,-117){\line(1,0){65}}
		\put(265,-35){\line(0,-1){82}}

		\put(255,-125){\line(6,1){159}}
		\put(414,-99){\line(0,-1){85}}
		\put(255,-125){\line(0,-1){59}}
		\put(255,-184){\line(1,0){159}}
	\end{picture}
	\caption{A picture for vertices in parts in $\Q$. The upper side includes curves and disks related to $\bigcup_{i=1}^k\Gamma_i'$-crossings. The lower side includes the corresponding vertices in $Y^*$ and parts in $\Q$.} \label{fig_crossing_disk_linking}
\end{figure}

Recall that every member of $\Y^*$ is free with respect to $\T'''$, and $Y^*$ is the union of those members. 
Now we define a partition $\Q$ of a subset of $Y^*$ as follows. (See Figure \ref{fig_crossing_disk_linking}.)
	\begin{itemize}		
		\item For each $\gamma \in \bigcup_{i=1}^k\Gamma_i'$ such that $\gamma$ does not contain any $\bigcup_{i=1}^k\Gamma'_i$-crossing, we define $Q_\gamma$ to be $\{\beta(a),\beta(b)\}$, where $a,b$ are the endpoints of $\gamma$. 
		\item For each $\bigcup_{i=1}^k\Gamma'_i$-crossing $x$, let $\gamma,\gamma'$ be the members of $\bigcup_{i=1}^k\Gamma_i'$ with $x \in \gamma \cap \gamma'$,
			\begin{itemize}
				\item for $t \in \{0,1\}$ and $\gamma'' \in \{\gamma,\gamma'\}$,
					\begin{itemize}
						\item define $v_{x,\gamma'',1-t}=\beta(\gamma''(t))$, if there exists no $\bigcup_{i=1}^k\Gamma_i'$-crossing between $\gamma''(t)$ and $x$ in $\gamma''$; 
						\item otherwise define $v_{x,\gamma'',1-t}=\eta(\Delta_{x,y,\gamma''})_{1-t}$, where $y$ is the $\bigcup_{i=1}^k\Gamma'_i$-crossing between $\gamma''(t)$ and $x$ contained in $\gamma''$ closest to $x$, and
					\end{itemize}
				\item define $Q_{x,0}=\{\eta(\Delta_x)_0,v_{x,\gamma,1},v_{x,\gamma',1}\}$ and $Q_{x,1}=\{\eta(\Delta_x)_1,v_{x,\gamma,0},v_{x,\gamma',0}\}$.
			\end{itemize}
		\item $\Q=\{Q_\gamma: \gamma$ does not contain any $\bigcup_{i=1}^k\Gamma'_i$-crossing$\} \cup \{Q_{x,0},Q_{x,1}: x$ is a $\bigcup_{i=1}^k\Gamma'_i$-crossing$\}$.
	\end{itemize}
By the definition of $\bigcup_{i=1}^k\Gamma'_i$, $\Q$ satisfies the topological feasibility condition.
Since $\lvert Y^* \rvert \leq 2k\sum_{(S,\Omega) \in \Se_2}r(\rho_S)+2\lvert E(H) \rvert \cdot \lvert V(H) \rvert q_{\lvert V(H) \rvert}+2k_0^*$, by Theorem \ref{linkage on surface} and Claim 6, there exists a set $\{G_Q: Q \in \Q\}$ of pairwise disjoint connected subdrawings of $G'''$ with $V(G_Q) \cap Y^* = Q$ for each $Q \in \Q$.
We may assume that each $G_Q$ is minimal in the sense that deleting any edge from $G_Q$ will make $G_Q$ disconnected.
So for each $(S,\Omega) \in \Se_1$, at most $\lvert \overline{\Omega} \rvert-1$ edges in $G'''$ given by $(S,\Omega)$ are contained in $\bigcup_{Q \in \Q}E(G_Q)$.
Since $\Se$ is effective with respect to $(\Se_1,\Se_2)$, $\bigcup_{Q \in \Q} G_Q$ together with $\bigcup_{i=1}^k\bigcup_{(S,\Omega) \in \Se_2}R_{S,i} \cup \bigcup_{x \in X}R_x$ shows that $G$ half-integrally packs $k$ $H$-subdivisions, where for each $v \in V(H)$, their branch vertices corresponding to $v$ are contained in $\{\pi_i(v): i \in [\ell]\}$.
Therefore, Statement 1 holds, a contradiction.
This proves the lemma.
\end{pf}

\begin{lemma} \label{put one more branch vertex into vortex}
For any connected graph $H$, surface $\Sigma$ and positive integers $k,\kappa,\rho$, there exist positive integers $\kappa^*=\kappa^*(H,\Sigma,k,\kappa,\rho),\rho^*=\rho^*(H,\Sigma,k,\kappa,\rho)$ and a constant function $\phi=\phi(H,\Sigma,k,\kappa,\rho)$ such that for every positive integer $\xi$, there exists an integer $\xi^*=\xi^*(H,\Sigma,k,\kappa,\rho,\xi)$ such that for every integer $\theta^*$, there exists an integer $\theta=\theta(H,\Sigma,k,\kappa,\rho,\xi,\theta^*)$ such that the following hold.
Let $G$ be a graph and $\T$ be a tangle in $G$. 
Let $\Se$ be a $\T$-central effective segregation with respect to a $(\kappa,\rho)$-witness $(\Se_1,\Se_2)$ such that there exists a proper $(\Sigma,\theta,\phi,\T)$-arrangement $\alpha$ with respect to $(\Se_1,\Se_2)$ in $\Sigma$.
Let $Z \subseteq V(G)$ with $\lvert Z \rvert \leq \xi$.
Let $\R=(R_v: v \in V(H))$ be a collection of subsets of $V(G)$.
If $t$ is an integer such that for every $\R$-compatible homeomorphic embedding $\pi$ from $H$ into $G$, either $\pi(V(H) \cup E(H)) \cap Z \neq \emptyset$, or $\lvert (\bigcup_{(S,\Omega) \in \Se_2}V(S)) \cap \pi(V(H)) \rvert \geq t$, then there exist $Z^* \subseteq V(G)$ with $\lvert Z^* \rvert \leq \xi^*$ and a $\T$-central effective segregation $\Se^*$ with respect to a $(\kappa^*,\rho^*)$-witness $(\Se_1^*,\Se_2^*)$ with a proper $(\Sigma,\theta^*,\phi,\T)$-arrangement with respect to $(\Se_1^*,\Se_2^*)$ in $\Sigma$ such that one of the following statements holds.
	\begin{enumerate}
		\item $G$ half-integrally packs $k$ $\R$-compatible subdivisions of $H$.
		\item There exist $(A,B) \in \T$ and an $\R$-compatible homeomorphic embedding $\pi^*$ from $H$ into $G$ such that $\pi^*(E(H)) \subseteq A$.
		\item For every $\R$-compatible homeomorphic embedding $\pi$ from $H$ into $G$, either $\pi(V(H) \cup E(H)) \cap Z^* \neq \emptyset$, or $\lvert (\bigcup_{(S,\Omega) \in \Se_2^*}V(S)) \cap \pi(V(H)) \rvert \geq t+1$.
	\end{enumerate}
\end{lemma}

\begin{pf}
Let $H$ be a connected graph, $\Sigma$ a surface, and $k,\kappa,\rho$ positive integers. 
We define the following.
	\begin{itemize}
		\item Let $\theta_{\ref{same imprint}}$, $\mu_{\ref{same imprint}}$, $\theta_{\ref{same imprint}}^*,\xi_{\ref{same imprint}},\lambda_{\ref{same imprint}},\eta_{\ref{same imprint}}$ be the functions $\theta_0^*,\mu^*$ and the numbers $\theta^*,\xi^*,\lambda^*,\eta^*$, respectively, mentioned in Lemma \ref{same imprint} by taking $H=H, k=k,\kappa=\kappa,\rho=\rho$ and $\Sigma=\Sigma$.
		\item Let $\theta_{\ref{gauge}}, \tau_{\ref{gauge}}, \rho_{\ref{gauge}},\lambda_{\ref{gauge}}$ be the numbers $\theta,\tau,\rho^*,\lambda^*$, respectively, mentioned in Lemma \ref{gauge} by taking $H=H,\kappa=\kappa,\rho=\rho,\theta_0=\theta_{\ref{same imprint}}$, $\mu=\mu_{\ref{same imprint}}$ and $\Sigma=\Sigma$.
		\item Define $\phi$ to be the constant function with $\phi(x)=\tau_{\ref{gauge}}$ for all integers $x$.	
		\item Let $M_{\ref{finitely many templates}}$ be the number $M$ mentioned in Lemma \ref{finitely many templates} by taking $\kappa=\kappa,\Sigma=\Sigma$, $(\Delta_1,...,\Delta_\kappa)$ to be any $\kappa$ oriented disks in $\Sigma$ with pairwise disjoint closure, $H=H$ and $\rho=\rho_{\ref{gauge}}$.
		\item Let $c_1= \lvert V(H) \rvert^{\lvert V(H) \rvert} \cdot \lvert E(H) \rvert^{2\lvert V(H) \rvert^2}$.
		\item Let $c=M_{\ref{finitely many templates}}^\kappa \cdot \lambda_{\ref{gauge}}^{2\kappa} \cdot c_1$.
		\item Define $\kappa^*,\rho^*$ to be the numbers $\kappa^*,\rho^*$, respectively, mentioned in Lemma \ref{sweeping balls into vortices} by taking $\kappa=\kappa,k=c \cdot \eta_{\ref{same imprint}}, \rho=\rho, \lambda=\lambda_{\ref{same imprint}}, \phi=\phi$.
	\end{itemize}
Let $\xi$ be a positive integer. 
Define $\xi^*=\xi+c \cdot \xi_{\ref{same imprint}}$.
Let $\theta^*$ be an integer. 
We define the following.
	\begin{itemize}
		\item Define $\kappa^*,\rho^*,\theta_{\ref{sweeping balls into vortices}}$ to be the numbers $\kappa^*,\rho^*,\theta$ mentioned in Lemma \ref{sweeping balls into vortices} by taking $\kappa=\kappa,k=c \cdot \eta_{\ref{same imprint}}, \rho=\rho, \lambda=\lambda_{\ref{same imprint}}, \phi=\phi$ and $\theta^*=\theta^*$.
		\item Define $\theta=\theta_{\ref{same imprint}}^*+\theta_{\ref{gauge}}+\theta_{\ref{sweeping balls into vortices}}+\theta^*+\xi^*$.
	\end{itemize}

Let $G$, $\T$, $\Se$, $\alpha$, $Z$, $\R=(R_v: v \in V(H))$ and $t$ be the ones mentioned in the statement of this lemma.

Suppose that this lemma does not hold.
Let $G'$ be the skeleton of $\alpha$ with respect to $(\Se_1,\Se_2)$.

By Lemma \ref{gauge}, there exists an $(H,\theta_{\ref{same imprint}},\mu_{\ref{same imprint}})$-gauge $\{\C_S: (S,\Omega) \in \Se_2\}$ with respect to $(\Se_1,\Se_2,\alpha)$ such that for each $(S,\Omega) \in \Se_2$, there are at most $\lambda_{\ref{gauge}}^2$ different pairs of cycles $(C_S,C_S')$ mentioned in Lemma \ref{gauge} among all homeomorphic embeddings from $H$ into $G$.
By Lemma \ref{finitely many templates}, for each $(S,\Omega) \in \Se_2$, there are at most $M_{\ref{finitely many templates}}$ different $\alpha(S,\Omega)$-templates with depth at most $\rho_{\ref{gauge}}$.
Furthermore, there are at most $c_1$ different sets $\{(H_x,\Omega_{H_x}): x \in X_\pi\}$ mentioned in the definition of an $(H,\theta_{\ref{same imprint}},\mu_{\ref{same imprint}})$-snapshots of $\pi$ with respect to $\{\C_S: (S,\Omega) \in \Se_2\}$ among all homeomorphic embeddings $\pi$ from $H$ to $G$.
We fix an ordering $\sigma$ of $\Se_2$.
So there are at most $c$ different $(H,\theta_{\ref{same imprint}},\mu_{\ref{same imprint}})$-snapshots of $\pi$ with respect to $\{\C_S: (S,\Omega) \in \Se_2\}$ and $\sigma$ among all homeomorphic embeddings $\pi$ from $H$ into $G$.
Therefore, the set of $\R$-compatible homeomorphic embeddings $\pi$ from $H$ into $G$ can be partitioned into $c$ (not necessarily non-empty) subsets $A_1,A_2,...,A_c$ such that for each $i \in [c]$, the $\R$-compatible homeomorphic embeddings in $A_i$ have the same $(H,\theta_{\ref{same imprint}},\mu_{\ref{same imprint}})$-snapshot with respect to $\{\C_S:(S,\Omega) \in \Se_2\}$ and $\sigma$. 

For each $i \in [c]$, applying Lemma \ref{same imprint} by taking $\{\pi_1,...,\pi_\ell\}$ to be $A_i$, since Statements 1 and 2 of this lemma do not hold, either 
	\begin{itemize}
		\item there exists $Z_i \subseteq V(G)$ with $\lvert Z_i \rvert \leq \xi_{\ref{same imprint}}$ such that $Z_i \cap \pi(E(H)) \neq \emptyset$ for each $\pi \in A_i$, or 
		\item there exist $W_i \subseteq V(G')$ with $\lvert W_i \rvert \leq \eta_{\ref{same imprint}}$ and a set $\{\Lambda_w: w \in W_i\}$, where each $\Lambda_w$ is a $\lambda_{\ref{same imprint}}$-zone around some $w \in W_i$ in $G'$ such that for each $\pi \in A_i$, $\bigcup_{(S,\Omega)\in \Se_1,\alpha(S,\Omega) \subseteq \bigcup_{w \in W_i}\overline{\Lambda_w}}V(S)$ contains $\pi(v_i)$ for some $v_i \in V(H)$ with $\pi(v_i) \not \in V(\bigcup_{(S,\Omega)\in\Se_2}S)$.
	\end{itemize}
Define $Z^*$ to be the union of $Z$ and $Z_i$ for all $i \in [c]$ in which $Z_i$ defined.
So $\lvert Z^* \rvert \leq \xi^*$.
Define $\W^*$ to be the union of the sets $\{\Lambda_w: w \in W_i\}$ over all numbers $i \in [c]$ in which $W_i$ is defined.
So $\lvert \W^* \rvert \leq c\cdot \eta_{\ref{same imprint}}$.

By Lemma \ref{sweeping balls into vortices}, there exists a $\T$-central segregation $\Se^*$ with a $(\kappa^*,\rho^*)$-witness $(\Se_1^*,\Se_2^*)$ such that $\Se_1^* \subseteq \Se_1$, $(\bigcup_{(S,\Omega) \in \Se_2}S) \cup (\bigcup_{(S,\Omega) \in \Se_1, \alpha(S,\Omega) \subseteq \bigcup_{\Lambda \in \W^*}\overline{\Lambda}}S) \subseteq \bigcup_{(S,\Omega) \in \Se_2^*}S$, and $\Se^*$ has a $(\Sigma, \theta^*,\phi,\T)$-arrangement with respect to $(\Se_1^*,\Se_2^*)$ in $\Sigma$.
Since $\Se_1^* \subseteq \Se_1$, $\Se^*$ is effective respect to $(\Se_1^*,\Se_2^*)$.
Furthermore, if $\pi$ is an $\R$-compatible homeomorphic embedding from $H$ into $G$ with $\pi(V(H) \cup E(H)) \cap Z^*=\emptyset$, then $\pi \in A_i$ for some $i \in [c]$ with $W_i$ defined, so some vertex $v \in V(H)$ with $\pi(v) \not \in V(\bigcup_{(S,\Omega)\in\Se_2}S)$ has the property that $\pi(v) \in \bigcup_{(S,\Omega)\in \Se_1,\alpha(S,\Omega)\subseteq \bigcup_{w \in W_i}\overline{\Lambda_w}}V(S) \subseteq \bigcup_{(S,\Omega) \in \Se_2^*}V(S)$.
Therefore, for every $\R$-compatible homeomorphic embedding $\pi$ from $H$ into $G$ with $\pi(V(H) \cup E(H)) \cap Z^* = \emptyset$, we have $\lvert (\bigcup_{(S,\Omega) \in \Se_2^*}V(S)) \cap \pi(V(H)) \rvert - \lvert (\bigcup_{(S,\Omega) \in \Se_2}V(S)) \cap \pi(V(H)) \rvert \geq 1$.
So Statement 3 of this lemma holds, a contradiction.
\end{pf}

\begin{lemma} \label{arrangement EP}
For any connected graph $H$, surface $\Sigma$ and positive integers $k,\kappa,\rho$, there exist positive integers $\theta,\xi$ and a constant function $\phi$ such that the following hold.
Let $G$ be a graph, $\T$ a tangle in $G$, and $\R=(R_v: v \in V(H))$ a collection of subsets of $V(G)$.
If there exists a $\T$-central effective segregation $\Se$ of $G$ with respect to a $(\kappa,\rho)$-witness $(\Se_1,\Se_2)$ such that there exists a proper $(\Sigma, \theta,\phi,\T)$-arrangement with respect to $(\Se_1,\Se_2)$ in $\Sigma$, then one of the following statements holds.
	\begin{enumerate}
		\item $G$ half-integrally packs $k$ $\R$-compatible subdivisions of $H$.
		\item There exists $Z \subseteq V(G)$ with $\lvert Z \rvert \leq \xi$ such that $\pi(V(H) \cup E(H)) \cap Z \neq \emptyset$ for every $\R$-compatible homeomorphic embedding from $H$ into $G$.
		\item There exist $(A,B) \in \T$ and an $\R$-compatible homeomorphic embedding $\pi^*$ from $H$ into $G$ such that $\pi^*(E(H)) \subseteq A$.
	\end{enumerate}
\end{lemma}

\begin{pf}
Let $H$ be a connected graph, $\Sigma$ be a surface, and $k,\kappa,\rho$ be positive integers.
Let $\kappa_0=\kappa$, $\rho_0=\rho$ and $\xi_0=0$.
For every $i \geq 1$, let $\kappa_i,\rho_i,\xi_i,\phi_i$ be the numbers $\kappa^*,\rho^*,\xi^*$ and function $\phi$, respectively, mentioned in Lemma \ref{put one more branch vertex into vortex} by taking $H=H, \Sigma=\Sigma,k=k,\kappa=\kappa_{i-1},\rho=\rho_{i-1},\xi=\xi_{i-1}$.
Define $\xi=\xi_{\lvert V(H) \rvert+1}$ and $\phi=\sum_{i=1}^{\lvert V(H) \rvert+1}\phi_{i}$.

Let $\theta_{\lvert V(H) \rvert+1}=1$.
For $i \in [\lvert V(H) \rvert+1]$, let $\theta_{i-1}$ be the number $\theta$ mentioned in Lemma \ref{put one more branch vertex into vortex} by taking $H=H,\Sigma=\Sigma,k=k,\kappa=\kappa_{i-1},\rho=\rho_{i-1},\xi=\xi_{i-1},\theta^*=\theta_i$.
Define $\theta=\theta_0$.

Let $G$, $\T$, $\R=(R_v: v \in V(H))$ and $\Se$ be the ones mentioned in this lemma.
Assume that Statements 1 and 3 of this lemma do not hold.
Let $Z_0=\emptyset$.
Note that $\Se$ has the property that for every $\R$-compatible homeomorphic embedding $\pi$ from $H$ into $G$, either $\pi(V(H) \cup E(H)) \cap Z_0 \neq \emptyset$ or $\lvert (\bigcup_{(S,\Omega) \in \Se_2}V(S)) \cap \pi(V(H)) \rvert \geq 0$.
Denote $\Se$ by $\Se^{(0)}$.
We can inductively apply Lemma \ref{put one more branch vertex into vortex} to $\Se^{(i-1)}$ for each $i \in [\lvert V(H) \rvert+1]$ to construct a set $Z_i \subseteq V(G)$ with $\lvert Z_i \rvert \leq \xi_i$ and a $\T$-central effective segregation $\Se^{(i)}$ with respect to a $(\kappa_i,\rho_i)$-witness $(\Se^{(i)}_1,\Se^{(i)}_2)$ with a proper $(\Sigma,\theta_i,\phi,\T)$-arrangement with respect to  $(\Se^{(i)}_1,\Se^{(i)}_2)$ in $\Sigma$ such that for every $\R$-compatible homeomorphic embedding $\pi$ from $H$ into $G$, either $\pi(V(H) \cup E(H)) \cap Z_i \neq \emptyset$ or $\lvert (\bigcup_{(S,\Omega) \in \Se_2^{(i)}}V(S)) \cap \pi(V(H)) \rvert \geq i$.

For every $\R$-compatible homeomorphic embedding $\pi$ from $H$ into $G$, since $\lvert \pi(V(H)) \rvert=\lvert V(H) \rvert$, we have $\pi(V(H) \cup E(H)) \cap Z_{\lvert V(H) \rvert+1} \neq \emptyset$.
This proves Statement 2 by taking $Z=Z_{\lvert V(H) \rvert+1}$.
\end{pf}

\section{Half-integral Erd\H{o}s-P\'{o}sa property} \label{sec: half-integral EP}

In this section, we will prove Theorem \ref{half EP} by combining results proved in earlier sections. 

The following is a restatement of a structure theorem for excluding a fixed graph as a minor proved in \cite{d}, which is a stronger form of a theorem in \cite{rs XVI}.

\begin{theorem}[{\cite[Theorem 7]{d}}] \label{D minor structure}
For every graph $L$, there exists a positive integer $\kappa$, such that for every nondecreasing function $\phi$, there exist integers $\theta,\xi,\rho$ such that if $G$ is a graph with a tangle $\T$ with order at least $\theta$ controlling no $L$-minor, then there exist $Z \subseteq V(G)$ with $\lvert Z \rvert \leq \xi$ and an effective $(\T-Z)$-central segregation $\Se$ of $G-Z$ with respect to a $(\kappa,\rho)$-witness $(\Se_1,\Se_2)$ such that $\Se$ has a proper $(\Sigma,\phi(p),\phi,\T-Z)$-arrangement with respect to $(\Se_1,\Se_2)$, where $\Sigma$ is a surface in which $L$ cannot be drawn and $p$ is the minimum such that every member of $\Se_2$ is a vortex of depth at most $p$.
\end{theorem}

\begin{lemma} \label{connected one side}
For any connected graph $H$ and positive integer $k$, there exist positive integers $\theta,\xi$ such that if $G$ is a graph with a tangle $\T$ of order at least $\theta$, and $\R=(R_v: v \in V(H))$ is a collection of subsets of $V(G)$, then one of the following statements holds.
	\begin{enumerate}
		\item $G$ half-integrally packs $k$ $\R$-compatible subdivisions of $H$.
		\item There exists $Z \subseteq V(G)$ with $\lvert Z \rvert \leq \xi$ such that $\pi(V(H) \cup E(H)) \cap Z \neq \emptyset$ for every $\R$-compatible homeomorphic embedding $\pi$ from $H$ into $G$.
		\item There exist $(A,B) \in \T$ and an $\R$-compatible homeomorphic embedding $\pi^*$ from $H$ into $G$ such that $\pi^*(E(H)) \subseteq A$.
	\end{enumerate}
\end{lemma}

\begin{pf}
Let $H$ be a connected graph, and let $k$ be a positive integer.
Let $\theta_{\ref{with clique minor}},t_{\ref{with clique minor}},\xi_{\ref{with clique minor}}$ be the numbers $\theta,t,\xi$ mentioned in Lemma \ref{with clique minor} by taking $H=H$ and $k=k$.
Let $\kappa_{\ref{D minor structure}}$ be the number $\kappa$ mentioned in Theorem \ref{D minor structure} by taking $L=K_{t_{\ref{with clique minor}}}$.
For each surface $\Sigma$ in which $K_{t_{\ref{with clique minor}}}$ cannot be drawn and for each positive integer $\rho$, let $\theta_{\ref{arrangement EP}}^{(\Sigma,\rho)},\xi_{\ref{arrangement EP}}^{(\Sigma,\rho)}$ and $\phi_{\ref{arrangement EP}}^{(\Sigma,\rho)}$ be the numbers $\theta,\xi$ and function $\phi$, respectively, mentioned in Lemma \ref{arrangement EP} by taking $H=H,\Sigma=\Sigma, k=k,\kappa=\kappa_{\ref{D minor structure}}$ and $\rho=\rho$.
For each positive integer $\rho$, let $\theta_{{\ref{arrangement EP}}}^{(\rho)}=\sum_{\Xi}\theta_{\ref{arrangement EP}}^{(\Xi,\rho)}$, $\xi_{{\ref{arrangement EP}}}^{(\rho)} = \sum_{\Xi}\xi_{\ref{arrangement EP}}^{(\Xi,\rho)}$ and $\phi_{\ref{arrangement EP}}^{(\rho)} = \sum_{\Xi}\phi_{\ref{arrangement EP}}^{(\Xi,\rho)}$, where these sums are over all surfaces $\Xi$ in which $K_{t_{\ref{with clique minor}}}$ cannot be drawn.
Let $\phi$ be the function with domain ${\mathbb Z}$ such that $\phi(x)=\theta_{\ref{arrangement EP}}^{(x)} + \xi_{\ref{arrangement EP}}^{(x)} + \phi_{\ref{arrangement EP}}^{(x)}(x)$ for every positive integer $x$, and $\phi(x)=0$ for every non-positive integer $x$.
Let $\theta_{\ref{D minor structure}},\xi_{\ref{D minor structure}},\rho_{\ref{D minor structure}}$ be the numbers $\theta,\xi,\rho$ mentioned in Theorem \ref{D minor structure} by taking $L=K_{t_{\ref{with clique minor}}}$ and $\phi=\phi$.
Define $\theta=\theta_{\ref{with clique minor}}+\theta_{\ref{D minor structure}}$ and $\xi=\xi_{\ref{with clique minor}}+\xi_{\ref{D minor structure}}+\sum_{p=1}^{\rho_{\ref{D minor structure}}}\xi_{\ref{arrangement EP}}^{(p)}$.

Let $G$ be a graph and let $\T$ be a tangle in $G$ with order at least $\theta$.
Let $\R=(R_v: v \in V(H))$ be a collection of subsets of $V(G)$.

If $\T$ controls a $K_{t_{\ref{with clique minor}}}$-minor, then one of the statements of this lemma holds by Lemma \ref{with clique minor}.
So we may assume that $\T$ does not control a $K_{t_{\ref{with clique minor}}}$-minor.
By Theorem \ref{D minor structure}, there exist $Z_0 \subseteq V(G)$ with $\lvert Z_0 \rvert \leq \xi_{\ref{D minor structure}}$ and a $(\T-Z_0)$-central effective segregation $\Se$ of $G-Z_0$ with respect to a $(\kappa_{\ref{D minor structure}},\rho_{\ref{D minor structure}})$-witness $(\Se_1,\Se_2)$ such that there exists a proper $(\Sigma,\phi(p),\phi,\T-Z_0)$-arrangement $\alpha$ with respect to $(\Se_1,\Se_2)$ in a surface $\Sigma$ in which $K_{t_{\ref{with clique minor}}}$ cannot be drawn, where $p$ is the minimum such that every member of $\Se_2$ is a $p$-vortex.
Hence $(\Se_1,\Se_2)$ is a $(\kappa_{\ref{D minor structure}},p)$-witness of $\Se$, and for all distinct members $(S_1,\Omega_1),(S_2,\Omega_2) \in \Se_2$, $m_{\T'}(\overline{\Omega_1},\overline{\Omega_2}) \geq \phi(p) \geq \phi_{\ref{arrangement EP}}^{(\Sigma,p)}(p)$, where $\T'$ is the tangle of the skeleton of $\alpha$ with respect to $(\Se_1,\Se_2)$ that is conformal with $\T-Z_0$ mentioned in the definition of a $(\Sigma,\phi(p),\phi,\T-Z_0)$-arrangement.
Note that $\phi(p) \geq \theta_{\ref{arrangement EP}}^{(p)} \geq \theta_{\ref{arrangement EP}}^{(\Sigma,p)}$.
So $\alpha$ is a $(\Sigma,\theta_{\ref{arrangement EP}}^{(\Sigma,p)},\phi_{\ref{arrangement EP}}^{(\Sigma,p)},\T-Z_0)$-arrangement of $\Se$ with respect to $(\Se_1,\Se_2)$ in $\Sigma$.
Therefore, one of the statements of this lemma holds by Lemma \ref{arrangement EP}.
\end{pf}

\bigskip

Now we are ready to prove Theorem \ref{half EP} for connected graphs $H$.
The following is a restatement.

\begin{lemma} \label{connected half EP}
For every connected graph $H$, there exists a function $f$ such that if $G$ is a graph and $\R=(\R_v: v \in V(G))$ is a collection of subsets of $V(G)$, then for every positive integer $k$, either $G$ half-integrally packs $k$ $\R$-compatible subdivisions of $H$, or there exists $Z \subseteq V(G)$ with $\lvert Z \rvert \leq f(k)$ such that $\pi(V(H) \cup E(H)) \cap Z \neq \emptyset$ for every $\R$-compatible homeomorphic embedding $\pi$ from $H$ into $G$.
\end{lemma}

\begin{pf}
Let $H$ be a connected graph.
For every positive integer $k$, let $\theta_k$ and $\xi_k$ be the numbers $\theta,k$ mentioned in Lemma \ref{connected one side} by taking $H=H$ and $k=k$.
Define $f$ to be the function with domain ${\mathbb N} \cup \{0\}$ such that $f(0)=0$, and for $i \geq 1$, define $f(i)=2f(i-1)+3\theta_i+\xi_i$.

Now we prove by induction on $k$ that either $G$ half-integrally packs $k$ $\R$-compatible subdivisions of $H$, or there exists $Z \subseteq V(G)$ with $\lvert Z \rvert \leq f(k)$ such that $\pi(V(H) \cup E(H)) \cap Z \neq \emptyset$ for every $\R$-compatible homeomorphic embedding $\pi$ from $H$ into $G$.
This statement clearly holds if $k=0$.
So we may assume that $k \geq 1$ and this statement holds for every nonnegative integer $k'$ with $k'<k$.

Suppose to the contrary that this statement does not hold.
So $G$ does not half-integrally packs $k$ $\R$-compatible subdivisions of $H$, and there does not exist a subset of $V(G)$ with size at most $f(k)$ intersecting all $\R$-compatible subdivisions of $H$ in $G$.

Suppose that $E(H)=\emptyset$.
Since $H$ is connected, $\lvert V(H) \rvert=1$. 
Since $G$ does not half-integrally pack $k$ $\R$-compatible subdivisions of $H$, the unique member $R$ of $\R$ has at most $k-1 \leq f(k)$ vertices.
But $R$	intersects all $\R$-compatible subdivisions of $H$ in $G$, a contradiction.

Hence $E(H) \neq \emptyset$.
Since $H$ is connected, $\pi(V(H)) \subseteq \pi(E(H))$.

\noindent{\bf Claim 1:} For every separation $(C,D)$ of $G$ of order less than $\theta_k$, exactly one of $C-V(D)$ and $D-V(C)$ contains an $\R$-compatible subdivision of $H$.

\noindent{\bf Proof of Claim 1:}
Let $(A,B)$ be a separation of $G$ of order less than $\theta_k$.
If neither $A-V(B)$ nor $B-V(A)$ contains any $\R$-compatible subdivision of $H$, then since $H$ is connected, $G-V(A \cap B)$ does not contain any $\R$-compatible subdivision of $H$, so $V(A \cap B)$ is a set with size at most $\theta_k \leq f(k)$ such that $\pi(V(H) \cup E(H)) \cap V(A \cap B) \neq \emptyset$ for every $\R$-compatible homeomorphic embedding $\pi$ from $H$ into $G$, a contradiction.
So at least one of $A-V(B)$ and $B-V(A)$ contains an $\R$-compatible subdivision of $H$.

Suppose that both $A-V(B)$ and $B-V(A)$ contain $\R$-compatible subdivisions of $H$.
Then $B$ (and $A$, respectively) does not half-integrally pack $k-1$ $\R$-compatible subdivisions of $H$, for otherwise $G$ half-integrally packs $k$ $\R$-compatible subdivisions of $H$.
By the induction hypothesis, there exist $Z_A,Z_B \subseteq V(G)$ with size at most $f(k-1)$ such that $\pi_A(E(H)) \cap Z_A \neq \emptyset \neq \pi_B(E(H)) \cap Z_B$ for any $\R$-compatible homeomorphic embedding $\pi_A$ from $H$ into $A$ and any $\R$-compatible homeomorphic embedding $\pi_B$ from $H$ into $B$.
Since $H$ is connected, $\pi(V(H) \cup E(H)) \cap (Z_A \cup Z_B \cup V(A \cap B)) \neq \emptyset$ for every $\R$-compatible homeomorphic embedding $\pi$ from $H$ into $G$.
But $\lvert Z_A \cup Z_B \cup V(A \cap B) \rvert \leq 2f(k-1)+\theta_k \leq f(k)$, a contradiction.
So the claim holds.
$\Box$

Define $\T$ to be the set of separations of $G$ of order less than $\theta_k$ such that a separation $(C,D)$ of $G$ of order less than $\theta_k$ belongs to $\T$ if and only if $D-V(C)$ contains an $\R$-compatible subdivision of $H$ but $C-V(D)$ does not.

\noindent{\bf Claim 2:} $\T$ is a tangle in $G$ of order $\theta_k$.

\noindent{\bf Proof of Claim 2:}
Claim 1 implies that $\T$ satisfies (T1).
It suffices to prove that $\T$ satisfies (T2) and (T3).

Suppose that $\T$ does not satisfy (T2).
So there exist $(A_i,B_i) \in \T$ for $i \in [3]$ such that $A_1 \cup A_2 \cup A_3=G$.
Let $Z_i=V(A_i \cap B_i)$ for $i \in [3]$.
We claim that $\pi(V(H) \cup E(H)) \cap (Z_1 \cup Z_2 \cup Z_3) \neq \emptyset$ for every $\R$-compatible homeomorphic embedding $\pi$ from $H$ into $G$.
If $\pi(V(H) \cup E(H)) \cap Z_i = \emptyset$ for some $i \in [3]$, then $\pi(E(H)) \subseteq B_i-V(A_i)$ by the definition of $\T$ and Claim 1, since $H$ is connected.
But $A_1 \cup A_2 \cup A_3=G$, so $\pi(V(H) \cup E(H)) \cap (Z_1 \cup Z_2 \cup Z_3) \neq \emptyset$. 
Hence $\pi(V(H) \cup E(H)) \cap (Z_1 \cup Z_2 \cup Z_3) \neq \emptyset$ for every $\R$-compatible homeomorphic embedding $\pi$ from $H$ into $G$.
However, $\lvert Z_1 \cup Z_2 \cup Z_3 \rvert \leq 3\theta_k \leq f(k)$, a contradiction.
So $\T$ satisfies (T2).

Suppose that $\T$ does not satisfy (T3).
Then there exists a separation $(C,D) \in \T$ such that $V(C)=V(G)$.
Since $C-V(D)$ does not contain an $\R$-compatible subdivision of $H$, $V(C \cap D)$ is a set with size at most $\theta_k \leq f(k)$ intersecting every $\R$-compatible subdivision of $H$ in $G$, a contradiction.
So $\T$ satisfies (T3) and is a tangle of order $\theta_k$.
$\Box$

By Lemma \ref{connected one side}, since $\xi_k \leq f(k)$, there exist $(A^*,B^*) \in \T$ and an $\R$-compatible homeomorphic embedding $\pi^*$ from $H$ into $G$ such that $\pi^*(E(H)) \subseteq A^*$.
Since $\pi(E(H))=\pi(V(H) \cup E(H))$, $B^*-V(A^*)$ does not half-integrally pack $k-1$ $\R$-compatible subdivisions of $H$.
By the induction hypothesis, there exists $Z' \subseteq V(B^*)-V(A^*)$ with $\lvert Z' \rvert \leq f(k-1)$ such that $\pi(V(H) \cup E(H)) \cap Z' \neq \emptyset$ for every $\R$-compatible homeomorphic embedding $\pi$ from $H$ into $B^*-V(A^*)$.
Since $H$ is connected and $A^*-V(B^*)$ does not contain any $\R$-compatible subdivision of $H$, $V(A^* \cap B^*) \cup Z'$ intersects every $\R$-compatible subdivision of $H$ in $G$.
But $\lvert V(A^* \cap B^*) \cup Z' \rvert \leq \theta_k+f(k-1) \leq f(k)$, a contradiction.
This proves the lemma.
\end{pf}

\bigskip

Now we are ready to prove Theorem \ref{half EP}.
The following is a restatement.

\begin{theorem} \label{final half EP restatement}
For every graph $H$, there exists a function $f$ such that for every graph $G$, collection $\R=(R_v: v \in V(H))$ of subsets of $V(G)$ and positive integer $k$, either $G$ half-integrally packs $k$ $\R$-compatible subdivisions of $H$, or there exists a set $Z \subseteq V(G)$ with $\lvert Z \rvert \leq f(k)$ such that $\pi(V(H) \cup E(H)) \cap Z \neq \emptyset$ for every $\R$-compatible homeomorphic embedding $\pi$ from $H$ into $G$.
\end{theorem}

\begin{pf}
Let $H$ be a graph.
Let $c$ be the number of components of $H$.
We shall prove this theorem by induction on $c$.
The theorem holds if $c=1$ by Lemma \ref{connected half EP}.
So we may assume that $c \geq 2$ and this theorem holds for every graph with less than $c$ components.

Define $\HH$ to be the family of the graphs that can be obtained from $H$ by applying one of the following operations.
	\begin{itemize}
		\item Add an edge incident with two vertices in different components.
		\item Subdivide an edge in a component of $H$ to obtain a new vertex $v$, and add an edge incident with $v$ and a vertex in another component.
		\item Pick two different components of $H$, subdivide one edge in each of them, and add an edge incident with these two new vertices.
	\end{itemize}
Note that $\lvert \HH \rvert \leq (\lvert V(H) \rvert + \lvert E(H) \rvert)^2$, and every graph in $\HH$ has less than $c$ components.
We may assume that $V(H') \supseteq V(H)$ for every $H' \in \HH$.
Define $\HH'$ to be the family of graphs, where each member of $\HH'$ can be obtained from $H$ by deleting $c'$ components for some $c' \geq 1$.
We may assume that $V(H') \subseteq V(H)$ for every $H' \in \HH'$.

By the induction hypothesis, for each $H' \in \HH \cup \HH'$, there exists a function $f_{H'}$ such that for any graph graph $G$, collection $\R=(R_v: v \in V(H'))$ of subsets of $V(G)$ and positive integer $k$, either $G$ half-integrally packs $k$ $\R$-compatible subdivisions of $H'$, or there exists a set $Z_{H'} \subseteq V(G)$ with $\lvert Z_{H'} \rvert \leq f_{H'}(k)$ such that $\pi(V(H') \cup E(H')) \cap Z_{H'} \neq \emptyset$ for every $\R$-compatible homeomorphic embedding $\pi$ from $H'$ into $G$.

Define $f$ to be the function such that $f(\cdot)=\sum_{H' \in \HH} f_{H'}(\cdot) + (kc)^c\sum_{H'\in \HH'}f_{H'}(\cdot)$. 

Suppose to the contrary that $f$ does not satisfy the conclusion of this theorem.
So there exist a graph $G$, a collection $\R=(R_v: v \in V(H))$ of subsets of $V(G)$ and a positive integer $k$ such that $G$ does not half-integrally pack $k$ $\R$-compatible subdivisions of $H$, and there does not exist a subset of $V(G)$ with size at most $f(k)$ intersecting every $\R$-compatible subdivision of $H$ in $G$.

For every $H' \in \HH \cup \HH'$, define $\R^{H'}$ to be the collection $(R^{H'}_v: v \in V(H'))$, where $R^{H'}_v=R_v$ if $v \in V(H') \cap V(H)$, and $R^{H'}_v=V(G)$ otherwise.
Since each graph $H' \in \HH$ contains a subdivision of $H$, if $G$ half-integrally packs $k$ $\R^{H'}$-compatible subdivisions of $H'$ for some $H' \in \HH$, then $G$ half-integrally packs $k$ $\R$-compatible subdivisions of $H$, a contradiction.
So for each $H' \in \HH$, there exists $Z_{H'} \subseteq V(G)$ with $\lvert Z_{H'} \rvert \leq f_{H'}(k)$ such that $G-Z_{H'}$ does not contain any $\R^{H'}$-compatible subdivision of $H'$.
Let $Z_0=\bigcup_{H' \in \HH} Z_{H'}$.
Note that $\lvert Z_0 \rvert \leq \sum_{H' \in \HH}f_{H'}(k) \leq f(k)$.
So there exists an $\R$-compatible subdivision of $H$ in $G-Z_0$.

Let $G'$ be an $\R$-compatible subdivision of $H$ in $G-Z_0$.
Suppose that some component $C$ of $G-Z_0$ contains two components $C_1,C_2$ of $G'$.
Then there exists a path $P$ in $C$ from $V(C_1)$ to $V(C_2)$ internally disjoint from $V(C_1) \cup V(C_2)$.
But $G' \cup P$ is an $\R^{H'}$-compatible subdivision of some graph $H' \in \HH$, a contradiction.
Hence for every $\R$-compatible subdivision of $H$ in $G-Z_0$, every component of $G-Z_0$ contains at most one component of this subdivision of $H$.

Let $G_1,G_2,...,G_t$ be the components of $G-Z_0$.
Let $H_1,H_2,...,H_c$ be the components of $H$.
Define a bipartite graph $W$ with a bipartition $(L,R)$, where $V(L)=\{u_{i,j}: i \in [c], j \in [k]\}$ and $V(R)=\{v_i: i \in [t]\}$, such that for all $i \in [c], j \in [k], \ell \in [t]$, $u_{i,j}$ is adjacent in $W$ to $v_\ell$ if and only if $G_\ell$ contains an $\R^{H_i}$-compatible subdivision of $H_i$.
Note that if $W$ has a matching saturating $L$, then for each $i \in [c]$, there exist a set of $k$ distinct components of $G-Z_0$, where each member of this set contains an $\R^{H_i}$-compatible subdivision of $H_i$, and the sets of $k$ components are disjoint for different $i$, so $G-Z_0$ contains $k$ disjoint $\R$-compatible subdivisions of $H$, a contradiction.
Hence $W$ does not have a matching saturating $L$.
By Hall's theorem, there exists $L' \subseteq L$ such that $\lvert \{v \in R: uv \in E(W), u \in L'\} \rvert < \lvert L' \rvert$.
We choose such $L'$ to be as large as possible.
The maximality of $L'$ implies that if $u_{i,j} \in L'$ for some $i \in [c]$ and $j \in [k]$, then $u_{i,j'} \in L'$ for every $j' \in [k]$.
Let $H_{L'}$ be the graph $\bigcup_{i \in [c], u_{i,1} \in L'}H_i$, and let $H_{L-L'}$ be the graph $H-V(H_{L'})$.

Let $R'=\{v \in R: uv \in E(W), u \in L'\}$.
Let $G_{R'}=\bigcup_{\ell \in [t],v_\ell \in R'}G_\ell$, and let $G_{R-R'}=(G-Z_0)-V(G_{R'})$.
Then $W-(L' \cup R')$ has a matching saturating $L-L'$ by the maximality of $L'$ and Hall's theorem.
Hence $G_{R-R'}$ contains $k$ disjoint $\R^{H_{L-L'}}$-compatible subdivisions of $H_{L-L'}$.
So $G_{R'}$ does not half-integrally pack $k$ $\R$-compatible subdivisions of $H_{L'}$.
Furthermore, by the definition of $R'$, if there exists a set $Z \subseteq V(G_{R'})$ such that $Z$ intersects every $\R^{H_{L'}}$-compatible subdivision of $H_{L'}$ in $G_{R'}$, then $Z$ intersects every $\R$-compatible subdivision of $H$ in $G-Z_0$.

Suppose that $L'$ is a proper subset of $L$.
Then $H_{L'}$ has less components than $H$ and belongs to $\HH'$, so there exists $Z' \subseteq V(G_{R'})$ with $\lvert Z' \rvert \leq f_{H_{L'}}(k)$ such that $Z'$ intersects every $\R$-compatible subdivision of $H_{L'}$ in $G-Z_0$.
Hence $Z_0 \cup Z'$ is a set with size at most $f(k)$ intersecting every $\R$-compatible subdivision of $H$ in $G$, a contradiction.

So $L'=L$.
Therefore, $G_{R'}$ is a graph with at most $\lvert L' \rvert-1 \leq kc-1$ components. 
Without loss of generality, let $R'=[\lvert R' \rvert]$.
For each injection $\iota: [c] \rightarrow [\lvert R' \rvert]$, there exists $i_\iota \in [c]$ such that $G_{\iota(i_\iota)}$ does not half-integrally pack $k$ $\R^{H_{i_\iota}}$-compatible subdivisions of $H_{i_\iota}$; otherwise $G_{R'}$ half-integrally packs $k$ $\R$-compatible subdivisions of $H$.
So for each injection $\iota: [c] \rightarrow [\lvert R' \rvert]$, there exists $Z_\iota \subseteq V(G_{\iota(i_\iota)})$ with $\lvert Z_\iota \rvert \leq f_{H_{i_\iota}}(k)$ such that $Z_\iota$ intersects every $\R^{H_{i_\iota}}$-compatible subdivision of $H_{i_\iota}$ in $G_{\iota(i_\iota)}$.
Define $Z$ to be the union of $Z_\iota$ over all injections $\iota$ from $[c]$ to $[\lvert R' \rvert]$.
Since there are at most $\lvert R' \rvert^c \leq (kc)^c$ such injections, and $H_{i_\iota} \in \HH'$ for each such injection $\iota$, $\lvert Z \rvert \leq (kc)^c \sum_{H' \in \HH'}f_{H'}(k)$.

Let $G''$ be any $\R$-compatible subdivision of $H$ in $G-Z_0$.
Since every component of $G-Z_0$ contains at most one component of $G''$, there exists an injection $\iota_{G''}: [c] \rightarrow [\lvert R' \rvert]$ such that for every $i \in [c]$, $G_{\iota_{G''}(i)}$ contains the component of $G''$ that is an $\R^{H_i}$-compatible subdivision of $H_i$.
Then $Z_{\iota_{G''}} \subseteq Z$ intersects $G''$.
Therefore, $Z$ intersects every $\R$-compatible subdivision of $H$ in $G-Z_0$.
In other words, $G-(Z_0 \cup Z)$ does not contain any $\R$-compatible subdivision of $H$.
Consequently, $Z_0 \cup Z$ intersects every $\R$-compatible subdivision of $H$ in $G$. 
Note that $\lvert Z_0 \rvert + \lvert Z \rvert \leq \sum_{H' \in \HH}f_{H'}(k)+(kc)^c\sum_{H' \in \HH'}f_{H'}(k) \leq f(k)$, a contradiction.
This proves the theorem.
\end{pf}

\bigskip

\noindent{\bf Acknowledgement:} The author thanks Sergey Norin for his encouragement for working on this paper. The author thanks the anonymous reviewer for very careful reading and many constructive suggestions.

\end{document}